\newif\ifconfver
\confverfalse      

\newif\ifcutshort      
\cutshorttrue

\newif\ifcutshortlvltwo  
\cutshortlvltwofalse

\ifconfver
\documentclass[10pt,twocolumn,twoside]{IEEEtran}
\else
\documentclass[11pt, draftclsnofoot, onecolumn]{IEEEtran}
\fi

\usepackage{blkarray}
\usepackage{bigdelim}
\usepackage[ruled]{algorithm2e}
\usepackage{enumerate}

\usepackage{pifont}
\usepackage[dvipsnames]{xcolor}
\usepackage{cite}
\usepackage{float}
\usepackage{threeparttable}

\usepackage{amsmath,amsfonts,bm,amssymb,epsfig,psfrag,ifthen,color}
\usepackage{booktabs}
\usepackage{multirow}
\usepackage{array}
\usepackage{calc}
\usepackage{graphicx}
\usepackage{psfrag}

\usepackage{stfloats}
\usepackage{url}
\usepackage[caption=false]{subfig}
\usepackage{longtable}
\usepackage{amsmath,amsfonts,amssymb,amsbsy,bm,paralist,theorem,ifthen,color}
\usepackage{mathtools}

\def\MatrixFont{\bf}
\def\VectorFont{\bf}

\newcommand{\mA}{{\MatrixFont A}}
\newcommand{\mB}{{\MatrixFont B}}

\newcommand{\mD}{{\MatrixFont D}}
\newcommand{\mE}{{\MatrixFont E}}

\newcommand{\mI}{{\MatrixFont I}}

\newcommand{\mO}{{\MatrixFont O}}

\newcommand{\mQ}{{\MatrixFont Q}}

\newcommand{\mW}{{\MatrixFont W}}

\newcommand{\va}{{\VectorFont a}}
\newcommand{\vbb}{{\VectorFont b}}

\newcommand{\vp}{{\VectorFont p}}

\newcommand{\vu}{{\VectorFont u}}
\newcommand{\vv}{{\VectorFont v}}

\newcommand{\vx}{{\VectorFont x}}
\newcommand{\vy}{{\VectorFont y}}
\newcommand{\vz}{{\VectorFont z}}

\newcommand{\vzero}{{\VectorFont 0}}
\newcommand{\vone}{{\VectorFont 1}}

\newcommand{\vlambda}{{\boldsymbol{\lambda}}}
\newcommand{\valpha}{{\boldsymbol{\alpha}}}

\newcommand{\vbeta}{{\boldsymbol{\beta}}}

\newcommand{\green}{\color{SeaGreen}}

\newtheorem{Lemma}{Lemma}

\newtheorem{theorem}{Theorem}
\newtheorem{Def}{Definition}

\newtheorem{Exa}{Example}
\newtheorem{Rmk}{Remark}

\ifconfver
\def\blue{\color{blue}}
\else

\def\blue{\color{blue}}
\fi
\def\red{\color{red}}

\begin{document}

    \bibliographystyle{IEEEtran}

    \title{ {Distributed Scaled Proximal ADMM Algorithms for Cooperative Localization in WSNs}}

    \ifconfver \else {\linespread{1.1} \rm \fi

        \author{\vspace{0.8cm}Mei~Zhang,
        	~Zhiguo~Wang, ~Feng~Yin, \IEEEmembership{Senior Member,~IEEE},
        	~and Xiaojing~Shen\\
 \thanks{Mei Zhang, Zhiguo Wang (corresponding author), and Xiaojing Shen with Department of Mathematics, Sichuan University, Chengdu, Sichuan 610064, China. E-mail: zhang$\_$mei@stu.scu.edu.cn, wangzhiguo@scu.edu.cn, shenxj@scu.edu.cn.}
 \thanks{Feng Yin is with the School of Science and Engineering, Chinese University of Hong Kong (Shenzhen), Shenzhen 518172, China, and also with Shenzhen Research Institute of Big Data, Shenzhen 518172, China (e-mail: yinfeng@cuhk.edu.cn).}
            }

        \maketitle

\vspace{-1.5cm}
\begin{center}
\today
\end{center}\vspace{0.5cm}

\begin{abstract}
Distributed cooperative localization in wireless networks is a challenging problem since it typically requires solving a large-scale nonconvex and nonsmooth optimization problem. In this paper, we reformulate the classic cooperative localization problem as a smooth and constrained nonconvex minimization problem while its loss function is separable over nodes. By utilizing the structure of the reformulation, we propose two novel scaled proximal alternating direction method of multipliers (SP-ADMM) algorithms, which can be implemented in a distributed manner. Compared with the classic semi-definite programming relaxation techniques, the proposed algorithms can provide more accurate position estimates with significantly lower computation complexity. The associated theoretical analysis shows that our algorithms {\blue globally converge to a KKT point} of the reformulated problem and a critical point of the original problem, with a favorable sublinear $\mathcal{O}\left(1/T\right)$ convergence rate, where $T$ is the iteration counter. Numerical experiments have consistently shown that the proposed SP-ADMM algorithms are superior to state-of-the-art methods in terms of localization accuracy and computational complexity across all tested scenarios, varying network size, number of anchors, average number of neighbors, and noise variance levels.
    \\\\
\end{abstract}


        \ifconfver \else \IEEEpeerreviewmaketitle} \fi

    \vspace{-0.5cm}

\section{Introduction}

Wireless sensor networks (WSNs) are widely used to deal with sensitive information in a variety of applications, including healthcare, military, Internet of Things, surveillance, and industrial \cite{singh2020handbook}. In the aforementioned applications, all collected information is meaningful only when the locations of the sensor nodes are accurately known. Therefore, localization is an enabling technique for WSNs. The cooperative localization problem aims to determine multiple sensor locations with the aid of a relatively small portion of anchors whose positions are precisely known and the relative noisy range measurements of adjacent nodes \cite{wymeersch2009cooperative}. The range measurements can be measured in a variety of ways, such as time-of-arrival (TOA) \cite{yin2015cooperative,pun2021local}, time-difference-of-arrival (TDOA) \cite{okello2011comparison}, angle-of-arrival (AOA) \cite{xu2017optimal} and received signal strength (RSS) \cite{yin2017received,jin2020bayesian}. In recent years, there has been a growing interest in estimating target positions through optimization techniques due to the faster response speed of position estimation by these methods and convergence guarantee.

Nonconvex and nonsmooth is the major difficulty of the least squares objective function in the maximum likelihood framework \cite{simonetto2014distributed}. It is hard to ﬁnd an optimal solution with low  computational complexity, and in consequence, most existing methods resort to developing approximate solutions through applying, for instance, the relaxation techniques, to the original nonconvex problem. A popular class of methods is based on semi-definite programming (SDP) relaxation \cite{biswas2006SDP,biswas2006SDPnoisy,wang2008further}, which can be solved by the interior-point algorithm. Although convex relaxation techniques guarantee convergence to a global minimum point, it is not necessarily a critical point of the original nonconvex formulation. In addition, the above convex relaxation methods are implemented in a centralized framework, since all measurements are collected and processed at a processing center. Centralized methods are vulnerable to the malfunction of any single node. There have been increasing efforts in developing distributed sensor network localization methods \cite{srirangarajan2008distributed,di2016next,jin2021exploiting}.

Distributed\footnote{The term distributed in this paper refers to an implementation that sensor nodes can locate themselves and neighboring nodes based on the local network information containing neighbors, without whole network data \cite{simonetto2014distributed,srirangarajan2008distributed}.} methods are able to avoid the major drawbacks of the centralized counterparts for large-scale networks, thus making them well-suited for wireless sensor network localization. Distributed methods not just have lower computational and communication complexity due to scalability but are also more robust to node failures \cite{buehrer2018collaborative,hong2017prox}. This inspired the use of distributed optimization approaches to solve the target positioning problem in large-scale networks. Compared with the centralized SDP methods, the work \cite{soares2015simple} proposed a more practical alternative. By leveraging the concept of convex envelope, the authors developed a method that is both scalable and simple to implement. They also analyzed the total number of iterations required by the algorithm to approach the optimal function value with a high probability for the convexified problem. Despite these advances, it still can be difficult to infer where the solutions are the critical points of the original nonconvex problem (discussed in \cite{piovesan2016cooperative,gur2020alternating}). In view of this, non-relaxed methods may be a workable viewpoint.

In this paper, we focus on the first-order method that solves the original nonconvex localization problem directly. In general, the design of this kind of algorithm first reformulates the problem by exploiting the internal structure of the problem and then develops an efficient optimization algorithm to solve the reformulated problem. For example, a nonconvex sequential greedy (NCSG) optimization algorithm was proposed in \cite{shi2010distributed}. It is also proved that this algorithm owns the convergence guarantee inherited from the non-linear Gauss-Seidel framework \cite{bertsekas1989parallel}. Notwithstanding, the limit point of the generated sequence converges to the KKT point of the reformulated problem, the connection with the original problem is yet not provided.

One efficient way is to reformulate the original localization problem as a two-block nonconvex optimization problem with linear equality constraints, and then a distributed ADMM method is proposed in \cite{erseghe2015distributed,luke2017simple} to solve the resulting nonconvex reformulation. In \cite{wang2021distributed,zhang2020proximal}, they show that the ADMM algorithm may be unstable when the objective function is nonconvex and nonsmooth. In order to obtain a favorable initial point for the method developed in \cite{erseghe2015distributed}, the authors in \cite{piovesan2016cooperative} proposed a hybrid ADMM (ADMM-H) method that contains two stages. In the first stage, a reliable starting point of the target problem is given by exploiting the convex relaxation technique introduced in \cite{soares2015simple}. In the second stage, the solution will be constantly updated by adopting the method developed in \cite{erseghe2015distributed}. While the simulations demonstrated faster convergence, there remain limitations both in the local minimum theoretical guarantee and practical implementation. This is discussed in detail by \cite{gur2020alternating} and our Section \ref{sec: simulation}.

Recently, for the single source localization problem, the authors in \cite{luke2017simple} adopted a simple variational representation of the Euclidean norm, then derived an equivalent smooth reformulation with ball constraint. Later on, it was extended for the multi-source localization problem in \cite{gur2020alternating}, and they proposed an alternating minimization (AM) method to solve the reformulated constrained smooth and nonconvex optimization problem. There are two versions of the AM method in \cite{gur2020alternating}, the fully centralized (AM-FC) and the fully distributed (AM-FD). Due to the sequential\footnote{In this paper, sequential means sensor nodes perform their update calculations in turn, that is, they have to wait for a part of neighbor nodes to complete the update before running their local update steps.} nature of the AM-FD method, it is impractical to apply it to large-scale networks. Therefore, the authors further propose a unifying AM (AM-U) algorithm to remedy this difficulty. The basic idea of AM-U is to divide the sensors in the network into several disjoint clusters and then apply the AM-FD method separately for each cluster. The convergence of the AM-FD can be guaranteed. Since the difficulty of nonconvex and nonsmooth, the convergence rate for the original problem is not provided.

In this work, we aim to exploit the advantages of \cite{erseghe2015distributed} and \cite{luke2017simple} to develop an efficient distributed algorithm with lower computational complexity for the nonconvex and nonsmooth localization problem and establish a theoretical guarantee of its performance.

Our contributions are summarized as follows.

\begin{table*}[!t]
	\resizebox{\textwidth}{!}{
		\begin{threeparttable}[b]
			\caption{Comparisons of Different Algorithms}
			\label{comparisons of different alg}
			\begin{tabular}
				[l]{c|ccccc|c|c|cc}
				\hline&&&&&&\multicolumn{3}{c}{Sensor $i$}&\\\cline{7-10}
				&Convex&First order&Step&Convergence&&Computational& Communication&Storage\\
				Algorithm&Relaxations&Method&Size&Rate&Parallelized&Complexity&Cost &Space\\
				\hline SDP \cite{biswas2006SDP}&\ding{52}&\ding{55}&-&\ding{55}&-&$O(n^3)$&$nN_i$&$\frac{1}{2}\sum_{i=1}^NN_i-|\mathcal{E}_{a}|^\divideontimes+n(N-m)$\\[1mm]
				SF \cite{soares2015simple}&\ding{52}&\ding{52}&fixed&\ding{55}\tnote{$^{\ast}$}&\ding{52}&$O(nN_i)$&$nN_i$&$2n+N_i+1$(parallel method)\\[1mm]
				AM-FD \cite{gur2020alternating}&\ding{55}&\ding{52}&fixed&\ding{55}&\ding{55}& $O(nN_i)$&$nN_i$&$n+nN_i+N_i+1$\\[1mm]
				ADMM-H \cite{piovesan2016cooperative}&hybrid\tnote{$^{\dagger}$}&\ding{55}&-&\ding{55}&\ding{52}&$O(n^3T_i+nN_i)$\tnote{$^{\star}$}&$2nN_i+N_i$&$2N_i(N_i+1)n^2+4nN_i+N_i+9$\\[1mm]
				\textbf{SP-ADMM (ours)}&\ding{55}&\ding{52}&fixed&$\mathcal{O}(1/T)$\tnote{$^{\ddagger}$}&\ding{52}&$O(nN_i)$&$2nN_i$& $4nN_i+N_i+3$ (Algorithm \ref{al1})\\
				\hline
			\end{tabular}
			\begin{tablenotes}
				\item [$\ast$] is due to the convergence analysis of the SF for the convex relaxation problem rather than the original nonconvex problem.
				\item [$\dagger$] Here, ``hybrid" means a two-stage algorithm, including a convex relaxation stage and a nonconvex stage.
				\item [$\star$] $T_i$ refers to the number of iterations required by the nonconvex Newton algorithm to converge.
				\item [$\ddagger$] $T$ refers to the total number of iterations of the algorithm.
				\item[$\divideontimes$] $|\mathcal{E}_{a}|$ denotes the number of edges in which both nodes are anchors. The value of storage space refers to the minimum total storage space required by the SDP per step, as described in reference \cite{ding2021optimal}.
			\end{tablenotes}
		\end{threeparttable}
	}
\end{table*}

\begin{itemize}	
	\item By introducing an auxiliary variable for the nonsmooth Euclidean norm, we reformulate the classic cooperative localization problem as a smooth and constrained nonconvex minimization problem, whose loss function is separable over nodes and has two block variables that eventually leads to a nice optimization structure.
	\item To exploit the nice structure of the new reformulation, we proposed a scaled proximal ADMM (SP-ADMM) algorithm, which is suitable for distributed computing and parallel implementation. Moreover, to further reduce the storage space at each node, we proposed Algorithm \ref{al1} that is a simplified version of Algorithm \ref{al2}. As shown in Table~\ref{comparisons of different alg}, the proposed algorithms enjoy lower computation complexity and storage space when compared with the existing SDP relaxation method and ADMM-H method, respectively.
	 \item By utilizing a novel potential function, {\blue we demonstrates that the global convergence
	 	of the sequence generated by the proposed algorithm, that is, the whole sequence converges to a unique KKT point of the reformulated problem and a critical point of the original problem.}
	 Remarkably, these algorithms also exhibit sublinear convergence, with a convergence rate of $\mathcal{O}\left(1/T\right)$, where $T$ represents the iteration counter. To the best of our knowledge, this is the first result that shows the sublinear convergence rate of a distributed ADMM algorithm for a nonconvex and nonsmooth localization problem.
	\item Numerical experiments conducted on a variety of networks have consistently demonstrated that the proposed algorithms outperform existing methods in terms of both localization accuracy and computational efficiency. The experimental data encompasses scenarios with different network sizes, number of anchors, average number of neighbor nodes, and variance levels of the measurement noise.
\end{itemize}

{\bf Notation:} We use lowercase bold letters to denote vectors, $\vone_{N_i}$ and $\vzero_{N_i}$ to represent $N_i$-dimensional column vector with all ones and all zeros, respectively. Capital bold letters represent matrices, specially, $\mI_{N_i}$ and $\mO_{N_i}$ denote $N_i\times N_i$-dimensional identity matrix and zero matrix, respectively. $\|\vx\|$ is the Euclidean norm of a real vector $\vx$, $\otimes$ denotes the Kronecker product. We use $\text{vec}\left(\vx_{i},i\in\mathcal{N}\right)$ to denote the concatenated vector of $\vx_{i}$ for all $i\in\mathcal{N}$ and $\textbf{D}\text{iag}\left(\vz\right)$ to denote the diagonal matrix with the coefficients of $\vz$ along the diagonal. The projection operator of set $\mathcal{B}^{N_i}$ is defined as:
\begin{displaymath}
\textmd{proj}_{\mathcal{B}^{N_i}}\left(\vu_i^t\right):=\arg\min_{\vu_i\in\mathcal{B}^{N_i}}\frac{1}{2}\|\vu_i-\vu_i^t\|^2.
\end{displaymath}
Lastly, let $f:\mathcal{C}\rightarrow(-\infty,+\infty]$ be a proper closed and convex function and $\mW$ be a positive semi-definite matrix. Then the scaled proximal operator of $f$ is given by
\begin{equation*}
\textmd{prox}_{f}^{\mW}\left(\vz\right):=\arg\min_{\vv\in\mathcal{C}} f\left(\vv\right)+\frac{1}{2}\|\vv-\vz\|^2_{\mW},
\end{equation*}
where $\|\cdot\|_{\mW}$ is the scaled norm induced by $\mW$, i.e., $\|\vz\|_{\mW}^{2}:=\langle\vz,\mW\vz\rangle$ for every $\vz$ in $\mathbb{R}^{\left(2N_i+1\right)n}$.

{\bf Synopsis:} Section \ref{sec: pro reformulate} introduces the reformulated smooth constrained nonconvex minimization problem. The proposed SP-ADMM algorithms are presented in Section \ref{sec: alg develop}. Section \ref{sec: conv analysis} presents the theoretical results of the convergence conditions and convergence rate of the SP-ADMM algorithm. The performance of the proposed SP-ADMM algorithm is illustrated in Section \ref{sec: simulation}, and the conclusion is given in Section \ref{sec:conclusion}.

\section{Problem Formulation}\label{sec: pro reformulate}

\subsection{Problem Statement}

The wireless sensor network that we consider throughout this paper is represented as an undirected and connected graph, $\mathcal{G}=(\mathcal{N},\mathcal{E})$, and the topology is assumed to be known. The node set $\mathcal{N}=\{1,2,\ldots,N\}$ consists of $N$ nodes, some of which are anchors with known true positions collected in the set  $\mathcal{A}=\{\va_{N-m+1},\ldots,\va_{N}\}\subset\mathcal{N}$. For each node $i\in\mathcal{N}$, we define $\mathcal{N}_i=\left\{j\mid\left(i,j\right)\in\mathcal{E}\right\}$ as the set of adjacent nodes to node $i$, and $N_i$ as the cardinality of $\mathcal{N}_i$. The true position of node $i$ is represented by $\vp_i\in\mathbb{R}^n$ for $i\in\mathcal{N}$, and the collection of all node positions is denoted by  $\vp=\text{vec}(\vp_i,i\in\mathcal{N})\in\mathbb{R}^{nN}$. The available noisy range measurement between node $i$ and its adjacent node $j\in\mathcal{N}_i$ is denoted as $d_{i,j}$, and we assume that $d_{i,j}=d_{j,i}$ following the convention in \cite{yin2015cooperative,soares2015simple}. Specifically, the noisy range measurement $d_{i,j}$ can be expressed as in \cite{yousefi2014cooperative,vaghefi2015cooperative},
\begin{equation}\label{m}
d_{i,j}=\|\vp_i-\vp_j\|+w_{i,j},~i\in\mathcal{N},\,j\in\mathcal{N}_i,
\end{equation}
where $w_{i,j}$ are the zero-mean, independent, and identically-distributed Gaussian measurement noise terms.

Using these notations, the maximum likelihood estimator, as our baseline, can be obtained through solving the following nonconvex constrained optimization problem (following \cite{erseghe2015distributed})
\begin{subequations}\label{1}
	\begin{align}
	\mathop{\arg\min}\limits_{\vp\in\mathbb{R}^{nN}}&~\sum_{i\in\mathcal{N}}\sum_{j\in\mathcal{N}_i}\frac{1}{2}\left(\left\|\vp_i-\vp_j\right\|-d_{i,j}\right)^{2}\\
	\text{subject to}&~\vp_k=\va_k,~\forall k\in\mathcal{A}.
	\end{align}
\end{subequations}

\subsection{Problem Reformulation}

\begin{figure}[!t]
	\centering
	\includegraphics[scale=0.9]{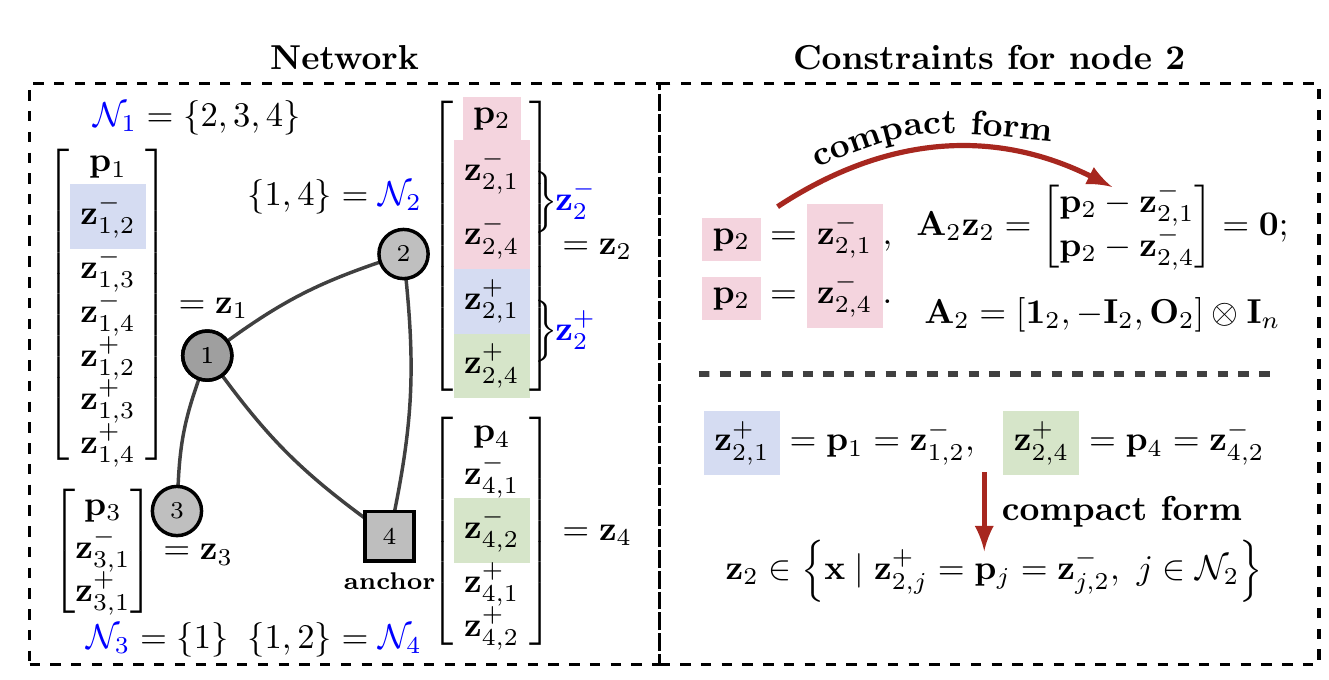}	
	\caption{An illustrating example of a connected wireless sensor network in 2-D space ($N=4,m=1$).}
	\label{gra}
\end{figure}

\noindent
Our first step is to derive an equivalent smooth and constrained reformulation of problem \eqref{1}, which provides the key insight toward algorithm design and its convergence of this paper. Notice that the objective function in problem \eqref{1} can be written explicitly as
\begin{equation}\label{obj1}
\sum_{i\in\mathcal{N}}\sum_{j\in\mathcal{N}_i}\Big[\frac{1}{2}\|\vp_i-\vp_j\|^2-d_{i,j}\underbrace{\|\vp_i-\vp_j\|}_{\text{nonsmooth}}+\frac{1}{2}d_{i,j}^2\Big].
\end{equation}
Obviously, there is a nonsmooth term in the objective function. Motivated by the recent works \cite{luke2017simple,gur2020alternating}, we apply the Cauchy-Schwartz inequality to obtain
\begin{equation}\label{pp}
\|\vp_i-\vp_j\|=\max_{\vu_{i,j}\in\mathcal{B}}~\vu_{i,j}^T(\vp_i-\vp_j),
\end{equation}
where $\mathcal{B}:=\{\vx\in\mathbb{R}^n\mid\|\vx\|\leq1\}$ is a unit ball in $\mathbb{R}^{n}$ with the center at the origin and $\vu_{i,j}\in\mathbb{R}^{n}$ is an auxiliary variable, then we have
\begin{equation}\label{pp1}
-\|\vp_i-\vp_j\|=\min_{\vu_{i,j}\in\mathcal{B}}~-\vu_{i,j}^T(\vp_i-\vp_j).
\end{equation}
Substituting \eqref{pp1} into \eqref{obj1}, problem \eqref{1} is rewritten as a minimization problem of a smooth function over a ball constraint set as follows
\begin{subequations}\label{11}
	\begin{align}
	\label{11_obj}&\mathop{\arg\min}\limits_{\vp,\vu}~\sum_{i\in\mathcal{N}}\sum_{j\in\mathcal{N}_i}\Big[\frac{1}{2}\|\vp_i-\vp_j\|^2-d_{i,j}\vu_{i,j}^T(\vp_i-\vp_j)\Big]\\
	&\text{subject to}~\vu_i\in\mathcal{B}^{N_i},~\forall i\in\mathcal{N},\\\label{anchor}
	&\qquad\qquad\vp_k=\va_k,~\forall k\in\mathcal{A}.
	\end{align}
\end{subequations}
where $\vu:=\text{vec}\left(\vu_i,i\in\mathcal{N}\right)$,  $\vu_i:=\text{vec}\left(\vu_{i,j},j\in\mathcal{N}_i\right)$, and $\mathcal{B}^{N_i}:=\{(\vx_1,\vx_2,\ldots,\vx_{N_i})\mid\vx_j\in\mathcal{B},j=1,\ldots,N_i\}$ is the Cartesian product of $N_i$ balls $\mathcal{B}$. To decouple the loss function of problem \eqref{11} over nodes, we introduce auxiliary variables that duplicate the positions of neighboring nodes. Specifically, we define $\vz_{i,j}^+\in\mathbb{R}^{n}$ as a copy of the position $\vp_j$ for each neighbor $j\in\mathcal{N}_i$, i.e.,
\begin{align}\label{p_j}
\vz_{i,j}^+:=\vp_j,~j\in\mathcal{N}_i.
\end{align}
Similarly, we define $\vz_{i,j}^{-}\in\mathbb{R}^{n}$ as a copy of the position $\vp_i$ assigned to neighbor $j\in\mathcal{N}_i$, i.e.,
\begin{align}\label{p_i}
\vz_{i,j}^{-}:=\vp_i,~j\in\mathcal{N}_i.
\end{align}
By collecting all variables associated with node $i$ together, we define a new variable $\vz_i$ as 
\begin{equation}\label{z}
\vz_i:=
\begin{bmatrix}\vp_{i}\\
\vz_i^-\\ \vz_i^+\end{bmatrix}\in\mathbb{R}^{(2N_i+1)n},~i\in\mathcal{N}.
\end{equation}
Here, $\vz_i^-:=\text{vec}\left(\vz_{i,j}^-,j\in\mathcal{N}_i\right)\in\mathbb{R}^{N_in}$ is the collection of $N_i$ copies of the position $\vp_i$, and $\vz_i^+:=\text{vec}\left(\vz_{i,j}^{+},~j\in\mathcal{N}_i\right)\in\mathbb{R}^{N_in}$  represents the collection of all replicas of $\vp_j$ from neighboring node $j\in\mathcal{N}_i$. Using these notations, the objective function
\eqref{11} can be written in a separable form as follows:
\begin{align}\label{obj_Qm}
\sum_{i\in\mathcal{N}}\Big[\underbrace{\frac{1}{2}\left\|\mQ_i\vz_i\right\|^{2}-\vu_i^T\mD_i\mQ_i\vz_i}_{F_i(\vz_i,\vu_i)}+\delta_{\mathcal{B}^{N_i}}\left(\vu_i\right)\Big],
\end{align}
where $\mD_i:=\textbf{D}\text{iag}\left(\text{vec}\left(d_{i,j},j\in\mathcal{N}_i\right)\right)\otimes\mI_{n}\in\mathbb{R}^{N_in\times N_in}$ is the measurement matrix of node $i$ and
\begin{align}\label{mQ}
\mQ_i:=\left[\vone_{N_i},\mO_{N_i},-\mI_{N_i}\right]\otimes\mI_{n}\in\mathbb{R}^{N_in\times(1+2N_i)n}.
\end{align}
By using \eqref{z}, we can rewrite \eqref{p_i} as a compact form
\begin{align}\label{mA}
\mA_i\vz_i=\vzero,
\end{align}
where $\mA_i:=[\vone_{N_i},-\mI_{N_i},\mO_{N_i}]\otimes\mI_{n}\in\mathbb{R}^{N_in\times(1+2N_i)n}$. Moreover, for a pair of connected sensors $i$ and $j$, we can deduce from equation \eqref{p_i} that $\vz_{j,i}^- = \vp_j$. Combining this equation with \eqref{p_j}, we obtain an additional constraint that applies to all connected sensors:
\begin{align}\label{mz}
\vz:=\text{vec}(\vz_i,i\in\mathcal{N})\in\mathcal{Z},\quad\mathcal{Z}:=\{\vz|\vz_{i,j}^{+}=\vp_j=\vz_{j,i}^{-},~\forall~i\in\mathcal{N},~j\in\mathcal{N}_i\}.
\end{align}
Using \eqref{z}, the set of anchors in \eqref{anchor} can be expressed as
\begin{align}\label{anchor_set}
\vz\in\mathcal{X}:=\{\vz\mid\mE_i\vz_i=\va_i,\forall~i\in\mathcal{A}\},\quad\mE_i:=\left[1,\vzero^T_{N_i},\vzero^T_{N_i}\right]\otimes\mI_{n}\in\mathbb{R}^{n\times(1+2N_i)n}.
\end{align}
Finally, using \eqref{obj_Qm}, \eqref{mA}, and \eqref{mz}, the nonconvex optimization problem \eqref{11} can be equivalently reformulated into the compact form
\begin{subequations}\label{f}
	\begin{align}\label{f_con0}
	\mathop{\arg\min}\limits_{\vz,\vu}&\sum_{i\in\mathcal{N}}\Big[ \underbrace{\frac{1}{2}\left\|\mQ_i\vz_i\right\|^{2}-\vu_i^T\mD_i\mQ_i\vz_i}_{F_i(\vz_i,\vu_i)}+\delta_{\mathcal{B}^{N_i}}\left(\vu_i\right)\Big]\\
	\text{subject to}&\,\,\label{f_con}\mA_i\vz_i=\vzero,i\in\mathcal{N},
	\\&\label{z_set_con}\,\,\vz\in\mathcal{X},\,\,\vz\in\mathcal{Z}.
	\end{align}
\end{subequations}
To show these notations clearly, we give an example of a connected wireless sensor network in 2-D space in Fig. \ref{gra}.
\begin{Exa}\label{example}
	Fig.\ref{gra} presents a wireless sensor network with 3 (inexact) positions of node $i=1,2,3$ and one anchor of node 4. If we only focus on node 2, on the one hand, since $\mathcal{N}_2=\{1,4\}$ and $\vz_{2,1}^-,\vz_{2,4}^-$ are the copies of $\vp_2$, then we have
	\begin{align*}
	\vp_2=\vz_{2,1}^-,\quad\vp_2=\vz_{2,4}^-.
	\end{align*}
	On the other hand, Fig. \ref{gra} presents that there exit the same partial elements in $\vz_2,\,\vz_1$, and $\vz_4$. Specifically, $\vp_1$ is the element of both $\vz_1$ and $\vz_2$; $\vz_2$ and $\vz_4$ also contain $\vp_{4}$. To enforce this trivial observation, we obtain the following constraint
	\begin{align*}
	\vz_{2,1}^{+}=\vp_1=\vz_{1,2}^-,\quad\vz_{2,4}^{+}=\vp_4=\vz_{4,2}^-,
	\end{align*}
	which is part of the constraint set $\mathcal{Z}$ deﬁned in \eqref{mz}.
\end{Exa}

\begin{Rmk}
	The objective function $F_i(\vz_i,\vu_i)$ in problem \eqref{f} is convex with respect to $\vz_i$ when $\vu_i$ is fixed, however, it is NOT a jointly convex function for $(\vz_i,\vu_i)$. In addition, one can observe that the objective function is separable but the linear constraint $\mathcal{Z}$ is NOT separable in $\vz$. In Section \ref{sec: alg develop}, we introduce a scaled proximal term to deal with this problem.
\end{Rmk}

\section{Proposed SP-ADMM algorithms}\label{sec: alg develop}

\noindent
In this section, we present the proposed SP-ADMM algorithms for solving the nonconvex and nonsmooth optimization problem derived in \eqref{f}. The proposed SP-ADMM algorithms are built upon the ADMM method \cite{luke2017simple} developed for large-scale nonconvex problems.

First, we introduce the augmented Lagrangian (AL) function of problem \eqref{f} as follows:
\begin{equation*}
\mathcal{L}(\vz,\vu,\boldsymbol{\lambda}):=\sum_{i\in\mathcal{N}}\mathcal{L}_i\left(\vz_i,\vu_i,\vlambda_i\right),
\end{equation*}
where
\begin{align}\label{L}
\mathcal{L}_i\left(\vz_i,\vu_i,\vlambda_i\right):=F_i\left(\vz_i,\vu_i\right)+\delta_{\mathcal{B}^{N_i}}\left(\vu_i\right)+\langle\vlambda_i,\mA_i\vz_i\rangle+\frac{c}{2}\|\mA_i\vz_i\|^{2},~i\in\mathcal{N},
\end{align}
and $\vlambda_i:=\text{vec}\left(\vlambda_{i,j},j\in\mathcal{N}_i\right)\in\mathbb{R}^{N_in}$ correspond to the Lagrangian multipliers, $c>0$ is a penalty coefficient. We apply the following scaled proximal ADMM updates \cite{boyd2011distributed} to obtain
\begin{align}\label{z_up}
\vz^{t+1}&=\mathop{\arg\min}\limits_{\substack{\vz\in\mathcal{Z}\\\vz\in\mathcal{X}}}\sum_{i\in\mathcal{N}}\mathcal{L}_i\left(\vz_i,\vu^t_i,\vlambda^t_i\right)\hspace{-0.03in}+\frac{c}{2}\|\vz_i-\vz_i^t\|^{2}_{\mB_i^T\mB_i},\\\label{u_up}
\vu^{t+1}&=\mathop{\arg\min}\limits_{\vu}\sum_{i\in\mathcal{N}}\mathcal{L}_i\left(\vz_i^{t+1},\vu_i,\vlambda^t_i\right)+\frac{\rho}{2}\|\vu_i-\vu_i^t\|^{2},\\\label{lambda_up}
\vlambda_i^{t+1}&=\vlambda_i^t+c\mA_i\vz_i^{t+1},\,\,\forall i\in\mathcal{N},
\end{align}
where $\rho>0$ is a penalty coefficient and $t$ denotes the iteration step.

The AL function $\mathcal{L}$ is separable in each of the variables $\left(\vz_i,\vu_i,\vlambda_i\right),\,i\in\mathcal{N}$ and $\mathcal{L}_i,\,i\in\mathcal{N}$ is convex as a function separately of $\left(\vz_i,\vu_i\right)$ when the others are fixed. These features are key to distributed implementation. In addition, we remark that the (scaled) proximal term is critical in both the implementation efficiency and the theoretical analysis. It is used to ensure the following attractive properties:
\begin{itemize}
	\item $\vz_i\rightarrow\mathcal{L}_i\left(\vz_i,\vu_i^t,\vlambda_i^t\right)+\frac{c}{2}\|\vz_i-\vz_i^t\|^{2}_{\mB_i^T\mB_i}$ is strongly convex;
	\item $\vu_i\rightarrow\mathcal{L}_i\left(\vz_i^{t+1},\vu_i,\vlambda^t_i\right)+\frac{\rho}{2}\|\vu_i-\vu_i^t\|^{2}$ is strongly convex;
	\item The sequence $\{\left(\vz_i^t,\vu_i^t\right)\}$ has lower computational complexity.
\end{itemize}

It should be noted that applying distinct penalty coefficients for $\vz_i$ and $\vu_i$ is necessary due to they control diﬀerent aspects of the optimization problem. The proximal term controlled by $\rho$ only aﬀects the update for $\vu$, whereas the AL term controlled by $c$ aﬀects the update for $\vz$ and requires the use of a dual multiplier $\vlambda$ to enforce the constraint. The choice of $c$ can have a crucial impact on the convergence behavior of the algorithm. Moreover, using diﬀerent penalty coefficients enables users to ﬁne-tune the algorithm’s parameters for optimal performance in practice presented in Section \ref{sec: simulation}.

To see the first point, we reformulate the objective function of problem \eqref{z_up}. Combining \eqref{f_con0} and \eqref{L}, and rearranging the above quadratic term, we have
\begin{align}\nonumber
\mathcal{L}_i\left(\vz_i,\vu_i^t,\vlambda_i^t\right)+\frac{c}{2}\|\vz_i-\vz_i^t\|^{2}_{\mB_i^T\mB_i}=&\frac{1}{2}\left\|\vz_i\right\|^{2}_{\mW_i}-\left\langle\mQ_i^T\mD_i\vu_i^t-\mA_i^T\vlambda_i^t+c\mB_i^T\mB_i\vz_i^t,\vz_i\right\rangle\\\label{L110}
&+\delta_{\mathcal{B}^{N_i}}\left(\vu_i^t\right)+\frac{c}{2}\|\vz_i^t\|^{2}_{\mB_i^T\mB_i},
\end{align}
where
\begin{equation}
\label{new2}\mW_i=\mQ_i^T\mQ_i+c\mA_i^T\mA_i+c\mB_i^T\mB_i.
\end{equation}
Suppose $\mB_i^T\mB_i$ is chosen such that $\mW_i$ is a positive definite matrix, then the objective function of subproblem \eqref{z_up} is strongly convex, which ensures that the subproblem \eqref{z_up} has a unique solution. By completing the square of \eqref{L110}, we get
\begin{align}\label{L111}
\mathcal{L}_i\left(\vz_i,\vu_i^t,\vlambda_i^t\right)+\frac{c}{2}\|\vz_i-\vz_i^t\|^{2}_{\mB_i^T\mB_i}=\frac{1}{2}\|\vz_i-\tilde{\vz}_i^{t+1}\|_{\mW_i}^2-\frac{1}{2}\|\tilde{\vz}_i^{t+1}\|^2_{\mW_i}+\delta_{\mathcal{B}^{N_i}}\left(\vu_i^t\right)+\frac{c}{2}\|\vz_i^t\|_{\mB_i^T\mB_i}^{2},
\end{align}
where
\begin{equation}\label{new1}
\tilde{\vz}_i^{t+1}=\mW_i^{-1}\left(\mQ_i^T\mD_i\vu_i^t-\mA_i^T\vlambda_i^t+c\mB_i^T\mB_i\vz_i^t\right).
\end{equation}
Substituting \eqref{L111} into \eqref{z_up} and omitting the constant terms, the optimization problem \eqref{z_up} is equivalent to
\begin{equation}\label{lem_zt1}
\vz^{t+1}=\mathop{\arg\min}\limits_{\vz\in\mathcal{X},\vz\in\mathcal{Z}}~\sum_{i\in\mathcal{N}}\frac{1}{2}\left\|\vz_{i}-\tilde{\vz}_{i}^{t+1}\right\|^{2}_{\mW_{i}}.
\end{equation}
Optimization problem \eqref{lem_zt1} presents two challenges. Firstly, obtaining an analytic solution is difficult due to the linear constraints of both $\mathcal{Z}$ and $\mathcal{X}$. Secondly, calculating the inverse matrix $\mW_i$ in \eqref{new1} to obtain $\tilde{\vz}_i^{t+1}$ can be computationally expensive. In the following, we demonstrate that an appropriate choice of $\mB_i^T\mB_i$ can overcome these challenges. To facilitate the calculation of the inverse of $\mW_i$, we set
\begin{equation}\label{cbb}
c\mB_i^T\mB_i=c|\mA_i^T\mA_i|+|\mQ_i^T\mQ_i|,
\end{equation}
where $|\cdot|$ takes element-wise absolute value of a matrix. From the definition of $\mA_i$ and $\mQ_i$ in \eqref{mA} and \eqref{mQ} respectively, we get
\begin{align}\label{AtA}
\mA_i^T\mA_i=\begin{bmatrix}
N_i&-\boldsymbol{1}_{N_i}^T&\boldsymbol{0}_{N_i}^T\\
-\boldsymbol{1}_{N_i}&\mI_{N_i}&\mO_{N_i}\\
\boldsymbol{0}_{N_i}&\mO_{N_i}&\mO_{N_i}
\end{bmatrix}\otimes\mI_{n},\\
\label{QtQ}
\mQ_i^T\mQ_i=\begin{bmatrix}
N_i&\vzero_{N_i}^T&-\boldsymbol{1}_{N_i}^T\\
\vzero_{N_i}&\mO_{N_i}&\mO_{N_i}\\
-\vone_{N_i}&\mO_{N_i}&\mI_{N_i}
\end{bmatrix}\otimes\mI_{n}.
\end{align}
Substituting \eqref{AtA}-\eqref{QtQ} into \eqref{cbb}, we obtain
\begin{equation}\label{b2}
c\mB_i^T\mB_i=
\begin{bmatrix}
\left(c+1\right)N_i&c\cdot\vone_{N_i}^T&\vone_{N_i}^T\\
c\cdot\vone_{N_i}&c\cdot\mI_{N_i}&\mO_{N_i}\\
\vone_{N_i}&\mO_{N_i}&\mI_{N_i}
\end{bmatrix}\otimes\mI_{n}.
\end{equation}
With the aid of \eqref{AtA}-\eqref{b2}, we can guarantee that $\mW_i$ defined in \eqref{new2} is a  positive definite diagonal matrix of the form
\begin{equation}\label{w}
\mW_i=2\cdot\textbf{D}\text{iag}\left(\left[\left(c+1\right)N_i,c\cdot\boldsymbol{1}_{N_i}^T,\boldsymbol{1}_{N_i}^T\right]\right)\otimes\mI_{n}.
\end{equation}
It makes $\tilde{\vz}_i^{t+1}$ in \eqref{new1} easy to compute. Moreover, by denoting $\tilde{\vz}_i:=[\left(\tilde{\vp}_{i}\right)^T,\left(\tilde{\vz}_i^-\right)^T,\left(\tilde{\vz}_i^+\right)^T]^T$, where its last two parts are
\begin{equation*}
\tilde{\vz}_i^-:=\text{vec}\left(\tilde{\vz}_{i,j}^{-},j\in\mathcal{N}_i\right),~\tilde{\vz}_i^+:=\text{vec}\left(\tilde{\vz}_{i,j}^{+},j\in\mathcal{N}_i\right),
\end{equation*}
then $\vz_i^{t+1}$ can be deduced using the following remark.
\begin{Rmk}\label{r_z_ztilde}
	Given $\mB_i$, $i\in\mathcal{N}$ in \eqref{cbb}, the analytic formula of the optimal solution $\vz_{i}^{t+1}$ in \eqref{lem_zt1} is derived as follows
	\begin{align}\label{ztt1}	\vp_{i}^{t+1}&=\tilde{\vp}_{i}^{t+1},\\\label{zt}
	\left(\vz_{i}^-\right)^{t+1}&=\text{\rm vec}\left(\frac{c}{c+1} (\tilde{\vz}_{i,j}^{-})^{t+1}+\frac{1}{c+1}(\tilde{\vz}_{j,i}^{+})^{t+1},~j\in\mathcal{N}_i\right),\\\label{ztt2}	\left(\vz_{i}^+\right)^{t+1}&=\text{\rm vec}\left(\frac{1}{c+1}(\tilde{\vz}_{i,j}^{+})^{t+1}+\frac{c}{c+1}(\tilde{\vz}_{j,i}^{-})^{t+1},~j\in\mathcal{N}_i\right).
	\end{align}
	The proof is shown in the Section \ref{proof_r_z_ztilde} in the supplement material. From \eqref{ztt1}-\eqref{ztt2}, one can see that $\vz_i^{t+1}$ can be obtained in a distributed fashion using only $\tilde{\vz}_j^{t+1}$, $j\in\mathcal{N}_i$ from its neighbors.
\end{Rmk}

To see the last point for $\vu_i$, by omitting the terms that are irrelevant with $\vu$, the optimization problem in \eqref{u_up} reduces to
\begin{equation*}
\mathop{\arg\min}\limits_{\vu}\sum_{i\in\mathcal{N}}\left(-\langle\mQ_i^T\mD_i\vu_i,\vz_i^{t+1}\rangle+\delta_{\mathcal{B}^{N_i}}\left(\vu_i\right)+\frac{\rho}{2}\|\vu_i-\vu_i^t\|^{2}\right).
\end{equation*}
The above objective function is separable in $\vu_i$ and further performing projection yields
\begin{equation}\label{u_proj}
\vu_i^{t+1}=\text{proj}_{\mathcal{B}^{N_i}}\left(\vu_i^t+\frac{1}{\rho}\mD_i\mQ_i\vz_i^{t+1}\right),
\end{equation}
which represents a projection onto the unit ball that can be reduced to
\begin{equation}
\label{u}
\vu_{i,j}^{t+1}=\frac{\tilde{\vu}_{i,j}^{t+1}}{\max\{1,\|\tilde{\vu}_{i,j}^{t+1}\|\}},\,j\in\mathcal{N}_i,
\end{equation}
where
\begin{equation*}	\text{vec}\left(\tilde{\vu}^{t+1}_{i,j},j\in\mathcal{N}_i\right)=\vu_i^t+\frac{1}{\rho}\mD_i\mQ_i\vz_i^{t+1}.
\end{equation*}

In summary, problem \eqref{z_up} and problem \eqref{u_up} can be equivalently rewritten as problem \eqref{lem_zt1} and problem \eqref{u_proj}, respectively, and the corresponding distributed solutions are given in \eqref{ztt1}-\eqref{ztt2} and \eqref{u}. An outline of the distributed procedure is described in Algorithm \ref{al2}.

\begin{algorithm}
	\caption{Distributed SP-ADMM algorithm}
	\KwIn{parameter $c,\rho,$ and initial values of $\vz^{0}_i, \vu^{0}_i,\vlambda^{0}_i=\vzero$, $i=1,\ldots,N$.}Let
	\vspace{-0.25in}
	\begin{equation*}\hspace{0.1in}
	\mW_i=2\cdot\textbf{D}\text{iag}\left(\left[\left(c+1\right)N_i,c\cdot\boldsymbol{1}_{N_i}^T,\boldsymbol{1}_{N_i}^T\right]\right)\otimes\mI_{n}.
	\end{equation*}
	\For{$t\in\{0,\ldots,T\}$}{
		\For{client $i=1,2,\ldots,N$ in parallel}{
			\begin{align}\label{ztilde}
			\tilde{\vz}_i^{t+1}&=\mW_i^{-1}\Big(\mQ_i^T\mD_i\vu_i^t-\mA_i^T\vlambda_i^t+c\mB_i^T\mB_i\vz_i^t\Big)\\\label{xtilde_a}
			\tilde{\vp}_{i}^{t+1}&=\va_i,\,\text{if}~i\in\mathcal{A}
			\end{align}}
		{\bf communication:}\\
		$\Rightarrow$ Broadcast $(\tilde{\vz}_{i,j}^{-})^{t+1},\,\,(\tilde{\vz}_{i,j}^{+})^{t+1}$ to region $j\in\mathcal{N}_i$,\\
		$\Leftarrow$ Receive $(\tilde{\vz}_{j,i}^{-})^{t+1},\,(\tilde{\vz}_{j,i}^{+})^{t+1}$ from region $j\in\mathcal{N}_i$.\\
		\For{client $i=1,2,\ldots,N$ in parallel}{
			\begin{align}
			&\text{Update}~\vz_i^{t+1}~\text{via}~\eqref{ztt1}-\eqref{ztt2}\\\nonumber
			&\text{Update}~\vu_i^{t+1}~\text{via}~\eqref{u}\\\label{du2}
			&\vlambda_i^{t+1}=\vlambda_i^t+c\,\mA_i\vz^{t+1}_i
			\end{align}	
	}}
	\label{al2}
\end{algorithm}

\begin{Rmk}
	We show that the SP-ADMM algorithm is closely related to the ADMM algorithm derived in \cite{luke2017simple} for single source localization. Specifically, let us set $\vp^t=\vp_{i}^t$, $\va_j=\vp_{j}^t$ for all $t\geq0,\,j\in\mathcal{N}_i=\mathcal{A}$ (single source localization). Then, by \eqref{ztt1} and the definition of $\mE_i$ in \eqref{anchor_set}, we obtain
	\begin{equation}\label{pz}
	\vp^{t+1}=\tilde{\vp}_{i}^{t+1}=\mE_i\tilde{\vz}_i^{t+1},~i\notin\mathcal{A}.
	\end{equation}
	Multiplying $\mE_i$ on both sides of \eqref{new1} and rearranging terms yields
	\begin{align}\label{pi0}
	\tilde{\vp}_i^{t+1}=\sum_{j\in\mathcal{N}_i}\frac{1}{2\left(c+1\right)N_i}\Big[d_{ij}\vu_{i,j}^t+\left(2c+1\right)\left(\vp_{i}^t-\vp_{j}^t\right)-\vlambda_{i,j}^t+\left(2c+1\right)\vp_{j}^t+(\vz_{i,j}^+)^t-c\left(\vp_{i}^t-(\vz_{i,j}^-)^t\right)\Big].
	\end{align}
	If $\vp_{i}^t=(\vz_{i,j}^-)^t$, which is $\mA_i\vz_i^t=\vzero,\,i\in\mathcal{N}$, then combining \eqref{pz} with \eqref{pi0}, replacing $\vp_{j},\,\vz_{i,j}^+$ with $\va_j$, and rearranging terms again we have
	\begin{align}\label{pi}
	\vp^{t+1}=\frac{1}{m}\sum_{j=1}^{m}\Big[\va_j+\frac{1}{2\left(c+1\right)m}\left(d_{j}\vu_{j}^t+\left(2c+1\right)\left(\vp^t-\va_j\right)-\vlambda_{j}^t\right)\Big].
	\end{align}
	Here we have replaced $d_{i,j},\vu_{i,j}^t,\vlambda_{i,j}^t$ by $d_{j},\vu_{j}^t,\vlambda_{j}^t$ since $i$ is a fixed value in the single source localization problem. Hence, if $\vz_i^t$ satisfies the linear equality constraint \eqref{f_con}, then equation \eqref{pi} is identical to the ADMM algorithm given in [\cite{luke2017simple}, Equation (3.8)] except for the constant coefficients of the variables. Therefore, the proposed SP-ADMM algorithm can be regarded as an extended version of the  ADMM algorithm in [\cite{luke2017simple} with extra capability to handle the multi-source localization problems.
	
	However, in the case of multi-source localization, the ADMM algorithm can be significantly more complex than in single-source localization due to coupled variables and increased computational complexity. Coupled variables make subproblems in ADMM entangled and thus the convergence of the algorithm is hard to prove. Moreover, the increased number of variables in multi-source localization increases computational complexity, leading to slower convergence rates and reduced localization accuracy. 
\end{Rmk}

In Algorithm \ref{al2}, we observed that in order to carry out the new round of iterations, node $i$ has to store the following values: $\vz_i^t,\vu_i^t,\vlambda_i^t,\mW_i,\mB_i^T\mB_i,\mD_i,\mQ_i,\mA_i$, $c,\rho,N_i$. If we assume that the storage unit occupied by any real number is one, then node $i$ requires $13n^2N_i^2+10n^{2}N_i+2n^2+4nN_i+n+3$ storage units in total. Note that although the matrices and parameters involved are fixed, the rest of the vectors are updated with each iteration. To reduce the storage space required by each sensor to run Algorithm \ref{al2}, we transform the update step as follows.

Let us start by giving an explicit formula for the parts of the update variables in Algorithm \ref{al2}. Using the form of $\mA_i,\mQ_i$, and $\mD_i$, defined in \eqref{mA}, \eqref{mQ}, and \eqref{obj_Qm}, respectively, we have
\begin{align}\nonumber
\mA_i^T\vlambda_i^t=&\Big[\sum_{j\in\mathcal{N}_{i}}\left(\vlambda_{i,j}^t\right)^{T},\text{vec}
\left(-\vlambda_{i,j}^t,j\in\mathcal{N}_i\right)^T,\vzero_{N_{i}}^{T}\Big]^{T},\\\label{QD_Alam}
\mQ_i^T\mD_{i}\vu_i^t=&\Big[\sum_{j\in\mathcal{N}_{i}}\left(d_{i,j}\vu_{i,j}^t\right)^{T},\vzero_{N_{i}}^{T},\text{vec}\left(-d_{i,j}\vu_{i,j}^t,j\in\mathcal{N}_i\right)^T\Big]^T.	
\end{align}
Also, it follows from the definition \eqref{b2} that
\begin{align}\nonumber
c\mB_i^T\mB_i\vz_i^t&=\Big[\Big((c+1)N_i\vp_i^t+c\sum_{j\in\mathcal{N}_{i}}(\vz_{i,j}^-)^{t}+\sum_{j\in\mathcal{N}_{i}}(\vz_{i,j}^+)^{t}\Big)^{T},\text{vec}\left(c\left(\vp_i^t+(\vz_{i,j}^-)^{t}\right),j\in\mathcal{N}_i\right)^T,\\\label{cBBzt}&\qquad\text{vec}\left(\vp_i^t+(\vz_{i,j}^+)^{t},j\in\mathcal{N}_i\right)^T\Big]^{T}.
\end{align}
By substituting \eqref{QD_Alam}-\eqref{cBBzt} into \eqref{ztilde} yields the following update
\begin{align}\label{xtildeii}
\tilde{\vp}_{i}^{t+1}=&\frac{1}{2(c+1)N_i}\sum_{j\in\mathcal{N}_i}\Big[d_{i,j}\vu_{i,j}^t-\vlambda_{i,j}^t+c\left(\vp_i^t+(\vz_{i,j}^-)^{t}\right)+\vp_i^t+(\vz_{i,j}^+)^{t}\Big],\\
\label{zt11}(\tilde{\vz}_{i,j}^-)^{t+1}=&\frac{1}{2c}\vlambda_{i,j}^t+\frac{1}{2}\left(\vp_i^t+(\vz_{i,j}^-)^{t}\right),\\ \label{zt12}(\tilde{\vz}_{i,j}^+)^{t+1}=&-\frac{1}{2}d_{i,j}\vu_{i,j}^t+\frac{1}{2}\left(\vp_i^t+(\vz_{i,j}^+)^{t}\right).
\end{align}
Using once again the form of $\mQ_i,\mD_i$ and $\mA_i$, it follows from the \eqref{u} and \eqref{du2} that
\begin{align}\nonumber
\tilde{\vu}_{i,j}^{t+1}=&\vu_{i,j}^t+\frac{d_{i,j}}{\rho}\left(\vp_i^{t+1}-(\vz_{i,j}^+)^{t+1}\right),\\
\label{u1}\vu_{i,j}^{t+1}=&\frac{1}{\max\{1,\|\tilde{\vu}_{i,j}^{t+1}\|\}}\tilde{\vu}_{i,j}^{t+1},\\
\label{lambda1}\vlambda_{i,j}^{t+1}=&\vlambda_{i,j}^t+c\left(\vp_i^{t+1}-(\vz_{i,j}^-)^{t+1}\right).
\end{align}

Next, we intend to remove $\tilde{\vz}_i^{t+1}$ to reduce the storage space and computation required by Algorithm \ref{al2}. Combining \eqref{ztt1} with \eqref{xtilde_a} and \eqref{xtildeii}, we have
\begin{equation}\label{z3}	
{\vp_{i}^{t+1}}=
\begin{cases}
\frac{\sum_{j\in\mathcal{N}_i}\left(2d_{i,j}\vu_{i,j}^t-2\vlambda_{i,j}^t+\valpha_{i,j}^t+\vbeta_{i,j}^t\right)}{2(c+1)N_i},&{\text{if}}\ i\notin\mathcal{A},\\
{\va_i,}&{\text{if}}\ {i\in \mathcal{A}.}
\end{cases}
\end{equation}
where
\begin{align}\label{ab0}
\valpha_{i,j}^t:&=\vlambda_{i,j}^t+c\left(\vp_i^t+(\vz_{i,j}^-)^{t}\right),\\\label{ab}
\vbeta_{i,j}^t:&=-d_{i,j}\vu_{i,j}^t+\vp_i^t+(\vz_{i,j}^+)^{t}.
\end{align}
By adding $\frac{1}{c}(\tilde{\vz}_{j,i}^+)^{t+1}$ to the both side of \eqref{zt11} and using \eqref{zt}, we get
\begin{equation}\label{zt1}
\frac{c+1}{c}(\vz_{i,j}^-)^{t+1}=\frac{1}{2c}\vlambda_{i,j}^t+\frac{1}{2}\left(\vp_i^t+(\vz_{i,j}^-)^{t}\right)+\frac{1}{c}(\tilde{\vz}_{j,i}^+)^{t+1}.
\end{equation}
Then, substituting \eqref{zt12} into \eqref{zt1} and using \eqref{ab0}-\eqref{ab}, it yields
\begin{equation}\label{zt2}
(\vz_{i,j}^-)^{t+1}=\frac{1}{2(c+1)}\left(\valpha_{i,j}^t+\vbeta_{j,i}^t\right).
\end{equation}
Similarly, we add $c(\tilde{\vz}_{j,i}^-)^{t+1}$ to \eqref{zt12}, then using \eqref{ztt2} and \eqref{ab0}-\eqref{ab} to obtain
\begin{equation}\label{zt3}
(\vz_{i,j}^+)^{t+1}=\frac{1}{2\left(c+1\right)}\left(\vbeta_{i,j}^t+\valpha_{j,i}^t\right).
\end{equation}
Applying \eqref{zt3} to \eqref{u1} and \eqref{ab}, we get
\begin{align}
\nonumber\tilde{\vu}_{i,j}^{t+1}=&\vu_{i,j}^t+\frac{d_{i,j}}{\rho}\vp_i^{t+1}-\frac{d_{i,j}}{2\rho\left(c+1\right)}\left(\vbeta_{i,j}^t+\valpha_{j,i}^t\right),\\
\label{u2}\vu_{i,j}^{t+1}=&\frac{1}{\max\{1,\|\tilde{\vu}_{i,j}^{t+1}\|\}}\tilde{\vu}_{i,j}^{t+1},\\\label{ab2}
\vbeta_{i,j}^{t+1}=&-d_{i,j}\vu_{i,j}^{t+1}+\vp_{i}^{t+1}+\frac{1}{2\left(c+1\right)}\left(\vbeta_{i,j}^t+\valpha_{j,i}^t\right).
\end{align}
Substituting \eqref{zt2} into \eqref{lambda1} and \eqref{ab0}, we have
\begin{align}\label{lambda2}
\vlambda_{i,j}^{t+1}&=\vlambda_{i,j}^t+c\vp_{i}^{t+1}-\frac{c}{2(c+1)}\left(\valpha_{i,j}^t+\vbeta_{j,i}^t\right),\\\label{ab1}
\valpha_{i,j}^{t+1}&=\vlambda_{i,j}^t+2c\vp_{i}^{t+1}.
\end{align}

\begin{algorithm}[t]
	\caption{Simplified SP-ADMM algorithm}
	\KwIn{parameter $c,\rho,$ and initial values of $\vz^{0}\in\mathcal{Z},\vu^{0}_{i,j},\vlambda^{0}_{i,j}=\vzero$, $i\in N,\,j\in\mathcal{N}_i$.}Compute $\{\left(\valpha_{i,j}^{0},\vbeta_{i,j}^{0}\right)\}_{i\in\mathcal{N},j\in\mathcal{N}_i}$ defined in \eqref{ab0}-\eqref{ab}\\
	\For{$t\in\{0,\ldots,T\}$}{
		{\bf communications: }\\
		$\Rightarrow$ Broadcast $\valpha_{i,j}^t,\vbeta_{i,j}^t$ to region $j\in\mathcal{N}_i$,\\
		$\Leftarrow$ Receive $\valpha_{j,i}^t,\vbeta_{j,i}^t$ from region $j\in\mathcal{N}_i$.\\
		\For{client $i=1,2,\ldots,N$ in parallel}{
			\begin{align*}
			\text{Update}~\vp_{i}^{t+1}~\text{via}~\eqref{z3}\\\text{Update}~\vu_{i,j}^{t+1}~\text{via}~\eqref{u2}\\\text{Update}~\vbeta_{i,j}^{t+1}~\text{via}~\eqref{ab2}\\\text{Update}~\valpha_{i,j}^{t+1}~\text{via}~\eqref{ab1}\\\text{Update}~\vlambda_{i,j}^{t+1}~\text{via}~\eqref{lambda2}
			\end{align*}
	} }
	\label{al1}
\end{algorithm}

The simplified version of the algorithm is summarized in Algorithm \ref{al1}. From the updating steps of Algorithm \ref{al1}, it can be seen that each node only needs to store the following values: $\valpha_{i,j}^t,\vbeta_{i,j}^t,\vu_{i,j}^t,\vlambda_{i,j}^t,d_{i,j},j\in\mathcal{N}_i$ and $c,\rho,N_i$. Compared with the storage space required by node $i$ of Algorithm \ref{al2}, Algorithm \ref{al1} demands only $4nN_i+N_i+3$ storage units, which is significantly reduced compared to that required by Algorithm \ref{al2}. We remark that the update \eqref{z3} and \eqref{u2}-\eqref{ab1} of Algorithm \ref{al1} should be interpreted as an improved version that optimizes the storage requirement of Algorithm \ref{al2}. Although Algorithm \ref{al1} oﬀers superior storage eﬃciency, the presentation of Algorithm \ref{al2} still has value for several reasons. Firstly, Algorithm \ref{al2} serves as a foundation for understanding the principles and techniques behind Algorithm \ref{al1}. By presenting Algorithm \ref{al2} ﬁrst, we can demonstrate the motivation and rationale for the modiﬁcations made in Algorithm \ref{al1}. Additionally, the convergence analysis for Algorithm \ref{al2} is more direct since Algorithm \ref{al1} streamlines the steps of Algorithm \ref{al2}.

\begin{Rmk}(\textbf{Complexity Comparison}) For each iteration of Algorithm \ref{al2}, the computation cost is derived by \eqref{ztilde} that involves the inversion of diagonal matrix $\mW_i\in\mathbb{R}^{n\left(2N_i+1\right)\times n\left(2N_i+1\right)}$, which is independent of the iteration index $t$, and thus needs to be computed only once. Besides, equation \eqref{u} requires calculating the Euclidean norm of  $\vu_{i}\in\mathbb{R}^{nN_i\times1}$. Hence, the computational complexity of the proposed SP-ADMM algorithm is $\mathcal{O}\left(nN_i\right)$, which is lower than SDP relaxation method $O(n^{3})$ and ADMM-H, but equal to SF and AM-FD. The communication cost per iteration per node is directly proportional to the number of scalar variables that sensors need to transmit to their neighboring nodes \cite{simonetto2014distributed}. Thus, the AM-FD, SDP, and SF algorithms have a communication cost of $nN_i$, while the proposed SP-ADMM requires sending two variables to its neighbor, resulting in a communication cost per iteration of $2nN_i$. It should be noted that the proposed SP-ADMM has lower communication cost than ADMM-H since the latter also needs to transmit a penalty parameter to its neighboring nodes at each iteration. The detailed comparisons between the proposed SP-ADMM algorithm with the other methods have been shown in Table~\ref{comparisons of different alg}.
\end{Rmk}

\section{Convergence Analysis}\label{sec: conv analysis}

\noindent
In this section, we present a proof for the convergence of the proposed Algorithm \ref{al2}, which consists of three main steps. First, we provide an upper bound \cite{fan2021spectrally} for the difference  $\mathcal{L}(\vz^{t+1},\vu^{t+1},\vlambda^{t+1})-\mathcal{L}(\vz^{t},\vu^{t},\vlambda^{t})$. Then, we construct a potential function $\varsigma^t$ based on the AL function $\mathcal{L}(\vz^{t},\vu^{t},\vlambda^{t})$, which monotonically decreases with each iteration and has a lower bound. {\blue Employing the Kurdyka–\L ojasiewicz property \cite{bolte2014proximal,guo2017convergence,bolte2018first,boct2020proximal,yashtini2022convergence}, we establish the global convergence of the sequence generated  by Algorithm \ref{al2}.} Next, we introduce the optimal gap function $\mathcal{F}(\vz,\vu,\vlambda)$ to show that the proposed Algorithm \ref{al2} converges to both a KKT point of problem \eqref{f} and a critical point of the original problem \eqref{1}. Lastly, we define the $\varepsilon$-solution and prove that the algorithm sequence converges with a sublinear rate.

First, we bound the update step of the AL function $\mathcal{L}\left(\vz,\vu,\vlambda\right)$ at each iteration.
\begin{Lemma}\label{LL}
	Suppose $c\mB_i^T\mB_i$ takes the form of \eqref{cbb}, and let $\{\left(\vz_i^t,\vu_i^t,\vlambda_i^t\right)\}$ be a sequence generated by Algorithm  \ref{al2}. Then, for all $t\geq1$, we have
	\begin{align}\nonumber
	&\mathcal{L}\left(\vz^{t+1},\vlambda^{t+1},\vu^{t+1}\right)-\mathcal{L}\left(\vz^t,\vlambda^t,\vu^t\right)\\\nonumber\leq&\sum_{i\in\mathcal{N}}\Big[-\frac{c}{2}\|\vz_i^{t+1}-\vz_i^t\|^{2}_{\mB_i^T\mB_i}-\frac{\rho}{2}\|\vu_i^{t+1}-\vu_i^t\|^{2}\\\nonumber
	&\qquad+\frac{3\left(N_{\max}+1\right)}{c}\|\mQ_i\left(\tilde{\vz}_i^{t+1}-\tilde{\vz}_i^t\right)-\mD_i\left(\vu_i^t-\vu_i^{t-1}\right)\|^{2}\\\nonumber
	&\qquad+3c\left(N_{\max}+1\right)\|\tilde{\vz}_i^{t+1}-\vz_i^{t+1}-\left(\tilde{\vz}_i^t-\vz_i^t\right)\|^{2}_{\mA_i^T\mA_i}\\\label{Lt1_Lt}
	&\qquad+3\left(1+c\right)\left(1+N_{\max}\right)\|\tilde{\vz}_i^{t+1}-\vz_i^t-\left(\tilde{\vz}_i^t-\vz_i^{t-1}\right)\|^{2}_{\mB_i^T\mB_i}\Big],
	\end{align}
	where $N_{\max}:=\max\{N_i,i\in\mathcal{N}\}$.
\end{Lemma}
\begin{IEEEproof}
	See the Section \ref{proof_LL} in the supplement material.
\end{IEEEproof}
Obviously, the right-hand-side (rhs) of \eqref{Lt1_Lt} is a sum of three positive terms, no matter how large $c$ and $\rho$ are, there is no guarantee that $\mathcal{L}\left(\vz,\vu,\vlambda\right)$ decreases at each iteration step. Thus, it cannot be used as a potential function. In search for an appropriate potential function, we add the constraint violation and the proximal term to the AL function, then a novel potential function is designed as follows
\begin{align}\label{poten}
\varsigma^t=\sum_{i\in\mathcal{N}}\frac{c}{2}\left[\kappa_{1}\|\mA_i\tilde{\vz}_i^t\|^{2}+\kappa_{2}\|\mA_i\vz_i^t\|^{2}+\frac{\rho}{{\blue2}c}\|\vu_i^t-\vu_i^{t-1}\|^{2}+\left(\kappa_{1}+\kappa_{2}\right)\|\vz_i^t-\vz_i^{t-1}\|^{2}_{\mB_i^T\mB_i}\right]+\mathcal{L}(\vz^t,\boldsymbol{\lambda}^t,\vu^t),
\end{align}
where $\kappa_{1},\kappa_{2},c,\rho>0$ are some positive constants, which can be determined by the following Lemma.

\begin{Lemma}\label{lem1-2}
	Suppose the sequence $\{\left(\vz_i^t,\vu_i^t,\vlambda_i^t\right)\}$ is generated by Algorithm \ref{al2} and $c\mB_i^T\mB_i$ takes the form of \eqref{cbb}. Then we have the following
	\begin{align}\nonumber
	\varsigma^{t+1}-\varsigma^t\leq&\sum_{i\in\mathcal{N}}\Big[-\frac{1}{2}\|\vz_i^{t+1}-\vz_i^t\|^{2}_{\mW_i}-\frac{\kappa_{1}-1}{2}\|\vz_i^{t+1}-\vz_i^t\|^{2}_{\mQ_i^T\mQ_i}-\frac{c\left(\kappa_{1}-1\right)}{2}\|\vz_i^{t+1}-\vz_i^t\|^{2}_{\mA_i^T\mA_i}\\\nonumber
	&{\blue-\frac{\rho}{4}\|\vu_i^{t+1}-\vu_i^{t}\|^{2}}-\frac{c\left(\kappa_{1}-6\left(N_{\max}+1\right)\right)}{2}\|\tilde{\vz}_i^{t+1}-\vz_i^{t+1}-\left(\tilde{\vz}_i^t-\vz_i^t\right)\|^{2}_{\mA_i^T\mA_i}\\\nonumber
	&-\frac{c\kappa_{1}-6\left(N_{\max}+1\right)}{2c}\|\mQ_i(\tilde{\vz}_i^{t+1}-\tilde{\vz}_i^t)-\mD_i\left(\vu_i^t-\vu_i^{t-1}\right)\|^{2}\\\nonumber
	&-\frac{c\kappa_{1}\hspace{-0.1cm}-6\left(1+c\right)\hspace{-0.1cm}\left(N_{\max}+1\right)}{2}\|\tilde{\vz}_i^{t+1}-\vz_i^t-(\tilde{\vz}_i^t-\vz_i^{t-1})\|^{2}_{\mB_i^T\mB_i}\\\nonumber
	&-\left(\frac{c\kappa_{2}\cdot\tilde{\tau}_{\min}}{2N_{\text{sum}}n\left(c+1\right)^{2}}-\frac{\left(N_{\max}+1\right)c\kappa_{1}}{2}\right)\|\tilde{\vz}_i^{t+1}-\vz_i^t\|^{2}\\\label{lem_var}
	&-\left(\frac{\rho}{{\blue4}}-d_{\max}^{2}\left(\kappa_{1}+\kappa_{2}\right)\right)\|\vu_i^t-\vu_i^{t-1}\|^{2}\Big],
	\end{align}	
	where $N_{\text{sum}}:=\sum_{i\in\mathcal{N}}N_i$ is the total number of neighboring nodes, $d_{\max}:=\max\{d_{i,j},i\in\mathcal{N},j\in\mathcal{N}_i\}$ is the maximum measurement distance and $\tilde{\tau}_{\min}:=\min\{\left(c+1\right)^{2}N_i^{2}+c^{2}N_i+N_i,i\in\mathcal{N}\}$.
\end{Lemma}
\begin{IEEEproof}
	See the Section \ref{apc} in the supplement material.
\end{IEEEproof}

From the above analysis, it is clear that as long as $\kappa_{1},\kappa_{2}$ and $\rho$ are sufficiently large, the rhs of \eqref{lem_var} is less than zero. As such, the potential function $\varsigma^t$ decreases at each iteration of SP-ADMM. Below, we derive the precise bounds for $\kappa_{1},\kappa_{2}$ and $\rho$. First, following the first five rows on the rhs of \eqref{lem_var}, a sufficient condition for $\kappa_{1}$ is given below (for any given $c>0$)
\begin{equation}\label{k1_all}
\kappa_{1}\geq 6\left(N_{\max}+1\right)\left(1+\frac{1}{c}\right).
\end{equation}
Second, for any given $c$ and $\kappa_{1}$, according to the sixth row on the rhs of \eqref{lem_var}, we need
\begin{equation}\label{k2_all}
\kappa_{2}\geq\frac{N_{\text{sum}}n\left(c+1\right)^{2}\left(N_{\max}+1\right)\kappa_{1}}{\tilde{\tau}_{\min}}.
\end{equation}
Finally, given $c,\kappa_{1},\kappa_{2}$, and based on the last row on the rhs of \eqref{lem_var}, parameter $\rho$ requires to satisfy
\begin{equation}\label{c_all}
\rho\geq {\blue4}d_{\max}^{2}\left(\kappa_{1}+\kappa_{2}\right).
\end{equation}
We conclude that if \eqref{k1_all}-\eqref{c_all} are satisfied, then the potential function $\varsigma^t$ will decrease at every iteration.

Now, we are ready to establish the main result. To this end, let us define the function $\mathcal{F}\left(\vz^t ,\vu^t,\vlambda^t\right)$ as the optimal gap of problem \eqref{f} given by

\begin{align}\label{F}
\mathcal{F}\left(\vz^t,\vu^t,\vlambda^t\right)=\sum_{i\in\mathcal{N}}\Big[\|\vz_i^t-\text{proj}_{\mathcal{X},\mathcal{Z}}\left(\vz_i^t-\left(\nabla_{\vz_i}
F_i\left(\vz_i^t,\vu_i^t\right)+\mA_i^T\vlambda_i^t\right)\right)\|^{2}+\|\mA_i\vz_i^t\|^{2}+\|\vu_i^t-\vu_i^{t-1}\|^{2}\Big].
\end{align}
Based on the function $\mathcal{F}\left(\vz^t,\vu^t,\vlambda^t\right)$ defined above, Lemma \ref{KKT_Critical} establishes its relationship with a KKT point of problem \eqref{f} and a critical point of problem \eqref{11}.
\begin{Lemma}\label{KKT_Critical}
	When $\mathcal{F}\left(\vz^t,\vu^t,\vlambda^t\right)=0$, then $(\vz^t,\vu^t,\vlambda^t)$ is a KKT solution of problem \eqref{f} satisfying:
	\begin{subequations}\label{KKT_cond}
		\begin{align}
		\label{KKT_cond1}\vz^t\in&\mathop{\arg\min}\limits_{\substack{\vz\in\mathcal{Z}\\\vz\in\mathcal{X}}}\sum_{i\in\mathcal{N}}F_i\left(\vz_i,\vu_i^t\right)+\langle\vlambda_i^t,\mA_i\vz_i\rangle,\\
		\label{KKT_cond2}\vzero\in&\nabla_{\vu_i}F_i\left(\vz_i^t,\vu_i^t\right)+\partial\delta_{\mathcal{B}^{N_i}}(\vu_{i}^t),\\
		\label{KKT_cond3}\vzero=&\mA_i\vz_i^t.
		\end{align}
	\end{subequations}
	Moreover, $(\vp^t,\vu^t)$ is a critical point of problem \eqref{11}, {\blue and $\vp^t$ is a critical point of the original problem \eqref{1}.}
\end{Lemma}
\begin{IEEEproof}
	See the Section \ref{proof_KKT_Critical} in the supplement material.
\end{IEEEproof}
Lemma \ref{KKT_Critical} suggests that we can demonstrate the sublinear convergence of Algorithm \ref{al2} to a critical point of the original problem \eqref{1} by proving the sublinear convergence of the sequence $\{(\vz^t,\vu^t,\vlambda^t)\}$ to a KKT stationary point of problem \eqref{f}. This convergence can be illustrated by invoking the following definition and lemma.

\begin{Def}
	Given problem \eqref{f}, we define $\left(\vz,\vu,\vlambda\right)$ as an $\varepsilon$-solution if $\|\mA_i\vz_i\|<\varepsilon$ and there exists a vector $\vv_1\in\nabla_{\vz_i}F_i(\vz_i,\vu_i)+\mA_i\vlambda_i+\partial\delta_{\mathcal{X}\times \mathcal{Z} }\left(\vz_i\right)$ and a vector  $\vv_2\in\nabla_{\vz_i}F_i(\vz_i,\vu_i)+\partial\delta_{\mathcal{B}^{N_i}}(\vu_i)$ such that $\|\vv_1\|\leq\varepsilon,\|\vv_2\|\leq\varepsilon$.
\end{Def}
\begin{Lemma}\label{sublinear-lem}
	Let parameters $c$ and $\rho$ satisfy (62)-(64), then for any $t>1$, there exists a constant $M>0$ such that we can find an $s\in\{1,2,...,t-1\}$ satisfying $(\vz_i^{s+1},\vu_i^{s+1},\vlambda_i^{s+1})$ being a $M/\sqrt{t-1}$-solution. This implies that an $\epsilon$-solution can be obtained within $M^2/\epsilon^2$ iterations.
\end{Lemma}
\begin{IEEEproof}
	See the Section \ref{proof_sublinear-lem} in the supplement material.
\end{IEEEproof}

\begin{theorem}\label{them2}
	For each nodes $i\in\mathcal{N}$, suppose the sequence $\left\{\left(\vz_i^t,\vu_i^t,\vlambda_i^t\right)\right\}$ is generated by Algorithm \ref{al2} and the conditions \eqref{k1_all}-\eqref{c_all} are satisfied by properly setting the parameters. Then we have
	\begin{itemize}
		\item {\bf Lower Bound:}
		\begin{equation*}
		\exists\,\,\underline{\varsigma}>-\infty\,\,s.t.\,\,\varsigma^t\geq\underline{\varsigma},\,\,\forall t>0.
		\end{equation*}
		\item {\bf Eventual Consensus:}
		\begin{align*}
		&\lim_{t\to\infty}\vz^{t+1}_i-\vz^t_i\rightarrow\mathbf{0},\quad\lim_{t\to\infty}\vlambda^{t+1}_i-\vlambda^t_i\rightarrow\mathbf{0},\\
		&\lim_{t\to \infty}\vu^{t+1}_i-\vu^t_i\rightarrow\mathbf{0},\quad\lim_{t\to \infty}\mA_i\vz^t_i\rightarrow\mathbf{0}.
		\end{align*}
		\item {\blue{\bf Global Convergence:} the sequence $\left\{\vy^t=\left(\vz^t,\vu^t,\vlambda^t\right)\right\}$ has a finite length, i.e.,  $\sum_{t=1}^{\infty}\|\vy^{t+1}-\vy^t\|<\infty,$ and it converges to a KKT point of problem \eqref{f}.}
		\item {\bf Convergence to Stationary Points\footnote{Given the linear constraint and nonsmooth term in the nonconvex optimization problem \eqref{f}, we refer to the points that satisfy the KKT condition as the stationary points (or KKT points) of the problem \eqref{f}, such as{\blue\cite{hong2016convergence,zhang2020proximal,kumar2016asynchronous,sun2019two,zhang2020improved}.}}:}
		$\mathcal{F}\left(\vz^t,\vu^t,\vlambda^t\right)\rightarrow\mathbf{0}$ and the iteration sequence $\left\{\left(\vz^t,\vu^t,\vlambda^t\right)\right\}$ converges to a KKT stationary point of problem \eqref{f}.
		\item {\bf Sublinear Convergence Rate:} For any given $\epsilon_{1}>0$, suppose that $\mathcal{F}(\vz^t,\vu^t,\vlambda^t)$ in \eqref{F} is less than $\epsilon_{1}$ for the first time in the  $T$-th iteration step, i.e.,
		\begin{equation*}
		T:=\arg\min_t\mathcal{F}(\vz^t,\vu^t,\vlambda^t)\leq\epsilon_{1}.
		\end{equation*}
		Then there exists a positive constant $\epsilon_{2}>0$ (the specific form is defined in the \eqref{epsilon2} of the supplement material) such that  $\epsilon_{1}\leq\frac{\epsilon_{2}}{T-1}$, which means the convergence rate of function $\mathcal{F}(\vz^t,\vu^t,\vlambda^t)$ is $\mathcal{O}(1/T)$.
	\end{itemize}
\end{theorem}
\begin{IEEEproof}
	See the Section \ref{proof_them} in the supplement material.
\end{IEEEproof}
To our knowledge, Theorem \ref{them2} is the first result that shows the $\mathcal{O}(1/T)$ convergence rate of the distributed ADMM algorithm for the nonconvex and nonsmooth localization problem. Since Algorithm \ref{al1} is a simplified version of Algorithm  \ref{al2} in Section \ref{sec: alg develop}, the convergence result also holds for Algorithm \ref{al1}. In the next section, numerical results will corroborate that the SP-ADMM algorithms can demonstrate more favorable convergence behavior than the benchmark methods.

\section{Numerical Results}\label{sec: simulation}

\noindent
We consider the two-dimensional sensor network localization problem and evaluate the performance of the proposed Algorithm \ref{al2} against several other state-of-the-art methods, including SDP \cite{biswas2006SDP}, SF \cite{soares2015simple}, AM-FD \cite{gur2020alternating} and  ADMM-H \cite{piovesan2016cooperative}. We remark that the proposed Algorithm \ref{al1} only simplifies the storage space of Algorithm \ref{al2}. Therefore, there is no difference in positioning accuracy between Algorithm \ref{al2} and Algorithm \ref{al1}. Note also that the first two methods employ convex relaxation, the third one solves the original nonconvex problem directly, and the last one is a hybrid convex/nonconvex solver. Following \cite{gur2020alternating}, we will also examine the performance of SP-ADMM-NAG50 and AM-FD-NAG50, where SP-ADMM-NAG50 means that we run 50 steps of the NAG method \cite{nesterov1983method} to obtain a good initial point, after that, the proposed SP-ADMM is used until convergence. For fairness, AM-FD-NAG50 uses the same initialization.

The criteria under which we compare the above algorithms are the practical running time counted as the maximal computation time among all parallel computing components and the averaged root mean squared error (RMSE) in a particular iteration, as described in \cite{erseghe2015distributed,piovesan2016cooperative},
\begin{align*}
\text{RMSE}(t)=\sqrt{\sum_{i\in(\mathcal{N}/\mathcal{A})}\frac{\|\vp^{t}_{i}-\vp_i\|^2}{N-m}},
\end{align*}
where $\mathcal{N}/\mathcal{A}$ represents the set of agents with unknown positions, which has $N-m$ elements,
$\vp_i$ is the true position of node $i$ and  $\vp^{t}_{i}$ is the estimated position of node $i$ in the $t$-th iteration of the algorithm. In addition, to check the convergence of the proposed algorithm as stated in Theorem \ref{them2}, let
\begin{align}\nonumber
S(t)&=\sum_{i\in\mathcal{N}}\|\nabla_{\vz_i}F_i\left(\vz_i^t,\vu_i^t\right)+\mA_i^T\vlambda_i^t\|^{2},\\\nonumber
U(t)&=\sum_{i\in\mathcal{N}}\|\vu_i^t-\vu_i^{t-1}\|^{2},\\
P(t)&=\sum_{i\in\mathcal{N}}\|\mA_i\vz_i^t\|^{2},
\end{align}
where $S(t)$ and $U(t)$ are defined as the stationarity gap while $P(t)$ is defined as the feasibility gap for problem \eqref{f}. The MATLAB code that implements the proposed SP-ADMM is available at: https://github.com/zm-stu/SP-ADMM.

\begin{table}[!t]
	\renewcommand{\arraystretch}{1.3}
	\caption{Benchmark Network and Algorithm Setup}
	\label{Benchmark parameter}
	\center
	\setlength{\tabcolsep}{3mm}{
		\begin{tabular}[l]{ccccc|cccccc|c}
			\hline
			&&&&&\multicolumn{7}{c}{Method Parameters} \\\cline{6-12}
			\multicolumn{5}{c|}{Network Parameters}&\multicolumn{6}{c|}{ ADMM-H }& AM-FD\\\hline	
			$N$&$m$&$C_{\text{range}}$&$\sigma_{\text{add}}$&$D_{\text{avg}}$&$\epsilon_{c}$&$\zeta_{c}$&$\tau_{c}$&$\theta_c$&$\delta_c$&$\lambda_{\max}$&$\vu^0$\\
			\hline
			500&10&0.3&0.02&14.15&0.002&0.25&0.008&0.98&1.01&$10^3$&$\vzero$\\\hline
			1000&20&0.1&0.007&11.09&0.003&0.07&0.002&0.98&1.01&$10^3$&$\vzero$\\
			\hline
	\end{tabular}}
\end{table}

\subsection{Benchmark Network}
\begin{figure}[!t]
	\centering
	\subfloat[AWGN]{\includegraphics[scale=0.33]{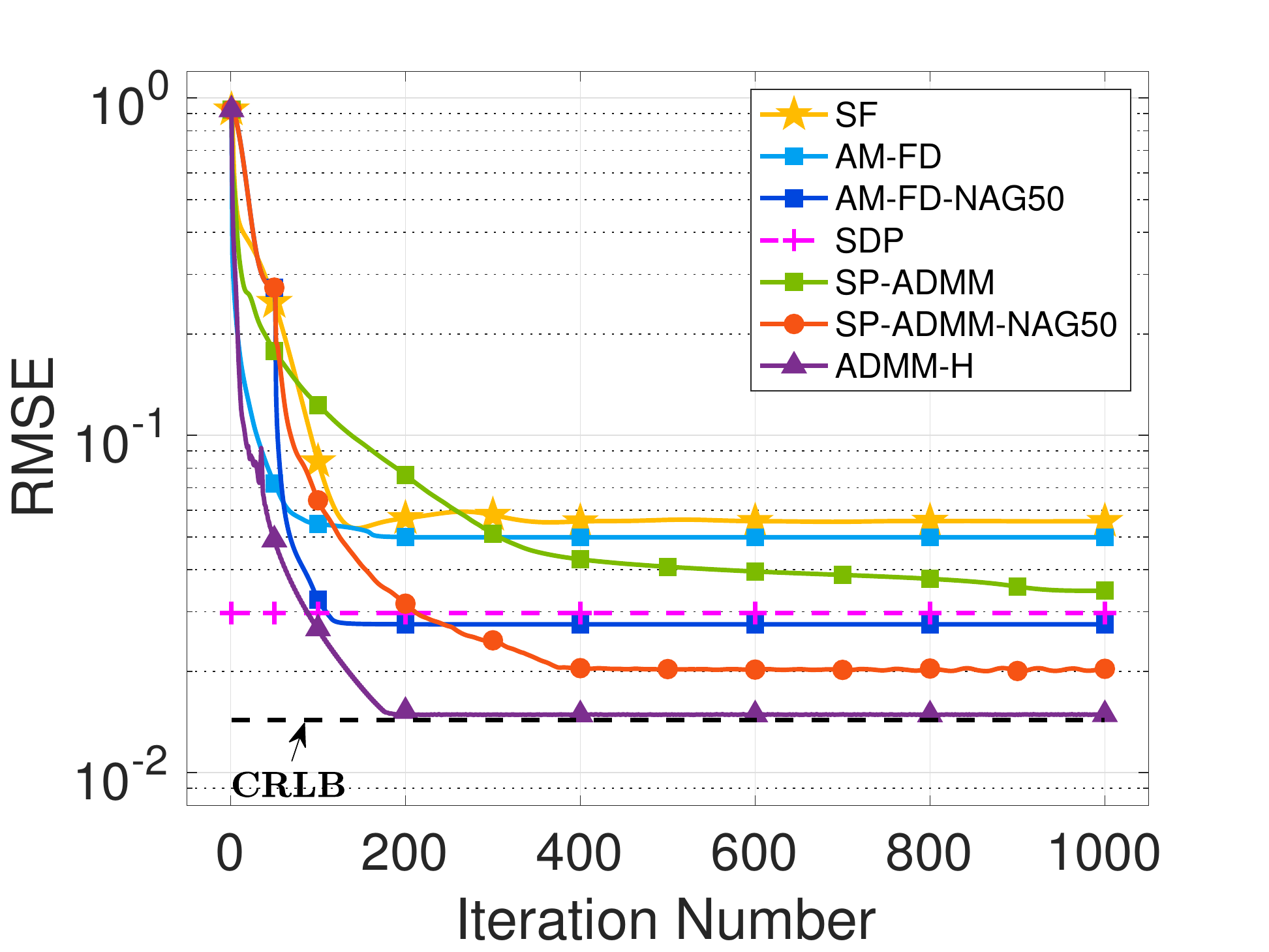}\label{500_a}}
	\subfloat[AWGN]{\includegraphics[scale=0.24]{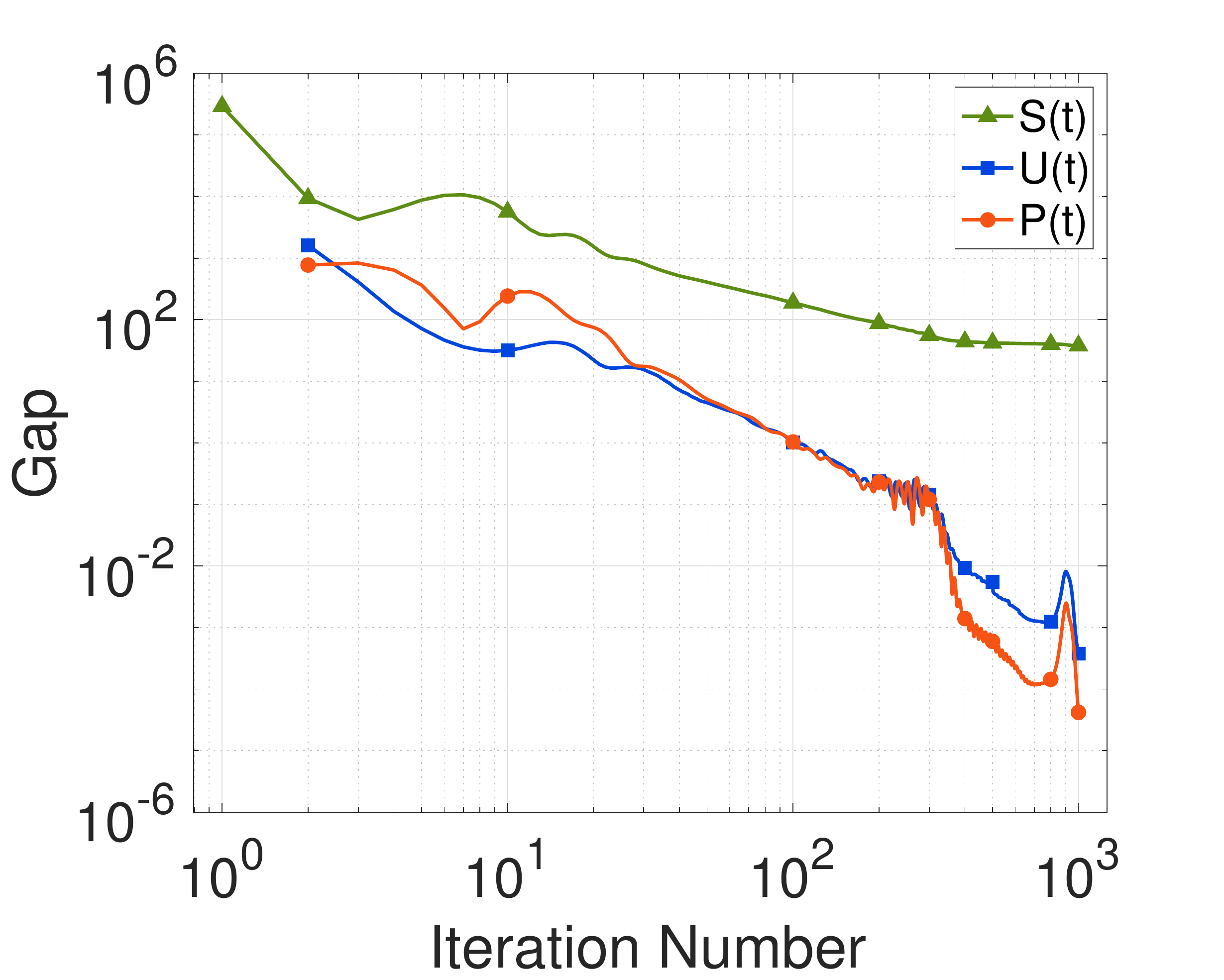}\label{500_b}}
	
	\subfloat[Range-dependent noise]{\includegraphics[scale=0.33]{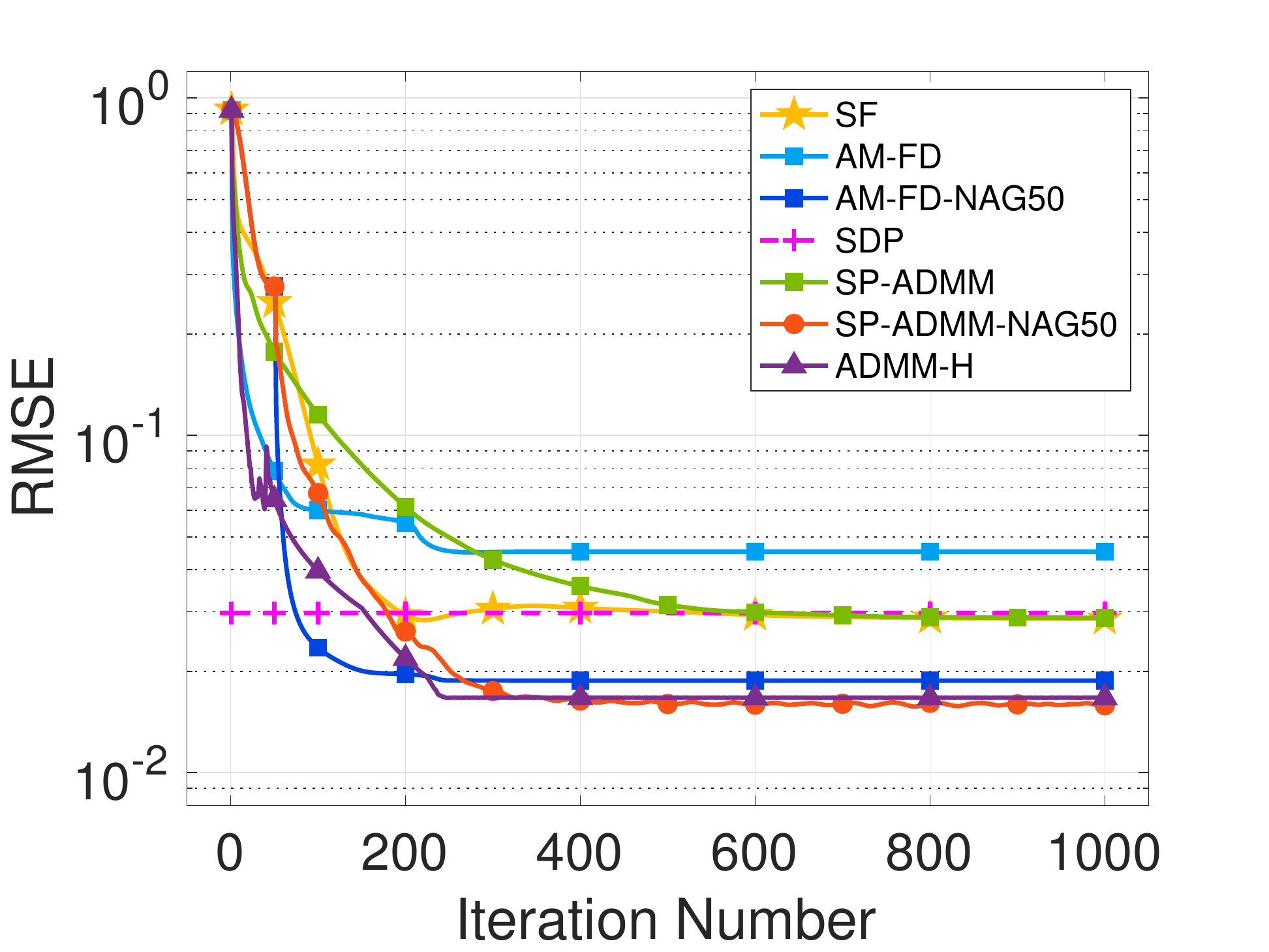}\label{500_c}}
	\subfloat[Range-dependent noise]{
		\includegraphics[scale=0.24]{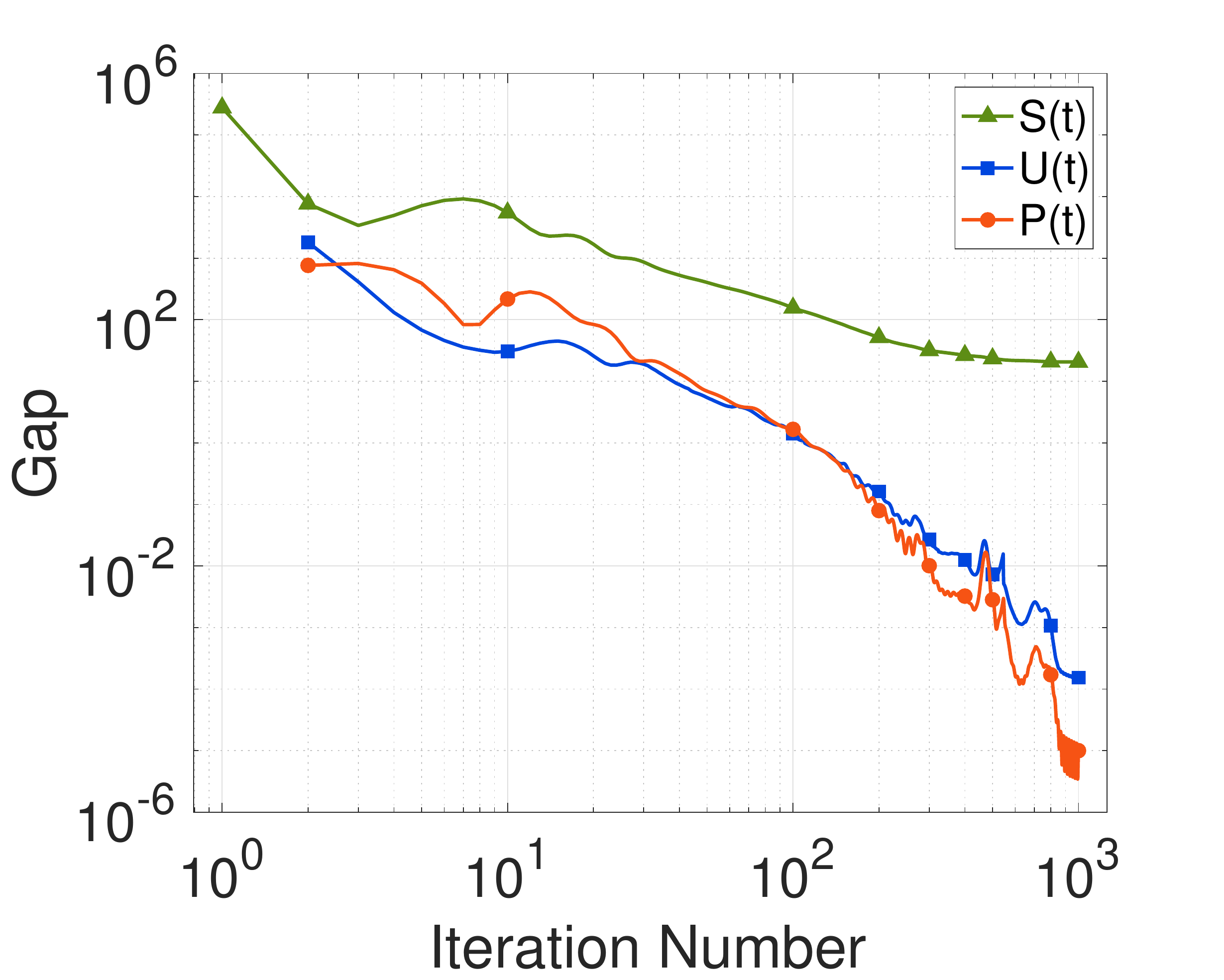}\label{500_d}}
	\caption{Performance with benchmark network ($N=500, m=10$). Measurement noise: AWGN (first row) and Range-dependent Gaussian noise $\sigma_{i,j}^{2}=\sigma_{\text{add}}\|\vp_i-\vp_j\|^{2}$ (second row). RMSE value (left column); feasibility gap and stationarity gap (right column).}
	\label{num_500}
\end{figure}

In the experiment, various methods are applied to the benchmark network data collected by the Stanford Computational Imaging Lab \cite{YY}. Table~\ref{Benchmark parameter} lists the details of the two networks and the specific parameters used for both the ADMM-H method and AM-FD method, which refer to the settings in \cite{piovesan2016cooperative} and \cite{gur2020alternating}. $C_{\text{range}}$ and $D_{\text{avg}}:=\frac{1}{N}\sum_{i\in\mathcal{N}}N_i$ in Table~\ref{Benchmark parameter} represent the communication range and the average number of neighbors, respectively. For the proposed algorithm, parameters are set to $\rho=c=0.11,\vu^0=\vzero$ for $N=500$ networks and $\rho=c=0.0197,\vu^0=\vzero$ for $N=1000$ networks. The initial point $\vz^{0}$ for both networks is selected from a uniform distribution, $\textbf{Unif}\left(-1,1\right)^{n\left(4\mid\mathcal{E}\mid+N\right)}$. Following \cite{jin2021exploiting}, we consider two different kinds of measurement noise: one is an additive white Gaussian noise (AWGN) with standard deviation $\sigma_{i,j}=\sigma_{\text{add}}$, while the other is the range dependent Gaussian noise, namely a zero mean Gaussian distribution with range dependent variance $\sigma_{i,j}^{2}=\sigma_{\text{add}}\|\vp_i-\vp_j\|^{2}$. The purpose of including the range-dependent Gaussian noise experiment is to evaluate the robustness of our proposed algorithm.

\begin{figure}[!t]
	\centering
	\subfloat[AWGN]{
		\includegraphics[scale=0.33]{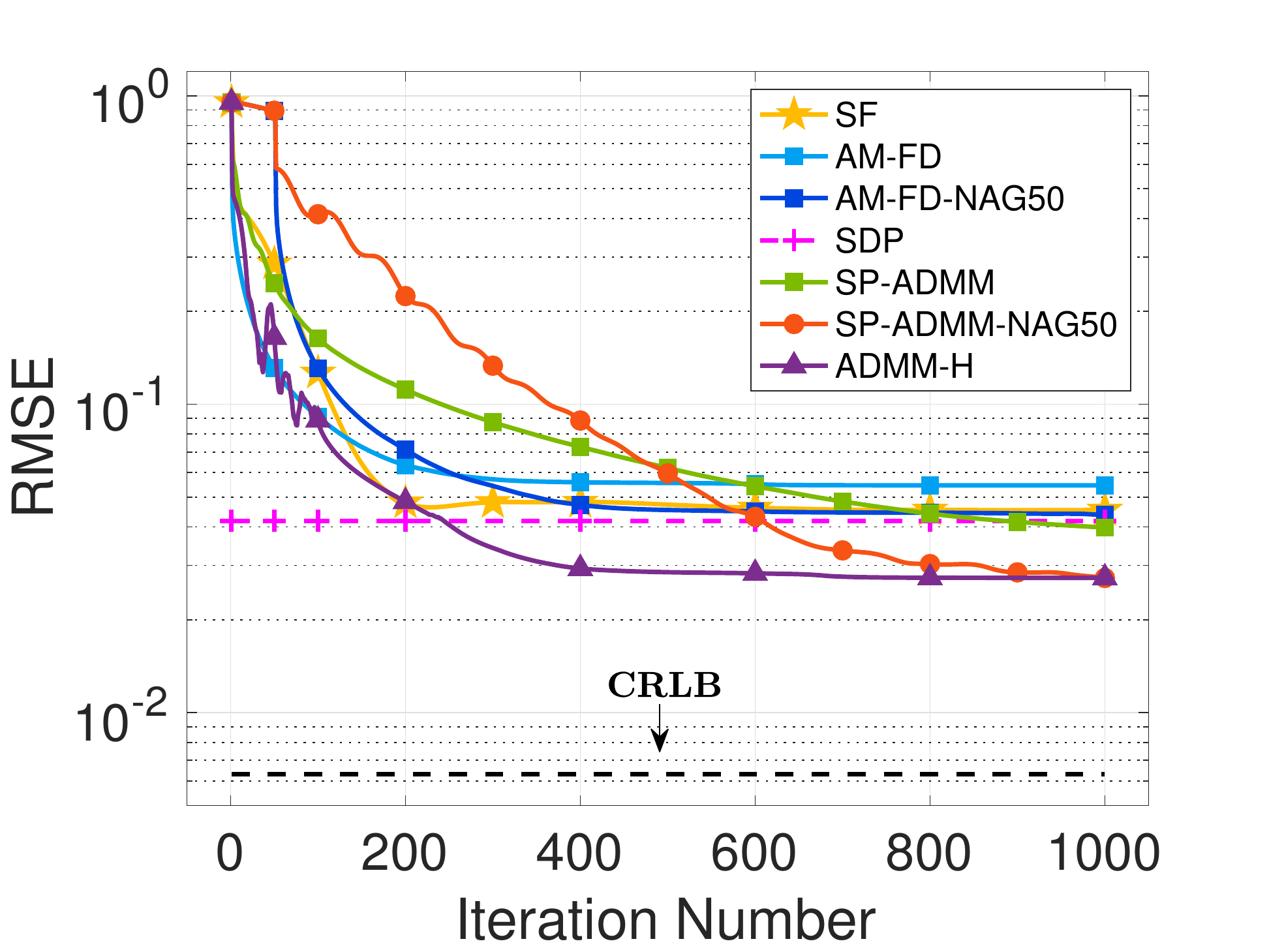}\label{1000_a}
	}
	\subfloat[AWGN]{
		\includegraphics[scale=0.24]{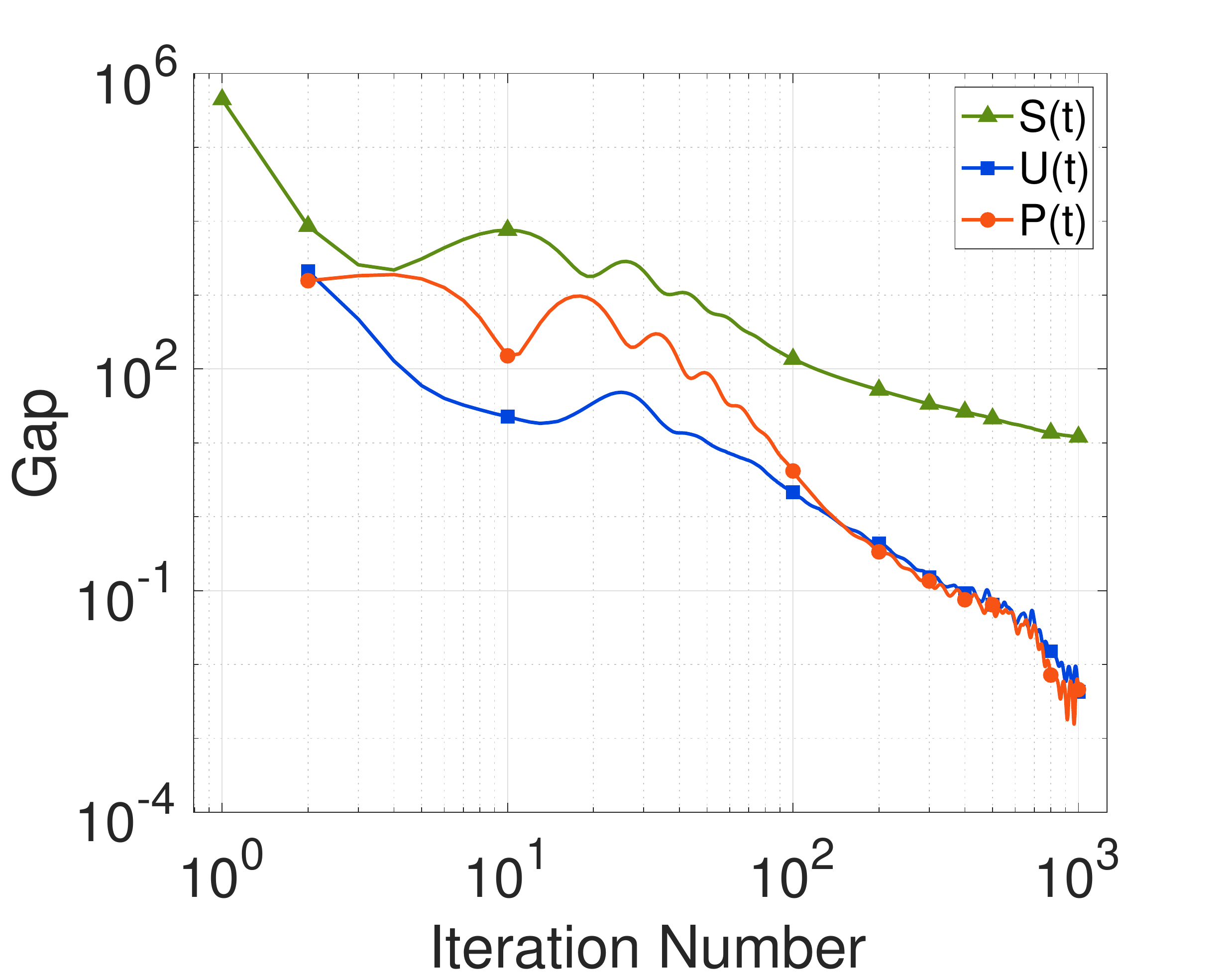}\label{1000_b}}
	\,
	\subfloat[Range-dependent noise]{
		\includegraphics[scale=0.33]{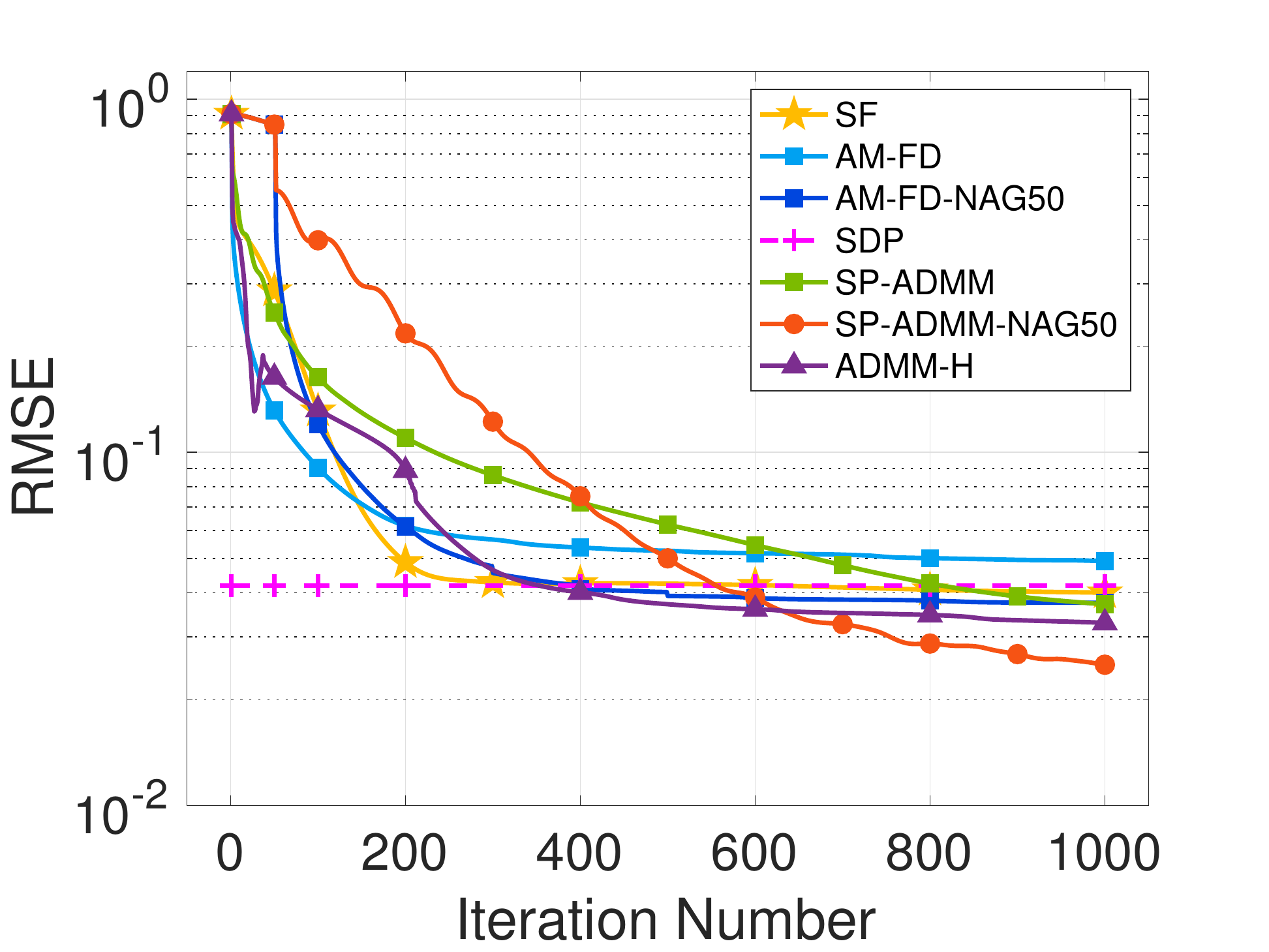}\label{1000_c}
	}
	\subfloat[Range-dependent noise]{
		\includegraphics[scale=0.24]{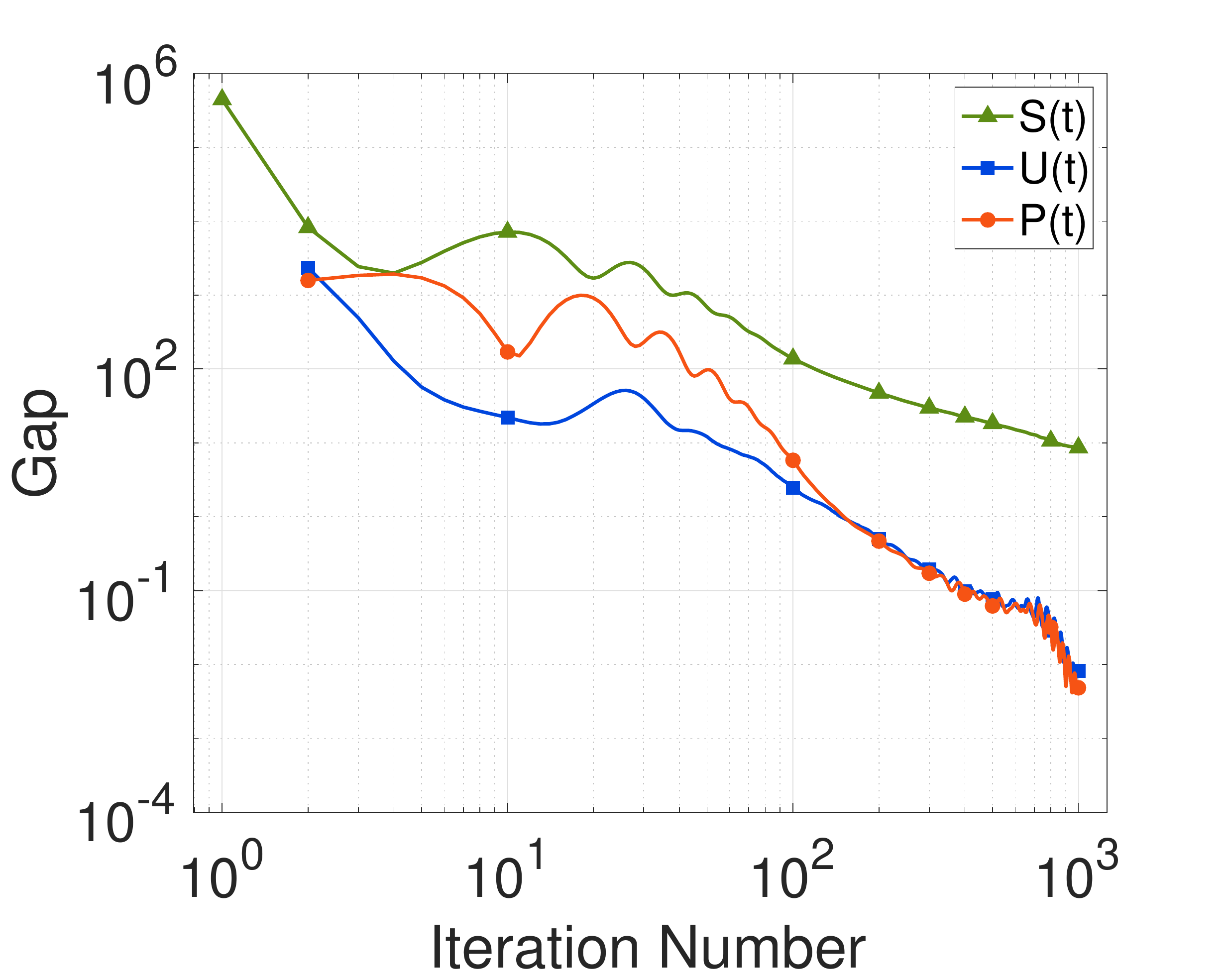}\label{1000_d}}
	\caption{Performance with benchmark network (N=1000, m=20). Measurement noise: AWGN (first row) and Range-dependent Gaussian noise $\sigma_{i,j}^{2}=\sigma_{\text{add}}\|\vp_i-\vp_j\|^{2}$ (second row). RMSE value (left column); feasibility gap and stationarity gap (right column).}
	\label{num_1000}
\end{figure}
\begin{table}[!t]
	\renewcommand{\arraystretch}{1.3}
	\center
	\begin{threeparttable}[b]
		\caption{Comparisons of Running Time}
		\label{benchmark running time}
		\scriptsize
		\setlength{\tabcolsep}{4mm}{
			\begin{tabular}[l]{c|c|c|c}
				\hline
				&&\multicolumn{2}{c}{Run time (seconds)}\\\cline{3-4}
				Algorithm&RMSE (AWGN)&Parallelized&Sequential(Per Step)\\\hline\hline
				\multicolumn{4}{c}{Benchmark Network ($N=500,\text{step}=1000,\text{CRLB}\approx 0.014$)}\\\hline
				SF&9.85e02&0.05401&1.24e03\,(1.2428)\\
				AM-FD&4.99e-02&-\tnote{$a$}&6.19e02\,(0.6193)\\
				AM-FD-NAG50&2.75e-02&-&6.75e02\,(0.6753)\\
				SDP&2.97e-02&-&1.98e02\,(6.1815)\\
				ADMM-H&1.49e-02&0.1518&3.50e03\,(3.4957)\\
				\textbf{SP-ADMM-NAG50}&2.03e-02&\textbf{0.01091}&6.64e02\,(0.6642)\\
				\textbf{SP-ADMM}&2.98e-02&\textbf{0.01784}&6.30e02\,(0.6298)\\\hline
				\multicolumn{4}{c}{Benchmark Network ($N=1000,\text{step}=1000,\text{CRLB}\approx0.005$)}\\\hline
				SF&7.41e-02&0.045&6.09e03\,(6.085)\\
				AM-FD&5.07e-02&-&1.43e02\,(0.1430)\\
				AM-FD-NAG50&4.17e-02&-&1.45e02\,(0.1455)\\
				SDP&4.19e-02&-&1.56e03\,(23.9601)\\
				ADMM-H&2.73e-02&0.152&1.06e04\,(10.5978)\\
				\textbf{SP-ADMM-NAG50}&2.94e-02&\textbf{0.034}&9.63e01\,(0.0965)\\
				\textbf{SP-ADMM}&3.73e-02&\textbf{0.033}&9.04e01\,(0.0905)\\
				\hline
		\end{tabular}}
		\begin{tablenotes}
			\item [$a$] These symbols ``-" in the table mean that the methods cannot be implemented in parallel.
		\end{tablenotes}
	\end{threeparttable}
\end{table}
\begin{figure}[!t]
	\centering
	\subfloat[$N=500,m=10$]{
		\includegraphics[scale=0.5]{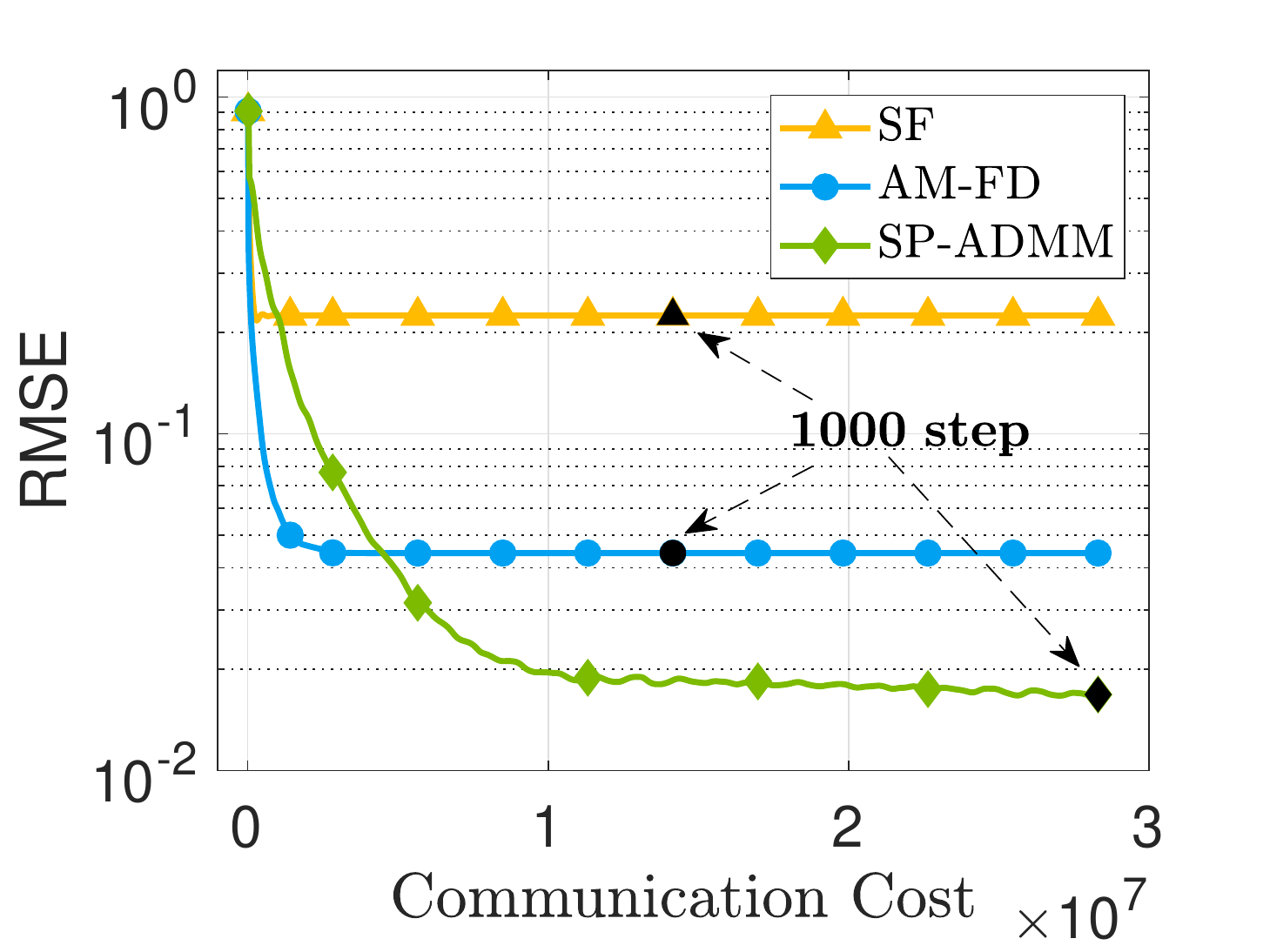}\label{comm_500}}
	\hfil\hspace{-0.65cm}
	\subfloat[$N=1000,m=20$]{
		\includegraphics[scale=0.5]{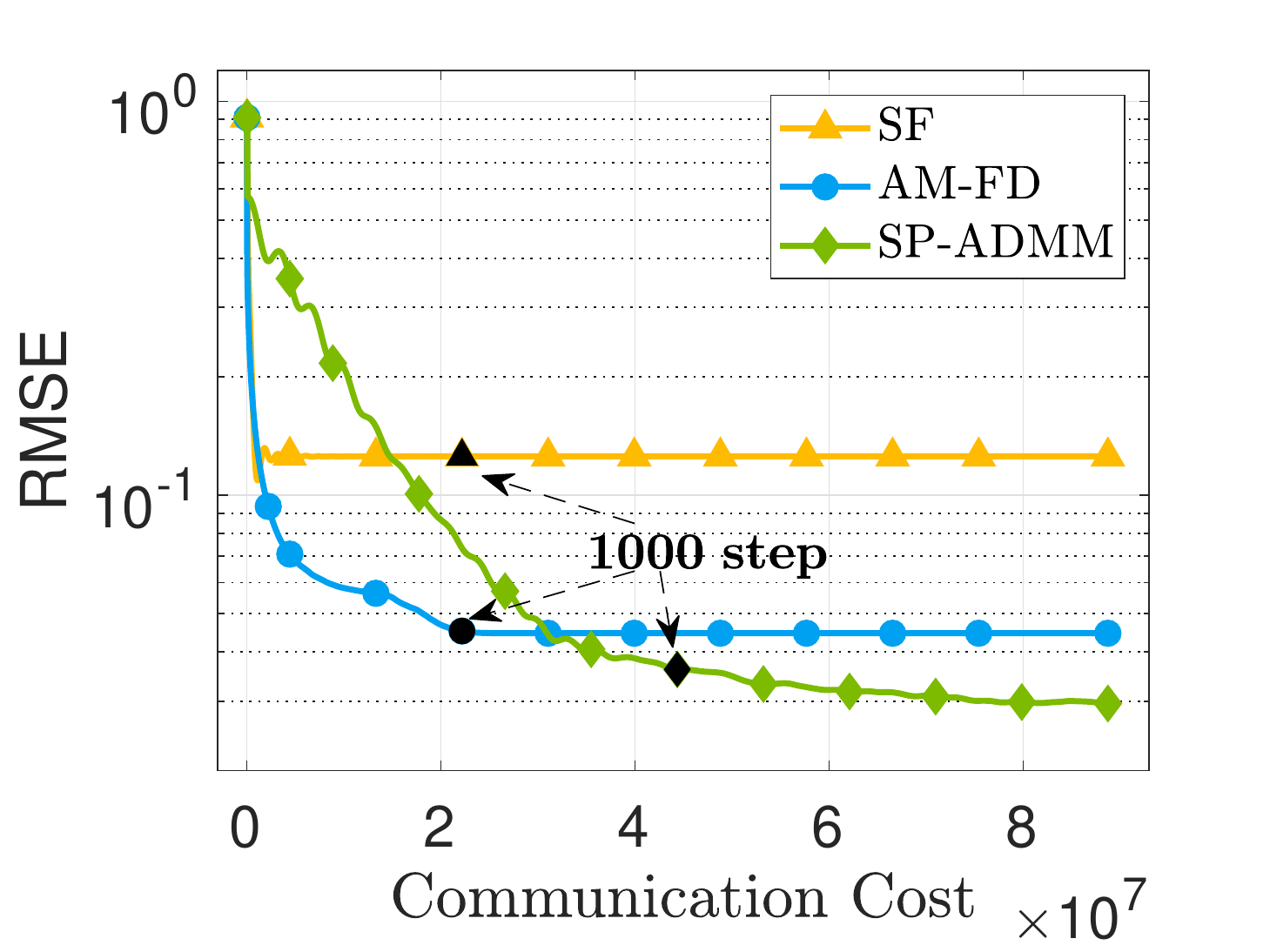}\label{comm_1000}}
	\caption{Comparison of RMSE and communication cost for different methods in the benchmark network with AWGN noise.}
	\label{comm_cost}
\end{figure}

\textbf{Convergence performance:}
Let us first examine the convergence
behaviors of the proposed algorithm, SDP \cite{biswas2006SDP}, SF \cite{soares2015simple}, AM-FD \cite{gur2020alternating}, and ADMM-H \cite{piovesan2016cooperative}. Fig. \ref{num_500}\subref{500_a} and Fig. \ref{num_500}\subref{500_c} display the RMSE versus the iteration number when we use AWGN and range-dependent Gaussian noise, respectively. As shown in Fig. \ref{num_500}\subref{500_a} that ADMM-H based on a two-stage approach performs the best and its RMSE approaches the Cramer–Rao Lower Bound (CRLB) \cite{patwari2003relative}, but it is less advantageous in terms of the running time (see Table~\ref{benchmark running time}). We note that the AM-FD-NAG50 and SP-ADMM-NAG50 show lower RMSE than both SF and SDP, illustrating that non-relaxed problems generally result in better location estimation than that of the relaxed problems. When we focus on the first-order methods, Fig. \ref{num_500}\subref{500_a} and Fig. \ref{num_500}\subref{500_c} demonstrate that the SP-ADMM-NAG50 algorithm achieves higher accuracy with fewer iterations than other methods for both AWGN and range-dependent Gaussian noise scenarios. Furthermore, we see that AM-FD-NAG50 outperforms AM-FD, which reflects the advantage of using the NAG method to provide an initial value that is beneficial for the nonconvex localization problem \eqref{1}.

In Fig. \ref{num_1000}, we plot the RMSE versus the iteration number for a large network with $(N=1000, m=20)$. Further comparing Fig. \ref{num_1000}\subref{1000_a} (Fig. \ref{num_1000}\subref{1000_c}) with Fig. \ref{num_500}\subref{500_a} (Fig. \ref{num_500}\subref{500_c}), we observe that the performance of all the methods is degraded. Intriguingly, Fig. \ref{num_1000}\subref{1000_c} shows that the proposed method achieves the best performance under range-dependent Gaussian noise. It implies that the proposed method may be more suitable for a large network in reality, possibly due to the noise and error elasticity of ADMM \cite{simonetto2014distributed}.

\textbf{Gap:} To further examine the convergence of the SP-ADMM algorithm, we plot in Fig. \ref{num_500}\subref{500_b} (Fig. \ref{num_500}\subref{500_d}) and Fig. \ref{num_1000}\subref{1000_b} (Fig. \ref{num_1000}\subref{1000_d}) the curves of the optimal gap of problem \eqref{f} versus the iteration number. Observe explicitly that the performance gap of our proposed SP-ADMM algorithm reduces as the iteration number increases. This result verifies the efficacy of the Theorem \ref{them2}, which states that the iterative sequence $\left\{\left(\vz^t,\vu^t,\vlambda^t\right)\right\}$ generated by Algorithm \ref{al2} converges to a KKT stationary point of problem \eqref{f} and the convergence rate is $\mathcal{O}\left(1/T\right)$.

\textbf{Running time:} Table~\ref{benchmark running time} shows the running time of the methods in Fig. \ref{num_500} and Fig. \ref{num_1000}. It can be seen that the proposed SP-ADMM algorithm is most efficient in computation time than other methods. Note that while SF has a competitive running time compared to ours, its RMSE is bigger than ours. In addition, although the SDP solver here uses SeDuMi \cite{sturm1999using}, it still consumed a lot of computational time, especially when the network size $N$ is large. SDP terminates when the norm of the constraint gap reaches the order of $10^{-7}$, while the other methods perform 1000 iterations. In contrast, AM-FD and AM-FD-NAG50 require modestly increased computation time as the network size grows. Nevertheless, compared to other parallel methods, AM-FD and AM-FD-NAG50 are not the most sensible choices for large-scale networks and networks with a large average number of neighboring nodes, causing extra delays by the sequential structure under distributed networks. As a result, the structure of parallel implementation would be more efficient. Lastly, we note that the application of the ADMM-H method to large-scale networks is impractical due to the unaffordable sequential running time. Moreover, it can be seen from Table~\ref{Benchmark parameter} that the ADMM-H has more parameters to adjust for different networks.

\textbf{Communication cost:} To evaluate the efficiency of various communication methods, we performed a comparison of the RMSE of different methods under benchmark networks as communication costs increased, as depicted in Fig. \ref{comm_cost}. Our results reveal that although the proposed SP-ADMM incurs a higher communication cost than the SF and AM-FD algorithms during the same iteration (i.e., 1000 steps), it achieves higher accuracy at the same communication cost once convergence has been reached. Additionally, Fig. \ref{comm_cost} also indicates that the AM-FD and SF methods were incapable of meeting higher accuracy requirements, even with sufficient communication resources.

\subsection{Synthesized Network}
\begin{figure}[t]
	\subfloat[$N=2000$]{
		\includegraphics[scale=0.3]{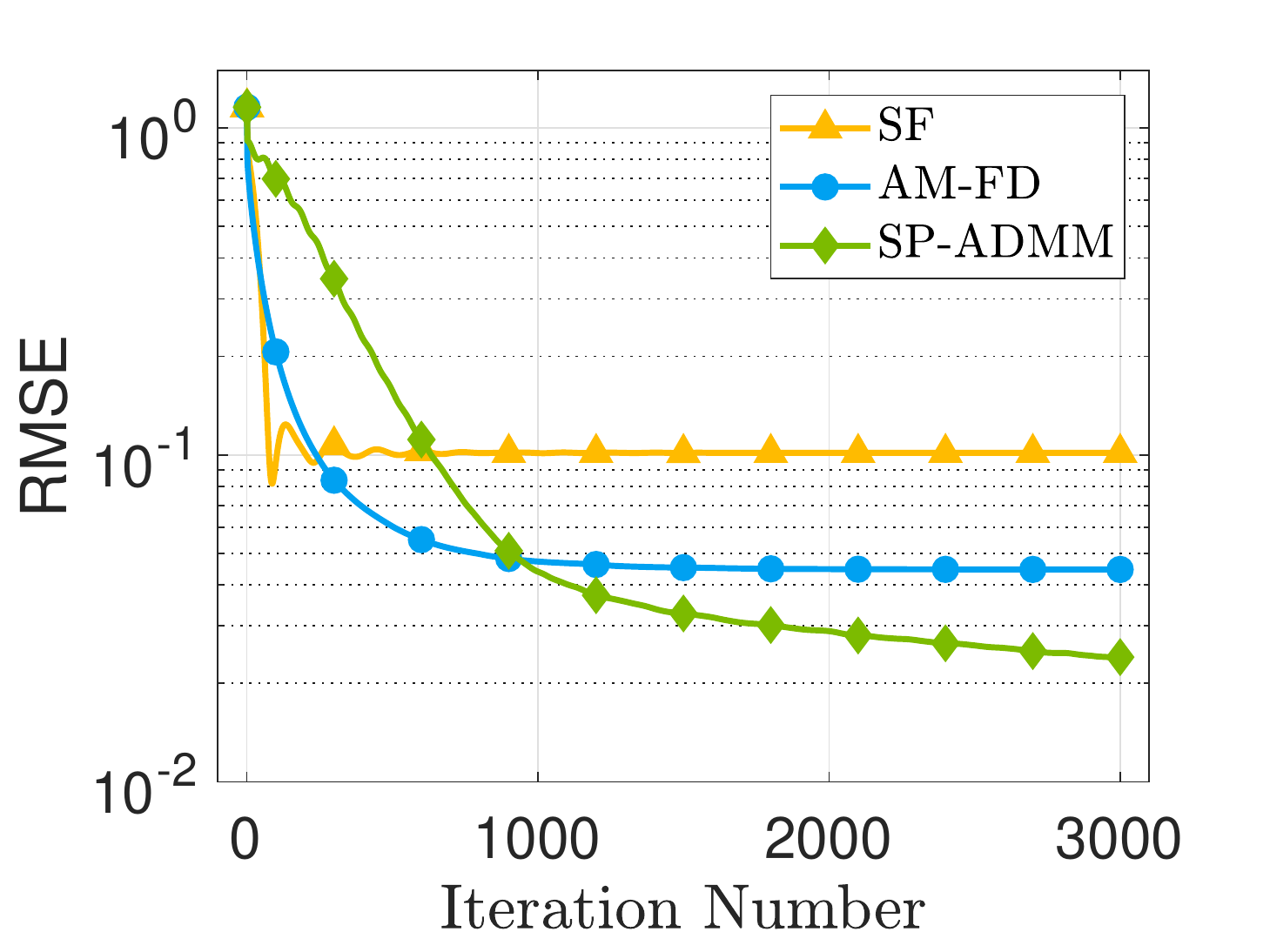}}
	\subfloat[$N=3000$]{
		\includegraphics[scale=0.3]{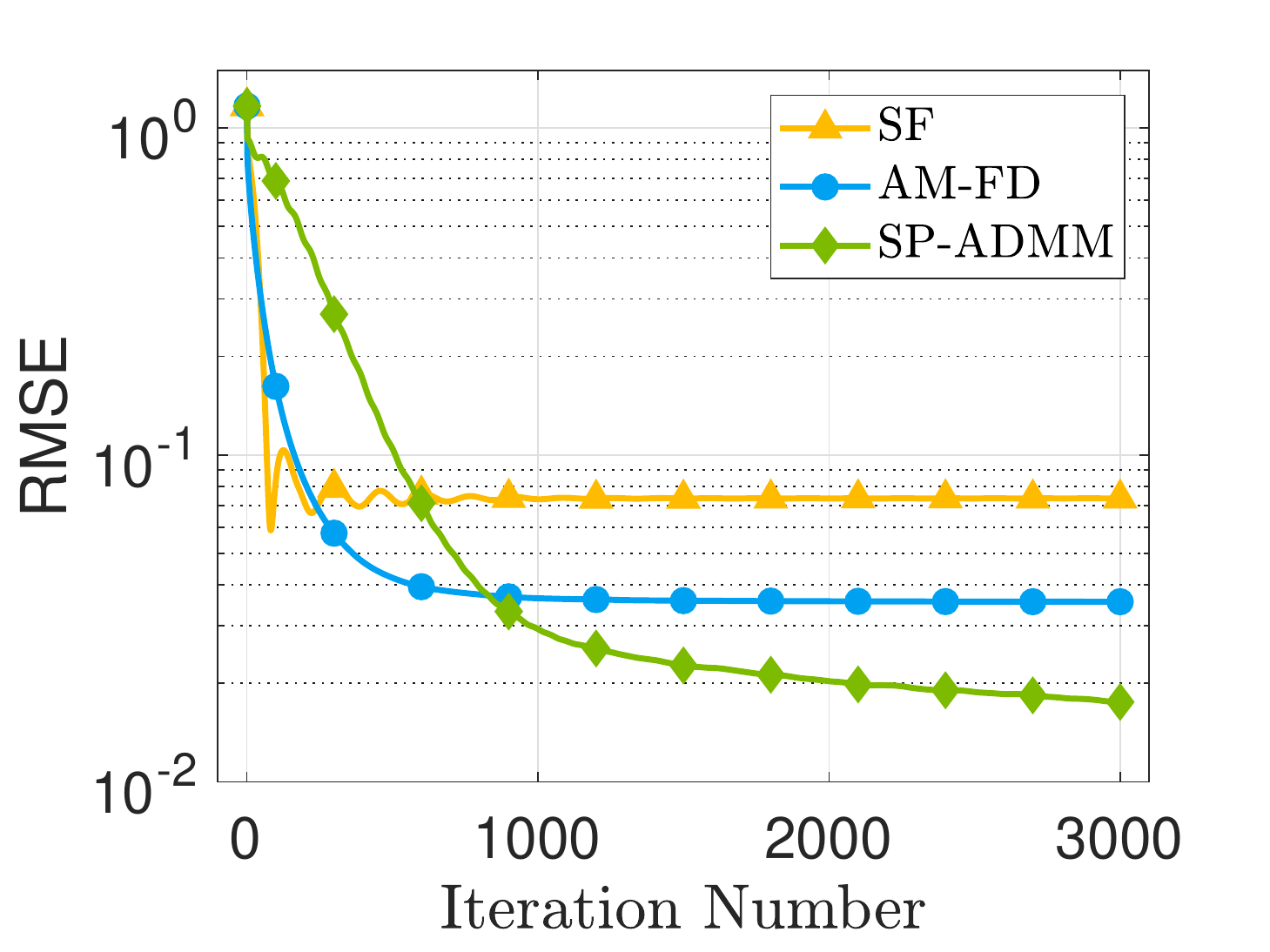}}
	\subfloat[$N=5000$]{
		\includegraphics[scale=0.3]{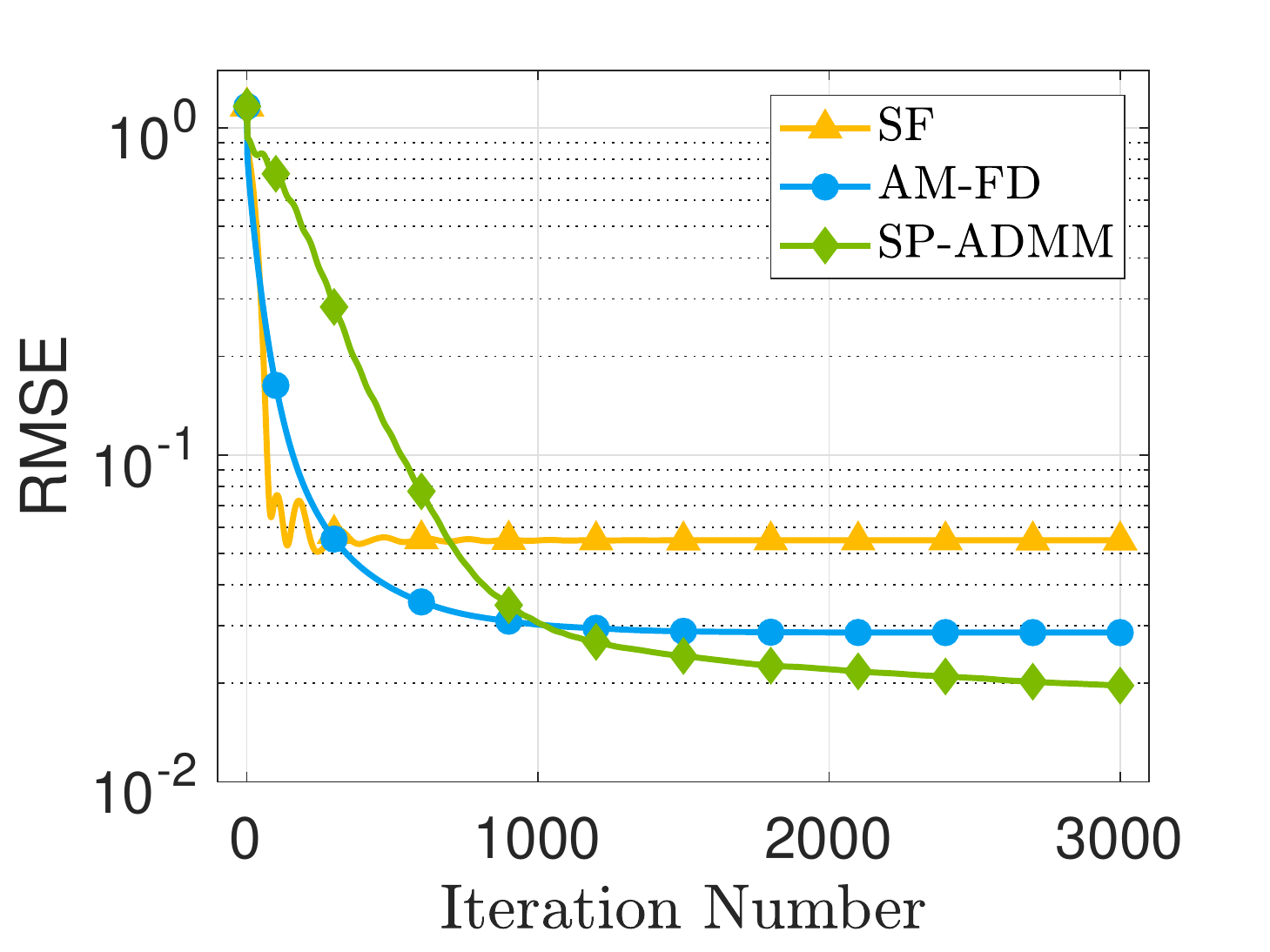}}
	\subfloat[$N=10000$]{
		\includegraphics[scale=0.3]{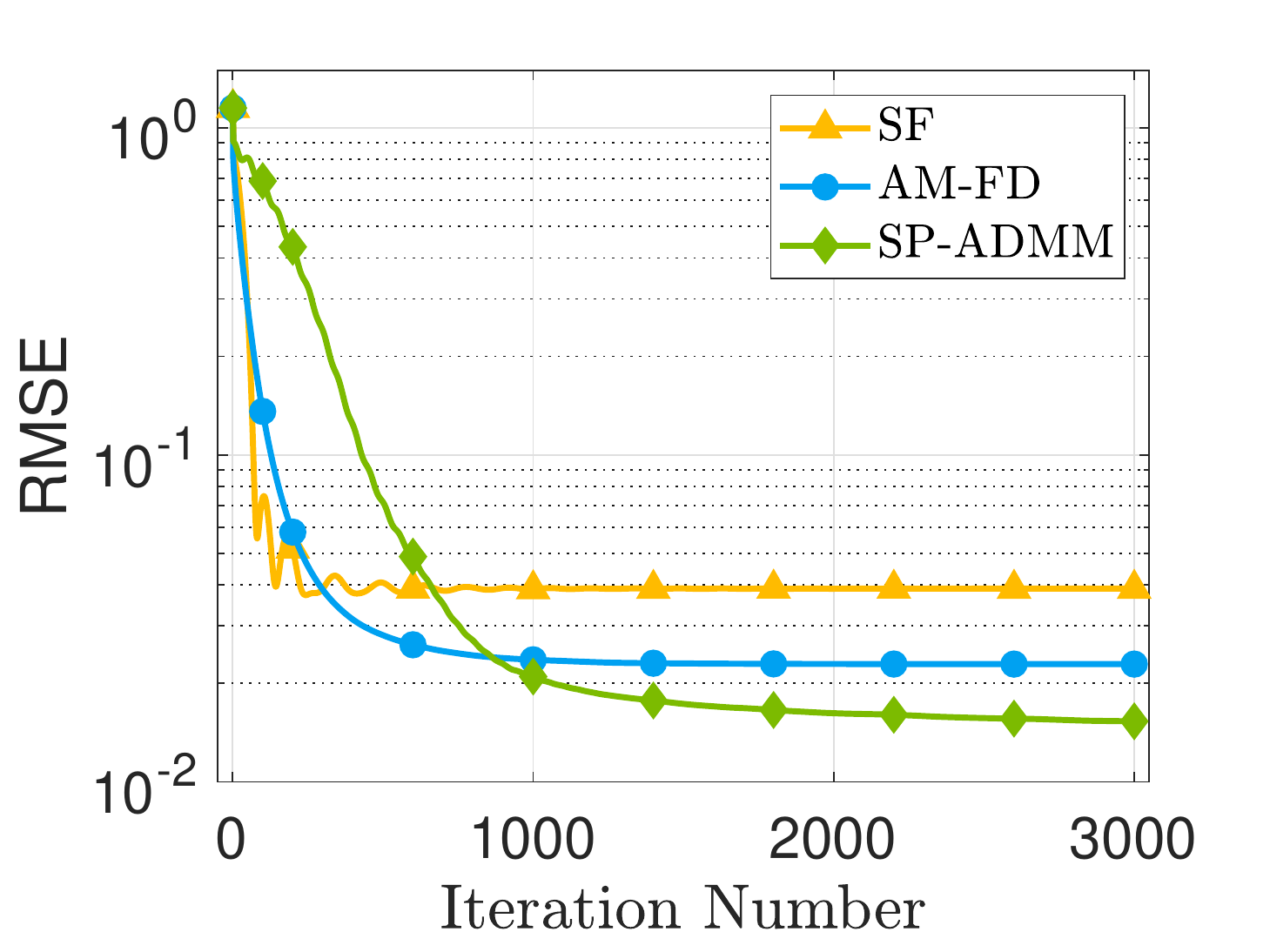}}
	
	\subfloat[$N=2000$]{
		\includegraphics[scale=0.3]{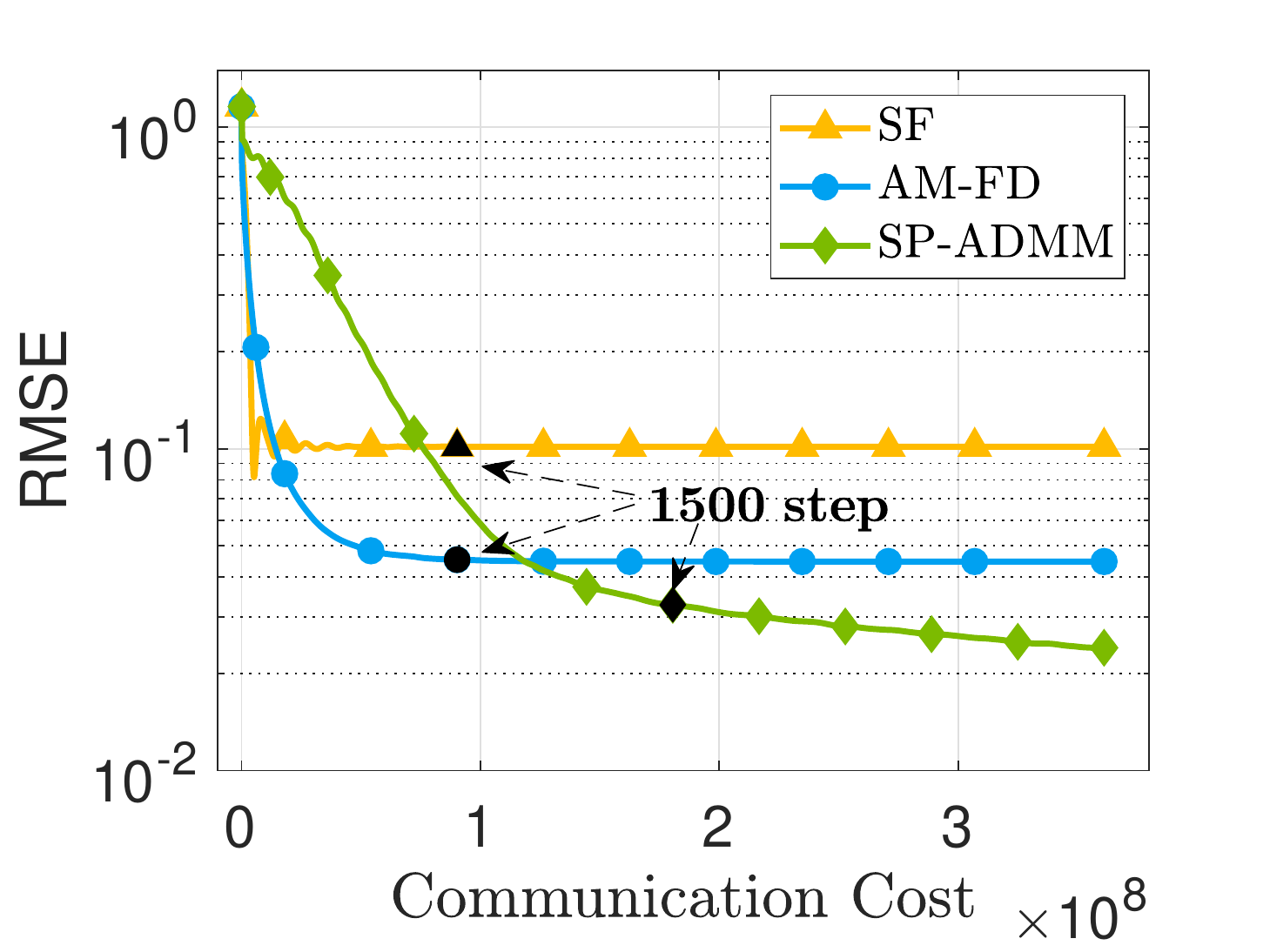}}
	\subfloat[$N=3000$]{
		\includegraphics[scale=0.3]{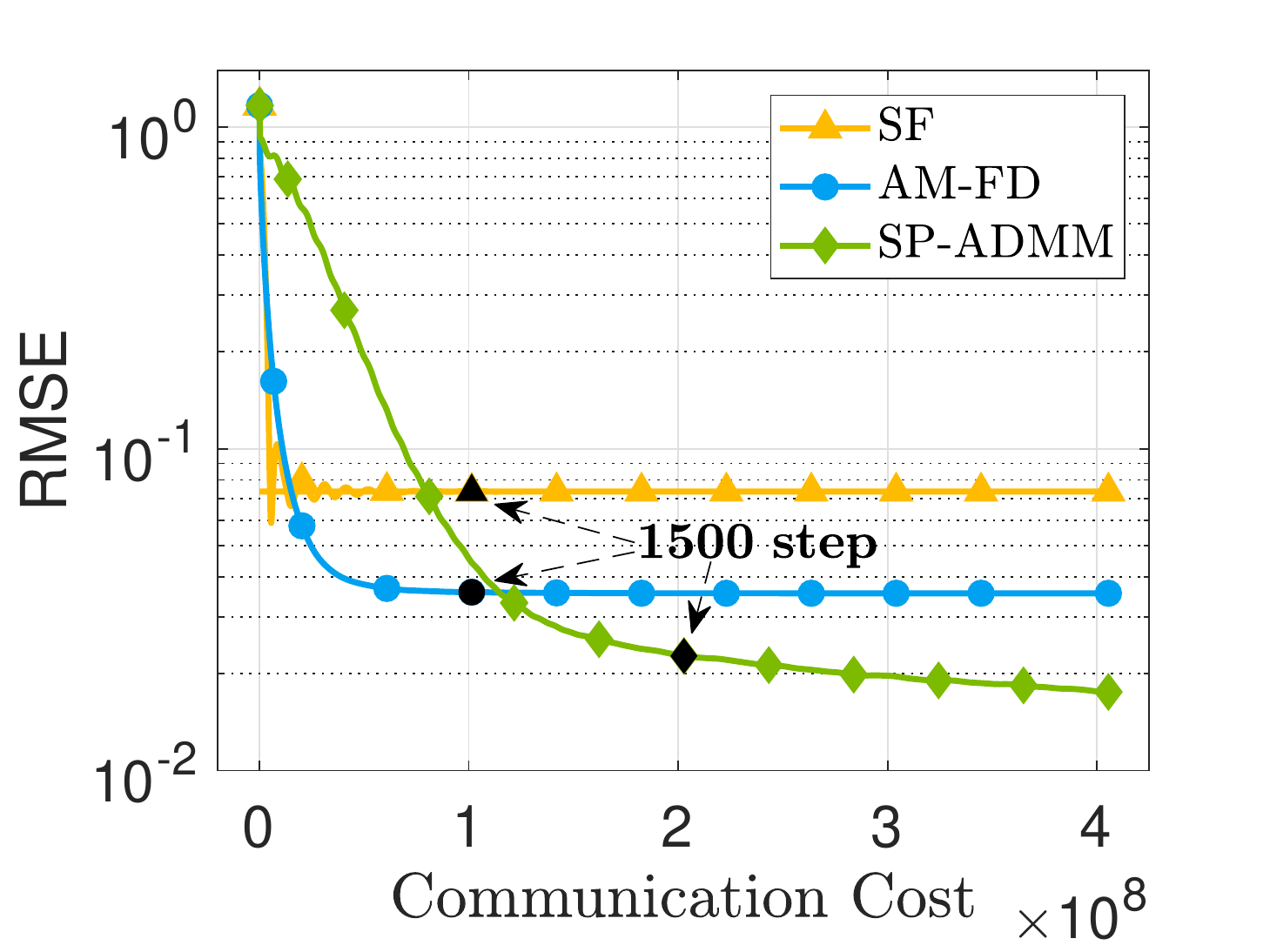}}
	\subfloat[$N=5000$]{
		\includegraphics[scale=0.3]{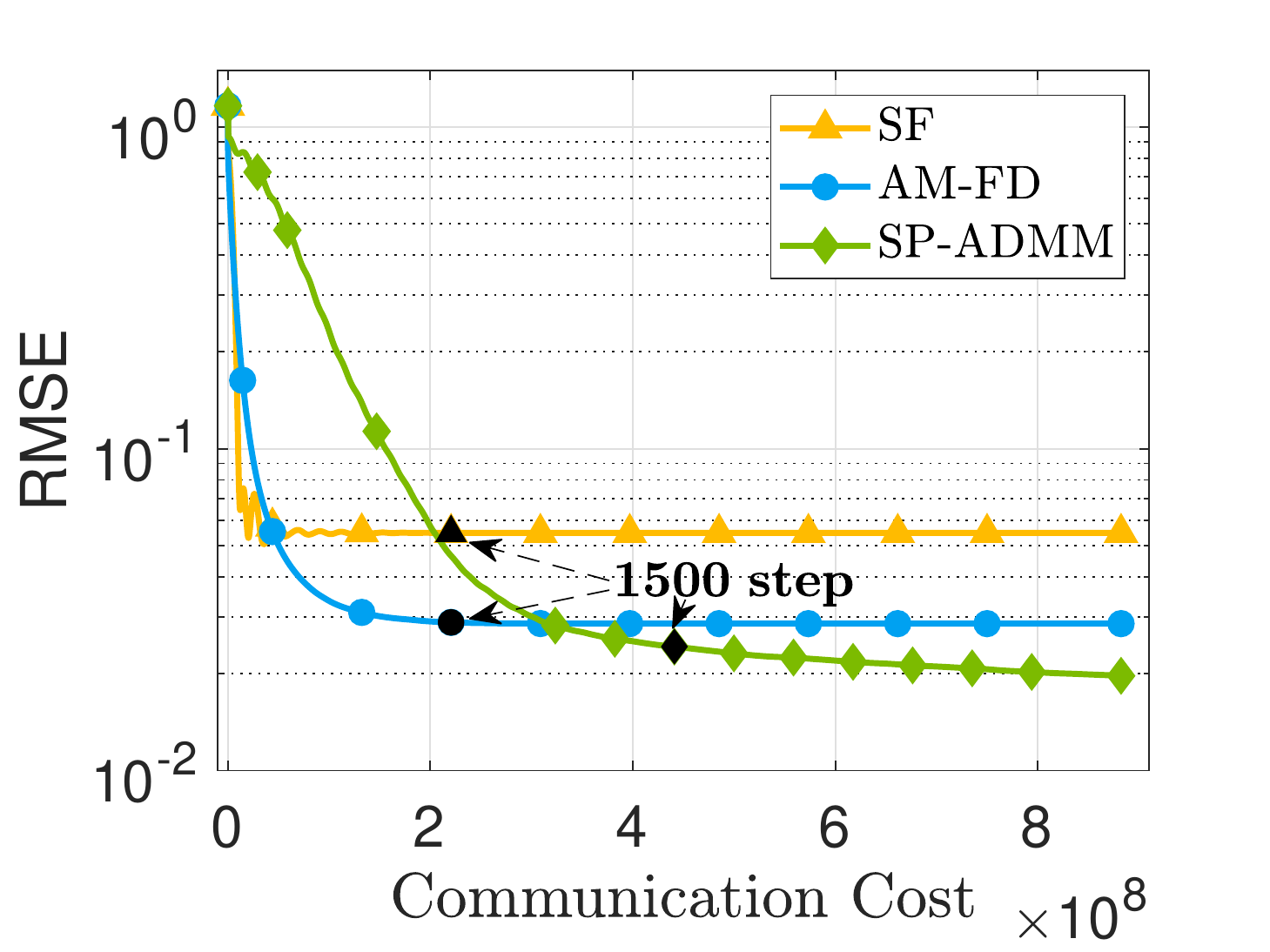}}
	\subfloat[$N=10000$]{
		\includegraphics[scale=0.3]{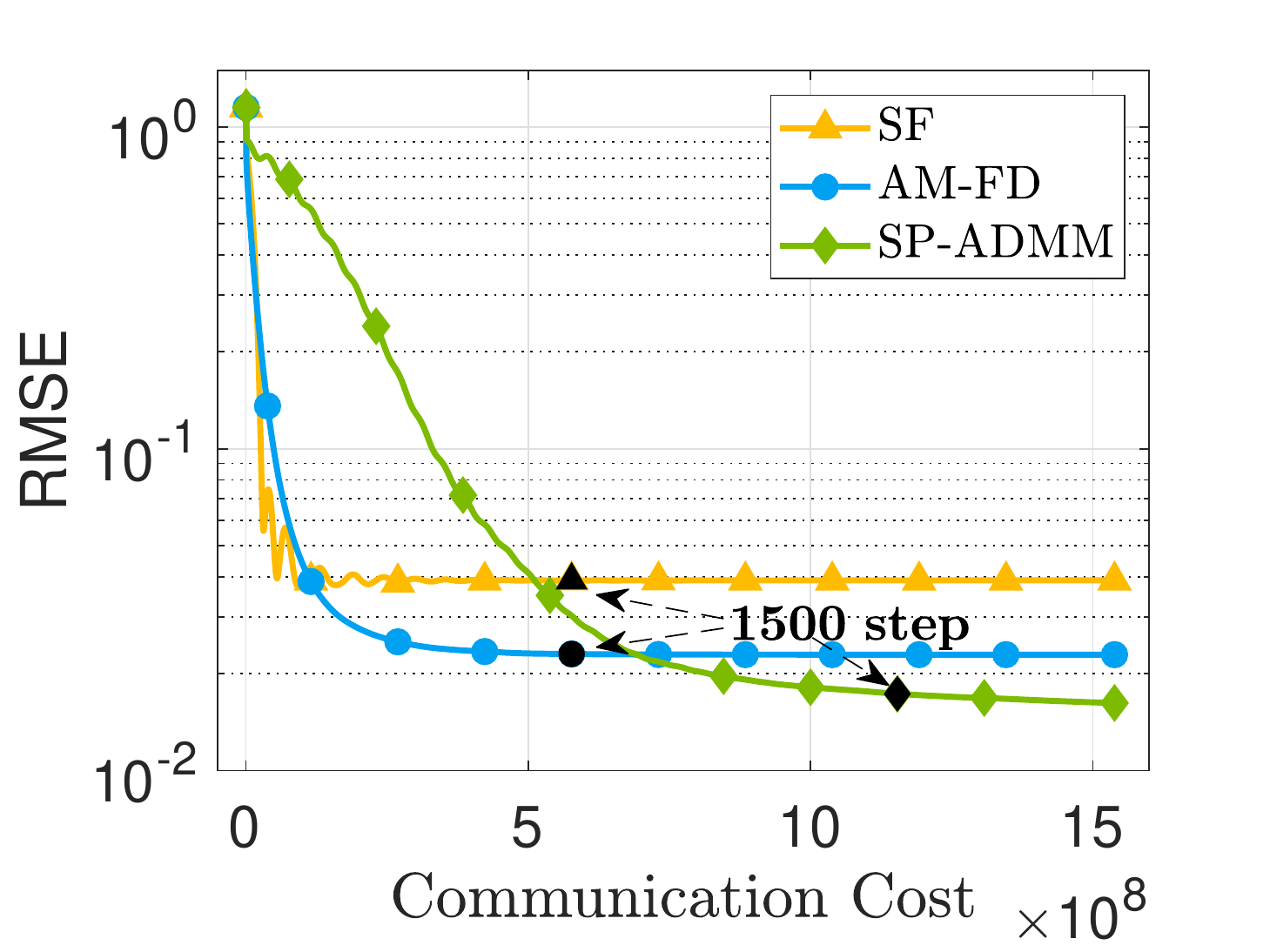}}
	\caption{Convergence performance of different methods with various sizes. (a)-(d) RMSE value vs. iteration number; (e)-(h) RMSE value vs. communication cost.}
	\label{change_N}
\end{figure}

In this subsection, we evaluate the performance of the proposed SP-ADMM algorithm in comparison to ADMM-H \cite{piovesan2016cooperative}, SDP \cite{biswas2006SDP}, SF \cite{soares2015simple}, and AM-FD \cite{gur2020alternating} under varying factors that could potentially impact localization accuracy.

\subsubsection{\textbf{Influence of the Number of Sensors ($N$)}}

We initiated our investigation by evaluating the performance of random networks with sizes $N$=2000, 3000, 5000, and 10000. Similar to the benchmark networks, the anchor number was set to 2$\%$ of the network size,  and the value of $\sigma_{\text{add}}$ of AWGN was chosen to be 7$\%$ of the communication range. The results of our experiments are presented in Fig. \ref{change_N}, where Fig. \ref{change_N}(a)-(d) depicts RMSE versus iteration steps and Fig. \ref{change_N}(e)-(h) illustrates RMSE versus communication loss. Similar to the results obtained for the benchmark networks shown in Fig. \ref{comm_cost}, our findings demonstrate that the proposed SP-ADMM algorithm incurs a higher communication cost compared to the SF and AM-FD algorithms during the same number of iterations (i.e., 1500 steps). However, once convergence has been achieved, the SP-ADMM algorithm provides higher accuracy at the same communication cost. Moreover, Fig. \ref{change_N}(e)-(h) show that even with sufficient communication resources, the AM-FD and SF methods are incapable of meeting higher accuracy requirements.

\subsubsection{\textbf{Influence of the Number of Anchors ($m$)}}
\begin{table}[!t]
	\renewcommand{\arraystretch}{1.3}
	\caption{Synthesized Network and Algorithm Setup}
	\label{synthesized parameter}
	\center
	\setlength{\tabcolsep}{3mm}{
		\begin{tabular}[l]{ccccc|cccccc|c}
			\hline
			&&&&&\multicolumn{7}{|c}{Method Parameters}\\\cline{6-12}
			\multicolumn{5}{c|}{Network Parameters}&\multicolumn{6}{|c|}{ADMM-H}&AM-FD\\\hline	
			$N$&$m$&$C_{\text{range}}$&$\sigma_{\text{add}}$&$D_{\text{avg}}$&$\epsilon_{c}$&$\zeta_{c}$&$\theta_c$&$\delta_c$&$\lambda_{\max}$&$\tau_{c}$&$\vu^0$\\
			\hline\hline\multicolumn{12}{c}{Random (Changing $m$)}\\\hline
			495&5&0.3&0.02&12.88&\multirow{3}{*}{0.002}&\multirow{3}{*}{0.25}&\multirow{3}{*}{0.008}&\multirow{3}{*}{0.98}&\multirow{3}{*}{1.01}&\multirow{3}{*}{$10^3$}&\multirow{3}{*}{$\vzero$}\\
			515&20&0.3&0.02&16.48&&&&&&&\\
			525&30&0.3&0.02&18.50&&&&&&&\\
			\hline
	\end{tabular}}
\end{table}

\begin{figure}[!t]
	\centering
	\subfloat[$m=5$]{
		\includegraphics[width=2.55in,height=2in]{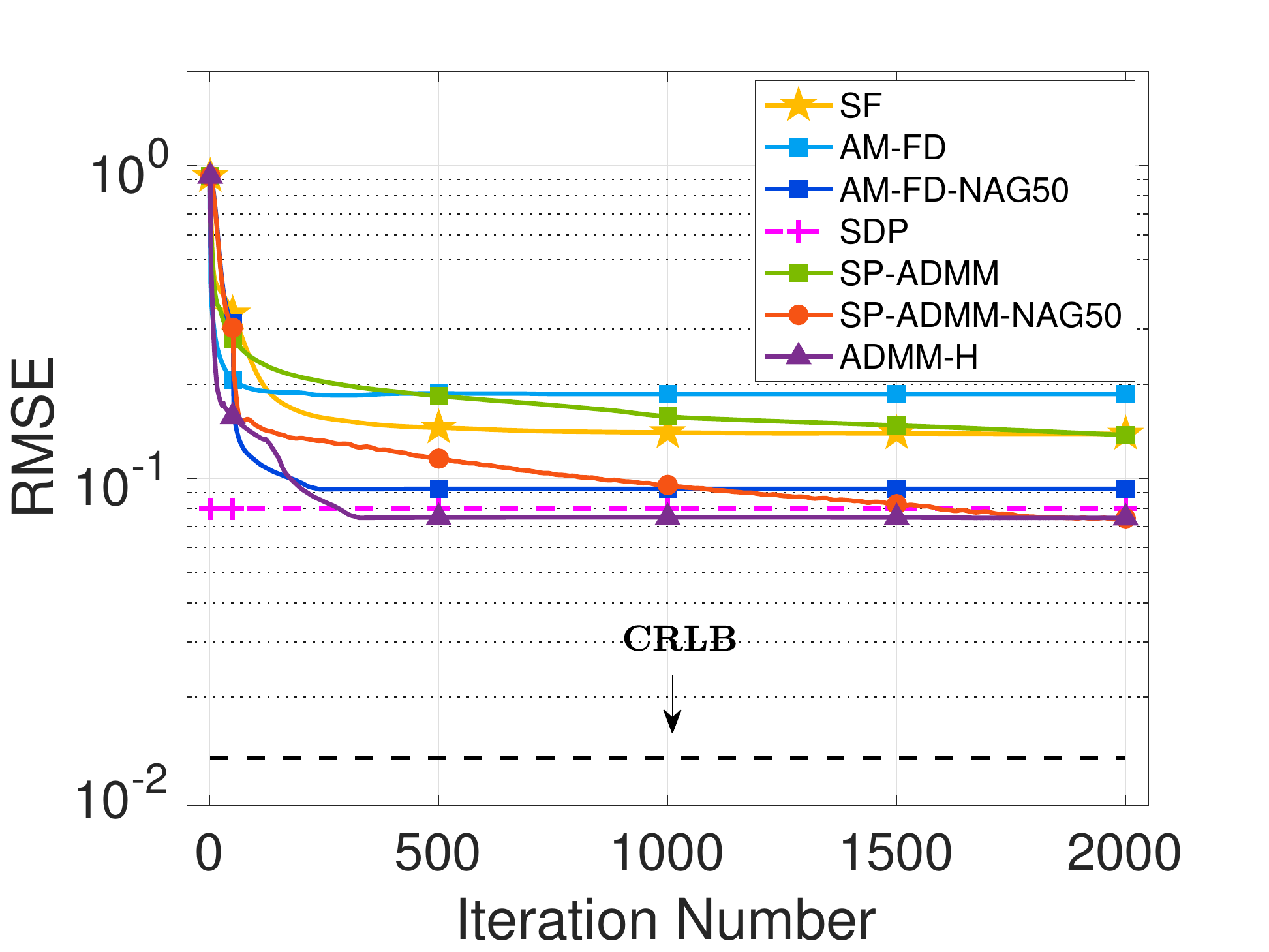}\label{m_5a}}
	\hfil\hspace{-0.65cm}
	\subfloat[$m=5$]{
		\includegraphics[width=2.55in,height=2in]{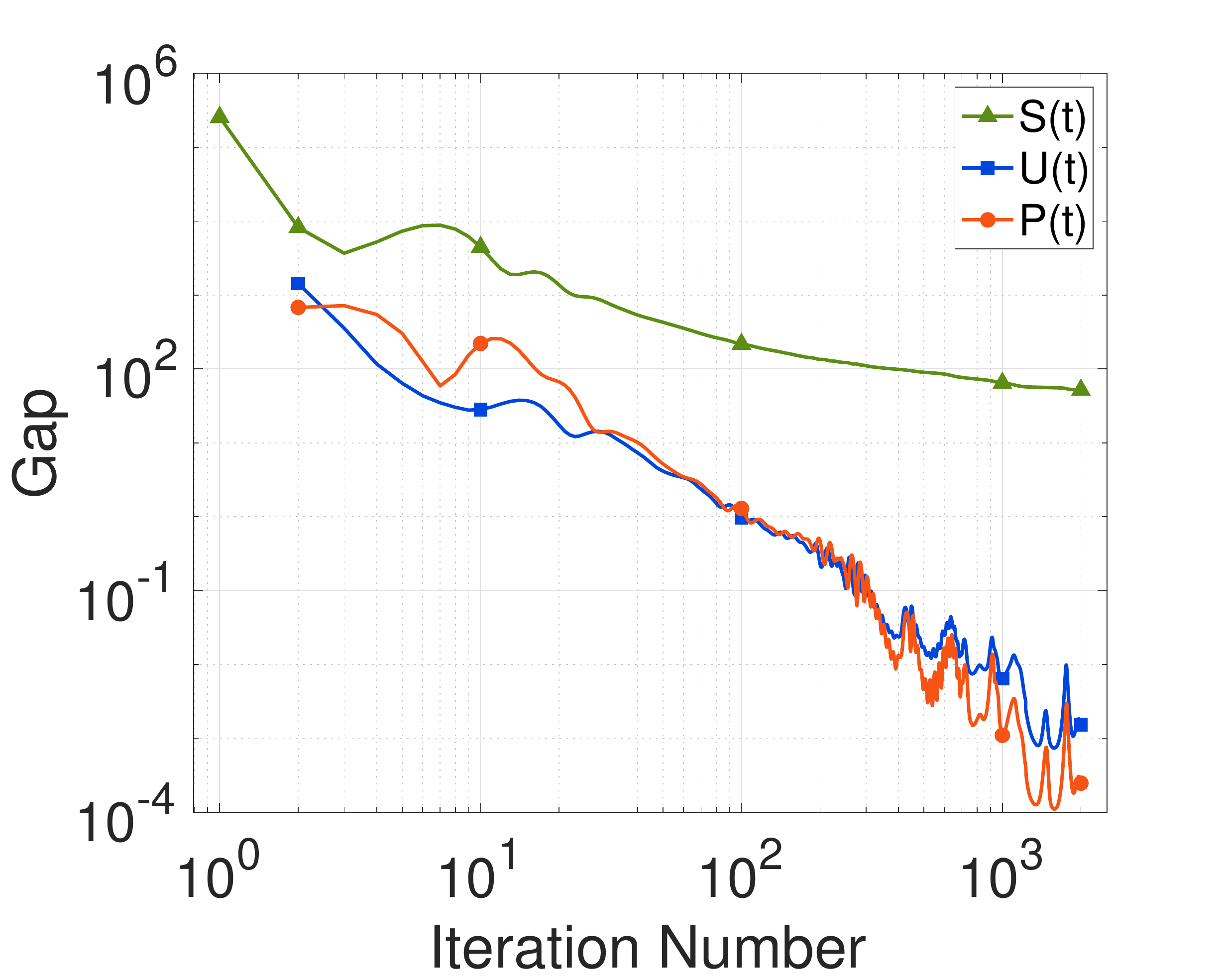}\label{m_5b}}
	\,
	\subfloat[$m=20$]{
		\includegraphics[width=2.55in,height=2in]{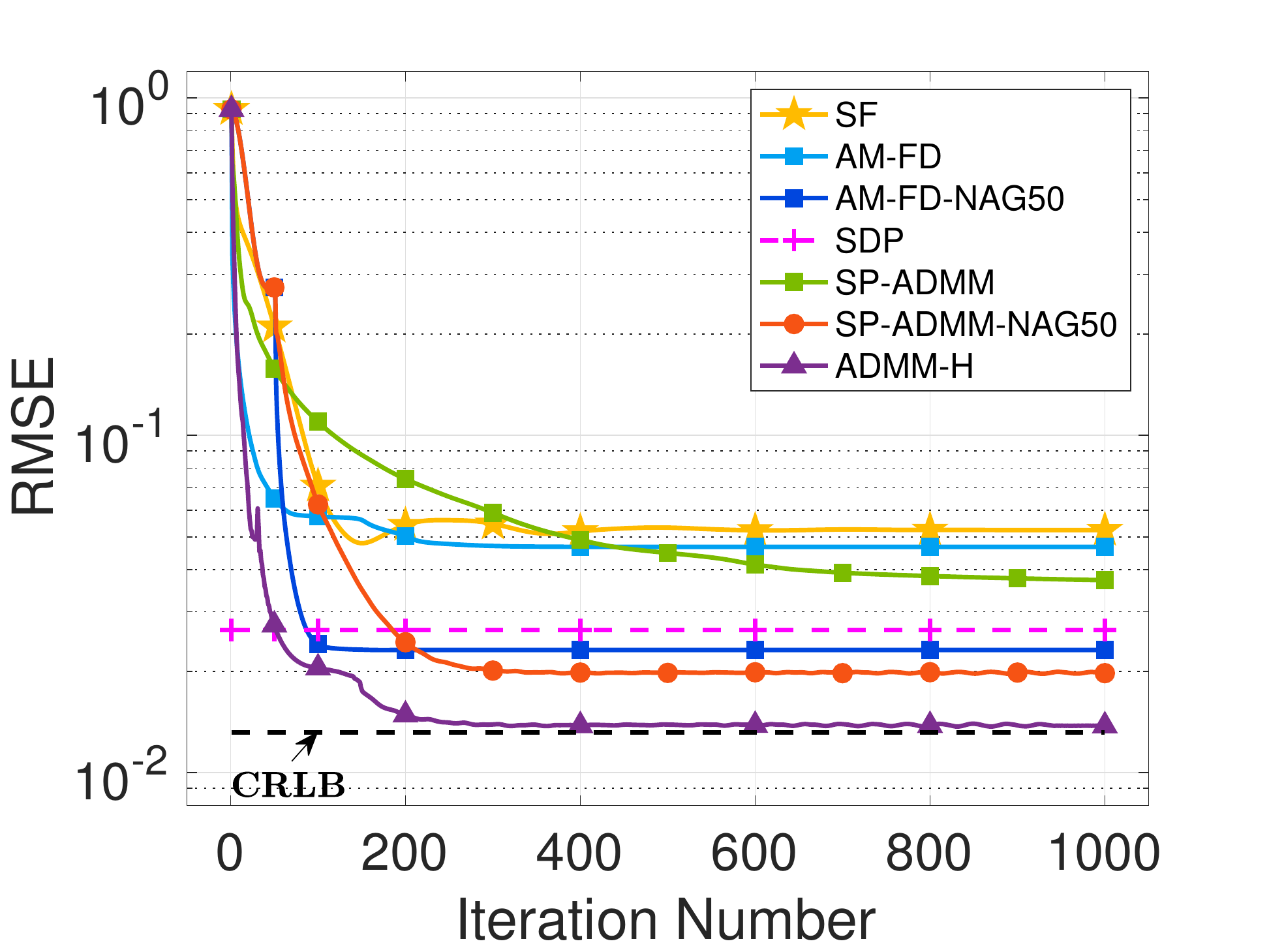}\label{m_20a}}
	\hfil\hspace{-0.65cm}
	\subfloat[$m=20$]{
		\includegraphics[width=2.55in,height=2in]{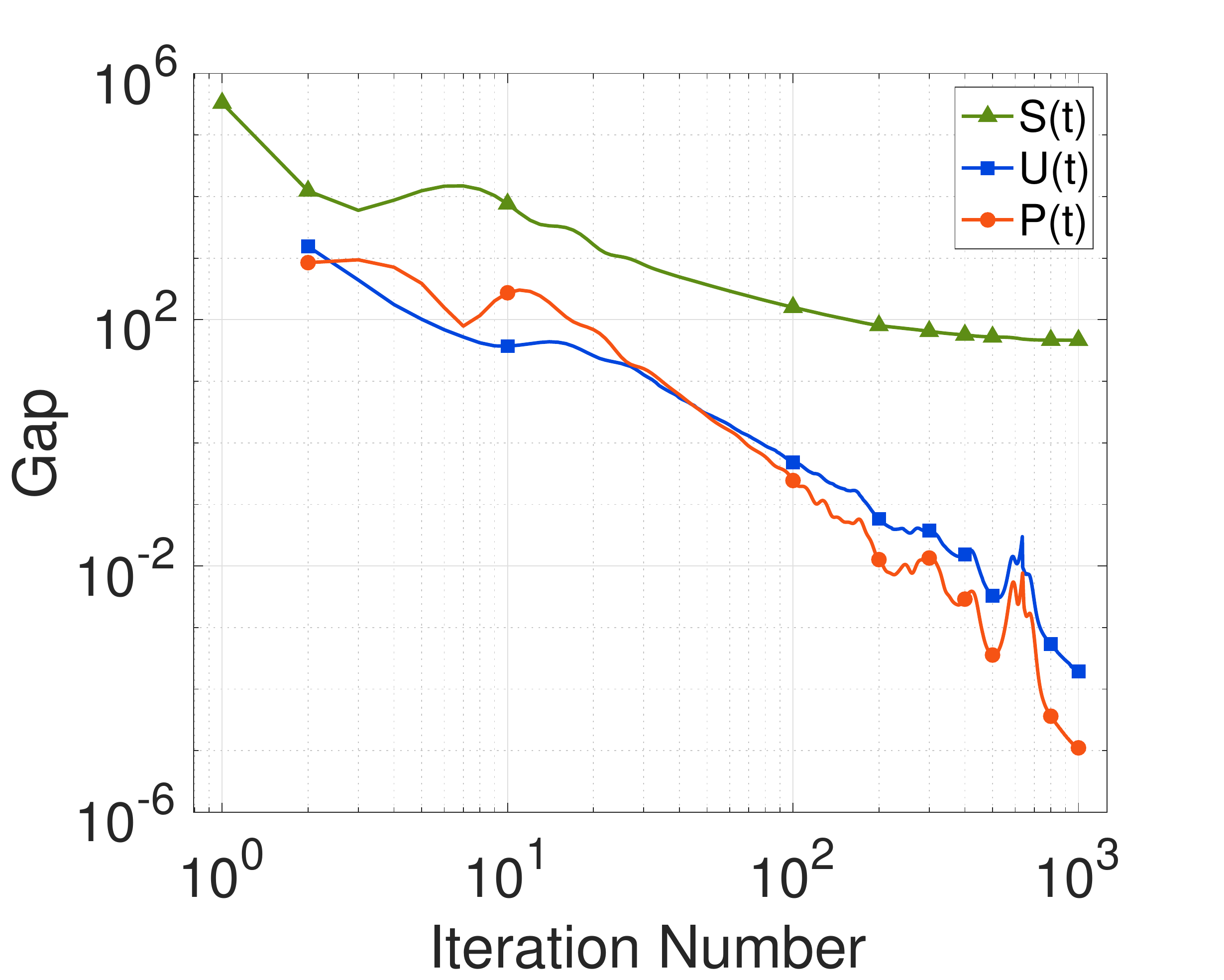}\label{m_20b}}
	\,
	\subfloat[$m=30$]{
		\includegraphics[width=2.55in,height=2in]{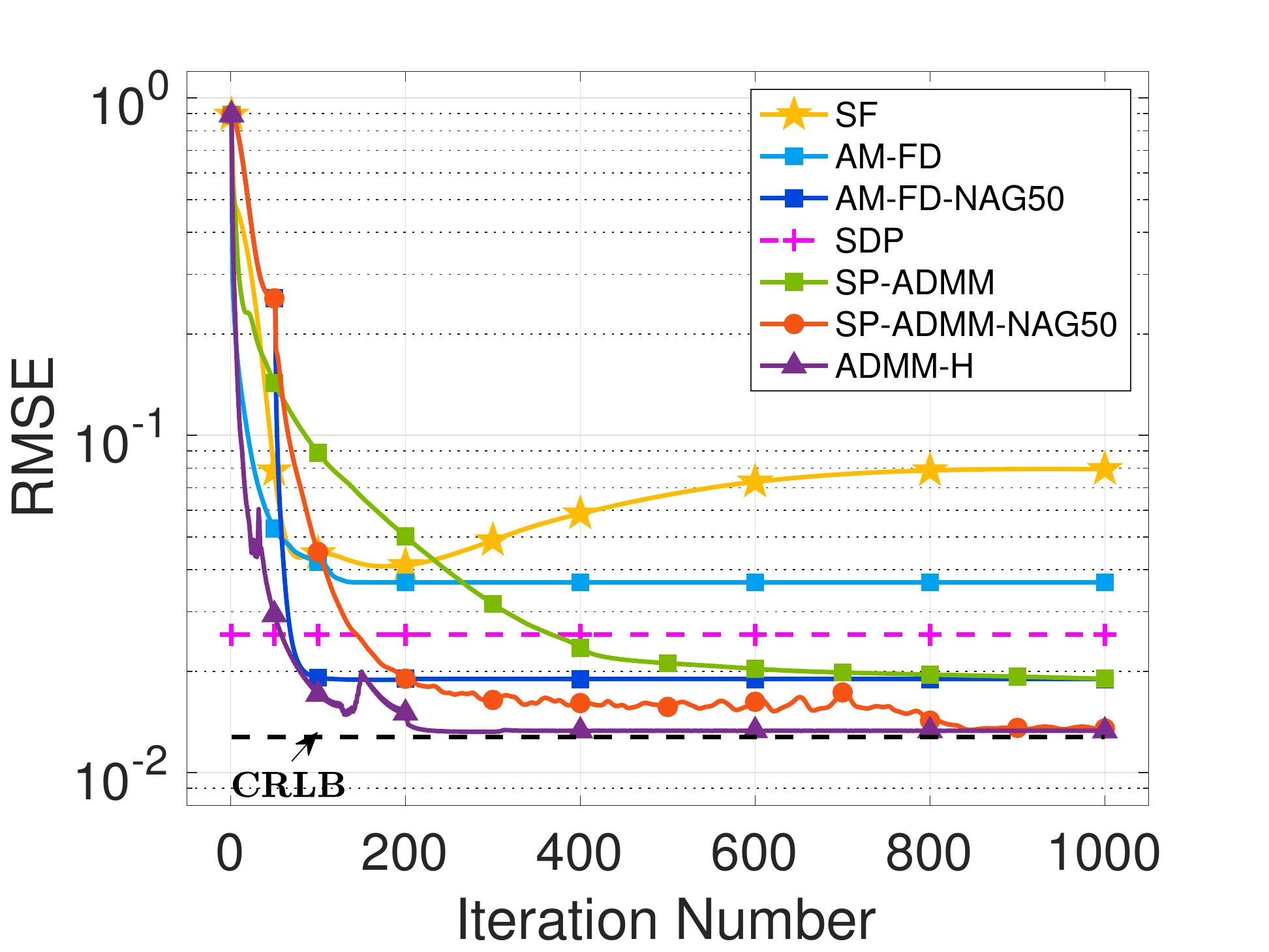}\label{m_30a}}
	\hfil\hspace{-0.65cm}
	\subfloat[$m=30$]{
		\includegraphics[width=2.55in,height=2in]{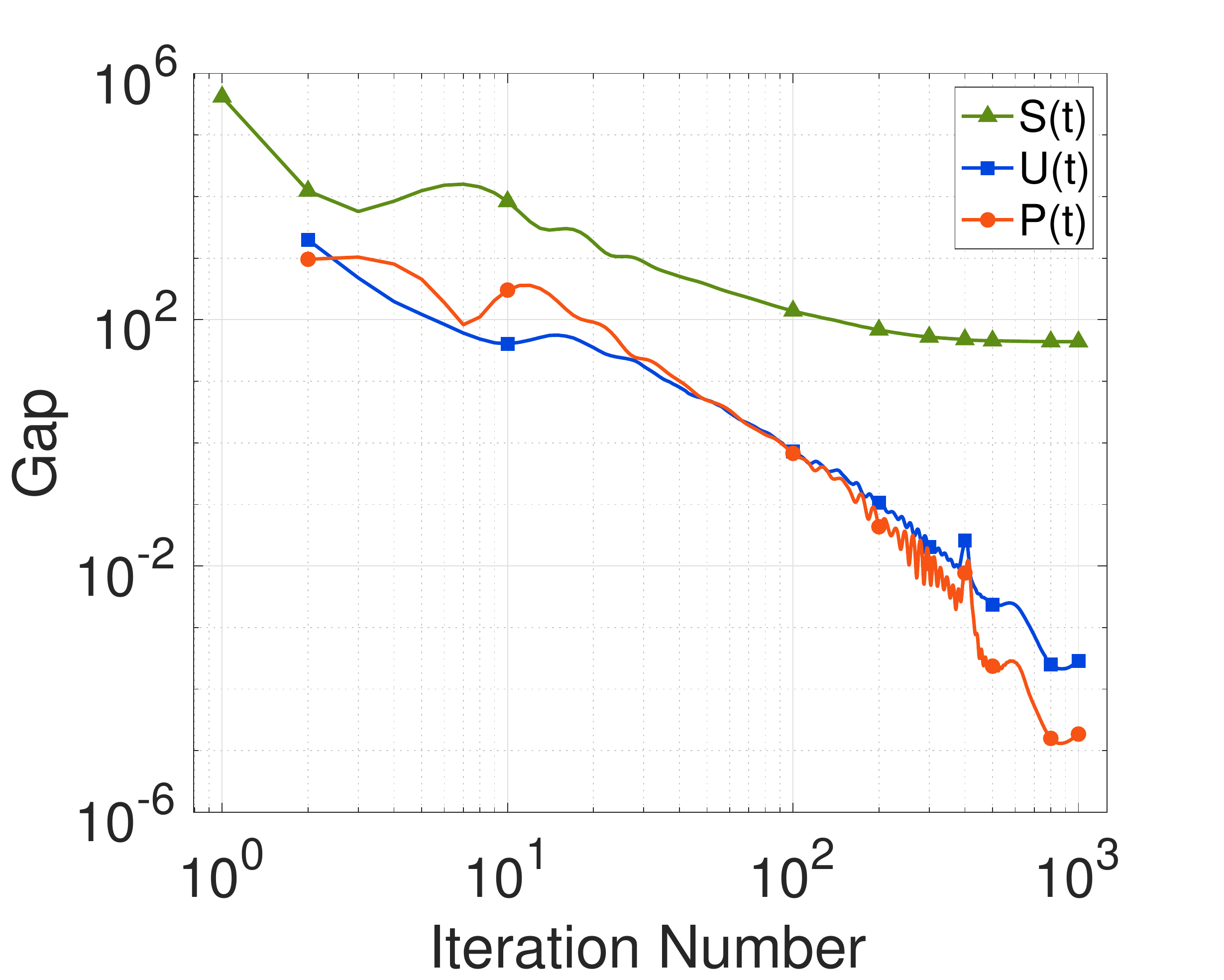}\label{m_30b}}
	\caption{Convergence performance under a different number of anchors with AWGN for the $N-m=490$ nodes network. RMSE value (left column); primal feasibility and stationarity gap (right column).}
	\label{change_m}
\end{figure}

\begin{table}[!t]
	\renewcommand{\arraystretch}{1.3}
	\center
	\begin{threeparttable}[b]
		\caption{Comparisons of Running Time}
		\label{synthesized running time}
		\setlength{\tabcolsep}{5mm}{
			\begin{tabular}[l]{c|c|c|c}
				\hline
				&&\multicolumn{2}{c}{Run time (seconds)}\\\cline{3-4}
				Algorithm&RMSE&Parallelized&Sequential(Per Step)\\\hline\hline
				\multicolumn{4}{c}{Random ($N=495,m=5,\text{step}=2000,\text{CRLB}\approx 0.016$)}\\\hline
				SF&4.64e-01&0.0402&1.97e02(0.0984)\\
				AM-FD&1.86e-01&-\tnote{$^b$}& 1.23e02\,(0.0613)\\
				AM-FD-NAG50&9.25e-02&-&1.24e02\,(0.0623)\\
				SDP&7.99e-02&-&1.59e02\,(4.9609)\\
				ADMM-H&7.47e-02&0.1412 &6.27e03\,(3.1338)\\
				\textbf{SP-ADMM-NAG50}&7.45e-02&\textbf{0.0256}&1.16e02\,(0.0581)\\
				\textbf{SP-ADMM}&1.38e-01&\textbf{0.0251}&1.11e02\,(0.0554)\\\hline
				\multicolumn{4}{c}{Random ($N=510,m=20,\text{step}=1000,\text{CRLB}\approx0.013$)}\\\hline
				SF&9.40e-02&0.04971&1.39e03\,(1.3411)\\
				AM-FD&4.67e-02&-&7.20e01\,(0.07202)\\
				AM-FD-NAG50&2.31e-02&-&7.66e01\,(0.0768)\\
				SDP&2.65e-02&-&3.24e02\,(10.4481)\\
				ADMM-H&1.38e-02&1.31&3.56e03\,(3.5575)\\
				\textbf{SP-ADMM-NAG50}&1.97e-02&\textbf{0.03232}&8.10e01\,(0.0881)\\
				\textbf{SP-ADMM}&3.72e-02&\textbf{0.0388}&8.29e01\,(0.0811)\\\hline
				\multicolumn{4}{c}{Random ($N=520,m=30,\text{step}=1000,\text{CRLB}\approx0.013$)}\\\hline
				SF&7.94e-02&0.0411&1.81e03\,(1.8142)\\
				AM-FD&3.67e-02&-&7.51e01(0.0752)\\
				AM-FD-NAG50&1.90e-02&-&8.20e01\,(0.0821)\\
				SDP&2.56e-02&-&5.57e02\,(18.5801)\\
				ADMM-H&1.33e-02&1.41&3.51e03\,(3.5150)\\
				\textbf{SP-ADMM-NAG50}&1.35e-02&\textbf{0.0343}&8.54e01\,(0.0854)\\
				\textbf{SP-ADMM}&1.90e-02&\textbf{0.0351}&8.51e01\,(0.0811)\\\hline
		\end{tabular}}
		\begin{tablenotes}
			\item [$b$] These symbols ``-" in the table mean that the methods cannot be implemented in parallel.
		\end{tablenotes}
	\end{threeparttable}
\end{table}

\noindent
The parameters for the random networks and ADMM-H method are summarized in Table~\ref{synthesized parameter}. The three networks of the experiment are formed by randomly removing 5 anchors, randomly generating 10 anchors, and randomly generating 20 anchors on the benchmark network ($N=500,m=10$). Fig. \ref{change_m} explores the influence of the number of anchor nodes on positioning accuracy. In all three networks, the parameters of our algorithm were set to $\rho=c=0.11,\vu^0=\vzero,$ and $\vz^0$ from a distribution  $\textbf{Unif}\left(-1,1\right)^{n\left(4\mid\mathcal{E}\mid+N\right)}$. As expected, the localization accuracy of most methods improves as the number of anchors increases. This is due to the anchors know their true locations, thus can provide more accurate position estimates of the neighboring nodes. In addition, we can see that except for the ADMM-H method, the SP-ADMM-NAG works better than other methods and gets closer to CRLB when increasing the number of anchors. A closer inspection shows that for the AM-FD (light blue line), which is also a nonconvex relaxation method, the proposed SP-ADMM (green line) performs better even though it converges not so fast at the beginning. AM-FD-NAG50 (dark blue line) and SP-ADMM-NAG50 (red line) also have such performance. The running time is shown in Table~\ref{synthesized running time}, and we find that the SP-ADMM is most computationally time efficient due to the parallel implementation. The results here are similar to those in Table~\ref{benchmark running time} and consistent with Table~\ref{comparisons of different alg}.

\subsubsection{\textbf{Influence of the Average Number of Neighbors ($D_{\text{avg}}$)}}
\begin{figure}[!t]
	\centering
	\subfloat[D$_{\text{avg}}=7.30$]{
		\includegraphics[scale=0.5]{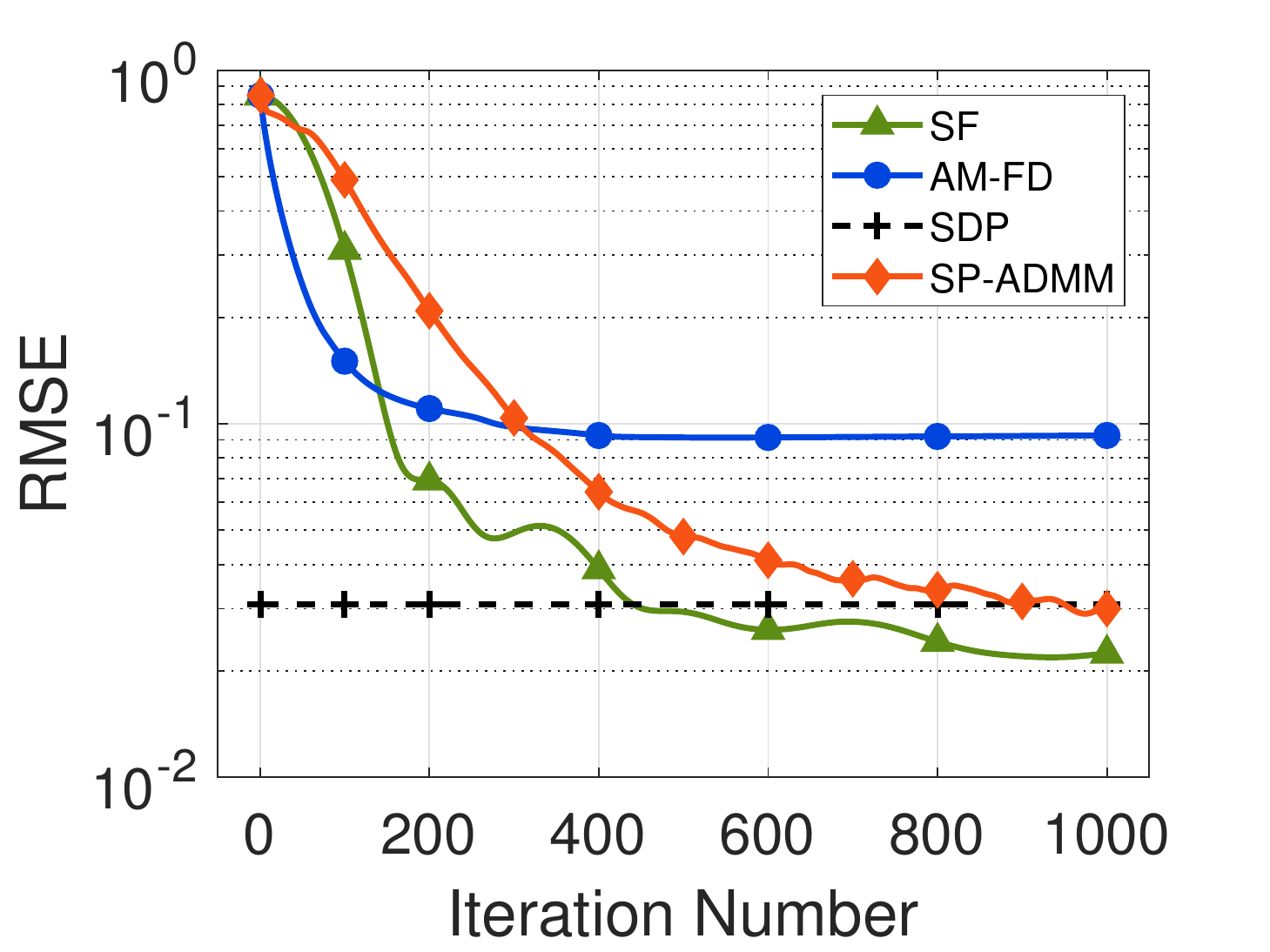}\label{D_avg_a}}\hfil\hspace{-0.2in}
	\subfloat[D$_{\text{avg}}=9.87$]{
		\includegraphics[scale=0.5]{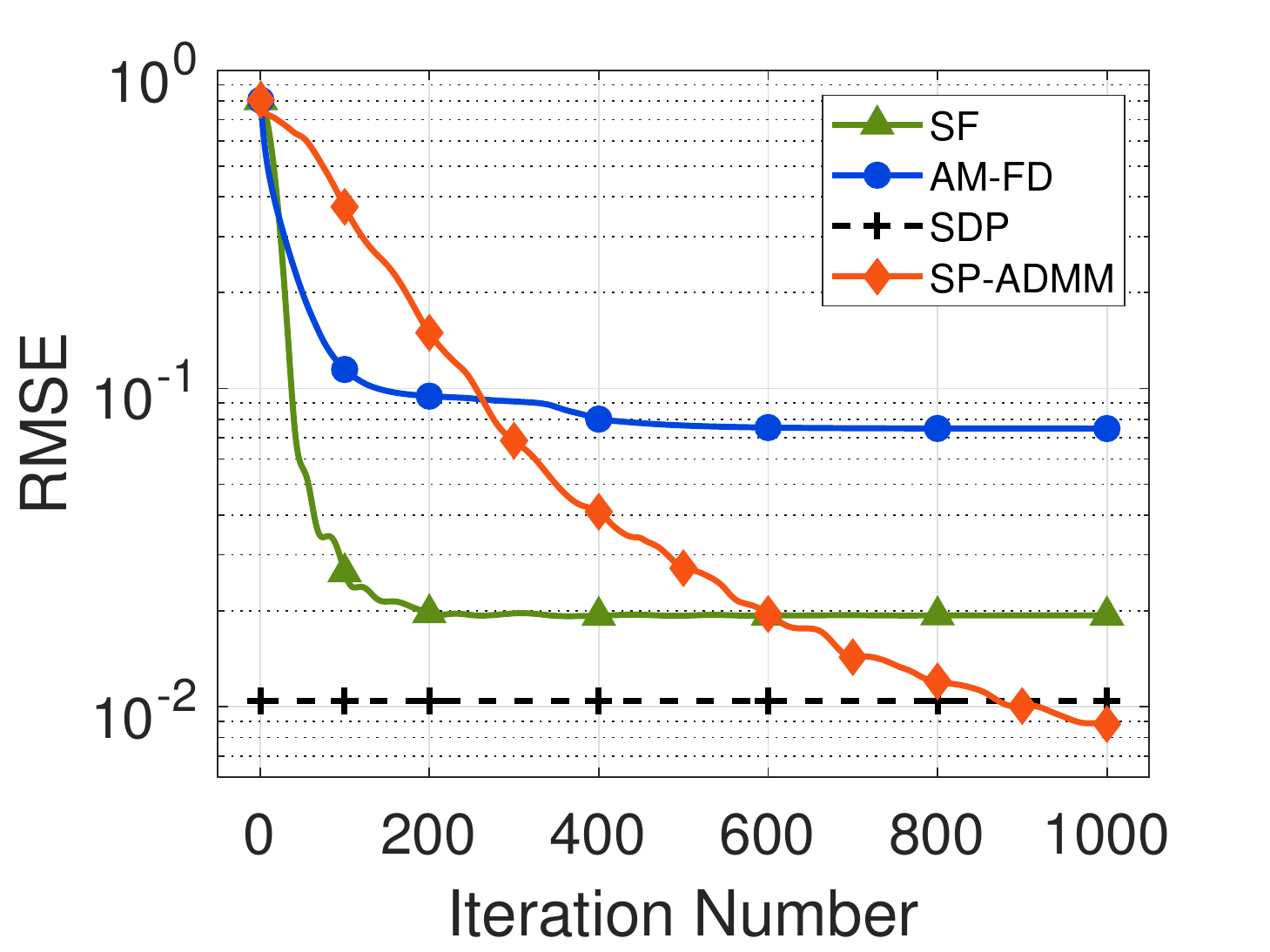}\label{D_avg_b}}\hfil\hspace{-0.2in}
	\,
	\subfloat[D$_{\text{avg}}=12.57$]{
		\includegraphics[scale=0.5]{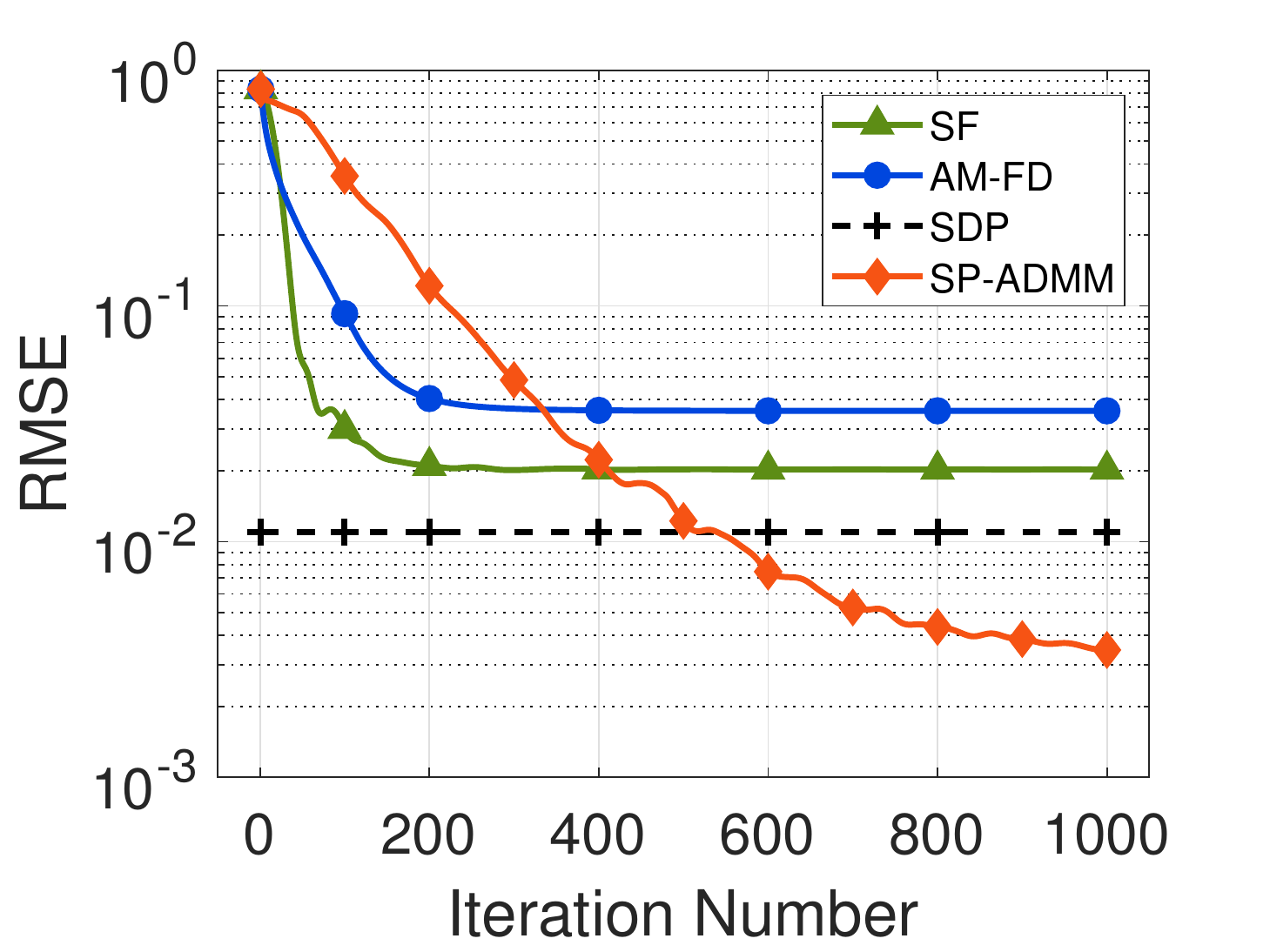}\label{D_avg_c}}\hfil\hspace{-0.2in}
	\subfloat[D$_{\text{avg}}=14.43$]{
		\includegraphics[scale=0.5]{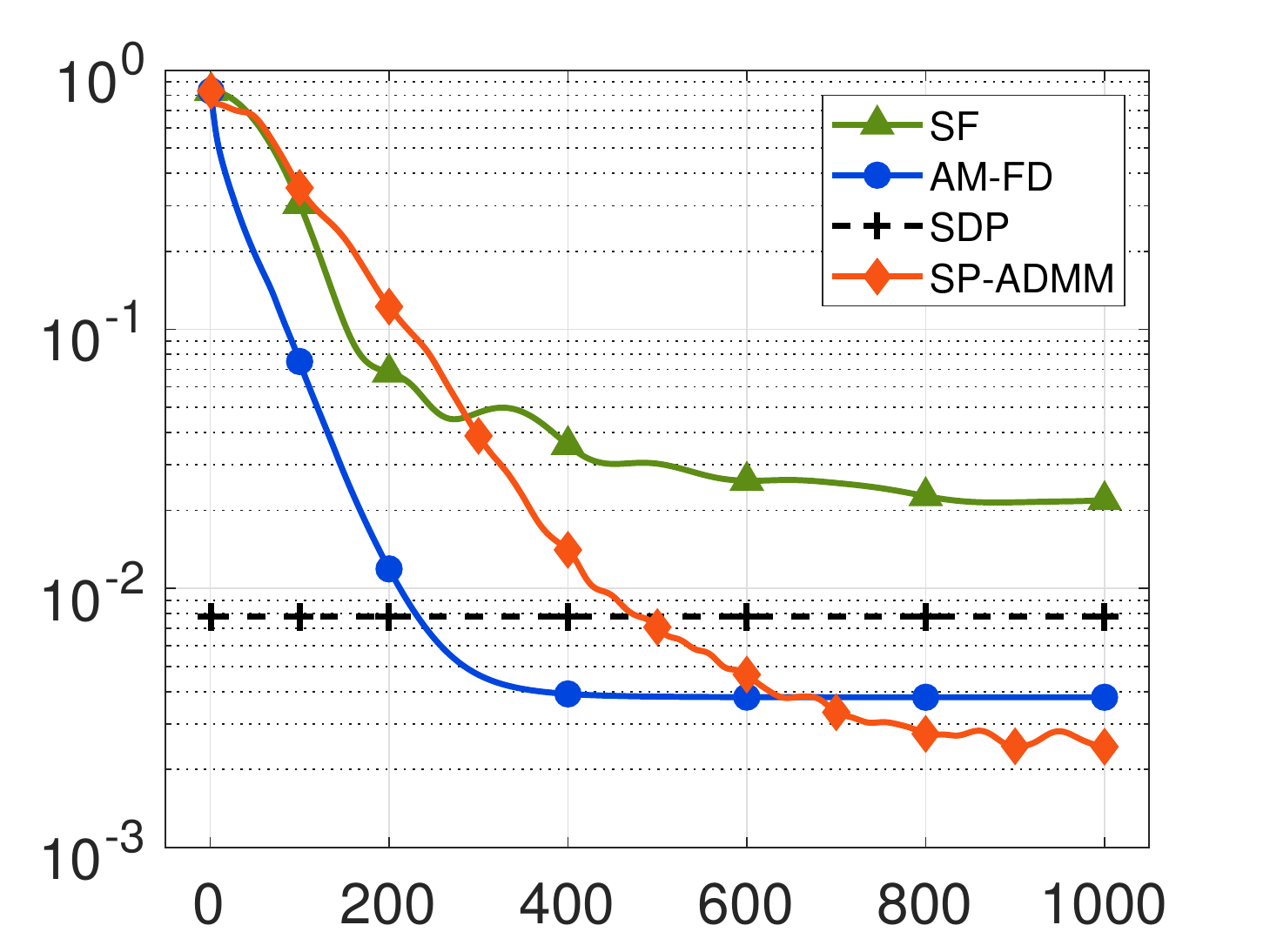}\label{D_avg_d}}\hfil\hspace{-0.2in}	
	\caption{Convergence performance under different average number of neighboring nodes (D$_{\text{avg}}$). Here, we adopt $\sigma_{\text{add}}=0.02$; $N=108,m=8$.}
	\label{change_Davg}
\end{figure}
\begin{figure}[!t]
	\centering
	\subfloat[SDP]{
		\includegraphics[scale=0.5]{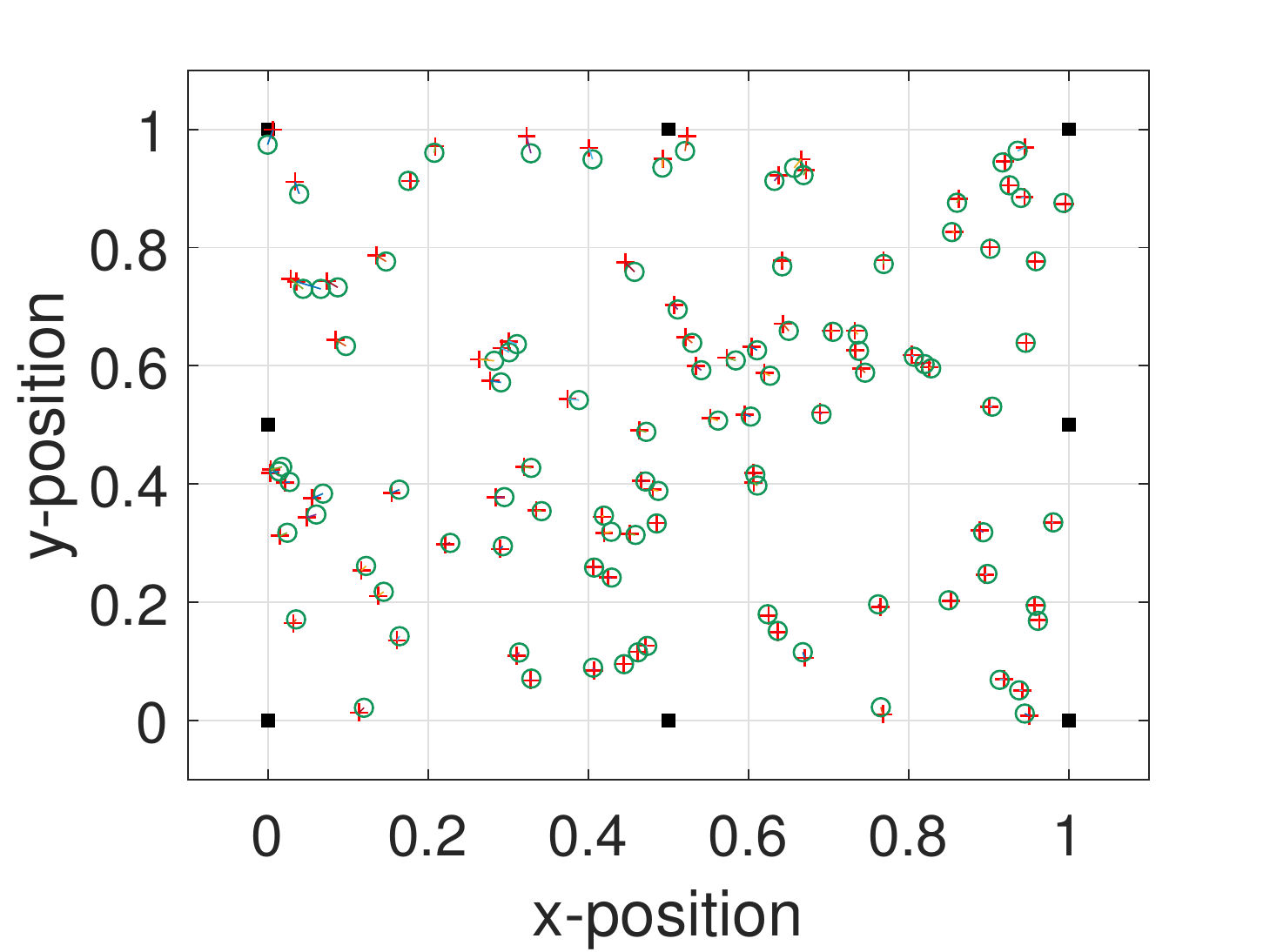}\label{sdp}
	}\hfil\hspace{-0.21in}
	\subfloat[SF]{
		\includegraphics[scale=0.5]{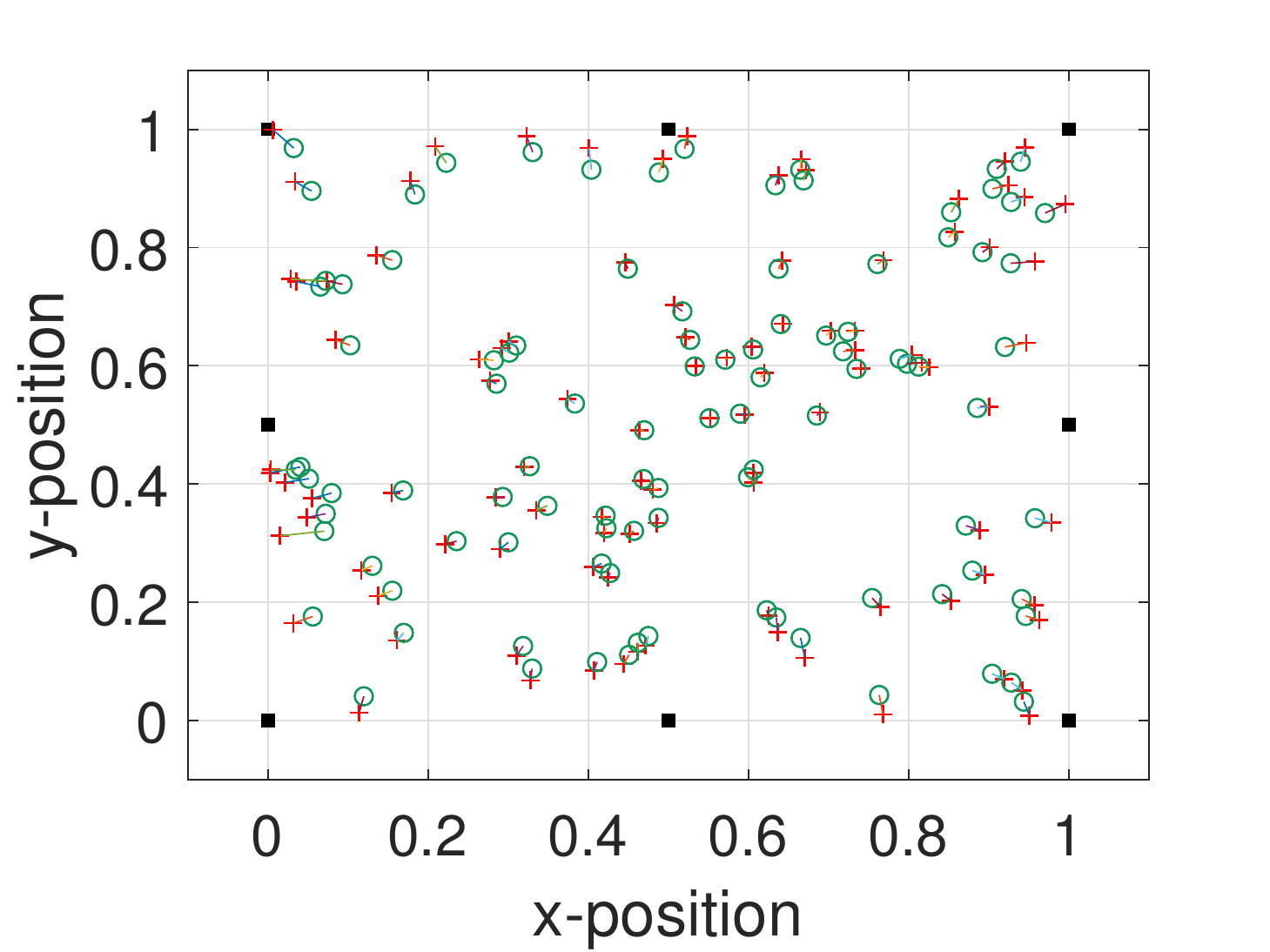}}\label{sf}
	\,
	\subfloat[AM-FD]{
		\includegraphics[scale=0.5]{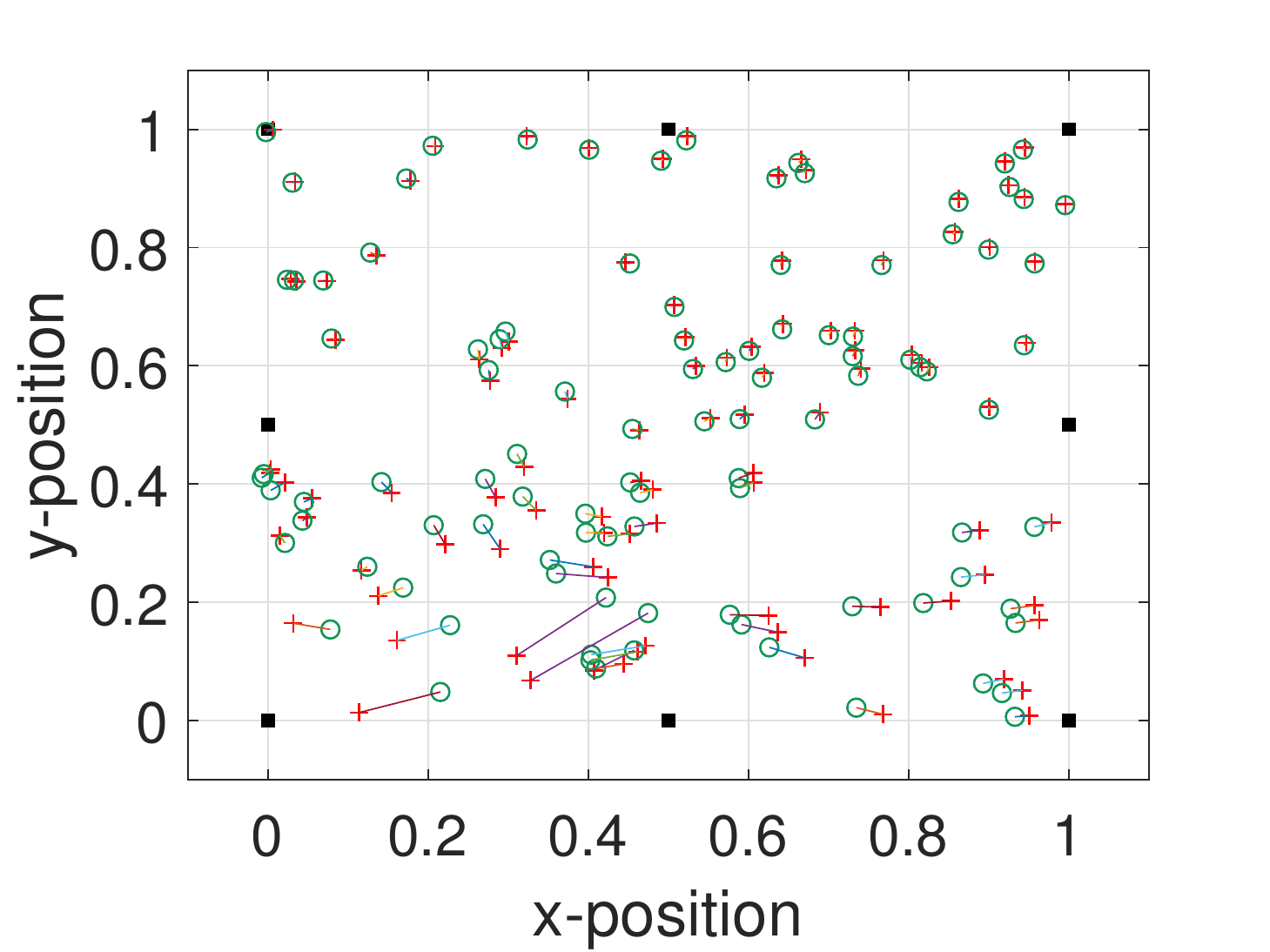}\label{am-fd}
	}\hfil\hspace{-0.21in}
	\subfloat[SP-ADMM]{
		\includegraphics[scale=0.5]{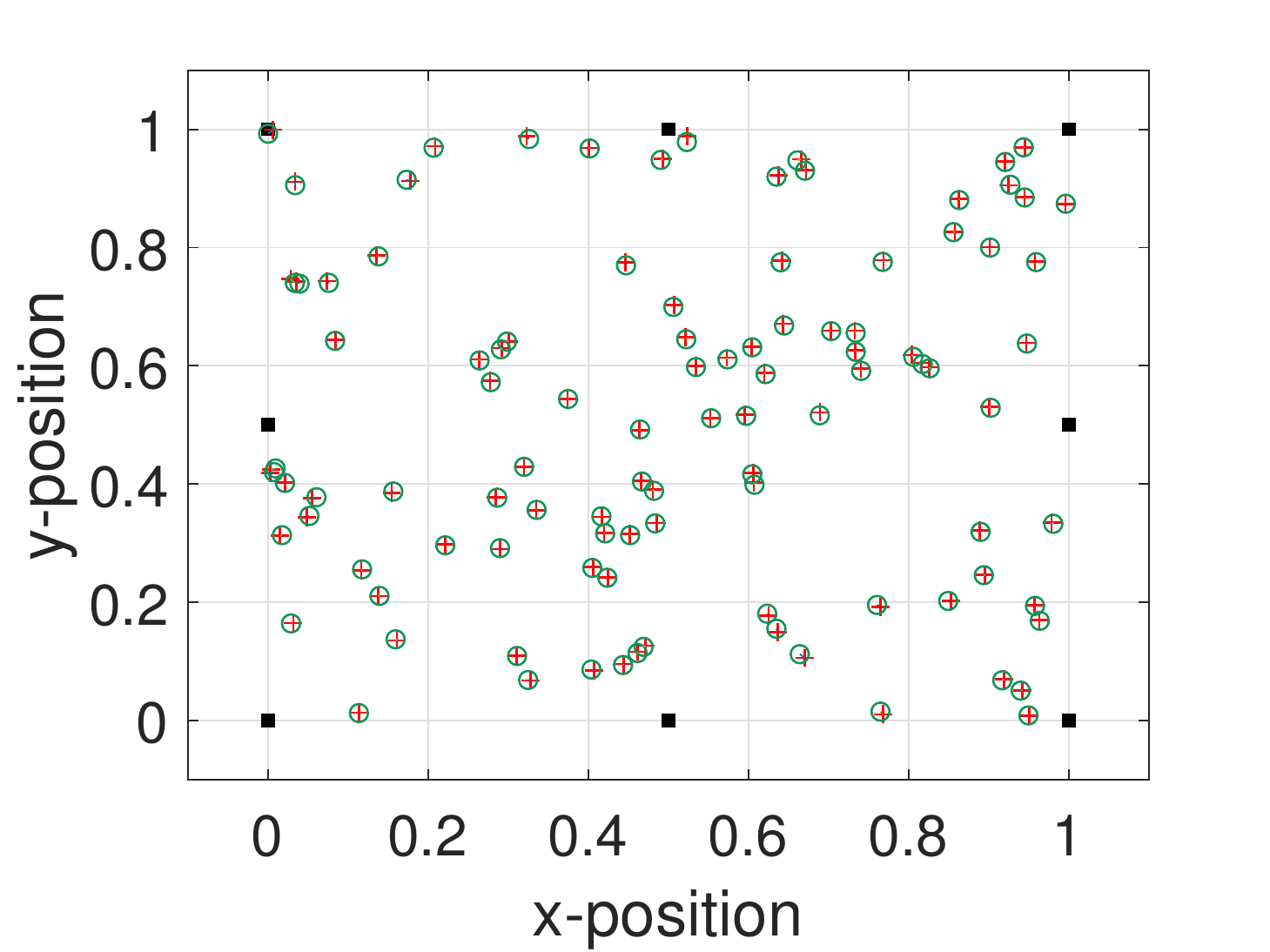}\label{sp-admm}}
	\caption{Position estimates obtained by the SDP, SF, AM-FD, and the proposed SP-ADMM algorithm in 2-D plane. Here, we adopt $\sigma_{\text{add}}=0.02,N=108,m=8$,  D$_{\text{avg}}=12.57$ (C$_{\text{range}}$ = 0.23). Here, the anchors are marked by $\blacksquare$, the true agent positions by ${\red +}$ and the estimated ones by ${\green \circ}$.}
	\label{position estimates}
\end{figure}

\noindent
Fig. \ref{change_Davg} investigates the influence of the average number of neighbors on the performance of four algorithms: the proposed SP-ADMM algorithm, SF, SDP, and AM-FD. The network $(N=108,m=8)$ randomly placed in the 2D $[0,1]\times[0,1]$ area. and we set the parameters $\rho=c=0.0265,\,\vz^0=\vzero$, and $\vu^0=0.5\cdot\vone$.  As shown in Fig. \ref{change_Davg}, the RMSE values of most algorithms decrease as the average number of neighbors increases. The reason may be that each node communicates with more neighbors at each iteration, which leads to the position estimated by different nodes for the same node reaching consistency faster. Notably, the SP-ADMM algorithm consistently achieves the lowest RMSE in all scenarios. Fig. \ref{position estimates} visualizes the position estimates obtained in Fig. \ref{change_Davg}\subref{D_avg_c}. It clearly shows that the positions estimated by our approach are very close to the true agents' locations.

\subsubsection{\textbf{Influence of the Measurement Noise Variance ($\sigma_{i,j}$)}}

\begin{figure}[!t]
	\centering
	\subfloat[$\sigma_{\text{add}}=0.01$]{
		\includegraphics[scale=0.5]{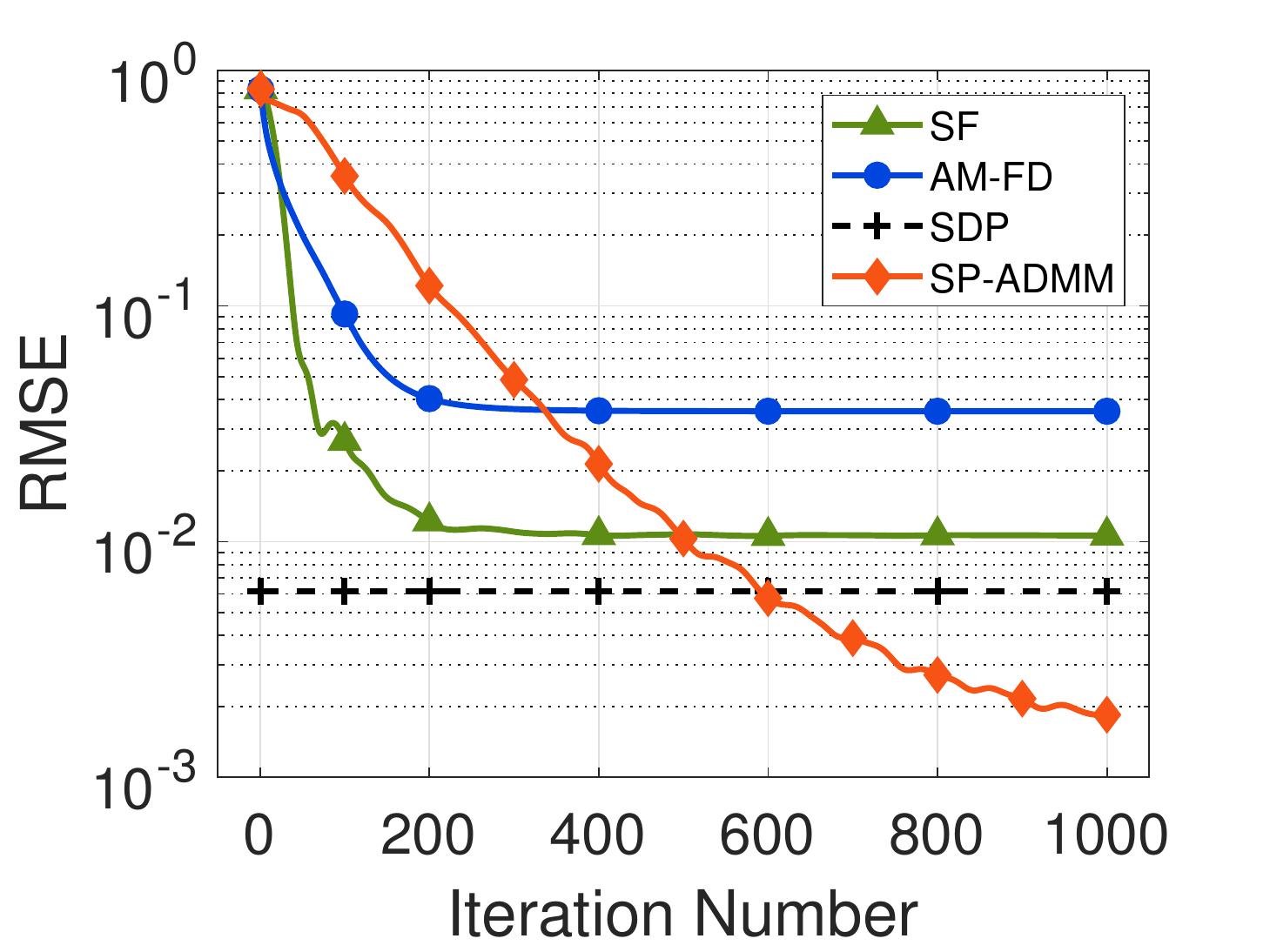}\label{sigma_a}
	}\hfil\hspace{-0.21in}
	\subfloat[$\sigma_{\text{add}}=0.03$]{
		\includegraphics[scale=0.5]{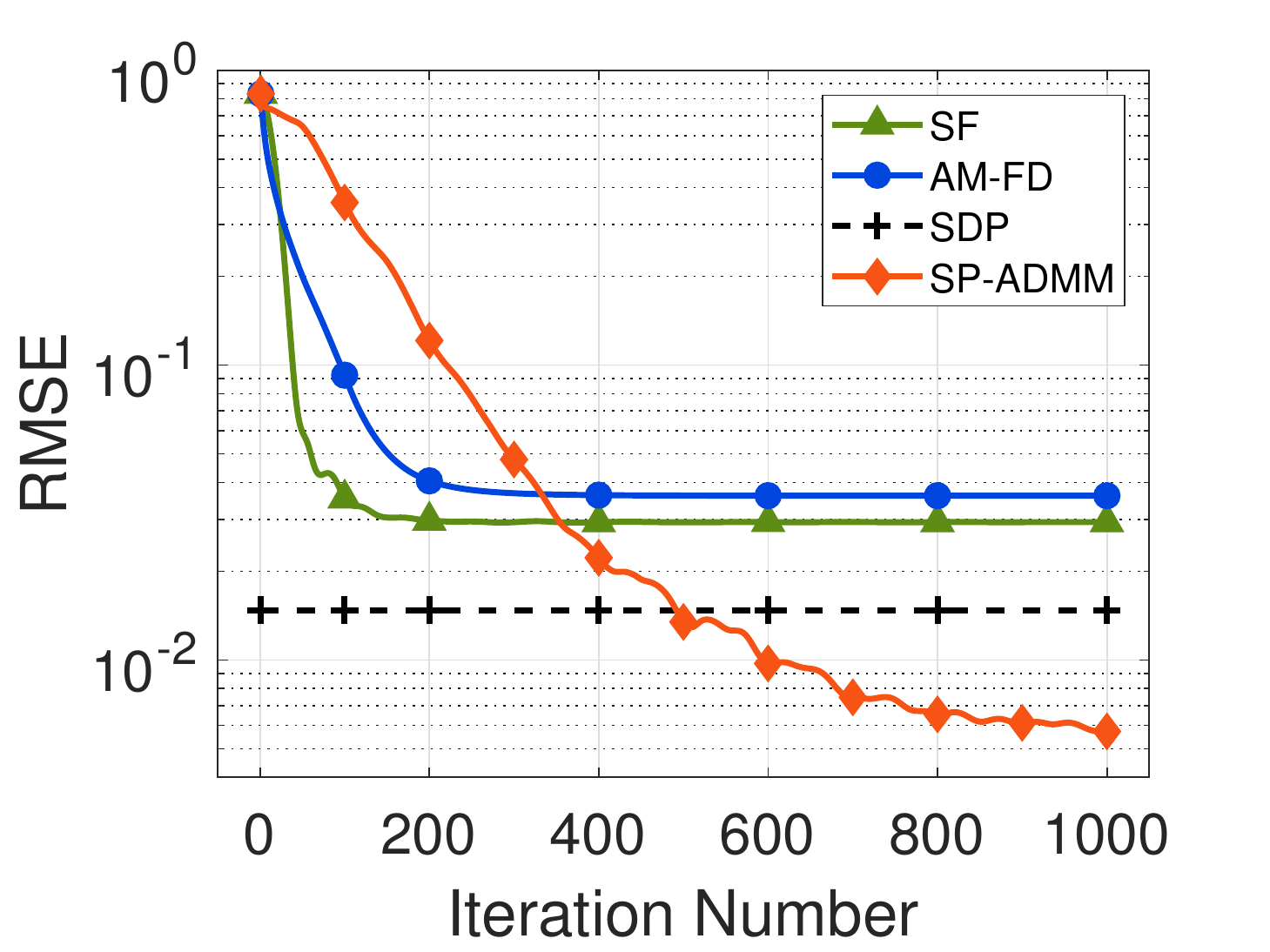}\label{sigma_b}}
	\,
	\subfloat[$\sigma_{\text{add}}=0.04$]{
		\includegraphics[scale=0.5]{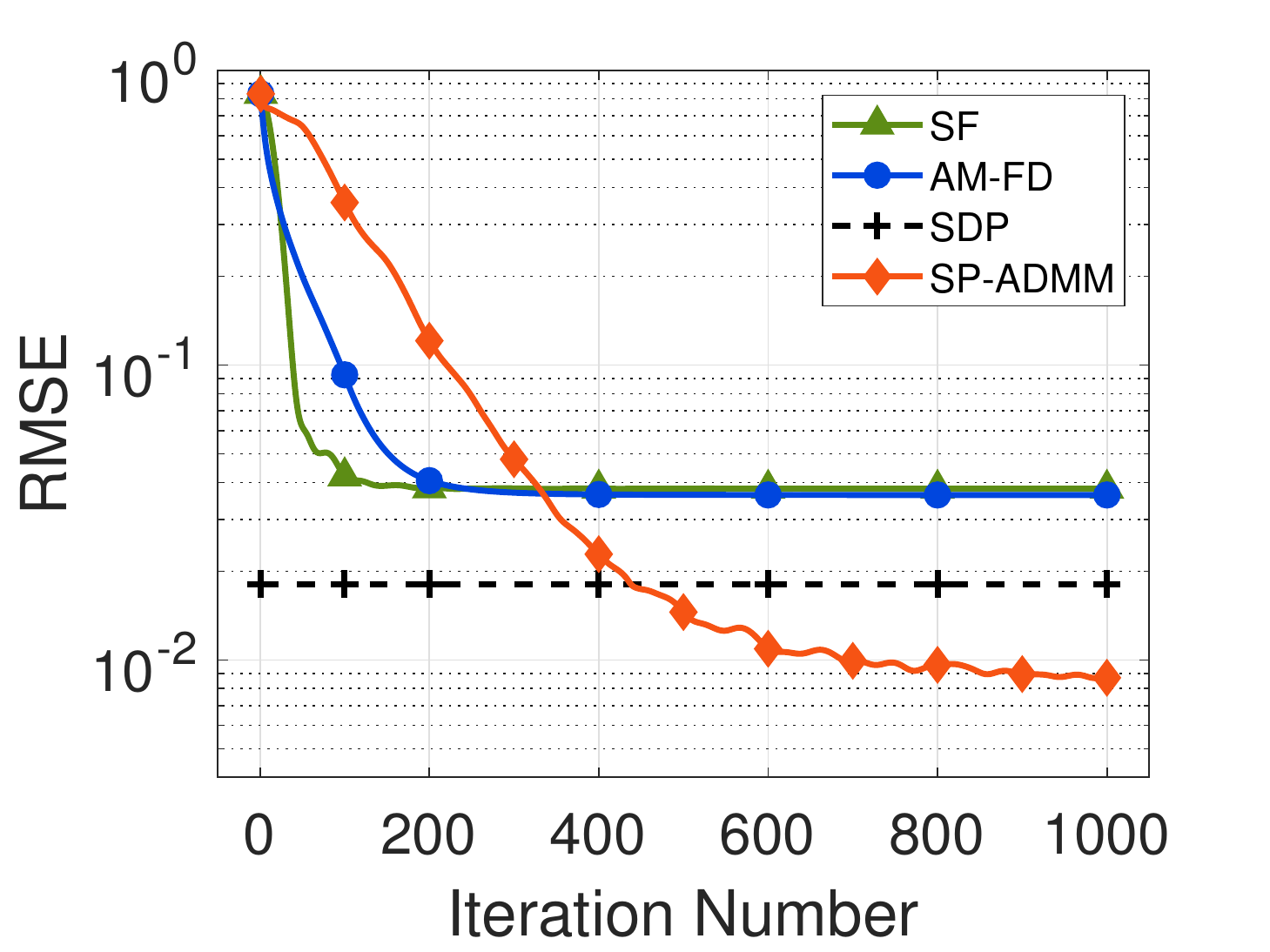}\label{sigma_c}
	}\hfil\hspace{-0.21in}
	\subfloat[$\sigma_{\text{add}}=0.05$]{
		\includegraphics[scale=0.5]{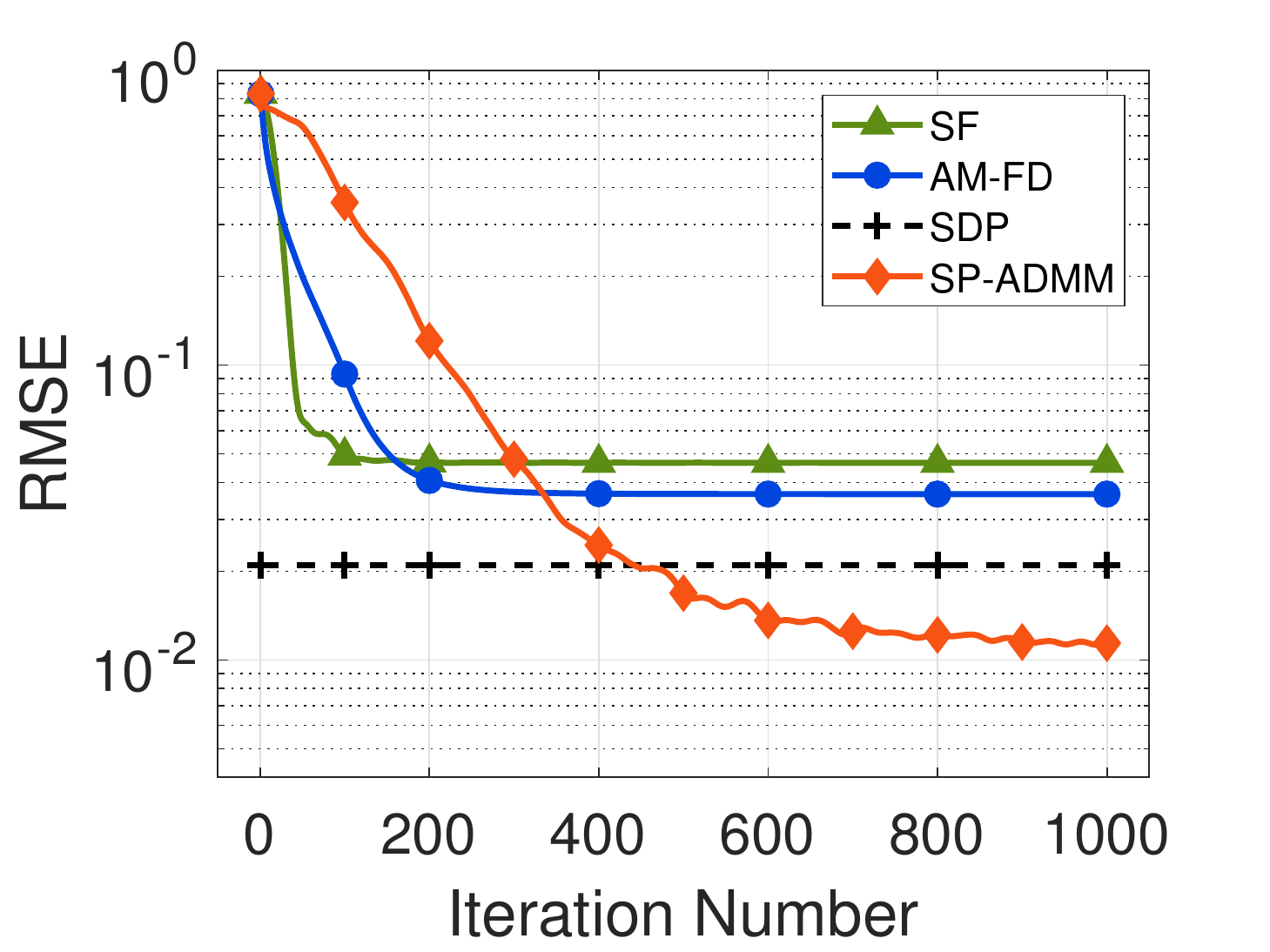}\label{sigma_d}}
	\caption{Convergence performance under different measurement noise variances. Here, we let $N=108,m=8$; communication range is set to 0.23.}
	\label{change_sigma}
\end{figure}
\begin{figure}[!t]
	\subfloat[synthesized network]{\includegraphics[scale=0.4]{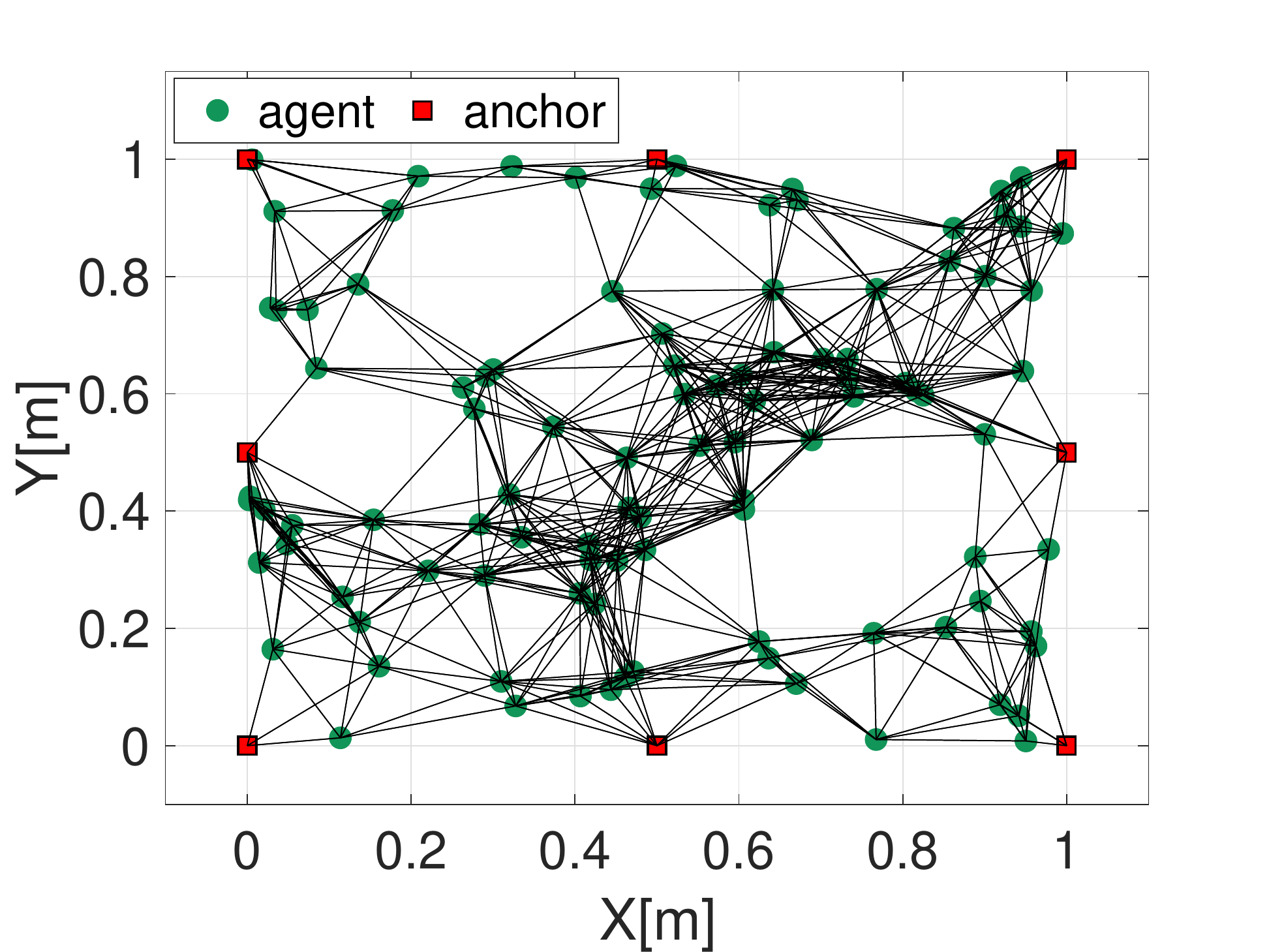}.pdf}
	\subfloat[performance with noise]{\includegraphics[scale=0.45]{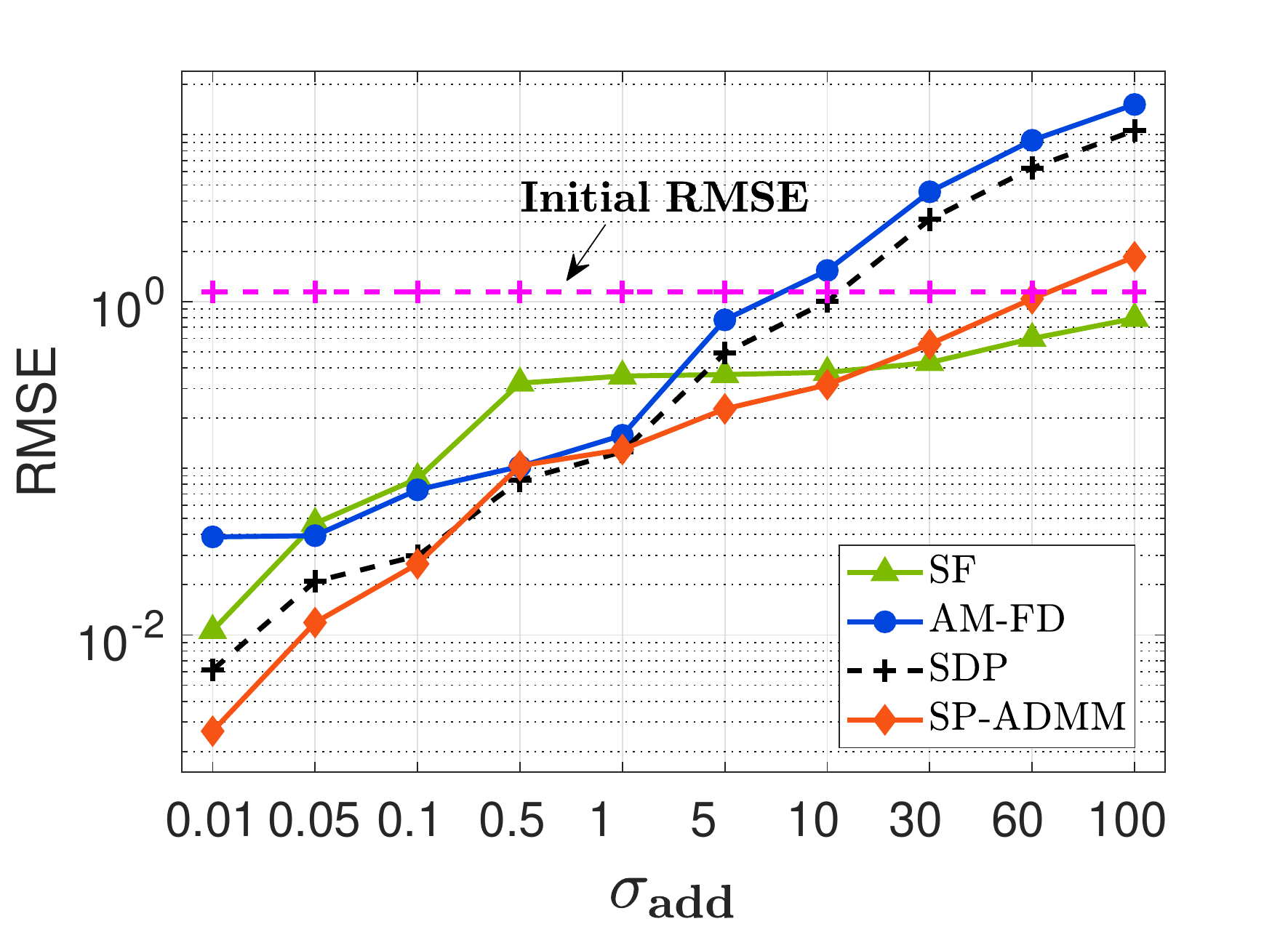}}
	\caption{Influence of range-dependent measurement noise $w_{i,j}$ on positioning accuracy across various methods. $w_{i,j}\sim \text{Gauss}(0,\sigma_{\text{add}}\|\vp_i-\vp_j\|^{2})$}\label{break2}
\end{figure}

\noindent
To test the robustness of the proposed algorithm against different measurement noise profiles, the localization scenarios are further extended to sensor networks with different measurement noise variances $\sigma_{i,j}^2=\sigma_{\text{add}}\|\vp_i-\vp_j\|^{2}$. The corresponding convergence results we obtained are depicted in Fig. \ref{change_sigma}. Both the network and our proposed algorithm used the same parameters as in Fig. \ref{change_Davg}\subref{D_avg_c}. Note that Fig. \ref{change_Davg}\subref{D_avg_c} shows the test results for $\sigma_{\text{add}}=0.02$, while we show some other configurations in Fig. \ref{change_sigma}. Here, we have similar observations as in Fig. \ref{change_Davg}. Although the performance of these methods demonstrates certain loss with the increase of the noise variance, the localization error of our proposed SP-ADMM algorithm is consistently lower in all scenarios. Moreover, one can see that the non-relaxed methods are relatively insensitive to measurement noise compared with the other convex relaxation methods. Perhaps the reason that the solution of the convex relaxation method is an approximation of the original problem.

We also conducted an experiment to investigate the influence of range-dependent noise on the positioning effectiveness of various methods in extreme cases. To this end, we randomly generated a set of 100 nodes and 8 anchors within a $[0,1]\times[0,1]$ square. The communication distance is set to 0.23, and the network topology is depicted in Fig. \ref{break2}(a). The distance measurements between sensors are subject to Gaussian-distributed noise with zero mean and variance $\sigma_{i,j}^2=\sigma_{\text{add}}\|\vp_i-\vp_j\|^2$. The initial positions are drawn from a uniform distribution $\mathbf{Unif}(-1,1)^{200}$,  and we calculated the corresponding Root-Mean-Squared-Error (RMSE) value as follows:
	\begin{align*}
	\text{RMSE}(0)=\sqrt{\frac{1}{N_{MC}}\frac{1}{100}\sum_{k=1}^{N_{MC}}\sum_{i\in(\mathcal{N}/\mathcal{A})}^{100}(\vp^0_{i}-\vp_i)^T(\vp^0_{i}-\vp_i)},
	\end{align*}
	where $N_{MC}$ denotes the count of independent  Monte Carlo trials, while $\vp^{0,k}_{i}$ symbolizes the initial position vector of the $i$-th node in the $k$-th Monte Carlo trial. For Fig. \ref{break2}, we set $N_{MC}=50$. When the RMSE of the algorithm is larger than the initial RMSE, it means that the algorithm has failed to improve the accuracy of the sensor positions and may have even made it worse. From Fig. \ref{break2}, one can observe that the proposed SP-ADMM failed to provide appropriate position estimates when the variance of the range-noise is bigger than $60$. Additionally, Fig. \ref{break2} shows that the RMSEs of all methods increase as the noise variance becomes larger. Nevertheless,  the proposed SP-ADMM shows the slowest performance degradation, indicating that the proposed method is more robust than AM-FD method and SDP method.

\subsection{Evaluation of the Proposed SP-ADMM Algorithm Versus Penalty Parameters}

\begin{figure}[!t]
	\centering
	\subfloat[fixed $\rho=0.11$]{
		\includegraphics[scale=0.4]{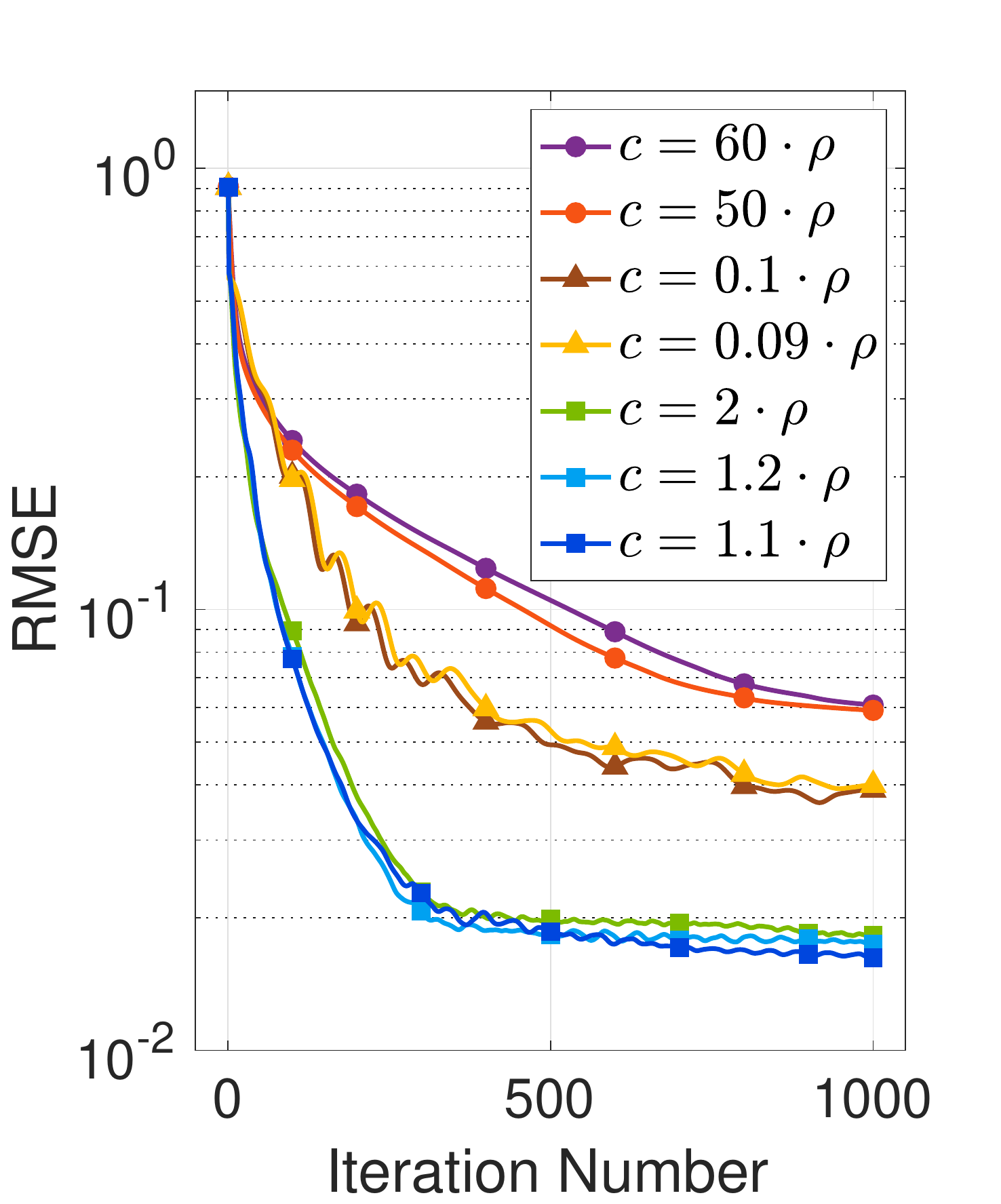}\label{fixed_rho_c}}
	\subfloat[fixed $c=0.11$]{
		\includegraphics[scale=0.4]{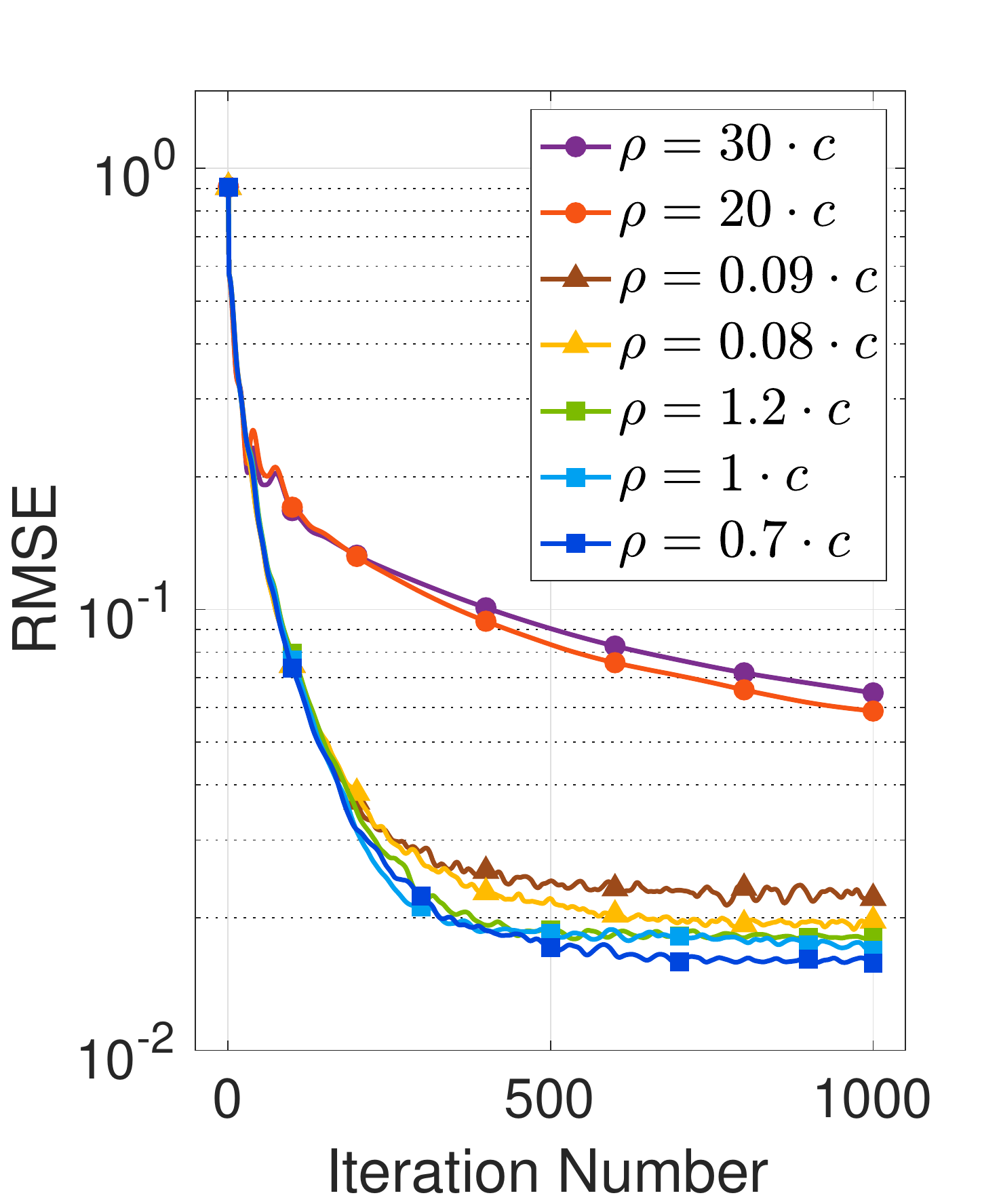}\label{fixed_c_rho}}\caption{RMSE curves  achieved by our proposed methods with different penalty coefficients under the benchmark network $(N=500, m=10)$ with AWGN noise.
	}\label{rho_c}
\end{figure}
\begin{figure}[!t]
	\centering		
	\subfloat[RMSE versus $\vu^0 (\vu^0=u\cdot\vone)$]{
		\includegraphics[scale=0.25]{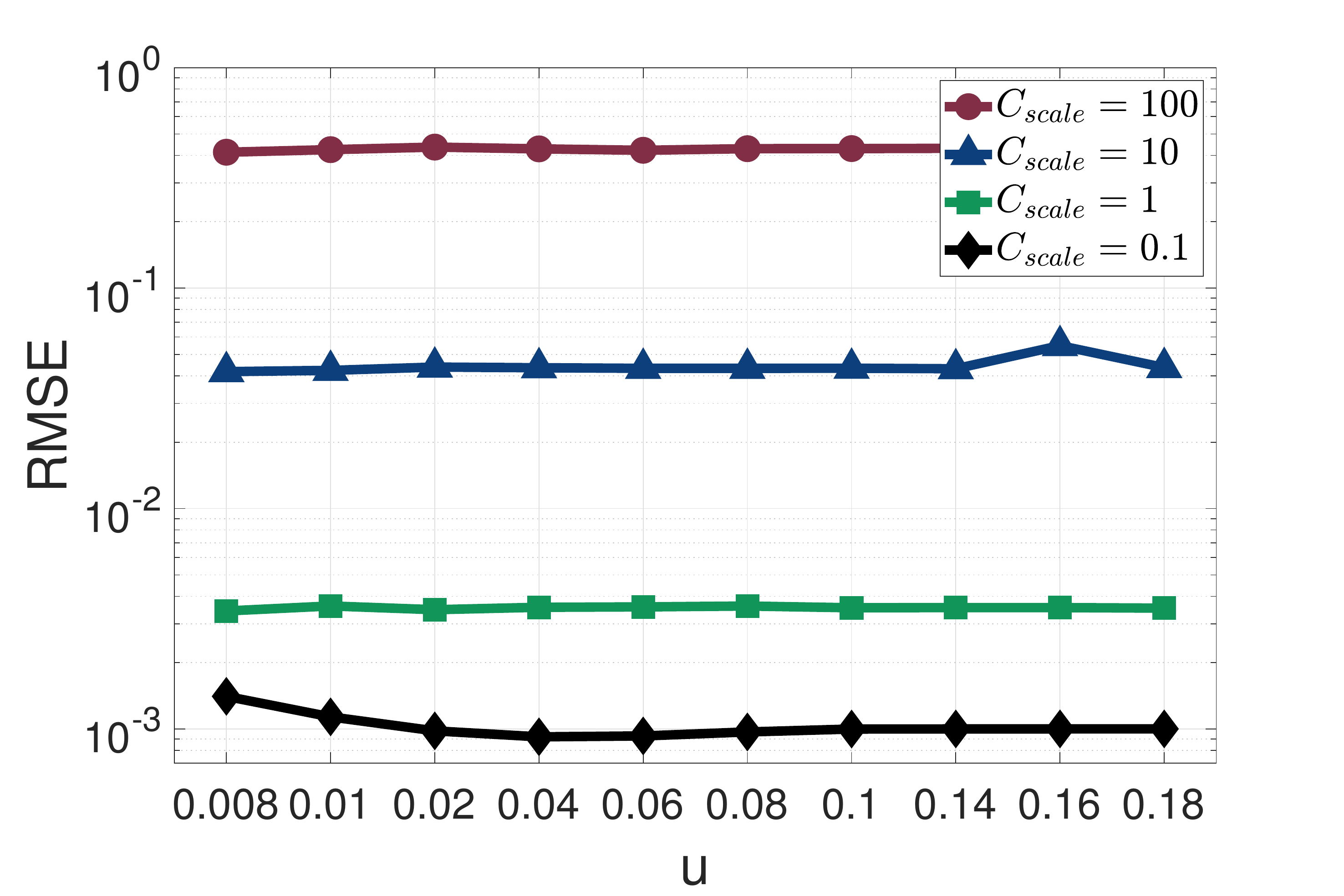}\label{u_ini}}
	\subfloat[RMSE versus $c$]{
		\includegraphics[scale=0.25]{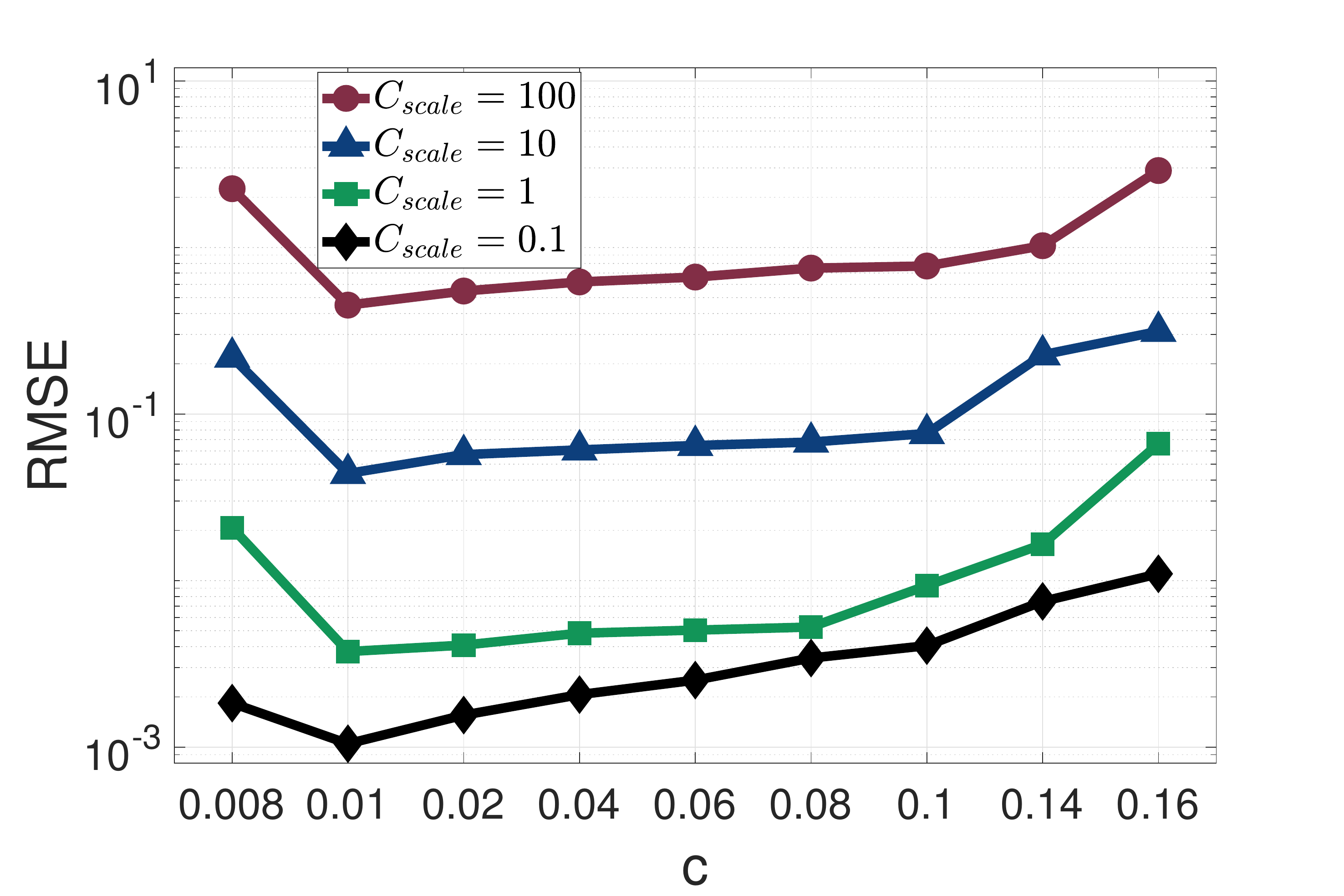}\label{c_scale}}
	\,
	\subfloat[RMSE versus $c$]{
		\includegraphics[scale=0.25]{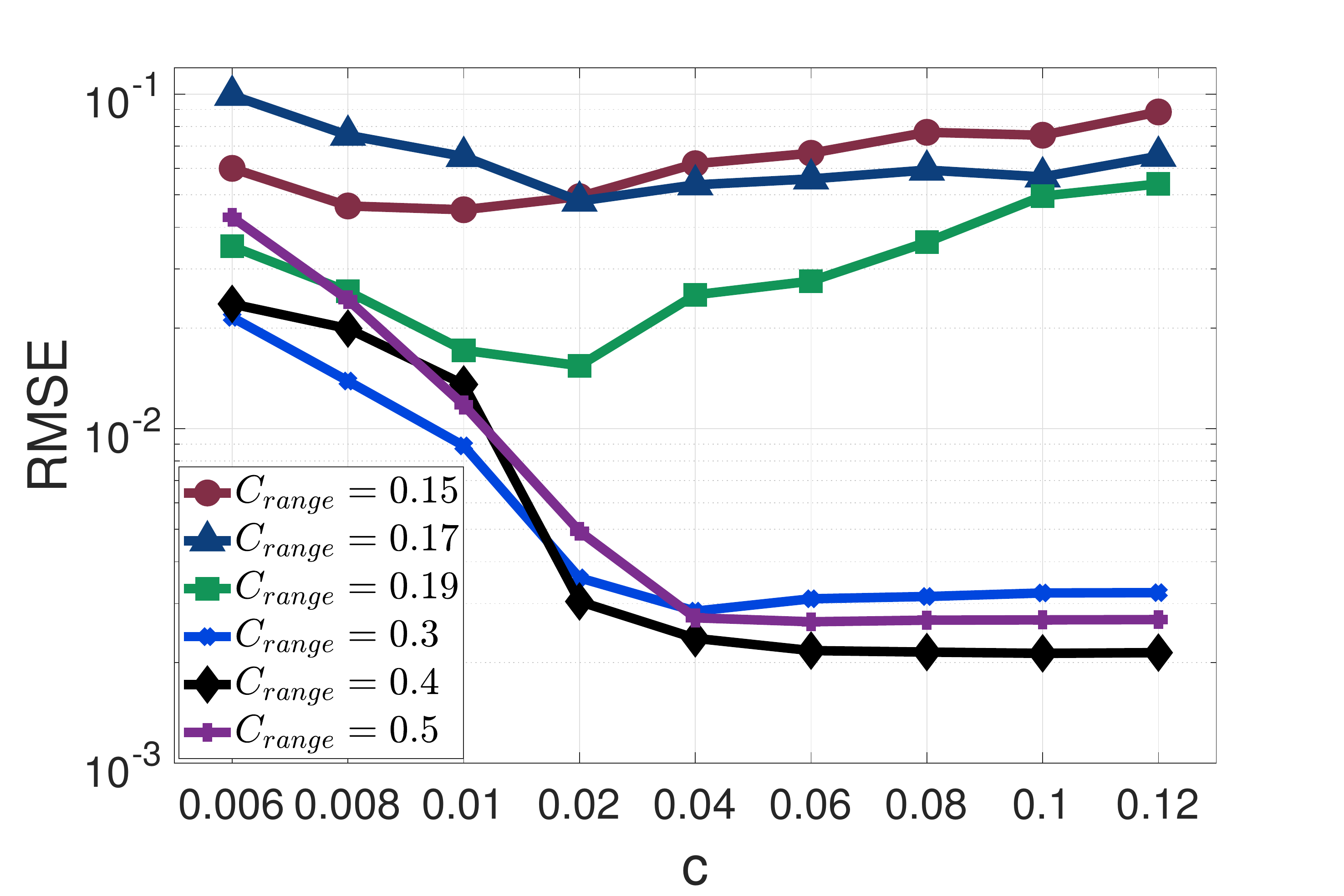}\label{c_range}}
	\subfloat[RMSE versus $N_{\max}$]{
		\includegraphics[scale=0.25]{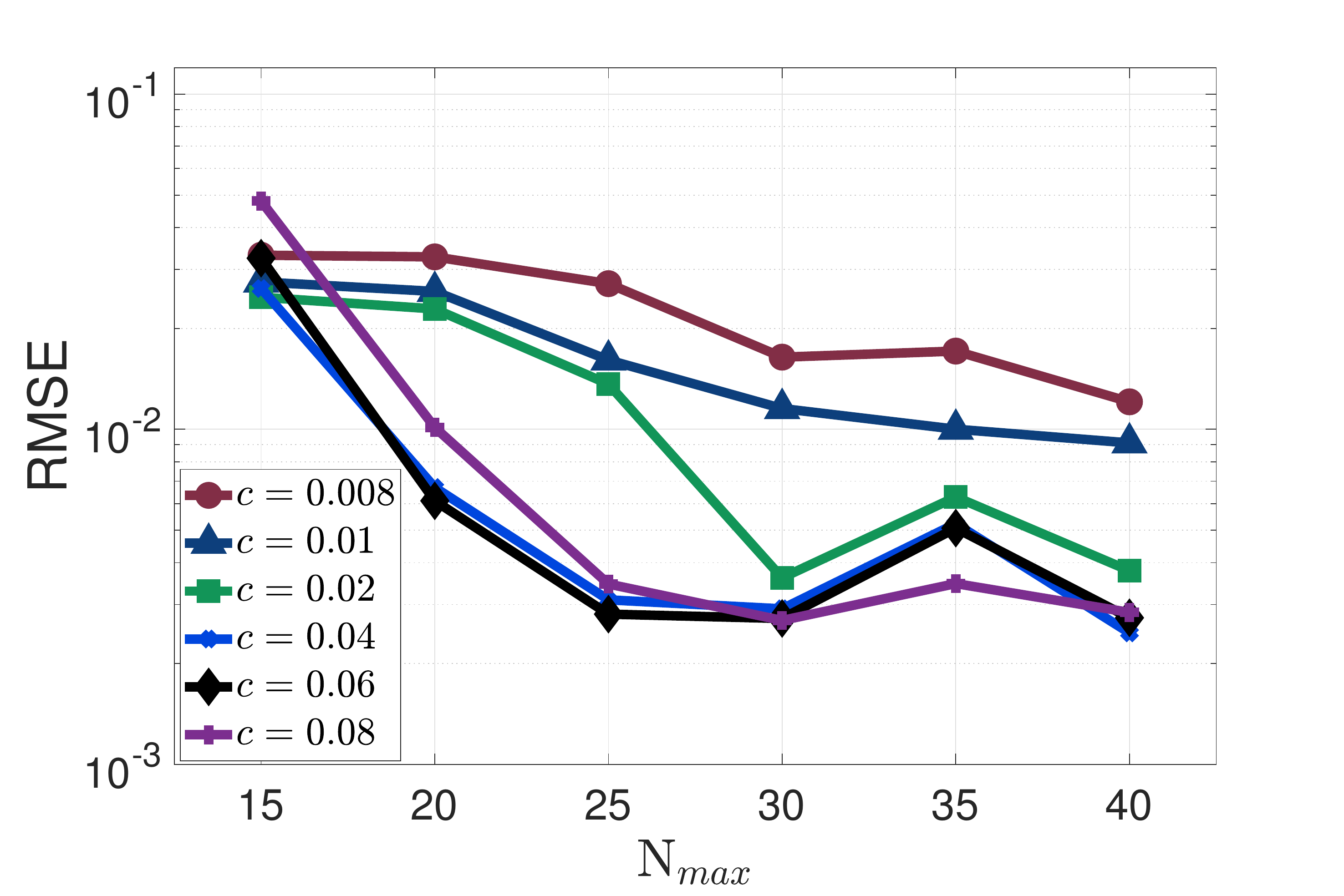}\label{c_n}}
	\caption{Localization accuracy of the SP-ADMM algorithm at iteration step $T= 1000$. Measurement noise variance $\sigma_{i,j}^2=0.02\|\vp_i-\vp_j\|^{2}$; $C_{\text{scale}}$ denotes the side length of the square area where the network is deployed; $\vu_i$ initialized as $\vu_i^0=u\cdot\vone_{N_i},\forall i\in\mathcal{N}$; $C_{\text{range}}$ refers to communication range;  $N_{\max}=\max\left\{N_i,i\in\mathcal{N}\right\}$. (a) penalty parameter $c$ is fixed at 0.0265; (b)(c)(d) initial parameter $u$ is fixed at 0.5.}
	\label{change_c}
\end{figure}

\noindent
We simulated the proposed SP-ADMM algorithm under different values of the penalty parameters $c,\rho$, as well as the initial value $\vu^{0}$. The different network parameters ($C_{\text{scale}}$, $C_{\text{range}},N_{\max}=\{N_i,i\in\mathcal{N}\}$) are also considered. Fig. \ref{change_c} summarizes the RMSE after 1000 iterations. Fig. \ref{rho_c} and Fig. \ref{change_c} not only show the sensitivity of the proposed SP-ADMM algorithm to parameters $\rho,c$, and $\vu^0$,  but also provide practical advice that can assist users, using the proposed SP-ADMM algorithm, to reduce the time in adjusting parameters $\rho,c,\vu^0$ on different networks.

\subsubsection{\textbf{The Selection Strategy of $\rho$ and $c$}}
\noindent
Theorem \ref{them2} guarantees that the proposed Algorithms will converge if $\rho$ is selected in accordance with \eqref{k1_all}-\eqref{c_all} for a given $c$. To corroborate this claim, we present in Fig. \ref{rho_c} the RMSE with iterations for different values of $\rho$, while maintaining $c=0.1$ under the benchmark network (N=500, m=10). Our experimental results demonstrate that, in practical implementation, the selection of proximal penalty coefficients impacts the convergence rate and final accuracy of the proposed algorithm. Specifically, for a benchmark network with $N=500$ and $m=10$, we observed that choosing a smaller value of $c$ or $\rho$ within an appropriate range resulted in the fastest convergence and highest accuracy. For instance, as depicted in Fig. \ref{rho_c}\subref{fixed_rho_c} (Fig. \ref{rho_c}\subref{fixed_c_rho}), when $c$ ($\rho$) is in the range of $[1.1\rho,2\rho]$ ($[0.7c,1.2c]$), the algorithm's performance improves with decreasing values of the parameter $c$. However, selecting a value of $c$ or $\rho$ that is too small or too large may significantly slow down the convergence rate and compromise the final accuracy, as evidenced by the yellow, brown, red, and purple lines in the two subfigures of Fig. \ref{rho_c}. The appropriate range of the parameter values is typically selected empirically based on the specific data used.

\subsubsection{\textbf{Localization Performance with $u$ and $c$ Selection}}

\noindent
Fig. \ref{change_c}\subref{u_ini} and Fig. \ref{change_c}\subref{c_scale} validate this statement. $C_{\text{scale}}$ denotes the side length of the square area where the network is deployed, for example, $C_{\text{scale}}=10$ indicates that the network deployed over $[0,10]\times[0,10]$ area. Fig. \ref{change_c}\subref{u_ini}-\subref{c_scale} employ the same sensor network as Figure , except for the variance of the measured noise $\sigma_{i,j}^2$, which is set to $\sigma_{i,j}^2=0.02\|\vp_i-\vp_j\|^2$. We can see that the curves at different $C_{\text{scale}}$ show tiny fluctuations for varying $u$ than in changing $c$. This is because the convergence rate is strongly dependent on the penalty parameter $c$: please see the supplementary material for the specific form of the convergence rate. These results indicate that the performance of the proposed SP-ADMM is not sensitive to initial parameter $\vu$ but sensitive to penalty parameter $c$. Besides, Fig. \ref{change_c}\subref{c_scale} shows that the appropriate interval value of $c$ for different $C_{\text{scale}}$ is consistent, and cannot be too large or too small. Therefore, parameter $c$ should be adjusted more carefully when the proposed SP-ADMM algorithm is applied to different networks.

\subsubsection{\textbf{The Selection Strategy of $c$ based on the Communication Range}}

\noindent
In this part, we focus on how the parameter $c$ should be chosen for networks with different communication ranges. The network composed of $N=108,m=8$ nodes located over $[0,1]\times[0,1]$ square area and the measurement noise with variance $\sigma_{i,j}^2=0.02\|\vp_i-\vp_j\|^{2}$. Fig. \ref{change_c}\subref{c_range} depicts the evolution of RMSE under networks at different communication ranges. Under the same number of iterative steps, it can be seen from Fig. \ref{change_c}\subref{c_range} that networks with large $C_{\text{range}}$ can obtain a lower RMSE by choosing a large $c$, networks with smaller $C_{\text{range}}$ can reach lower RMSE by selecting smaller $c$.

These results are consistent with \eqref{k1_all}-\eqref{c_all} derived from Lemma \ref{lem1-2}. In detail, the maximum measurement distance $d_{\max}$, the maximum number of adjacent nodes $N_{\max}$ and the total number of adjacent nodes $N_{\text{sum}}$ all increase as $C_{\text{range}}$ increases. Therefore, according to \eqref{k1_all}-\eqref{c_all}, $\rho$ needs to be increased to ensure the convergence of the proposed algorithm. Since we set $\rho=c$ in all experiments, it is reasonable for $c$ to larger. Moreover, we also obtain similar results for the benchmark network, such as $\rho=c=0.11$ ($C_{\text{range}}=0.3,N=500$) is greater than $\rho=c=0.0197$ ($C_{\text{range}}=0.1,N=1000$). These results serve as the rule-of-thumb for selecting the parameter $c$ of the proposed SP-ADMM algorithm when it is applied to different networks.

\subsubsection{\textbf{Convergence rate with} $N_{\max}$}
\noindent
Fig. \ref{change_c}\subref{c_n} displays the RMSE versus the maximum number of adjacent nodes $N_{\max}$ for different values of $c$. The sensor network graph is the same as Fig. \ref{change_c}\subref{c_range} at $C_\text{range}=0.25$. As seen in Fig. \ref{change_c}\subref{c_n}, a lower RMSE is obtained at a larger $N_{\max}$, which means that for the same penalty $c$, the proposed algorithm converges faster on a network with a larger $N_{\max}$. This is indeed consistent with our convergence rate of Theorem \ref{them2}. It can be seen from the supplement that the coefficient $\epsilon_{2}$ of the sublinear convergence rate is inversely proportional to $N_{\max}$ when $c,\rho,\kappa_{1}$ and $\kappa_{2}$ are fixed.
\subsection{Summary}
\noindent
The proposed algorithm has been evaluated against SDP, SF, AM-FD, and ADMM-H methods in terms of RMSE, running time, and communication cost across a range of network scenarios. These scenarios include different number of anchors and sensor nodes, average number of neighboring nodes, and different variance levels of the measurement noise under both benchmark and synthesized network configurations. Furthermore, we have included an empirical investigation on the impact of the penalty parameters $\rho$ and $c$ for our proposed SP-ADMM algorithm.
	
The results indicated that the proposed algorithm surpasses the SDP, SF, and AM-FD methods in terms of localization accuracy, while exhibiting lower computational complexity and communication cost than the ADMM-H method. Although the communication cost of the proposed SP-ADMM is higher than that of the SF and AM-FD algorithms per iteration, it provides higher accuracy at the same communication cost once convergence has been reached. Overall, the proposed SP-ADMM is a promising solution for cooperative localization problem in wireless sensor networks.

{\blue One possible future direction is to explore selective communication strategy among nodes to reduce the communication burden. This involves nodes randomly or sequentially choosing a subset of neighboring nodes to communicate with during each iteration. Additionally, addressing outliers in realistic measurement data and fusing other types of measurements are also important areas for future research.}

\section{Conclusion}\label{sec:conclusion}
In this paper, we have proposed novel distributed parallel SP-ADMM algorithms for wireless sensor network localization in Gaussian measurement noise. Our proposed algorithms directly tackle the challenging nonconvex and nonsmooth problem without resorting to convex relaxation. We have shown that Algorithm \ref{al2} has a low computational complexity (in Table~\ref{comparisons of different alg}) and proved that the algorithm converges to {\blue a critical point of the original problem (in Lemma \ref{KKT_Critical}, \ref{sublinear-lem},} and Theorem \ref{them2}) at rate $\mathcal{O}\left(1/T\right)$. Moreover, we have proposed Algorithm \ref{al1}, which is an improved version of Algorithm \ref{al2} in terms of storage space. Simulation and experimental results not only have shown that the proposed SP-ADMM algorithm is robust and performs well in the different average number of neighboring nodes, measurement noise variance, coverage area, and the number of anchors and non-anchors, but also provided suggestions of the parameters for migrating the proposed algorithm to different networks.

\vspace{-0.1cm}

\bibliography{refs}

\begin{thebibliography}{10}
\providecommand{\url}[1]{#1}
\csname url@samestyle\endcsname
\providecommand{\newblock}{\relax}
\providecommand{\bibinfo}[2]{#2}
\providecommand{\BIBentrySTDinterwordspacing}{\spaceskip=0pt\relax}
\providecommand{\BIBentryALTinterwordstretchfactor}{4}
\providecommand{\BIBentryALTinterwordspacing}{\spaceskip=\fontdimen2\font plus
\BIBentryALTinterwordstretchfactor\fontdimen3\font minus
  \fontdimen4\font\relax}
\providecommand{\BIBforeignlanguage}[2]{{%
\expandafter\ifx\csname l@#1\endcsname\relax
\typeout{** WARNING: IEEEtran.bst: No hyphenation pattern has been}%
\typeout{** loaded for the language `#1'. Using the pattern for}%
\typeout{** the default language instead.}%
\else
\language=\csname l@#1\endcsname
\fi
#2}}
\providecommand{\BIBdecl}{\relax}
\BIBdecl

\bibitem{singh2020handbook}
P.~K. Singh, B.~K. Bhargava, M.~Paprzycki, N.~C. Kaushal, and W.-C. Hong,
  \emph{Handbook of wireless sensor networks: issues and challenges in current
  Scenario's}.\hskip 1em plus 0.5em minus 0.4em\relax Springer, 2020, vol.
  1132.

\bibitem{wymeersch2009cooperative}
H.~Wymeersch, J.~Lien, and M.~Z. Win, ``Cooperative localization in wireless
  networks,'' \emph{Proceedings of the IEEE}, vol.~97, no.~2, pp. 427--450,
  2009.

\bibitem{yin2015cooperative}
F.~Yin, C.~Fritsche, D.~Jin, F.~Gustafsson, and A.~M. Zoubir, ``Cooperative
  localization in {WSN}s using {G}aussian mixture modeling: Distributed {ECM}
  algorithms,'' \emph{IEEE Transactions on Signal Processing}, vol.~63, no.~6,
  pp. 1448--1463, 2015.

\bibitem{pun2021local}
Y.-M. Pun and A.~M.-C. So, ``Local strong convexity of source localization and
  error bound for target tracking under time-of-arrival measurements,''
  \emph{IEEE Transactions on Signal Processing}, vol.~70, pp. 190--201, 2021.

\bibitem{okello2011comparison}
N.~Okello, F.~Fletcher, D.~Musicki, and B.~Ristic, ``Comparison of recursive
  algorithms for emitter localisation using {TDOA} measurements from a pair of
  {UAV}s,'' \emph{IEEE Transactions on Aerospace and Electronic Systems},
  vol.~47, no.~3, pp. 1723--1732, 2011.

\bibitem{xu2017optimal}
S.~Xu and K.~Do{\u{g}}an{\c{c}}ay, ``Optimal sensor placement for {3-D}
  angle-of-arrival target localization,'' \emph{IEEE Transactions on Aerospace
  and Electronic Systems}, vol.~53, no.~3, pp. 1196--1211, 2017.

\bibitem{yin2017received}
F.~Yin, Y.~Zhao, F.~Gunnarsson, and F.~Gustafsson, ``Received-signal-strength
  threshold optimization using {G}aussian processes,'' \emph{IEEE Transactions
  on Signal Processing}, vol.~65, no.~8, pp. 2164--2177, 2017.

\bibitem{jin2020bayesian}
D.~Jin, F.~Yin, C.~Fritsche, F.~Gustafsson, and A.~M. Zoubir, ``Bayesian
  cooperative localization using received signal strength with unknown path
  loss exponent: Message passing approaches,'' \emph{IEEE Transactions on
  Signal Processing}, vol.~68, pp. 1120--1135, 2020.

\bibitem{simonetto2014distributed}
A.~Simonetto and G.~Leus, ``Distributed maximum likelihood sensor network
  localization,'' \emph{IEEE Transactions on Signal Processing}, vol.~62,
  no.~6, pp. 1424--1437, 2014.

\bibitem{biswas2006SDP}
P.~Biswas, T.-C. Lian, T.-C. Wang, and Y.~Ye, ``Semidefinite programming based
  algorithms for sensor network localization,'' \emph{ACM Transactions on
  Sensor Networks (TOSN)}, vol.~2, no.~2, pp. 188--220, 2006.

\bibitem{biswas2006SDPnoisy}
P.~Biswas, T.-C. Liang, K.-C. Toh, Y.~Ye, and T.-C. Wang, ``Semidefinite
  programming approaches for sensor network localization with noisy distance
  measurements,'' \emph{IEEE Transactions on Automation Science and
  Engineering}, vol.~3, no.~4, pp. 360--371, 2006.

\bibitem{wang2008further}
Z.~Wang, S.~Zheng, Y.~Ye, and S.~Boyd, ``Further relaxations of the
  semidefinite programming approach to sensor network localization,''
  \emph{SIAM Journal on Optimization}, vol.~19, no.~2, pp. 655--673, 2008.

\bibitem{srirangarajan2008distributed}
S.~Srirangarajan, A.~H. Tewfik, and Z.-Q. Luo, ``Distributed sensor network
  localization using {SOCP} relaxation,'' \emph{IEEE Transactions on Wireless
  Communications}, vol.~7, no.~12, pp. 4886--4895, 2008.

\bibitem{di2016next}
P.~Di~Lorenzo and G.~Scutari, ``{NEXT}: In-network nonconvex optimization,''
  \emph{IEEE Transactions on Signal and Information Processing over Networks},
  vol.~2, no.~2, pp. 120--136, 2016.

\bibitem{jin2021exploiting}
D.~Jin, F.~Yin, A.~M. Zoubir, and H.~C. So, ``Exploiting sparsity of ranging
  biases for {NLOS} mitigation,'' \emph{IEEE Transactions on Signal
  Processing}, vol.~69, pp. 3782--3795, 2021.

\bibitem{buehrer2018collaborative}
R.~M. Buehrer, H.~Wymeersch, and R.~M. Vaghefi, ``Collaborative sensor network
  localization: Algorithms and practical issues,'' \emph{Proceedings of the
  IEEE}, vol. 106, no.~6, pp. 1089--1114, 2018.

\bibitem{hong2017prox}
M.~Hong, D.~Hajinezhad, and M.-M. Zhao, ``{Prox-PDA}: The proximal primal-dual
  algorithm for fast distributed nonconvex optimization and learning over
  networks,'' in \emph{International Conference on Machine Learning}.\hskip 1em
  plus 0.5em minus 0.4em\relax PMLR, 2017, pp. 1529--1538.

\bibitem{soares2015simple}
C.~Soares, J.~Xavier, and J.~Gomes, ``Simple and fast convex relaxation method
  for cooperative localization in sensor networks using range measurements,''
  \emph{IEEE Transactions on Signal Processing}, vol.~63, no.~17, pp.
  4532--4543, 2015.

\bibitem{piovesan2016cooperative}
N.~Piovesan and T.~Erseghe, ``Cooperative localization in {WSN}s: A hybrid
  convex/nonconvex solution,'' \emph{IEEE Transactions on Signal and
  Information Processing over Networks}, vol.~4, no.~1, pp. 162--172, 2016.

\bibitem{gur2020alternating}
E.~Gur, S.~Sabach, and S.~Shtern, ``Alternating minimization based first-order
  method for the wireless sensor network localization problem,'' \emph{IEEE
  Transactions on Signal Processing}, vol.~68, pp. 6418--6431, 2020.

\bibitem{shi2010distributed}
Q.~Shi, C.~He, H.~Chen, and L.~Jiang, ``Distributed wireless sensor network
  localization via sequential greedy optimization algorithm,'' \emph{IEEE
  Transactions on Signal Processing}, vol.~58, no.~6, pp. 3328--3340, 2010.

\bibitem{bertsekas1989parallel}
D.~Bertsekas, P.~Tsitsiklis, and N.~John, ``Parallel and distributed
  computation: Numeral methods,'' 1989.

\bibitem{erseghe2015distributed}
T.~Erseghe, ``A distributed and maximum-likelihood sensor network localization
  algorithm based upon a nonconvex problem formulation,'' \emph{IEEE
  Transactions on Signal and Information Processing over Networks}, vol.~1,
  no.~4, pp. 247--258, 2015.

\bibitem{luke2017simple}
D.~R. Luke, S.~Sabach, M.~Teboulle, and K.~Zatlawey, ``A simple globally
  convergent algorithm for the nonsmooth nonconvex single source localization
  problem,'' \emph{Journal of Global Optimization}, vol.~69, no.~4, pp.
  889--909, 2017.

\bibitem{wang2021distributed}
Z.~Wang, J.~Zhang, T.-H. Chang, J.~Li, and Z.-Q. Luo, ``Distributed stochastic
  consensus optimization with momentum for nonconvex nonsmooth problems,''
  \emph{IEEE Transactions on Signal Processing}, vol.~69, pp. 4486--4501, 2021.

\bibitem{zhang2020proximal}
J.~Zhang and Z.-Q. Luo, ``A proximal alternating direction method of multiplier
  for linearly constrained nonconvex minimization,'' \emph{SIAM Journal on
  Optimization}, vol.~30, no.~3, pp. 2272--2302, 2020.

\bibitem{ding2021optimal}
L.~Ding, A.~Yurtsever, V.~Cevher, J.~A. Tropp, and M.~Udell, ``An
  optimal-storage approach to semidefinite programming using approximate
  complementarity,'' \emph{SIAM Journal on Optimization}, vol.~31, no.~4, pp.
  2695--2725, 2021.

\bibitem{yousefi2014cooperative}
S.~Yousefi, X.-W. Chang, and B.~Champagne, ``Cooperative localization of mobile
  nodes in {NLOS},'' in \emph{2014 IEEE 25th Annual International Symposium on
  Personal, Indoor, and Mobile Radio Communication}.\hskip 1em plus 0.5em minus
  0.4em\relax IEEE, 2014, pp. 275--279.

\bibitem{vaghefi2015cooperative}
R.~M. Vaghefi and R.~M. Buehrer, ``Cooperative localization in {NLOS}
  environments using semidefinite programming,'' \emph{IEEE Communications
  Letters}, vol.~19, no.~8, pp. 1382--1385, 2015.

\bibitem{boyd2011distributed}
S.~Boyd, N.~Parikh, E.~Chu, B.~Peleato, J.~Eckstein \emph{et~al.},
  ``Distributed optimization and statistical learning via the alternating
  direction method of multipliers,'' \emph{Foundations and
  Trends{\textregistered} in Machine learning}, vol.~3, no.~1, pp. 1--122,
  2011.

\bibitem{fan2021spectrally}
W.~Fan, J.~Liang, G.~Lu, X.~Fan, and H.~C. So, ``Spectrally-agile waveform
  design for wideband mimo radar transmit beampattern synthesis via
  majorization-admm,'' \emph{IEEE Transactions on Signal Processing}, vol.~69,
  pp. 1563--1578, 2021.

\bibitem{bolte2014proximal}
J.~Bolte, S.~Sabach, and M.~Teboulle, ``Proximal alternating linearized
  minimization for nonconvex and nonsmooth problems,'' \emph{Mathematical
  Programming}, vol. 146, no. 1-2, pp. 459--494, 2014.

\bibitem{guo2017convergence}
K.~Guo, D.~Han, and T.-T. Wu, ``Convergence of alternating direction method for
  minimizing sum of two nonconvex functions with linear constraints,''
  \emph{International Journal of Computer Mathematics}, vol.~94, no.~8, pp.
  1653--1669, 2017.

\bibitem{bolte2018first}
J.~Bolte, S.~Sabach, M.~Teboulle, and Y.~Vaisbourd, ``First order methods
  beyond convexity and lipschitz gradient continuity with applications to
  quadratic inverse problems,'' \emph{SIAM Journal on Optimization}, vol.~28,
  no.~3, pp. 2131--2151, 2018.

\bibitem{boct2020proximal}
R.~I. Bo{\c{t}} and D.-K. Nguyen, ``The proximal alternating direction method
  of multipliers in the nonconvex setting: convergence analysis and rates,''
  \emph{Mathematics of Operations Research}, vol.~45, no.~2, pp. 682--712,
  2020.

\bibitem{yashtini2022convergence}
M.~Yashtini, ``Convergence and rate analysis of a proximal linearized admm for
  nonconvex nonsmooth optimization,'' \emph{Journal of Global Optimization},
  vol.~84, no.~4, pp. 913--939, 2022.

\bibitem{hong2016convergence}
M.~Hong, Z.-Q. Luo, and M.~Razaviyayn, ``Convergence analysis of alternating
  direction method of multipliers for a family of nonconvex problems,''
  \emph{SIAM Journal on Optimization}, vol.~26, no.~1, pp. 337--364, 2016.

\bibitem{kumar2016asynchronous}
S.~Kumar, R.~Jain, and K.~Rajawat, ``Asynchronous optimization over
  heterogeneous networks via consensus admm,'' \emph{IEEE Transactions on
  Signal and Information Processing over Networks}, vol.~3, no.~1, pp.
  114--129, 2016.

\bibitem{sun2019two}
K.~Sun and X.~A. Sun, ``A two-level distributed algorithm for general
  constrained non-convex optimization with global convergence,'' \emph{arXiv
  preprint arXiv:1902.07654}, 2019.

\bibitem{zhang2020improved}
X.~Zhang, J.~Ma, Z.~Cheng, S.~Huang, C.~W. de~Silva, and T.~H. Lee, ``Improved
  hierarchical admm for nonconvex cooperative distributed model predictive
  control,'' \emph{arXiv preprint arXiv:2011.00463}, 2020.

\bibitem{nesterov1983method}
Y.~E. Nesterov, ``A method for solving the convex programming problem with
  convergence rate ${O}(1/k^2)$,'' in \emph{Dokl. akad. nauk Sssr}, vol. 269,
  1983, pp. 543--547.

\bibitem{YY}
\BIBentryALTinterwordspacing
Y.~Ye, ``Computational optimization laboratory.stanford university.'' [Online].
  Available: \url{https://web.stanford.edu/~yyye/Col}
\BIBentrySTDinterwordspacing

\bibitem{patwari2003relative}
N.~Patwari, A.~O. Hero, M.~Perkins, N.~S. Correal, and R.~J. O'dea, ``Relative
  location estimation in wireless sensor networks,'' \emph{IEEE Transactions on
  signal processing}, vol.~51, no.~8, pp. 2137--2148, 2003.

\bibitem{sturm1999using}
J.~F. Sturm, ``{Using SeDuMi 1.02, a MATLAB} toolbox for optimization over
  symmetric cones,'' \emph{Optimization methods and software}, vol.~11, no.
  1-4, pp. 625--653, 1999.

\bibitem{beck2017first}
A.~Beck, \emph{First-order methods in optimization}.\hskip 1em plus 0.5em minus
  0.4em\relax Society for Industrial and Applied Mathematics, 2017.

\bibitem{Hong2017}
M.~Hong and Z.-Q. Luo, ``On the linear convergence of the alternating direction
  method of multipliers,'' \emph{Mathematical Programming}, vol. 162, no.~1,
  pp. 165--199, 2017.

\end{thebibliography}

\newpage

\begin{center}
	\LARGE
	\bf Supplementary Materials of Manuscript "Distributed Scaled Proximal ADMM Algorithms for Cooperative Localization in WSNs"
	
\end{center}

\begin{center}
	\large
	Mei Zhang, Zhiguo Wang, Feng Yin, and Xiaojing Shen	
	
\end{center}

\makeatletter
\setcounter{section}{0}
\renewcommand\thesection{\arabic{section}}
\renewcommand{\section}{
		\@startsection{section}
		{1}
		{0pt}
		{-3.5ex plus -1ex minus -0.2ex}
		{2.3ex plus 0.2ex }
		{\normalfont\Large\bfseries}
}
\makeatother

\setcounter{equation}{0}
	\renewcommand{\theequation}{S.\arabic{equation}}

	\section{Proof of the Remark \ref{r_z_ztilde}}\label{proof_r_z_ztilde}
	
	Recall that
	\begin{equation*}
	\mathcal{Z}=\left\{\vz|\vz_{i,j}^{+}=\vp_j=\vz_{j,i}^{-},~\forall~i\in\mathcal{N},~j\in\mathcal{N}_i\right\},
	\end{equation*}
	and $\mW_i$ is a diagonal matrix defined in \eqref{w}. Substituting them into \eqref{lem_zt1}, then the optimization \eqref{lem_zt1} is rewritten as
	\begin{align}\label{x_pro}
	\min_{\vz\in\mathcal{X}}\sum_{i\in\mathcal{N}}\Big[\left(c+1\right)N_i\left\|\vp_{i}-\tilde{\vp}_{i}^{t+1}\right\|^{2}+\sum_{j\in\mathcal{N}_i}c\big\|\vz_{i,j}^{-}-(\tilde{\vz}_{i,j}^{-})^{t+1}\big\|^{2}+\sum_{j\in\mathcal{N}_i}\big\|\vz_{i,j}^{+}-(\tilde{\vz}_{i,j}^{+})^{t+1}\big\|^{2}\Big].
	\end{align}
	From the definition of $\mathcal{X}$ in \eqref{anchor_set} and \eqref{x_pro}, we immediately obtain the optimal solution for $\vp_{i}^{t+1}$ as follows
	\begin{equation}\label{xii1}
	\vp_{i}^{t+1}=\left\{\begin{array}{ll}
	\tilde{\vp}_{i,i}^{t+1},&\text{for}~i\notin\mathcal{A},\\
	\va_i,&\text{for}~i\in\mathcal{A}.\end{array}\right.
	\end{equation}
	Combining \eqref{xii1} with \eqref{xtilde_a} yields
	\begin{equation}\label{xii}
	\vp_{i}^{t+1}=\tilde{\vp}_{i}^{t+1},\forall i\in\mathcal{N}.
	\end{equation}
	The terms in problem \eqref{x_pro} related to $\vz_{i,j}^-$ are
	\begin{equation*}
	c\left\|\vz_{i,j}^--(\tilde{\vz}_{i,j}^{-})^{t+1}\right\|^{2}+\left\|\vz_{i,j}^--(\tilde{\vz}_{j,i}^{+})^{t+1}\right\|^{2},\,\forall i\in\mathcal{N},j\in\mathcal{N}_i.
	\end{equation*}
	By the optimality condition, we have
	\begin{equation*}
	c\left((\vz_{i,j}^-)^{t+1}-(\tilde{\vz}_{i,j}^{-})^{t+1}\right)+\left((\vz_{i,j}^-)^{t+1}-(\tilde{\vz}_{j,i}^{+})^{t+1}\right)=\vzero.
	\end{equation*}
	After rearranging the terms, we obtain
	\begin{equation}\label{xji}
(\vz_{i,j}^-)^{t+1}=\frac{1}{c+1}\left(c\cdot(\tilde{\vz}_{i,j}^{-})^{t+1}+(\tilde{\vz}_{j,i}^{+})^{t+1}\right).
	\end{equation}
	Similarly, for $i\in\mathcal{N},j\in\mathcal{N}_i$, we can obtain the closed form solution for $\vz_{i,j}^+$ as follows
	\begin{equation}\label{xij}
	(\vz_{i,j}^+)^{t+1}=\frac{1}{c+1}\left(\tilde{\vz}_{i,j}^{+})^{t+1}+c\cdot(\tilde{\vz}_{j,i}^{-})^{t+1}\right),\, \forall i\in\mathcal{N},j\in\mathcal{N}_i.
	\end{equation}
	This completes the proof of Remark \ref{r_z_ztilde} by combining \eqref{z} and \eqref{xii}-\eqref{xij}.

	\section{Proof of the Lemma \ref{LL}}\label{proof_LL}
	
	\noindent
	In \eqref{L111}, we take $\vz_i=\vz_i^{t+1}$ and $\vz_i=\vz_i^t$, it yields that
	\begin{align}\label{LLL}
	\mathcal{L}_i\left(\vz_i^{t+1},\vu_i^t,\vlambda_i^t\right)+\frac{c}{2}\|\vz_i^{t+1}-\vz_i^t\|^{2}_{\mB_i^T\mB_i}-\mathcal{L}_i\left(\vz_i^t,\vu_i^t,\vlambda_i^t\right)=\frac{1}{2}\|\vz_i^{t+1}-\tilde{\vz}_i^{t+1}\|_{\mW_i}^2-\frac{1}{2}\|\vz_i^t-\tilde{\vz}_i^{t+1}\|_{\mW_i}^2.
	\end{align}
	Since $\vz_i^{t+1}$ is the optimal solution from Remark \ref{r_z_ztilde} and $\vz_i^t\in \mathcal{Z}\cap\mathcal{X}$, $\vz_i^{t+1}\in\mathcal{Z}\cap\mathcal{X}$, then we have
	\begin{equation}\label{zzw}
	\frac{1}{2}\|\vz_i^{t+1}-\tilde{\vz}_i^{t+1}\|^{2}_{\mW_i}\leq\frac{1}{2}\|\vz_i^t-\tilde{\vz}_i^{t+1}\|^{2}_{\mW_i}.
	\end{equation}
	Substituting \eqref{zzw} into \eqref{LLL}, we obtain
	\begin{equation}\label{LLL1}
	\mathcal{L}_i\left(\vz_i^{t+1},\vu_i^t,\vlambda_i^t\right)-\mathcal{L}_i\left(\vz_i^t,\vu_i^t,\vlambda_i^t\right)\leq-\frac{c}{2}\|\vz_i^{t+1}-\vz_i^t\|^{2}_{\mB_i^T\mB_i}.
	\end{equation}
	Using the $\vu$ update in \eqref{u_up} with the same technique, we can obtain the following inequality
	\begin{equation}\label{uu}
	\mathcal{L}_i\left(\vz_i^{t+1},\vu_i^{t+1},\vlambda_i^t\right)-\mathcal{L}_i\left(\vz_i^{t+1},\vu_i^t,\vlambda_i^t\right)\leq-\frac{\rho}{2}\|\vu_i^{t+1}-\vu_i^t\|^{2}.
	\end{equation}
	Recall that $\vlambda_i^{t+1}=\vlambda_i^t+c\mA_i\vz_i^{t+1}$ in \eqref{du2}, we have the trivial equality
	\begin{equation}\label{apa_lem2}
	\mathcal{L}_i\left(\vz_i^{t+1},\vu_i^{t+1},\vlambda_i^{t+1}\right)-\mathcal{L}_i\left(\vz_i^{t+1},\vu_i^{t+1},\vlambda_i^t\right)=\frac{1}{c}\|\vlambda_i^{t+1}-\vlambda_i^t\|^{2}.
	\end{equation}
	Applying \eqref{new1}, we have
	\begin{equation}\label{anew}
	\mA_i^T\vlambda_i^t=\mQ_i^T\mD_i\vu_i^t-\mW_i\tilde{\vz}_i^{t+1}+c\mB_i^T\mB_i\vz_i^t.
	\end{equation}
	Both sides plus the term $c\mA_i^T\mA_i\vz_i^{t+1}$, by \eqref{du2}, it yields
	\begin{align}\nonumber
	\mA_i^T\vlambda_i^{t+1}&=\mQ_i^T\mD_i\vu_i^t-\mW_i\tilde{\vz}_i^{t+1}+c\mB_i^T\mB_i\vz_i^t+c\mA_i^T\mA_i\vz_i^{t+1}\\\label{anew1}
	&=\mQ_i^T\mD_i\vu_i^t-\mQ_i^T\mQ_i\tilde{\vz}_i^{t+1}-c\mA_i^T\mA_i(\tilde{\vz}_i^{t+1}-\vz_i^{t+1})-c\mB_i^T\mB_i(\tilde{\vz}_i^{t+1}-\vz_i^t),
	\end{align}
	where the last equality dues to the definition of $\mW_i$ in \eqref{new2}. Let $\widetilde{\sigma}_{min}$ denote the smallest non-zero eigenvalue of $\mA_i^T\mA_i$, we have
	\begin{equation*}      \widetilde{\sigma}_{min}\|\vlambda_i^{t+1}-\vlambda_i^t\|^{2}\leq\|\mA_i^T(\vlambda_i^{t+1}-\vlambda_i^t)\|^{2}.
	\end{equation*}
	From the definition of $\mA_i^T\mA_i$ in \eqref{AtA}, it derives $\widetilde{\sigma}_{min}=1$. The above inequality combined with \eqref{anew1} implies that
	\begin{align}\nonumber
	&\|\vlambda_i^{t+1}-\vlambda_i^t\|^{2}\\\nonumber\leq
	&\|\mQ_i^T\mD_i\left(\vu_i^t-\vu_i^{t-1}\right)-\mQ_i^T\mQ_i(\tilde{\vz}_i^{t+1}-\tilde{\vz}_i^t)-c\mA_i^T\mA_i\left[\tilde{\vz}_i^{t+1}-\vz_i^{t+1}-\left(\tilde{\vz}_i^t-\vz_i^t\right)\right]\\\nonumber
	&-c\mB_i^T\mB_i\left[\tilde{\vz}_i^{t+1}-\vz_i^t-\left(\tilde{\vz}_i^t-\vz_i^{t-1}\right)\right]\|^{2}\\\nonumber\leq &3\left(N_i+1\right)\|\mQ_i\left(\tilde{\vz}_i^{t+1}-\tilde{\vz}_i^t\right)-\mD_i\left(\vu_i^t-\vu_i^{t-1}\right)\|^{2}+3c^{2}\left(N_i+1\right)\|\tilde{\vz}_i^{t+1}-\vz_i^{t+1}-\left(\tilde{\vz}_i^t-\vz_i^t\right)\|^{2}_{\mA_i^T\mA_i}\\\label{alamd1}
	&+3c^{2}\|\mB_i^T\mB_i\|\|\tilde{\vz}_i^{t+1}-\vz_i^t-\left(\tilde{\vz}_i^t-\vz_i^{t-1}\right)\|^{2}_{\mB_i^T\mB_i},
	\end{align}
	where the last inequality holds because of the triangle inequality and $\|\mQ_i^T\mQ_i\|=\|\mA_i^T\mA_i\|=N_i+1$ from \eqref{AtA}-\eqref{QtQ}, $\|\mA_i^T\mA_i\|$ denotes the spectral norm of a matrix $\mA_i^T\mA_i$. Along with the form of the matrix $c\mB_i^T\mB_i$ in \eqref{b2}, we have
	\begin{equation}\label{evib}
	\|\mB_i^T\mB_i\|\leq\frac{\left(1+c\right)\left(1+N_{\max}\right)}{c},
	\end{equation}
	where $N_{\max}:=\max\{N_i,i\in\mathcal{N}\}$. Substituting \eqref{evib} into \eqref{alamd1} and the resulting applied to \eqref{apa_lem2}, we have
	\begin{align}\nonumber	&\mathcal{L}_i\left(\vz_i^{t+1},\vu^{t+1},\vlambda_i^{t+1}\right)-\mathcal{L}\left(\vz_i^{t+1},\vu_i^{t+1},\vlambda_i^t\right)\\\nonumber\leq
	&\frac{3\left(N_{\max}+1\right)}{c}\|\mQ_i\left(\tilde{\vz}_i^{t+1}-\tilde{\vz}_i^t\right)-\mD_i\left(\vu_i^t-\vu_i^{t-1}\right)\|^{2}+3c\left(N_{\max}+1\right)\|\tilde{\vz}_i^{t+1}-\vz_i^{t+1}-\left(\tilde{\vz}_i^t-\vz_i^t\right)\|^{2}_{\mA_i^T\mA_i}\\\label{zzu}
	&+3\left(1+c\right)\left(1+N_{\max}\right)\|\tilde{\vz}_i^{t+1}-\vz_i^t-\left(\tilde{\vz}_i^t-\vz_i^{t-1}\right)\|^{2}_{\mB_i^T\mB_i}.
	\end{align}
	Combing \eqref{LLL1}, \eqref{uu}, \eqref{zzu}, moreover, taking summation over $i\in \mathcal{N}$, we obtain the final result \eqref{Lt1_Lt}.

	\section{Proof of the Lemma \ref{lem1-2}}\label{apc}

	\noindent
	The following Lemma shows the descent of the desired term.
	
	\begin{Lemma}\label{lem2-2}
		Suppose $c\mB_i^T\mB_i$ takes the form of \eqref{cbb}, then the following is true for Algorithm \ref{al2}
		\begin{align*}
		&\sum_{i\in\mathcal{N}}\frac{c}{2}\|\mA_i\tilde{\vz}_i^{t+1}\|^{2}+\frac{c}{2}\|\vz_i^{t+1}-\vz_i^t\|^{2}_{\mB_i^T\mB_i}\\\leq
		&\sum_{i\in\mathcal{N}}\Big[\frac{c}{2}\|\mA_i\tilde{\vz}_i^t\|^{2}+\frac{c}{2}\|\vz_i^t-\vz_i^{t-1}\|^{2}_{\mB_i^T\mB_i}-\frac{c}{2}\|\vz_i^{t+1}-\vz_i^t\|^{2}_{\mA_i^T\mA_i}-\frac{1}{2}\|\vz_i^{t+1}-\vz_i^t\|^{2}_{\mQ_i^T\mQ_i}\\
		&\qquad+\frac{d_{\max}^{2}}{2}\|\vu_i^t-\vu_i^{t-1}\|^{2}-\frac{1}{2}\|\mQ_i\left(\tilde{\vz}_i^{t+1}-\tilde{\vz}_i^t\right)-\mD_i\left(\vu_i^t-\vu_i^{t-1}\right)\|^{2}\\
		&\qquad-\frac{c}{2}\|\tilde{\vz}_i^{t+1}-\vz_i^{t+1}-(\tilde{\vz}_i^t-\vz_i^t)\|^{2}_{\mA_i^T\mA_i}+\frac{c\left(N_i+1\right)}{2}\|\vz_i^t-\tilde{\vz}_i^{t+1}\|^{2}\\
		&\qquad-\frac{c}{2}\|\tilde{\vz}_i^{t+1}-\vz_i^t-(\tilde{\vz}_i^t-\vz_i^{t-1})\|^{2}_{\mB_i^T\mB_i}\Big],
		\end{align*}
		where $d_{\max}:=\max\{d_{ij},i\in\mathcal{N},j\in\mathcal{N}_i\}$.
	\end{Lemma}
	\begin{IEEEproof}		
		According to \eqref{anew1}, we get
		\begin{align}\nonumber
		\langle\mA_i^T(\vlambda_i^{t+1}-\vlambda_i^t),\tilde{\vz}_i^{t+1}-\tilde{\vz}_i^t\rangle=&\left\langle-\mW_i\left(\tilde{\vz}_i^{t+1}-\tilde{\vz}_i^t\right)+\mQ_i^T\mD_i\left(\vu_i^t-\vu_i^{t-1}\right),\tilde{\vz}_i^{t+1}-\tilde{\vz}_i^t\right\rangle\\\label{g1}
		&+\left\langle c\mA_i^T\mA_i\left(\vz_i^{t+1}-\vz_i^t\right)+c\mB_i^T\mB_i\left(\vz_i^t-\vz_i^{t-1}\right),\tilde{\vz}_i^{t+1}-\tilde{\vz}_i^t\right\rangle.
		\end{align}
		
		Let us bound the left-hand side (lhs) and the rhs of \eqref{g1} separately. First, the lhs of \eqref{g1} can be expressed as
		\begin{equation*}
		\langle\mA_i^T(\vlambda_i^{t+1}-\vlambda_i^t),\tilde{\vz}_i^{t+1}-\tilde{\vz}_i^t\rangle=\left\langle c\mA_i^T\mA_i\vz_i^{t+1},\tilde{\vz}_i^{t+1}-\tilde{\vz}_i^t\right\rangle,
		\end{equation*}
		where the equality dues to dual update in \eqref{du2}. Furthermore,
		\begin{align}\nonumber
		&\left\langle c\mA_i^T\mA_i\vz_i^{t+1},\tilde{\vz}_i^{t+1}-\tilde{\vz}_i^t\right\rangle\\\nonumber=
		&\left\langle c\mA_i^T\mA_i\left(\vz_i^{t+1}-\vz_i^t\right),\tilde{\vz}_i^{t+1}-\tilde{\vz}_i^t\right\rangle-\left\langle c\mA_i^T\mA_i\vz_i^t,\tilde{\vz}_i^{t+1}-\tilde{\vz}_i^t\right\rangle\\\nonumber=
		&-\frac{c}{2}\|\vz_i^{t+1}-\vz_i^t-\left(\tilde{\vz}_i^{t+1}-\tilde{\vz}_i^t\right)\|_{\mA_i^T\mA_i}^{2}+\frac{c}{2}\|\tilde{\vz}_i^{t+1}-\tilde{\vz}_i^t\|^{2}_{\mA_i^T\mA_i}+\frac{c}{2}\|\vz_i^{t+1}-\vz_i^t\|^{2}_{\mA_i^T\mA_i}\\\label{ATlambda1}
		&+\frac{c}{2}\|\mA_i\left(\vz_i^t-\tilde{\vz}_i^t\right)\|^{2}-\frac{c}{2}\|\mA_i\left(\vz_i^t-\tilde{\vz}_i^{t+1}\right)\|^{2}+\frac{c}{2}\|\mA_i\tilde{\vz}_i^{t+1}\|^{2}-\frac{c}{2}\|\mA_i\tilde{\vz}_i^t\|^{2},
		\end{align}
		where the final equality dues to $\langle\va,\vbb\rangle=-\frac{1}{2}\|\va-\vbb\|^2+\frac{1}{2}\|\va\|^2+\frac{1}{2}\|\vbb\|^2$ for any $\va,\vbb\in\mathbb{R}^{\left(2N_i+1\right)n}$. And from the compatibility of the norm that
		\begin{equation}\label{cAzz}
		\frac{c}{2}\|\mA_i\left(\vz_i^t-\tilde{\vz}_i^{t+1}\right)\|^{2}\leq\frac{c\left(N_i+1\right)}{2}\|\vz_i^t-\tilde{\vz}_i^{t+1}\|^{2},
		\end{equation}
		where $\|\mA_i^T\mA_i\|=N_i+1$ due to \eqref{AtA}.
		
		Second, we have rewritten the first row on the rhs of the \eqref{g1}:
		\begin{align}\nonumber
		&\langle-\mW_i(\tilde{\vz}_i^{t+1}-\tilde{\vz}_i^t)+\mQ_i^T\mD_i(\vu_i^t-\vu_i^{t-1}),\tilde{\vz}_i^{t+1}-\tilde{\vz}_i^t\rangle\\\nonumber=
		&-\|\tilde{\vz}_i^{t+1}-\tilde{\vz}_i^t\|^{2}_{\mW_i}-\frac{1}{2}\|\mQ_i\left(\tilde{\vz}_i^{t+1}-\tilde{\vz}_i^t\right)-\mD_i\left(\vu_i^t-\vu_i^{t-1}\right)\|^{2}\\\label{QQzzQSuu}
		&+\frac{1}{2}\|\mQ_i\left(\tilde{\vz}_i^{t+1}-\tilde{\vz}_i^t\right)\|^{2}+\frac{1}{2}\|\mD_i(\vu_i^t-\vu_i^{t-1})\|^{2}.
		\end{align}
		For the second row on the rhs of the \eqref{g1}, we have
		\begin{align}\nonumber
		&\langle c\mA_i^T\mA_i\left(\vz_i^{t+1}-\vz_i^t\right)+c\mB_i^T\mB_i\left(\vz_i^t-\vz_i^{t-1}\right),\tilde{\vz}_i^{t+1}-\tilde{\vz}_i^t\rangle\\\nonumber=
		&\frac{c}{2}\|\mA_i\left(\tilde{\vz}_i^{t+1}-\tilde{\vz}_i^t\right)\|^{2}+\frac{c}{2}\|\mA_i\left(\vz_i^{t+1}-\vz_i^t\right)\|^{2}-\frac{c}{2}\|\tilde{\vz}_i^{t+1}-\tilde{\vz}_i^t-\left(\vz_i^{t+1}-\vz_i^t\right)\|^{2}_{\mA_i^T\mA_i}\\\label{cAAzz}
		&+\frac{c}{2}\|\tilde{\vz}_i^{t+1}-\tilde{\vz}_i^t\|^{2}_{\mB_i^T\mB_i}+\frac{c}{2}\|\vz_i^t-\vz_i^{t-1}\|^{2}_{\mB_i^T\mB_i}-\frac{c}{2}\|\tilde{\vz}_i^{t+1}-\vz_i^t-(\tilde{\vz}_i^t-\vz_i^{t-1})\|^{2}_{\mB_i^T\mB_i}.
		\end{align}
		Recall $\mW_i=\mQ_i^T\mQ_i+c\mA_i^T\mA_i+c\mB_i^T\mB_i$. Combing \eqref{QQzzQSuu} and \eqref{cAAzz}, then the rhs of \eqref{g1} can be expressed as
		\begin{align}\nonumber
		&\left\langle-\mW_i\left(\tilde{\vz}_i^{t+1}-\tilde{\vz}_i^t\right)+\mQ_i^T\mD_i\left(\vu_i^t-\vu_i^{t-1}\right),\tilde{\vz}_i^{t+1}-\tilde{\vz}_i^t\right\rangle\\\nonumber
		&+\left\langle c\mA_i^T\mA_i\left(\vz_i^{t+1}-\vz_i^t\right)+c\mB_i^T\mB_i\left(\vz_i^t-\vz_i^{t-1}\right),\tilde{\vz}_i^{t+1}-\tilde{\vz}_i^t\right\rangle\\\nonumber=
		&-\frac{1}{2}\|\mQ_i\left(\tilde{\vz}_i^{t+1}-\tilde{\vz}_i^t\right)-\mD_i\left(\vu_i^t-\vu_i^{t-1}\right)\|^{2}-\frac{c}{2}\|\tilde{\vz}_i^{t+1}-\tilde{\vz}_i^t-\left(\vz_i^{t+1}-\vz_i^t\right)\|^{2}_{\mA_i^T\mA_i}\\\nonumber
		&-\frac{1}{2}\|\tilde{\vz}_i^{t+1}-\tilde{\vz}_i^t\|^{2}_{\mW_i}-\frac{c}{2}\|\tilde{\vz}_i^{t+1}-\vz_i^t-(\tilde{\vz}_i^t-\vz_i^{t-1})\|^{2}_{\mB_i^T\mB_i}+\frac{c}{2}\|\mA_i\left(\vz_i^{t+1}-\vz_i^t\right)\|^{2}\\\label{right_hand}&+\frac{c}{2}\|\vz_i^t-\vz_i^{t-1}\|^{2}_{\mB_i^T\mB_i}+\frac{1}{2}\|\mD_i(\vu_i^t-\vu_i^{t-1})\|^{2}.
		\end{align}
		From the optimality conditions for the strongly convex optimization problem \eqref{lem_zt1}, for any $\vz_i\in\mathcal{Z}\cap\mathcal{X}$ we have
		\begin{align*}
		\sum_{i\in\mathcal{N}}\langle\mW_i(\vz_i^{t+1}-\tilde{\vz}_i^{t+1}),\vz_i^{t+1}-\vz_i\rangle&\leq0,\\
		\sum_{i\in\mathcal{N}}\langle\mW_i(\vz_i^t-\tilde{\vz}_i^t),\vz_i^t-\vz_i\rangle&\leq0.
		\end{align*}
		Plugging $\vz_i=\vz_i^t$ into the first inequality and $\vz_i=\vz_i^{t+1}$ into the second, adding the resulting inequalities, we have
		\begin{equation*}
		\sum_{i\in\mathcal{N}}\langle\mW_i\left(\tilde{\vz}_i^{t+1}-\tilde{\vz}_i^t+\vz_i^t-\vz_i^{t+1}\right),\vz_i^t-\vz_i^{t+1}\rangle\leq 0.
		\end{equation*}
		Rearranging the above inequality, we have
		\begin{equation}\label{wzztilde}
		\sum_{i\in\mathcal{N}}\|\vz_i^t-\vz_i^{t+1}\|^{2}_{\mW_i}\leq\sum_{i\in\mathcal{N}}\langle\mW_i\left(\tilde{\vz}_i^{t+1}-\tilde{\vz}_i^t\right),\vz_i^{t+1}-\vz_i^t\rangle.
		\end{equation}
		Using $\langle\va,\vbb\rangle=-\frac{1}{2}\|\va-\vbb\|^{2}+\frac{1}{2}\|\va\|^{2}+\frac{1}{2}\|\vbb\|^{2}$, then \eqref{wzztilde} becomes
		\begin{align}\label{zww}
		\sum_{i\in\mathcal{N}}\frac{1}{2}\|\vz_i^t-\vz_i^{t+1}\|^{2}_{\mW_i}\leq\sum_{i\in\mathcal{N}}\frac{1}{2}\|\tilde{\vz}_i^{t+1}-\tilde{\vz}_i^t\|^{2}_{\mW_i}-\frac{1}{2}\|\tilde{\vz}_i^{t+1}-\tilde{\vz}_i^t-\left(\vz_i^{t+1}-\vz_i^t\right)\|^{2}.
		\end{align}
		Substituting \eqref{cAzz} \eqref{zww} into \eqref{ATlambda1} \eqref{right_hand}, respectively, againg using $\mW_i=\mQ_i^T\mQ_i+c\mA_i^T\mA_i+c\mB_i^T\mB_i$ and combining \eqref{g1}, we obtain the result.
	\end{IEEEproof}
	
	Note that it remains to bound the term $\|\vz_i^t-\tilde{\vz}_i^{t+1}\|^2$ in Lemma \ref{lem2-2}. Next, we establish a simple Lemma to ensure that the term decreases.
	
	\begin{Lemma}\label{lem2-4}
		Let $\left\{\left(\vz_i^t,\vu_i^t,\vlambda_i^t\right)\right\}$ be the sequence generated by Algorithm \ref{al2}, $c\mB_i^T\mB_i$ takes the form of \eqref{cbb}. Then we have
		\begin{align*}
		&\sum_{i\in\mathcal{N}}\frac{c}{2}\|\mA_i\vz_i^{t+1}\|^{2}+\frac{c}{2}\|\vz_i^{t+1}-\vz_i^t\|^{2}_{\mB_i^T\mB_i}\\\leq
		&\sum_{i\in\mathcal{N}}\left[\frac{c}{2}\|\mA_i\vz_i^t\|^{2}+\frac{c}{2}\|\vz_i^t-\vz_i^{t-1}\|^{2}_{\mB_i^T\mB_i}+\frac{d_{\max}^{2}}{2}\|\vu_i^t-\vu_i^{t-1}\|^{2}-\frac{c\cdot\tilde{\tau}_{\min}}{2N_{\text{sum}}n\left(c+1\right)^{2}}\|\vz_i^t-\tilde{\vz}_i^{t+1}\|^{2} \right],
		\end{align*}
		where $\tilde{\tau}_{\min}:=\min\{\left(c+1\right)^{2}N_i^{2}+c^{2}N_i+N_i,i\in\mathcal{N}\}$ and $N_{\text{sum}}:=\sum_{i\in\mathcal{N}}N_i$.
	\end{Lemma}
	
	\begin{IEEEproof}
		Using \eqref{anew1} again, we can get
		\begin{align}\nonumber	\langle\mA_i^T(\vlambda_i^{t+1}-\vlambda_i^t),\vz_i^{t+1}-\vz_i^t\rangle=
		&\left\langle-\mW_i\left(\tilde{\vz}_i^{t+1}-\tilde{\vz}_i^t\right),\vz_i^{t+1}-\vz_i^t\right\rangle+c\|\mA_i\left(\vz_i^{t+1}-\vz_i^t\right)\|^{2}\\\label{g12}
		&+\left\langle\mQ_i^T\mD_i\left(\vu_i^t-\vu_i^{t-1}\right)+c\mB_i^T\mB_i\left(\vz_i^t-\vz_i^{t-1}\right),\vz_i^{t+1}-\vz_i^t\right\rangle.
		\end{align}
		First, the lhs of \eqref{g12} can be expressed as
		\begin{align}\nonumber	\langle\mA_i^T(\vlambda_i^{t+1}-\vlambda_i^t),\vz_i^{t+1}-\vz_i^t\rangle=
		&\nonumber\langle c\mA_i\vz_i^{t+1},\mA_i\vz_i^{t+1}-\mA_i\vz_i^t\rangle\\=
		&\label{lambdazz}\frac{c}{2}\|\mA_i\vz_i^{t+1}\|^{2}-\frac{c}{2}\|\mA_i\vz_i^t\|^{2}+\frac{c}{2}\|\mA_i\left(\vz_i^{t+1}-\vz_i^t\right)\|^{2}.
		\end{align}
		where the first equality dues to \eqref{du2}, and the second equality dues to $\langle\va,\va-\vbb\rangle=\frac{1}{2}\|\va-\vbb\|^2+\frac{1}{2}\|\va\|^2-\frac{1}{2}\|\vbb\|^2$ for any $\va,\vbb\in\mathbb{R}^{\left(2N_i+1\right)n}$.
		
		Next, we deal with the rhs of the \eqref{g12}. By \eqref{wzztilde} yields
		\begin{equation}\label{zw}
		\sum_{i\in\mathcal{N}}\left\langle-\mW_i\left(\tilde{\vz}_i^{t+1}-\tilde{\vz}_i^t\right),\vz_i^{t+1}-\vz_i^t\right\rangle\leq-\sum_{i\in\mathcal{N}}\|\vz_i^{t+1}-\vz_i^t\|^{2}_{\mW_i}.
		\end{equation}
		And using the Cauchy–Schwarz inequality, we have
		\begin{align}\nonumber	   &\left\langle\mQ_i^T\mD_i\left(\vu_i^t-\vu_i^{t-1}\right)+c\mB_i^T\mB_i\left(\vz_i^t-\vz_i^{t-1}\right),\vz_i^{t+1}-\vz_i^t\right\rangle\\\label{g2}\leq
		&\frac{1}{2}\|\mD_i\left(\vu_i^t-\vu_i^{t-1}\right)\|^{2}+\frac{1}{2}\|\mQ_i\left(\vz_i^{t+1}-\vz_i^t\right)\|^{2}+\frac{c}{2}\|\vz_i^t-\vz_i^{t-1}\|^{2}_{\mB_i^T\mB_i}+\frac{c}{2}\|\vz_i^{t+1}-\vz_i^t\|^{2}_{\mB_i^T\mB_i}.
		\end{align}
		Combing \eqref{g12}-\eqref{g2} and using $\mW_i=c\mB_i^T\mB_i+c\mA_i^T\mA_i+\mQ_i^T\mQ_i$, we have
		\begin{align}\nonumber
		&\sum_{i\in\mathcal{N}}\frac{c}{2}\|\mA_i\vz_i^{t+1}\|^{2}+\frac{c}{2}\|\vz_i^{t+1}-\vz_i^t\|^{2}_{\mB_i^T\mB_i}\\\label{cazb}\leq
		&\sum_{i\in\mathcal{N}}\Big[\frac{c}{2}\|\mA_i\vz_i^t\|^{2}+\frac{c}{2}\|\vz_i^t-\vz_i^{t-1}\|^{2}_{\mB_i^T\mB_i}+\frac{1}{2}\|\mD_i\left(\vu_i^t-\vu_i^{t-1}\right)\|^{2}-\frac{c}{2}\|\mA_i\left(\vz_i^{t+1}-\vz_i^t\right)\|^{2}\Big].
		\end{align}
		Hence, we will use the dual residual$\sum_{i\in\mathcal{N}}\|\mA_i\left(\vz_i^{t+1}-\vz_i^t\right)\|^{2}$ to bound $\sum_{i\in\mathcal{N}}\|\vz_i^t-\tilde{\vz}_i^{t+1}\|$. Recall the definition of $\vz_i$ in \eqref{z} and $\mA_i$ in \eqref{mA}, we have
		\begin{equation*}
		\mA_i\left(\vz_i^{t+1}-\vz_i^t\right)=\left[\vp_{i}^{t+1}-\vp_{i}^t-\left((\vz_{i,j}^{-})^{t+1}-(\vz_{i,j}^{-})^t\right)\right]_{j\in\mathcal{N}_i},
		\end{equation*}
		then it yields
		\begin{align}\label{1Azz}
		\vone_{N_in}^T\cdot\mA_i\left(\vz_i^{t+1}-\vz_i^t\right)=\sum_{j\in\mathcal{N}_i}\left[\vp_{i}^{t+1}-\vp_{i}^t-\left((\vz_{i,j}^{-})^{t+1}-(\vz_{i,j}^{-})^t\right)\right].
		\end{align}
		Substituting \eqref{ztt1}-\eqref{zt} into \eqref{1Azz}, rearranging the terms, and summing the above relation over $i\in\mathcal{N}$, we obtain
		\begin{align}\nonumber
		\sum_{i\in\mathcal{N}}\vone_{N_in}^T\cdot\mA_i\left(\vz_i^{t+1}-\vz_i^t\right)=
		&\sum_{i\in\mathcal{N}}\sum_{j\in\mathcal{N}_i}\Big[\frac{c}{c+1}\left(\tilde{\vp}_{i}^{t+1}-\vp_{i}^t-\left((\tilde{\vz}_{i,j}^{-})^{t+1}-(\vz_{i,j}^{-})^t\right)\right)\\\label{1Az}
		&\qquad\qquad+\frac{1}{c+1}\left(\tilde{\vp}_{i}^{t+1}-\vp_{i}^t-\left((\tilde{\vz}_{j,i}^{+})^{t+1}-(\vz_{i,j}^{-})^t\right)\right)\Big].
		\end{align}
		For the first row on the rhs of the \eqref{1Az}, using \eqref{mA} again, we have
		\begin{align}\label{2Az}
		\sum_{i\in\mathcal{N}}\sum_{j\in\mathcal{N}_i}\frac{c}{c+1}\left(\tilde{\vp}_{i}^{t+1}-\vp_{i}^t-\left((\tilde{\vz}_{i,j}^{-})^{t+1}-(\vz_{i,j}^{-})^t\right)\right)=
		\sum_{i\in\mathcal{N}}\vone_{N_in}^T\cdot\frac{c}{c+1}\mA_i\left(\tilde{\vz}_i^{t+1}-\vz_i^t\right).
		\end{align}
		For the second row on the rhs of the \eqref{1Az}, since $\vz_i^{t+1}$ is an optimal solution of \eqref{lem_zt1} for all $t\geq1$, we have $\vz_i^t\in\mathcal{Z}$, then we can deduce that
		\begin{equation}\label{z_const}
			(\vz_{i,j}^{-})^t=(\vz_{j,i}^{+})^t,~\forall~i\in\mathcal{N},~j\in\mathcal{N}_i.
		\end{equation}
		Hence, together with \eqref{z_const}, the second row on the rhs of the \eqref{1Az} can be expressed as
		\begin{align}\nonumber
		&\sum_{i\in\mathcal{N}}\sum_{j\in\mathcal{N}_i}\frac{1}{c+1}\left(\tilde{\vp}_{i}^{t+1}-\vp_{i}^t-\left((\tilde{\vz}_{j,i}^{+})^{t+1}-(\vz_{i,j}^{-})^t\right)\right)\\\nonumber=
		&\sum_{i\in\mathcal{N}}\sum_{j\in\mathcal{N}_i}\frac{1}{c+1}\left(\tilde{\vp}_{i}^{t+1}-\vp_{i}^t-\left((\tilde{\vz}_{j,i}^{+})^{t+1}-(\vz_{j,i}^{+})^t\right)\right)\\\label{3Az}=
		&\sum_{i\in\mathcal{N}}\vone_{N_in}^T\cdot\frac{1}{c+1}\mQ_i\left(\tilde{\vz}_i^{t+1}-\vz_i^t\right),
		\end{align}
		where the last step is due to \eqref{mQ} and rearranges the terms. Substituting \eqref{2Az} and \eqref{3Az} into \eqref{1Az}, we can get
		\begin{align}\label{4Az}
		\Big\|\sum_{i\in\mathcal{N}}\vone_{N_in}^T\cdot\mA_i\left(\vz_i^{t+1}-\vz_i^t\right)\Big\|^{2}=\frac{1}{\left(c+1\right)^{2}}\Big\|\sum_{i\in\mathcal{N}}\vone_{N_in}^T\cdot\left[c\mA_i+\mQ_i\right]\left(\tilde{\vz}_i^{t+1}-\vz_i^t\right)\Big\|^{2}.
		\end{align}
		Let $\mA:=\text{diag}\left(\mA_i,i\in\mathcal{N}\right)$ denote the block diagonal matrix whose diagonal coefficients correspond to the $\mA_i,i\in\mathcal{N}$ and $\mQ:=\text{diag}\left(\mQ_i,i\in\mathcal{N}\right)$, $\tilde{\vz}^{t+1}:=\text{vec}\left(\tilde{\vz}_i^{t+1},i\in\mathcal{N}\right)$, then \eqref{4Az} can be rewritten as a compact form:
		\begin{align}\label{enA}
		\|\vone_{N_{\text{sum}}n}^T\cdot\mA\left(\vz^{t+1}-\vz^t\right)\|^{2}=\frac{1}{\left(c+1\right)^{2}}\|\vone_{N_{\text{sum}}n}^T\cdot\left[c\mA+\mQ\right]\left(\tilde{\vz}^{t+1}-\vz^t\right)\|^{2},
		\end{align}
		where $N_{\text{sum}}:=\sum_{i\in\mathcal{N}}N_i$. Due to the Cauchy–Schwarz inequality, we can upper bound the lhs of \eqref{enA}
		\begin{equation}\label{EAzz}	N_{\text{sum}}n\cdot\|\mA\left(\vz^{t+1}-\vz^t\right)\|^{2}\geq\|\vone_{N_{\text{sum}}n}^T\cdot\mA\left(\vz^{t+1}-\vz^t\right)\|^{2}.
		\end{equation}
		For the rhs of \eqref{enA}, we have
		\begin{equation}\label{lower_az}		\|\vone_{N_{\text{sum}}n}^T\cdot\left[c\mA+\mQ\right]\left(\tilde{\vz}^{t+1}-\vz^t\right)\|^{2}\geq\tilde{\tau}_{\min}\|\tilde{\vz}^{t+1}-\vz^t\|^{2},
		\end{equation}
		where $\tilde{\tau}_{\min}$ is defined as the smallest non-zero eigenvalue of $\left[\vone_{N_{\text{sum}}n}^T\cdot\left(c\mA+\mQ\right)\right]^T\left[\vone_{N_{\text{sum}}n}^T\cdot\left(c\mA+\mQ\right)\right]$.
		From the definition of $\mA_i$ and $\mQ_i$ in \eqref{mA} and \eqref{mQ}, it derives $\tilde{\tau}_{\min}=\min\{\left(c+1\right)^{2}N_i^{2}+c^{2}N_i+N_i,i\in\mathcal{N}\}$.
		Substituting \eqref{EAzz} and \eqref{lower_az} into \eqref{enA}, we obtain
		\begin{equation}\label{Azz1}	\frac{c}{2}\|\mA\left(\vz^{t+1}-\vz^t\right)\|^{2}\geq\frac{c\cdot\tilde{\tau}_{\min}}{2N_{\text{sum}}n\left(c+1\right)^{2}}\|\tilde{\vz}^{t+1}-\vz^t\|^{2},
		\end{equation}
		which can be rewritten as a summation over network nodes:
		\begin{equation}\label{Azz}
		\sum_{i\in\mathcal{N}}\frac{c}{2}\|\mA_i\left(\vz^{t+1}_i-\vz_i^t\right)\|^{2}\geq\sum_{i\in\mathcal{N}}\frac{c\cdot\tilde{\tau}_{\min}}{2N_{\text{sum}}n\left(c+1\right)^{2}}\|\tilde{\vz}_i^{t+1}-\vz_i^t\|^{2}.
		\end{equation}	
		Finally, we substitute \eqref{Azz} to \eqref{cazb}, and the proof is complete.
	\end{IEEEproof}	
	Using Lemma \ref{LL}, Lemma \ref{lem2-2}, Lemma \ref{lem2-4} and the deﬁnition of the potential function $\varsigma^t$ in \eqref{poten}, we can get the desired result of \eqref{lem_var}.

	\section{Proof of the Lemma \ref{KKT_Critical}}\label{proof_KKT_Critical}
	If $\mathcal{F}\left(\vz^{t},\vu^{t},\vlambda^{t}\right)=0$, it implies that, for any node $i\in\mathcal{N}$,
	\begin{align}
	\label{z_proj}\vz_i^t-\text{proj}_{\mathcal{X},\mathcal{Z}}\left(\vz_i^t-\left(\nabla_{\vz_i}
	F_i\left(\vz_i^t,\vu_i^t\right)+\mA_i^T\vlambda_i^t\right)\right)&=0,\\
	\label{Az}\mA_i\vz_i^t&=0,\\
	\label{u_proj1}\vu_i^t-\vu_i^{t-1}&=0.
	\end{align}
	According to the definition of the projection operator, we have
	\begin{align}\label{optimal}
	\text{proj}_{\mathcal{X},\mathcal{Z}}\left(\vz_i^t-\left(\nabla_{\vz_i}
	F_i\left(\vz_i^t,\vu_i^t\right)+\mA_i^T\vlambda_i^t\right)\right)
	=\mathop{\arg\min}\limits_{\substack{\vz\in\mathcal{Z}\\\vz\in\mathcal{X}}}\left\|\vz_i-\vz_i^t+\left(\nabla_{\vz_i}
	F_i\left(\vz_i^t,\vu_i^t\right)+\mA_i^T\vlambda_i^t\right)\right\|^2,
	\end{align}
	Using \eqref{z_proj}, we deduce that $\vz_i^t$ is an optimal solution for problem \eqref{optimal}. Therefore, based on the first order necessary condition of \eqref{optimal}, we obtain
	\begin{align}\label{first_necessary}
	\langle\nabla_{\vz_i}
	F_i\left(\vz_i^t,\vu_i^t\right)+\mA_i^T\vlambda_i^t,\vx-\vz_i^t\rangle\geq\vzero,~\forall~\vx\in\mathcal{X},\vx\in\mathcal{Z}.
	\end{align}
	Since $F_i\left(\vz_i,\vu_i\right)+\langle\vlambda_i,\mA_i\vz_i\rangle$ is convex over $\mathcal{X}$ and $\mathcal{Z}$ for any fixed $\vu_i$, then the above expression \eqref{first_necessary} implies that $\vz_i^t$ is also sufficient for $\vz_i^t$ to minimize $F_i\left(\vz_i,\vu_i\right)+\langle\vlambda_i,\mA_i\vz_i\rangle$ over
	$\mathcal{X}$ and $\mathcal{Z}$, i.e.,
	\begin{align*} \vz_i^t\in\mathop{\arg\min}\limits_{\substack{\vz\in\mathcal{Z}\\\vz\in\mathcal{X}}}F_i\left(\vz_i,\vu_i^t\right)+\langle\vlambda_i^t,\mA_i\vz_i\rangle.
	\end{align*}
	Thanks to the separable property of the objective function about node $i$, we know
	$\vz^t$ shall satisfy the KKT condition
	\begin{align}\label{KKT_z}
	\vz^t\in\mathop{\arg\min}\limits_{\substack{\vz\in\mathcal{Z}\\\vz\in\mathcal{X}}}\sum_{i\in\mathcal{N}}F_i\left(\vz_i,\vu_i^t\right)+\langle\vlambda_i^t,\mA_i\vz_i\rangle.
	\end{align}
	Similarly, the primal variable $\vu_i$ is updated as follows
	\begin{align*}
	\vu_i^t = \text{proj}_{\mathcal{B}^{N_i}}(\vu_i^{t-1}-\nabla_{\vu_i}F_i\left(\vz_i^t,\vu_i^t\right)).
	\end{align*}
	Combing it with \eqref{u_proj1},  we have
	\begin{align*}
	0=\|\vu_i^{t-1}-\vu_i^t\|^2=\|\vu_i^{t-1}-\text{proj}_{\mathcal{B}^{N_i}}(\vu_i^{t-1}-\nabla_{\vu_i}F_i\left(\vz_i^t,\vu_i^t\right))\|^2,
	\end{align*}
	Therefore $\vu_{i}^{t-1}$ is the optimal solution of the following problem
	\begin{align}\label{proj_u}
	\min_{\vu_i\in\mathcal{B}^{N_i}}~\|\vu_i-\vu_i^{t-1}+\nabla_{\vu_i}F_i\left(\vz_i^t,\vu_i^t\right)\|^{2}.
	\end{align}
	Based on the optimality condition of problem \eqref{proj_u} and equation \eqref{u_proj1}, we have
	\begin{align}
	0\leq\langle\nabla_{\vu_i}F_i\left(\vz_i^t,\vu_i^t\right),\vu_i-\vu_i^t+\vu_i^t-\vu_i^{t-1}\rangle
	=\langle\nabla_{\vu_i}F_i\left(\vz_i^t,\vu_i^t\right),\vu_i-\vu_i^t\rangle,~\forall~\vu_i\in\mathcal{B}^{N_i}.
	\end{align}
	Since $F_i\left(\vz_i,\vu_i\right)$ is convex for any fixed $\vz_i$, it implies that
	\begin{align}\label{uF}
	\vu_i^t\in\mathop{\arg\min}\limits_{\vu_i\in\mathcal{B}^{N_i}}~F_i\left(\vz_i^t,\vu_i^t\right),
	\end{align}
	which is also written as
	\begin{align}\label{KKT_u}
	\vzero\in\nabla_{\vu_i}F_i\left(\vz_i^t,\vu_i^t\right)+\partial\delta_{\mathcal{B}^{N_i}}(\vu_{i}^t).
	\end{align}
	Combing \eqref{KKT_z}, \eqref{KKT_u} and \eqref{Az}, one can observe that $\left(\vz^{t},\vu^{t},\boldsymbol{\lambda}^{t}\right)$ satisfy the KKT condition \eqref{KKT_cond}.
	
	Next, we show that if $(\vz^t,\vu^t,\vlambda^t)$ is a KKT solution of Problem \eqref{f}, which satisfies \eqref{KKT_cond}, then $(\vp^t,\vu^t)$ is a critical point of the nonconvex problem \eqref{11} satisfying:
	\begin{subequations}\label{Critical1}
		\begin{align}
		\label{Critical_z1}
		\vzero&=\sum_{j\in\mathcal{N}_i}\left(\vp_i^t-\vp_j^t-d_{i,j}\vu_{i,j}^t+\vp_i^t-\vp_j^t+d_{j,i}\vu_{j,i}^t\right),~\forall~i\in(\mathcal{N}/\mathcal{A}),\\\label{Critical_u11}\vzero&\in-d_{i,j}\left(\vp_i^t-\vp_j^t\right)+\partial\delta_{\mathcal{B}}\left(\vu_{i,j}^t\right),~\forall~i\in\mathcal{N},~j\in\mathcal{N}_i,\\\label{Critical_u21}\vp_j^t&=\va_j,~\forall~j\in\mathcal{A}.
		\end{align}
	\end{subequations}

	Firstly, since $\vz^t$ is an optimal solution of \eqref{KKT_cond1}, we have $\vz^t\in\mathcal{X}$, thus \eqref{Critical_u21} is proved.
	
	Secondly, we prove $(\vp^t,\vu^t)$ satisfies \eqref{Critical_u11}. Since $\vz^t$ is an optimal solution of \eqref{KKT_cond1}, we have $\vz^t\in\mathcal{Z}$. By combining this with the KKT condition \eqref{KKT_cond3}, we can deduce that
	\begin{align}\label{zij+zij-}
	(\vz_{i,j}^{+})^t=\vp_j^t=(\vz_{j,i}^-)^t,~\forall~i\in\mathcal{N},~\forall~j\in\mathcal{N}_i.
	\end{align}
	By using the definition of $F_i(\vz_i^t,\vu^t_i)$ and \eqref{zij+zij-}, we obtain
	\begin{align}\label{nabuF}
	\nabla_{\vu_i}F_i(\vz_i^t,\vu_i^t)=-\mD_i\mQ_i\vz_i^t=\text{vec}\left(-d_{i,j}(\vp_i^t-(\vz_{i,j}^+)^t)~j\in\mathcal{N}_i\right)=\text{vec}(-d_{i,j}(\vp_i^t-\vp_j^t),~j\in\mathcal{N}_i),
	\end{align}
	where
	\begin{align}\nonumber
	\mD_i&=\textbf{D}\text{iag}\left(\text{vec}\left(d_{i,j},j\in\mathcal{N}_i\right)\right)\otimes\mI_{n},~\mQ_i=\left[\vone_{N_i},\mO_{N_i},-\mI_{N_i}\right]\otimes\mI_{n},\\\label{DQz}\vz_i^t&=\left[(\vp_i^t)^T,\text{vec}\left((\vz_{i,j}^-)^t,j\in\mathcal{N}_i\right)^T,\text{vec}\left((\vz_{i,j}^+)^t,j\in\mathcal{N}_i\right)^T\right]^T.
	\end{align}
	By combining \eqref{nabuF} with KKT condition \eqref{KKT_cond2}, we have that $(\vz_i^t,\vu_{i,j}^t)$ satisfies the critical condition \eqref{Critical_u11}.
	
	Finally, let us prove \eqref{Critical_z1}. Recall that
	\begin{align}\label{FXZ}
	F_i\left(\vz_i,\vu_i\right)&=\frac{1}{2}\left\|\mQ_i\vz_i\right\|^{2}-\vu_i^T\mD_i\mQ_i\vz_i,\\\nonumber\mathcal{X}&=\{\vz\mid\mE_i\vz_i=\va_i,\forall i\in\mathcal{A}\},~\mE_i=\left[1,\vzero^T_{N_i},\vzero^T_{N_i}\right]\otimes\mI_{n},~\mathcal{Z}=\{\vz_{i,j}^{+}=\vp_j=\vz_{j,i}^{-},~\forall~i\in\mathcal{N},j\in\mathcal{N}_i\}.
	\end{align}
	Substituting \eqref{DQz} and \eqref{FXZ} into \eqref{KKT_cond1}, then the optimization problem \eqref{KKT_cond1} is rewritten as
	\begin{align}\nonumber
	&\mathop{\arg\min}\limits_{\substack{\vz\in\mathcal{Z}\\\vz\in\mathcal{X}}}\sum_{i\in\mathcal{N}}F_i\left(\vz_i,\vu_i\right)+\langle\vlambda_i,\mA_i\vz_i\rangle\\\nonumber=&\mathop{\arg\min}\limits_{\substack{\vz_{i,j}^+=\vz_{j,i}^-\\\vp_i=\va_i,i\in\mathcal{A}}}~\sum_{i\in\mathcal{N}}\sum_{j\in\mathcal{N}_i}\left(\frac{1}{2}\|\vp_i-\vz_{i,j}^+\|^2-d_{i,j}\vu_{i,j}^T(\vp_i-\vz_{i,j}^+)+\vlambda_{i,j}^T(\vp_i-\vz_{i,j}^-)\right)\\\nonumber=&\mathop{\arg\min}\limits_{\vp_i,\vz_{i,j}^-}~\sum_{i\in(\mathcal{N}/\mathcal{A})}\sum_{j\in\mathcal{N}_i}\left(\frac{1}{2}\|\vp_i-\vz_{j,i}^-\|^2-d_{i,j}\vu_{i,j}^T(\vp_i-\vz_{j,i}^-)+\vlambda_{i,j}^T(\vp_i-\vz_{i,j}^-)\right)\\\label{equ}&\qquad\qquad+\sum_{i\in\mathcal{A}}\sum_{j\in\mathcal{N}_i}\left(\frac{1}{2}\|\va_i-\vz_{j,i}^-\|^2-d_{i,j}\vu_{i,j}^T(\va_i-\vz_{j,i}^-)+\vlambda_{i,j}^T(\va_i-\vz_{i,j}^-)\right).
	\end{align}
	where the last equality holds by substituting the constraints $\vz_{i,j}^+=\vz_{j,i}^-,\forall i\in\mathcal{N},j\in\mathcal{N}_i$ and $\vp_i=\va_i,i\in\mathcal{A}$ into the objective function of the first equality. Using the KKT condition \eqref{KKT_cond1}, we know that $\vz^t$ is an optimal solution of problem \eqref{equ}. By using the optimality condition of \eqref{equ} about variables $\vp_i$ and $\vz_{i,j}^-$, we have that
	\begin{align*}
	\vzero=&\sum_{j\in\mathcal{N}_i}\left(\vp_i^t-(\vz_{j,i}^-)^t-d_{i,j}\vu_{i,j}^t+\vlambda_{i,j}^t\right),~\forall~i\in(\mathcal{N}/\mathcal{A}),\\\vzero=&\sum_{j\in\mathcal{N}_i}\left((\vz_{i,j}^-)^t-\vp_j^t+d_{j,i}\vu_{j,i}^t-\vlambda_{i,j}^t\right),~\forall~i\in(\mathcal{N}/\mathcal{A}),
	\end{align*}
	Adding the above equalities together, we get
	\begin{align}\label{Criz}
	\vzero=\sum_{j\in\mathcal{N}_i}\left(\vp_i^t-(\vz_{j,i}^-)^t-d_{i,j}\vu_{i,j}^t+(\vz_{i,j}^-)^t-\vp_j^t+d_{j,i}\vu_{j,i}^t\right),~\forall~i\in(\mathcal{N}/\mathcal{A}),
	\end{align}
	By using the KKT condition \eqref{KKT_cond3}, we have $(\vz_{i,j}^-)^t=\vp_i^t$, $(\vz_{j,i}^-)^t=\vp_j^t$. Then \eqref{Criz} can be rewritten as
	\begin{align}\label{Criz1}
	\vzero=\sum_{j\in\mathcal{N}_i}\left(\vp_i^t-\vp_j^t-d_{i,j}\vu_{i,j}^t+\vp_i^t-\vp_j^t+d_{j,i}\vu_{j,i}^t\right),~\forall~i\in(\mathcal{N}/\mathcal{A}).
	\end{align}
	which shows that $\vz^t$ satisfy \eqref{Critical_z1}.
	
	We have thus established that if $(\vz^t,\vu^t,\vlambda^t)$ is a KKT solution of problem \eqref{f}, which satisfies \eqref{KKT_cond}, then $(\vp^t,\vu^t)$ is a critical point of problem \eqref{11}. {\blue Furthermore, we can demonstrate below that $\vp^t$ also qualifies as a critical point of the original problem \eqref{1}.
		
	Recall that the function $F_i\left(\vz_i,\vu_i\right)$ is deﬁned as
		\begin{align}\label{FF}
		F_i\left(\vz_i,\vu_i\right)=\frac{1}{2}\left\|\mQ_i\vz_i\right\|^{2}-\vu_i^T\mD_i\mQ_i\vz_i=\sum_{j\in\mathcal{N}_i}\frac{1}{2}\Big[\|\vp_i-\vp_j\|^{2}-d_{i,j}\vu_{i,j}^{T}\left(\vp_i-\vp_j\right)\Big],~\forall~i\in\mathcal{N}.
		\end{align}
		Therefore, combining \eqref{uF} with \eqref{FF} yields that $\vu_{i,j}^{t}$ is the optimal solution of the following problem
		\begin{align*}
		\mathop{\min}\limits_{\vu_{i,j}\in\mathcal{B}}-d_{i,j}\vu_{i,j}^{T}\left(\vp_i^t-\vp_j^t\right).
		\end{align*}
		From this, we can derive that
		\begin{align}\label{uuU}
		\vu_{j,i}^{t}=-\vu_{i,j}^{t},\,\forall~i\in\mathcal{N},j\in\mathcal{N}_i.
		\end{align}
		Using \eqref{uuU} and \eqref{Critical1}, we obtain that
		\begin{align*}
		\vp_i^t-\vp_j^t\in\partial\delta_{\mathcal{B}}\left(\vu_{i,j}^t\right)~\text{for all}~i\in\mathcal{N},~j\in\mathcal{N}_i.
		\end{align*}
		Since $\delta_{\mathcal{B}}\left(\cdot\right)$ is a proper, lower semicontinuous, and convex function,  it follows from [\cite{beck2017first} Theorem 4.20, p. 104] that
		\begin{align*}
		\vu_{i,j}^t\in\partial\delta^\ast_{\mathcal{B}}\left(\vp^t_i-\vp_j^t\right),~\text{for all}~i\in\mathcal{N},~j\in\mathcal{N}_i.
		\end{align*}
		Here, $\delta^\ast_{\mathcal{B}}$ represents the Fenchel conjugate of $\delta_{\mathcal{B}}$. From [\cite{beck2017first}, Example 2.31, p. 28], we also find that $\delta^\ast_{\mathcal{B}}\left(\cdot\right)=\|\cdot\|$. Hence $\vu_{i,j}^t\in\partial\|\cdot\|(\vp_i^t-\vp_j^t)$ for all $i\in\mathcal{N}$ and $j\in\mathcal{N}_i$. Substituting this fact, along with \eqref{uuU} and \eqref{Critical_u21}, into \eqref{Critical_z1}, we can get
		\begin{subequations}\label{cri_ori}
			\begin{align}
			\vzero&=\sum_{j\in\mathcal{N}_i}(\vp_i^t-\vp_j^t)-d_{i,j}\vu_{i,j}^t,\quad\forall i\in(\mathcal{N}/\mathcal{A}),\\\vu_{i,j}^t&\in\partial\|\cdot\|(\vp_i^t-\vp_j^t),\quad\forall i\in\mathcal{N},j\in\mathcal{N}_i,\\\vp_j^t&=\va_j,~\forall~j\in\mathcal{A}.
			\end{align}
		\end{subequations}
		Therefore, from \eqref{cri_ori} we know that $\vp^t$ is a critical point of the original nonconvex problem \eqref{1}.}

	\section{Proof of the Lemma \ref{sublinear-lem}}\label{proof_sublinear-lem}
	According to equations \eqref{k1_all}-\eqref{c_all}, it can be observed that the value of $\varsigma^t$ decreases after each iteration. Combining this observation with the first part of Theorem \ref{them2}, it can be concluded that for any $t>1$, there exists an $s\in{1,2,\ldots,t-1}$ such that $\varsigma^s-\varsigma^{s+1}\leq\frac{\varsigma^0-\underline{\varsigma}}{t-1}$. By utilizing equation \eqref{lem_var}, we can let
	\begin{align*}
	K=(\varsigma^0-\underline{\varsigma})\cdot\max&\Big\{\frac{1}{2},\frac{\kappa_1-1}{2},\frac{c(\kappa_1-1)}{2},{\blue\frac{\rho}{4}},\frac{c(\kappa_1-6(N_{\max}+1))}{2},\frac{c\kappa_1-6(N_{\max}+1)}{2c},\\&\quad\frac{c\kappa_{2}\cdot\tilde{\tau}_{\min}}{2N_{\text{sum}}n\left(c+1\right)^{2}}-\frac{\left(N_{\max}+1\right)c\kappa_{1}}{2},\frac{\rho}{{\blue4}}-d_{\max}^2(\kappa_1+\kappa_2)\Big\}.
	\end{align*}
	and we obtain
	\begin{align}\nonumber
	&\|\vz^{s+1}-\vz_i^s\|^{2}<\frac{K}{t-1},\quad\|\vu_{i}^{{\blue s+1}}-\vu_i^{{\blue s}}\|^{2}<\frac{K}{t-1},\quad\|\tilde{\vz}_i^{s+1}-\vz_i^{s+1}-\left(\tilde{\vz}_i^s-\vz_i^s\right)\|^{2}_{\mA_i^T\mA_i}<\frac{K}{t-1},\\\label{bound}&\|\mQ_i(\tilde{\vz}_i^{s+1}-\tilde{\vz}_i^s)-\mD_i\left(\vu_i^s-\vu_i^{s-1}\right)\|^{2}<\frac{K}{t-1},\quad\|\tilde{\vz}_i^{s+1}-\vz_i^s-(\tilde{\vz}_i^s-\vz_i^{s-1})\|^{2}_{\mB_i^T\mB_i}<\frac{K}{t-1}.
	\end{align}
	From problem \eqref{z_up}, we have
	\begin{align*}
	\vz_i^{s+1}=\mathop{\arg\min}\limits_{\vz}F_i\left(\vz_i,\vu_i^s\right)+\langle\vlambda_i^s,\mA_i\vz_i\rangle+\frac{c}{2}\|\vz_i-\vz_i^s\|^2_{\mB_i^T\mB_i}+\delta_{\mathcal{X}\times \mathcal{Z} }(\vz_i).
	\end{align*}
	Then by applying the optimality condition, we can get
	\begin{align*}
	\vzero\in\nabla_{\vz_i}F_i\left(\vz_i^{s+1},\vu_i^s\right)+\mA_i^T\vlambda_i^s+c\mB_i^T\mB_i(\vz_i^{s+1}-\vz_i^s)+\partial\delta_{\mathcal{X}\times \mathcal{Z}}(\vz_i^{s+1}).
	\end{align*}
	Letting
	\begin{align*}
	\vv_1=\nabla_{\vz_i}F_i\left(\vz_i^{s+1},\vu_i^{s+1}\right)+\mA_i^T\vlambda_i^{s+1}-\nabla_{\vz_i}F_i\left(\vz_i^{s+1},\vu_i^s\right)-\mA_i^T\vlambda_i^s-c\mB_i^T\mB_i(\vz_i^{s+1}-\vz_i^s),
	\end{align*}
	then we have
	\begin{align}\nonumber
	\|\vv_1\|=&\|-\mQ_i^T\mD_i(\vu_i^{s+1}-\vu_i^{s})+\mA_i^T(\vlambda_i^{s+1}-\vlambda_i^s)+c\mB_i^T\mB_i(\vz_i^{s+1}-\vz_i^{s})\|\\\label{v_norm}\leq&d_{\max}\sqrt{N_{i}+1}\|\vu_i^{s+1}-\vu_i^{s}\|+\sqrt{N_{i}+1}\|\vlambda_i^{s+1}-\vlambda_i^s\|+(1+c)(1+N_{\max})\|\vz_i^{s+1}-\vz_i^s\|,
	\end{align}
	where the last inequality holds due to $\|\mQ_i^T\mQ_i\|=\|\mA_i^T\mA_i\|=N_i+1$ and \eqref{evib}. Substituting \eqref{bound} and \eqref{alamd1} into \eqref{v_norm}, we have
	\begin{align*}
	\|\vv_1\|&\leq\frac{\sqrt{K}}{\sqrt{t-1}}\Big(d_{\max}\sqrt{N_i+1}+\sqrt{N_i+1}(\sqrt{N_i+1}+c\sqrt{N_i+1}+(1+c)(1+N_{\max}))\\&\qquad\qquad\quad+(1+c)(1+N_{\max})\Big)\\&\leq\frac{\sqrt{M_1}}{\sqrt{t-1}},
	\end{align*}
	where $M_1:=K(d_{\max}\sqrt{N_i+1}+\sqrt{N_i+1}(\sqrt{N_i+1}+c\sqrt{N_i+1}+(1+c)(1+N_{\max}))+(1+c)(1+N_{\max}))^2>0$.
	Recall the $\vu$ update step, we have
	\begin{align*}
	\vu_i^{s+1}=\mathop{\arg\min}\limits_{\vu}F_i(\vz_i^{s+1},\vu_i)+\frac{\rho}{2}\|\vu_i-\vu_i^s\|^{2}+\delta_{\mathcal{B}^{N_i}}(\vu_i),
	\end{align*}
	By the optimality condition,
	\begin{align*}
	\vzero\in\nabla_{\vu_i}F_i(\vz_i^{s+1},\vu_i^{s+1})+\rho(\vu_i^{s+1}-\vu_i^s)+\partial\delta_{\mathcal{B}^{N_i}}(\vu_i^{s+1}).
	\end{align*}
	Letting
	\begin{align}
	\vv_2=\nabla_{\vu_i}F_i(\vz_i^{s+1},\vu_i^{s})-\nabla_{\vu_i}F_i(\vz_i^{s+1},\vu_i^{s+1})-\rho(\vu_i^{s+1}-\vu_i^s)=-\rho(\vu_i^{s+1}-\vu_i^s),
	\end{align}
	then we have
	\begin{align}
	\vv_2\in\nabla_{\vu_i}F_i(\vz_i^{s+1},\vu_i^{s+1})+\partial\delta_{\mathcal{B}^{N_i}}(\vu_i^{s+1}).
	\end{align}
	By the inequality \eqref{bound}, we get
	\begin{align}
	\|\vv_2\|\leq\frac{\sqrt{K}}{\sqrt{t-1}}\cdot\rho=\frac{M_2}{\sqrt{t}},
	\end{align}
	where $M_2:=K\rho^2>0$. Then the Lemma 1 holds for $\vv_1,\vv_2$ and $M:=\max\{M_1,M_2\}$ and $(\vz^t,\vu^t,\vlambda^t)$ is a $M/\sqrt{t-1}$-solution. It indicates that the proposed algorithm converges to a KKT stationary point of problem \eqref{f} at a sublinear rate.

	\section{Proof of the Theorem \ref{them2}}\label{proof_them}
	
	\noindent
	\begin{Lemma}\label{var_low_bound}
		Suppose $c,\rho,\kappa_{1},\kappa_{2}$ are chosen according to \eqref{k1_all}-\eqref{c_all}. Then the following statement holds true
		\begin{equation*}
		\exists\,\,\underline{\varsigma}>-\infty\,\,s.t.\,\,\varsigma^t\geq\underline{\varsigma},\,\,\forall t>0.
		\end{equation*}
	\end{Lemma}
	
	\begin{IEEEproof}
		By using the update step \eqref{du2} of the dual variable $\vlambda_i^{t+1}$, we have
		\begin{align}\nonumber
		&\mathcal{L}(\vz^{t+1},\boldsymbol{\lambda}^{t+1},\vu^{t+1})\\\nonumber=
		&\sum_{i\in\mathcal{N}}F_i(\vz_i^{t+1},\vu_i^{t+1})+\langle\boldsymbol{\lambda}_i^{t+1},\mA_i\vz_i^{t+1}\rangle
		+\frac{c}{2}\|\mA_i\vz_i^{t+1}\|^{2}\\\nonumber=
		&\sum_{i\in\mathcal{N}}F_i(\vz_i^{t+1},\vu_i^{t+1})+\frac{1}{c}\langle\boldsymbol{\lambda}_i^{t+1},\boldsymbol{\lambda}_i^{t+1}-\boldsymbol{\lambda}_i^t\rangle+\frac{c}{2}\|\mA_i\vz_i^{t+1}\|^{2}\\\label{lzul}=
		&\sum_{i\in\mathcal{N}}\Big[F_i(\vz_i^{t+1},\vu_i^{t+1})+\frac{1}{2c}(\|\boldsymbol{\lambda}_i^{t+1}\|^{2}-\|\boldsymbol{\lambda}_i^t\|^{2}+\|\boldsymbol{\lambda}_i^{t+1}-\boldsymbol{\lambda}_i^t\|^{2})+\frac{c}{2}\|\mA_i\vz_i^{t+1}\|^{2}\Big].
		\end{align}
		Adding \eqref{lzul} from $t=1$ to $T$, we obtain
		\begin{align}\nonumber
		&\sum_{t=1}^T\mathcal{L}(\vz^{t+1},\boldsymbol{\lambda}^{t+1},\vu^{t+1})\\\label{Llowbound}=
		&\sum_{t=1}^T\sum_{i\in\mathcal{N}}\left[F_i\left(\vz_i^{t+1},\vu_i^{t+1}\right)+\frac{c}{2}\|\mA_i\vz_i^{t+1}\|^{2}+\frac{1}{2c}\left(\|\boldsymbol{\lambda}_i^{t+1}\|^{2}-\|\boldsymbol{\lambda}_i^t\|^{2}+\|\boldsymbol{\lambda}_i^{t+1}-\boldsymbol{\lambda}_i^t\|^{2}\right)\right].
		\end{align}
		From the deﬁnition of $F_i$ it follows that
		\begin{align}\nonumber
		\sum_{i\in\mathcal{N}}F_i\left(\vz_i^{t+1},\vu_i^{t+1}\right)=
		&\sum_{i\in\mathcal{N}}\frac{1}{2}\left\|\mQ_i\vz_i^{t+1}\right\|^{2}-\left(\vu_i^{t+1}\right)^t\mD_i\mQ_i\vz_i^{t+1}\\\nonumber=
		&\sum_{i\in\mathcal{N}}\frac{1}{2}\left\|\mQ_i\vz_i^{t+1}-\mD_i\vu_i^{t+1}\right\|^{2}-\frac{1}{2}\left\|\mD_i\vu_i^{t+1}\right\|^{2}\\\label{Flowbound}\geq
		&\sum_{i\in\mathcal{N}}-\frac{d_{\max}^{2}}{2}\left\|\vu_i^{t+1}\right\|^{2}>-\infty,\,\forall t>0,
		\end{align}
		where the last inequality is due to $\vu_i^t\in\mathcal{B}^{N_i},\forall i\in\mathcal{N},t>0$. Therefore, according to the definition \eqref{poten} of potential function, \eqref{Llowbound} and \eqref{Flowbound}, it shows that
		\begin{equation}\label{3}
		\sum_{t=1}^T\varsigma^t>-\infty,\forall T>0.
		\end{equation}
		When $c,\rho,\kappa_{1},\kappa_{2}$ satisfy \eqref{k1_all}-\eqref{c_all}, Lemma \ref{lem1-2} show that the potential function $\varsigma^t$ decreases at eatch iteration of Algorithm \ref{al2}. Hence, by \eqref{3} we can infer that $\varsigma^t\geq\underline{\varsigma}>-\infty$, $\forall t>0$.
	\end{IEEEproof}

	\begin{Lemma}\label{lem2-6}
		Suppose the parameters $c,\rho,\kappa_{1}$ and $\kappa_{2}$ satisfy \eqref{k1_all}-\eqref{c_all}, then the iterative sequence $\{\left(\vz_i^t,\vu_i^t,\vlambda_i^t\right)\}$ of Algorithm \ref{al2} satisfies
		\begin{align*}
		&\lim_{t\to\infty}\vz_i^{t+1}-\vz_i^t\rightarrow\mathbf{0},\,\,\lim_{t\to\infty}\mA_i\vz_i^t\rightarrow\mathbf{0},\\
		&\lim_{t\to\infty}\boldsymbol{\lambda}_i^{t+1}-\boldsymbol{\lambda}_i^t\rightarrow\mathbf{0},\lim_{t\to\infty}\vu_i^{t+1}-\vu_i^t\rightarrow\mathbf{0},\forall i\in\mathcal{N}.
		\end{align*}
	\end{Lemma}
	
	\begin{IEEEproof}
		Since $\mW_i$ is a positive definite diagonal matrix, combining \eqref{lem_var} and Lemma \ref{var_low_bound}, we have
		\begin{equation}
		\lim_{t\to\infty}\vz_i^{t+1}-\vz_i^t\rightarrow\vzero.
		\end{equation}	
		Using the update step of $\vlambda_i^{t+1}$ (cf.\eqref{du2}), we have from \eqref{apa_lem2}, \eqref{alamd1} and Lemma \ref{var_low_bound} that
		\begin{equation}\label{az0}
		\lim_{t\to\infty}\mA_i\vz_i^{t+1}\rightarrow\vzero,\quad\lim_{t\to\infty}\vlambda_i^{t+1}-\vlambda_i^t\rightarrow\mathbf{0}.
		\end{equation}
		By the inequality \eqref{lem_var} and Lemma \ref{var_low_bound}, we further obtain
		\begin{equation}
		\lim_{t\to\infty}\vu_i^{t+1}-\vu_i^t\rightarrow\mathbf{0}.
		\end{equation}
	\end{IEEEproof}
	
{\blue 
	\begin{Lemma}[Sufficient decrease condition]\label{lem-decrease}
		Suppose the sequence $\{\left(\vz^t,\vu^t,\vlambda^t\right)\}_{t\geq1}$ is generated by Algorithm \ref{al2}, $c\mB_i^T\mB_i$ takes the form of \eqref{cbb}, and the conditions \eqref{k1_all}-\eqref{c_all} are satisfied. Then we have
		\begin{align*}
		\varsigma^{t+1}-\varsigma^t\leq\sum_{i\in\mathcal{N}}\Big[-\min\{c,1\}\cdot\|\vz_i^{t+1}-\vz_i^t\|^{2}-\frac{\rho}{4}\|\vu_i^{t+1}-\vu_i^{t}\|^{2}-C_0\|\vlambda_i^{t+1}-\vlambda_i^t\|^{2}-C_1\|\tilde{\vz}_i^{t+1}-\vz_i^t\|^{2}\Big],
		\end{align*}	
		where $C_0:=\frac{\min\{c\kappa_{1}-6\left(N_{\max}+1\right),\,c\left(c\kappa_{1}-6\left(1+c\right)\left(N_{\max}+1\right)\right)\}}{\max\{3(N_i+1),\,3c(1+c)(1+N_{\max})\}}$ and  $C_1:=\frac{c\kappa_{2}\cdot\tilde{\tau}_{\min}}{2N_{\text{sum}}n\left(c+1\right)^{2}}-\frac{\left(N_{\max}+1\right)c\kappa_{1}}{2}$. $d_{\max}:=\max\{d_{i,j},i\in\mathcal{N},j\in\mathcal{N}_i\}$, 
		$N_{\max}:=\max\{N_i,i\in\mathcal{N}\}$, $N_{\text{sum}}:=\sum_{i\in\mathcal{N}}N_i$ is the total number of neighboring nodes and $\tilde{\tau}_{\min}:=\min\{\left(c+1\right)^{2}N_i^{2}+c^{2}N_i+N_i,i\in\mathcal{N}\}$. 
	\end{Lemma}
	\begin{IEEEproof} Substituting \eqref{evib} into \eqref{alamd1}, we have
\begin{align}\nonumber
\|\vlambda_i^{t+1}-\vlambda_i^t\|^{2}\leq&\max\{3(N_i+1),\,3c(1+c)(1+N_{\max})\}\cdot\Big(\|\mQ_i\left(\tilde{\vz}_i^{t+1}-\tilde{\vz}_i^t\right)-\mD_i\left(\vu_i^t-\vu_i^{t-1}\right)\|^{2}\\\label{lambdalambda}&+\|\tilde{\vz}_i^{t+1}-\vz_i^{t+1}-\left(\tilde{\vz}_i^t-\vz_i^t\right)\|^{2}_{\mA_i^T\mA_i}+\|\tilde{\vz}_i^{t+1}-\vz_i^t-\left(\tilde{\vz}_i^t-\vz_i^{t-1}\right)\|^{2}_{\mB_i^T\mB_i}\Big),\forall~i\in\mathcal{N}.
\end{align}
Combining \eqref{lambdalambda}, \eqref{w} with \eqref{lem_var}, we can derive the result.
\end{IEEEproof}

	\begin{Lemma}[Bounded sequence]\label{lem-bounded} Suppose that conditions \eqref{k1_all}-\eqref{c_all} are satisfied, the graph obtained by the network is connected, and there is at least one anchor sensor. Then the sequence $\{\left(\vz^t,\vu^t,\vlambda^t,\tilde{\vz}^t\right)\}_{t\geq1}$ generated by Algorithm \ref{al2} is bounded.
	\end{Lemma}
\begin{IEEEproof}
		We establish the boundedness of $\{\vu^t\}_{t\geq1}, \{\vz^t\}_{t\geq1},\{\vlambda^t\}_{t\geq1}$, and $\{\tilde{\vz}^t\}_{t\geq1}$ sequentially:
	\begin{enumerate}
		\item It is evident from Algorithm \ref{al2} that  $\{\vu^t\}_{t\geq1}\subset\mathcal{B}^{NN_i}$, where $\mathcal{B}$ represents the unit ball constraint. Hence, the sequence $\{\vu^t\}_{t\geq1}$ is bounded.
		\item Using the update rule of $\vlambda_i^{t+1}$, we obtain
		\begin{align}\label{lambda}	c\sum_{t=1}^{\infty}\sum_{i\in\mathcal{N}}\|\mA_i\vz_i^{t+1}\|^2=\sum_{t=1}^{\infty}\sum_{i\in\mathcal{N}}\|\vlambda_i^{t+1}-\vlambda_i^t\|^{2}\leq\frac{1}{C_0}\sum_{t=1}^{\infty}\left(\varsigma^t-\varsigma^{t+1}\right)\leq \frac{1}{C_0}\left(\varsigma^1-\underline{\varsigma}\right)<\infty.
		\end{align}
		where the first inequality is due to Lemma \ref{lem-decrease} and the second inequality comes from Lemma \ref{var_low_bound}. Consequently, $\{\|\mA_i\vz_i^{t}\|\}_{t\geq1}$ is bounded for all $i\in\mathcal{N}$. Similarly, by using Lemma \ref{lem-decrease} and Lemma \ref{var_low_bound} again, we can show that $\left\{\|\tilde{\vz}_i^{t+1}-\vz_i^t\|\right\}_{t\geq1}$ is also bounded for all $i\in\mathcal{N}$. Recall that
		\begin{align}\label{AQ}
		\mQ_i=\left[\vone_{N_i},\mO_{N_i},-\mI_{N_i}\right]\otimes\mI_{n},\quad\mA_i=\left[\vone_{N_i},-\mI_{N_i},\mO_{N_i}\right]\otimes\mI_{n},\quad \forall i\in\mathcal{N},
		\end{align}
		then with simple algebraic manipulation, we can derive
		\begin{align}\label{iAQ}
		\left[\vzero_{N_i},\mO_{N_i},-\mI_{N_i}\right]\otimes\mI_{n}\cdot\mQ_i^T=\mI_{N_in},\quad\left[\vzero_{N_i},\mO_{N_i},-\mI_{N_i}\right]\otimes\mI_{n}\cdot\mA_i^T=\mO_{N_in},\quad \forall i\in\mathcal{N}.
		\end{align}
		In addition, substituting \eqref{w} into the update step of $\tilde{\vz}_i^{t+1}$, we have
		\begin{align}\label{ALambda}
		\mQ_i^T\mQ_i\vz_i^t=-\mW_i\left(\tilde{\vz}_i^{t+1}-\vz_i^{t}\right)+\mQ_i^T\mD_i\vu_i^t-c\mA_i^T\mA_i\vz_i^t-c\mA_i^T\vlambda_i^t,\quad \forall i\in\mathcal{N}.
		\end{align}
		Multiplying both sides of \eqref{ALambda} by $\left[\vzero_{N_i},\mO_{N_i},-\mI_{N_i}\right]\otimes\mI_{n}$ and taking the norm, we obtain
		\begin{align*}
		\|\mQ_i\vz_i^t\|=\|\left[\vzero_{N_i},\mO_{N_i},-\mI_{N_i}\right]\otimes\mI_{n}\cdot\mW_i\left(\tilde{\vz}_i^{t+1}-\vz_i^{t}\right)+\mD_i\vu_i^t\|\leq2\|\tilde{\vz}_i^{t+1}-\vz_i^{t}\|+d_{\max}\|\vu_i^t\|,
		\end{align*}
		where $d_{\max}=\max\{d_{i,j},i\in\mathcal{N},j\in\mathcal{N}_i\}$. The first equality is due to \eqref{iAQ}, while the second inequality follows the triangle inequality and \eqref{w}. Since the right-hand side (rhs) of the above inequality is bounded, it follows that $\{\|\mQ_i\vz_i^t\|\}_{t\geq1}$ is bounded for all $i\in\mathcal{N}$. 
		Recall the definition of $\vz_i=[\vp_i^T,(\vz_i^-)^T,(\vz_i^+)^T]^T$ and $\mQ_i, \mA_j$ in \eqref{AQ}, we have
		\begin{align}\label{eqn:QA}
		\mQ_i\vz_i^t=\text{vec}\left(\vp_i^t-(\vz_{i,j}^+)^t,j\in\mathcal{N}_i\right),
		\mA_j\vz_j^t=\text{vec}\left(\vp_j^t-(\vz_{j,i}^-)^t,i\in\mathcal{N}_j\right).
		\end{align}
		By using  the triangle inequality, for all $t\geq0$, we have
		\begin{align}\label{eqn:vzp}
		\|(\vz_{i,j}^+)^t\|\leq\|\vp_i^t\|+\|\vp_i^t-(\vz_{i,j}^+)^t\|,\|\vp_j^t\|\leq\|(\vz_{j,i}^-)^t\|+\|\vp_j^t-(\vz_{j,i}^-)^t\|.
		\end{align}
		Recall the update step of $\vz_i^t$,  we have
		\begin{align}\label{AzQz}
		(\vz_{i,j}^+)^t=(\vz_{j,i}^-)^t,~\forall~i\in\mathcal{N},j\in\mathcal{N}_i,\quad\text{and}\quad\vp^{t}_i=\va_i,~\forall~i\in\mathcal{A}\subseteq\mathcal{N}.
		\end{align}
		Since we have at least one anchor, combining \eqref{eqn:vzp} with \eqref{AzQz}, we can get
		\begin{align}\label{pj}
		\|\vp_j^t\|\leq\|\va_i\|+\|\vp_i^t-(\vz_{i,j}^+)^t\|+\|\vp_j^t-(\vz_{j,i}^-)^t\|,\quad\forall~j\in\mathcal{N}_i,~i\in\mathcal{A}.
		\end{align}
		Therefore, the sequence $\{\|\vp_j^t\|\}_{t\geq1}$ is bounded for all $j\in\mathcal{N}_i,i\in\mathcal{A}$ due to the rhs of \eqref{pj} being bounded, which is a consequence of the boundedness of $\{\|\mQ_i\vz_i^t\|\}_{t\geq1}$ and $\{\|\mA_i\vz_i^t\|\}_{t\geq1}$ for all $i\in\mathcal{N}$ (see \eqref{eqn:QA}). Using the definition of $\vz_j^t$ and the triangle inequality, we have
		\begin{align}\label{zj}
		\|\vz_j^t\|\leq\sum_{i\in\mathcal{N}_j}3\|\vp_j^t\|+\|\vp_j^t-(\vz_{j,i}^-)^t\|+\|\vp_j^t-(\vz_{j,i}^+)^t\|,\quad\forall~j\in\mathcal{N}.
		\end{align}
		Since the rhs of \eqref{zj} is bounded, then the sequence $\{\|\vz_j^t\|\}_{t\geq1}$ is bounded for all $j\in\mathcal{N}_i,i\in\mathcal{A}$. By applying a similar argument as in the derivation of   \eqref{eqn:vzp}-\eqref{zj}, we can show that $\vz^t_i$ is bounded in a neighborhood of the nodes $j\in\mathcal{N}_i,i\in\mathcal{A}$. Given that the graph is connected, we can thus conclude that the sequence $\{\vz_i^t\}_{t\geq1}$ is bounded for all $i\in\mathcal{N}$.
		\item Rearranging the terms of \eqref{ALambda} and taking the norm yields
		\begin{align}\label{aalambda}
		\|\vlambda_i^t\|\leq\|\mA_i^T\vlambda_i^t\|\leq\|\mQ_i^T\mQ_i\vz_i^t\|+\|\mW_i\left(\tilde{\vz}_i^{t+1}-\vz_i^{t}\right)\|+\|\mQ_i^T\mD_i\vu_i^t\|+\|c\mA_i^T\mA_i\vz_i^t\|,
		\end{align}
		where the first inequality holds due to the smallest eigenvalue of $\mA_i\mA_i^T$ being 1 from \eqref{AQ}, and the second inequality follows from the triangle inequality.  Since the rhs of \eqref{aalambda} is bounded, it follows that  $\{\vlambda_i^t\}_{t\geq1}$ is bounded for all $i\in\mathcal{N}$.
		\item Since $\|\tilde{\vz}_i^{t+1}\|\leq\|\tilde{\vz}_i^{t+1}-\vz_i^t\|+\|\vz_i^t\|$, and the rhs is bounded, it follows that the sequence $\{\tilde{\vz}_i^t\}_{t\geq1}$ is bounded for all $i\in\mathcal{N}$.
	\end{enumerate}
\end{IEEEproof}
	
\begin{Lemma}[Subgradient bound]\label{lem-partial-bound}
Let $\{\vy^t=\left(\vz^t,\vu^t,\vlambda^t\right)\}_{t\geq1}$ be a sequence generated by Algorithm \ref{al2}, $\kappa_1,\kappa_{2},c$, and $\rho$ are parameters used in the potential function $\varsigma^t$. Then we have
\begin{align*}
\vv^{t+1}:=\left(\vv_{\vz^{t+1}},\vv_{\vu^{t+1}},\vv_{\vlambda^{t+1}},\vv_{\tilde{\vz}^{t+1}},\vv_{\vu^t},\vv_{\vz^t}\right)\in\partial\varsigma^{t+1},\quad\forall~t\geq 0,
\end{align*}
where
\begin{align}\nonumber
\vv_{\vz^{t+1}_i}:=&\mW_i\left(\vz_i^{t}-\tilde{\vz}_i^{t+1}\right)+\mQ_i^T\mD_i\left(\vu_{i}^t-\vu_{i}^{t+1}\right)+(\kappa_{2}+1)\mA_i^T\left(\vlambda_i^{t+1}-\vlambda_i^{t}\right)\\\nonumber&+\left(c\mA_i^T\mA_i+\mQ_i^T\mQ_i+c\left(\kappa_{1}+\kappa_{2}\right)\mB_i^T\mB_i\right)\left(\vz_i^{t+1}-\vz_i^{t}\right),\\\nonumber\vv_{\vu^{t+1}_i}:=&\frac{\rho}{2}\left(\vu_{i}^{t}-\vu_{i}^{t+1}\right),~\vv_{\vlambda^{t+1}_i}:=\frac{1}{c}\left(\vlambda_i^{t+1}-\vlambda_i^{t}\right),~\vv_{\tilde{\vz}^{t+1}_i}:=\kappa_1\mA_i^T\left(c\mA_i\left(\tilde{\vz}_i^{t+1}-\vz_i^{t+1}\right)+\vlambda_i^{t+1}-\vlambda_i^{t}\right),\\\label{partial_poten}\vv_{\vu^{t}_i}:=&\frac{\rho}{2}\left(\vu_{i}^{t}-\vu_{i}^{t+1}\right),~\vv_{\vz^{t}_i}:=c\left(\kappa_{1}+\kappa_{2}\right)\mB_i^T\mB_i\left(\vz_{i}^{t}-\vz_{i}^{t+1}\right),\quad i\in\mathcal{N}.
\end{align}
Moreover, for every $t\geq0$, it holds that
\begin{align}\nonumber \|\vv^{t+1}\|\leq&\sum_{i\in\mathcal{N}}\alpha_i\|\vz_i^{t+1}-\vz_i^{t}\|+\sum_{i\in\mathcal{N}}\left(\rho+d_{\max}\sqrt{N_i+1}\right)\|\vu_i^{t+1}-\vu_i^{t}\|+\sum_{i\in\mathcal{N}}\beta_i\|\vlambda_i^{t+1}-\vlambda_i^{t}\|\\\label{e}\leq&C_2\|\vy^{k+1}-\vy^k\|,
\end{align}
where
\begin{align}\nonumber
\alpha_i:=&\left(c+1\right)\left(N_i+1\right)+2\left(\kappa_{1}+\kappa_{2}\right)\left(1+c\right)\left(1+N_{\max}\right)+\frac{2N_i\left(c+1\right)^2\sqrt{nNN_{\text{sum}}(N_i+1)}}{\sqrt{\tilde{\tau}_{\min}}},\\\label{alphabe}\beta_i:=&\sqrt{N_i+1}\left(\kappa_{1}+1+\kappa_{2}\right)+\frac{1}{c},\quad C_2:=\sqrt{3N}\cdot\max\{\alpha_i,\rho+d_{\max}\sqrt{N_i+1},\beta_i,i\in\mathcal{N}\}.
\end{align}
\end{Lemma}
\begin{IEEEproof}
	By taking partial derivatives of $\varsigma^{t+1}$ with respect to $\vz_i^{t+1},\vu_i^{t+1},\vlambda_i^{t+1},\tilde{\vz}_i^{t+1},\vu_i^{t},\vz_i^{t}$ and using the update step of $\vlambda_i^{t+1}$, we obtain
	\begin{align}\nonumber
	\nabla_{\vz_i^{t+1}}\varsigma^{t+1}&=\nabla_{\vz_i}\mathcal{L}_i\left(\vz_i^{t+1},\vu_i^{t+1},\vlambda_i^{t+1}\right)+\kappa_{2}\mA_i^T\left(\vlambda_i^{t+1}-\vlambda_i^t\right)+c\left(\kappa_{1}+\kappa_{2}\right)\mB_i^T\mB_i\left(\vz_i^{t+1}-\vz_i^{t}\right),\\\nonumber\partial_{\vu_i^{t+1}}\varsigma^{t+1}&=\partial_{\vu_i}\mathcal{L}_i\left(\vz_i^{t+1},\vu_i^{t+1},\vlambda_i^{t+1}\right)+\frac{\rho}{2}\left(\vu_i^{t+1}-\vu_i^{t}\right),\\\nonumber\nabla_{\vlambda_i^{t+1}}\varsigma^{t+1}&=\nabla_{\vlambda_i}\mathcal{L}_i\left(\vz_i^{t+1},\vu_i^{t+1},\vlambda_i^{t+1}\right),\\\label{nablaL}\nabla_{\tilde{\vz}_i^{t+1}}\varsigma^{t+1}&=\vv_{\tilde{\vz}^{t+1}_i},\quad\nabla_{\vu_i^{t}}\varsigma^{t+1}=\vv_{\vu^{t}_i},\quad\nabla_{\vz_i^{t}}\varsigma^{t+1}=\vv_{\vz^{t}_i},
	\end{align}
	where $\vv_{\tilde{\vz}^{t+1}_i},\vv_{\vu^t_i}$, and $\vv_{\vz^t_i}$ are defined in \eqref{partial_poten}.
	Using the update step of $\tilde{\vz}_i^{t+1}$, we have
	\begin{align}\nonumber
	\nabla_{\vz_i}\mathcal{L}_i\left(\vz_i^{t+1},\vu_i^{t+1},\vlambda_i^{t+1}\right)=&\mW_i\left(\vz_i^{t}-\tilde{\vz}_i^{t+1}\right)+\mQ_i^T\mD_i\left(\vu_i^{t}-\vu_i^{t+1}\right)+\mA_i^T\left(\vlambda_i^{t+1}-\vlambda_i^{t}\right)\\\label{nablzz}&+\left(c\mA_i^T\mA_i+\mQ_i^T\mQ_i\right)\left(\vz_i^{t+1}-\vz_i^{t}\right).
	\end{align}
	Combining \eqref{nablzz} with \eqref{nablaL} yields $\vv_{\vz^{t+1}_i}=\nabla_{\vz_i^{t+1}}\varsigma^{t+1}$. By the optimality condition of $\vu$ subproblem \eqref{u_up}, we have
	\begin{align}\label{nabluu}
	\vzero\in\partial_{\vu_i}\mathcal{L}_i\left(\vz_i^{t+1},\vu_i^{t+1},\vlambda_i^{t}\right)+\rho\left(\vu_i^{t+1}-\vu_i^{t}\right)=\partial_{\vu_i}\mathcal{L}_i\left(\vz_i^{t+1},\vu_i^{t+1},\vlambda_i^{t+1}\right)+\rho\left(\vu_i^{t+1}-\vu_i^{t}\right).
	\end{align}
	Combining \eqref{nabluu} with \eqref{nablaL} yields
	$\vv_{\vu^{t+1}_i}\in\partial_{\vu_i^{t+1}}\varsigma^{t+1}$. Applying the update formula of $\vlambda_i^{t+1}$, yields
	\begin{align*}
	\vv_{\vlambda^{t+1}_i}=\frac{1}{c}\left(\vlambda_i^{t+1}-\vlambda_i^{t}\right)=\mA_i\vz_i^{t+1}=\nabla_{\vlambda_i}\mathcal{L}_i\left(\vz_i^{t+1},\vu_i^{t+1},\vlambda_i^{t+1}\right)=\nabla_{\vlambda_i^{t+1}}\varsigma^{t+1}.
	\end{align*}
	By using the expressions for $\vv^{t+1}$ and the triangle inequality, we obtain the following bound
	\begin{align}\label{v1}
	\|\vv^{t+1}\|\leq&\|\mW_i\|\|\tilde{\vz}_i^{t+1}-\vz_i^{t}\|+\left(2c(\kappa_{1}+\kappa_{2})\|\mB_i^T\mB_i\|+c\|\mA_i^T\mA_i\|+\|\mQ_i^T\mQ_i\|\right)\|\vz_i^{t+1}-\vz_i^{t}\|\\\nonumber&+\left(\sigma_{\max}^{1/2}\left(\mQ_i\mQ_i^T\right)\|\mD_i\|+\rho\right)\|\vu_i^{t+1}-\vu_i^{t}\|+\Big(\sigma^{1/2}_{\max}\left(\mA_i\mA_i^T\right)\left(\kappa_{2}+1+\kappa_1\right)+\frac{1}{c}\Big)\|\vlambda_i^{t+1}-\vlambda_i^{t}\|,
	\end{align}
	where $\sigma_{\max}\left(\mQ_i\mQ_i^T\right)$ and $\sigma_{\max}\left(\mA_i\mA_i^T\right)$ denote the largest eigenvalue of $\mQ_i\mQ_i^T$ and $\mA_i\mA_i^T$, respectively. By taking the square root of both sides of the inequality \eqref{Azz1}, and combining the result using the fact that $\|\vx\|_2\leq\|\vx\|_{1}\leq \sqrt{N}\|\vx\|_2,\,\forall~\vx\in\mathbb{R}^N$, we obtain
	\begin{align}\label{zz}
	\sum_{i\in\mathcal{N}}\|\tilde{\vz}_i^{t+1}-\vz_i^t\|\leq\frac{\left(c+1\right)\sqrt{nNN_{\text{sum}}}}{\sqrt{\tilde{\tau}_{\min}}}\sum_{i\in\mathcal{N}}\|\mA_i\left(\vz^{t+1}_i-\vz_i^t\right)\|.
	\end{align}
	Furthermore, using \eqref{w} and \eqref{AQ}, we have
	\begin{align}\label{cab}
	\|\mW_i\|=2(c+1)N_i,\sigma_{\max}\left(\mA_i\mA_i^T\right)=\sigma_{\max}\left(\mQ_i\mQ_i^T\right)=N_i+1.
	\end{align}
	Substituting \eqref{evib}, \eqref{zz}, and \eqref{cab} into \eqref{v1}, we get
	\begin{align}\label{first_result}
	\|\vv^{t+1}\|\leq\sum_{i\in\mathcal{N}}\alpha_i\|\vz_i^{t+1}-\vz_i^{t}\|+\sum_{i\in\mathcal{N}}\left(\rho+d_{\max}\sqrt{N_i+1}\right)\|\vu_i^{t+1}-\vu_i^{t}\|+\sum_{i\in\mathcal{N}}\beta_i\|\vlambda_i^{t+1}-\vlambda_i^{t}\|,
	\end{align}
	where $\alpha_i$ and $\beta_i$ are defined in \eqref{alphabe}. By the fact that $\|\vx\|_{1}\leq \sqrt{3N}\|\vx\|_2,\,\forall~\vx\in\mathbb{R}^{3N}$, we obtain
	\begin{align*}
	\|\vv^{t+1}\|\leq&C_2\|\vy^{t+1}-\vy^{t}\|,
	\end{align*}
	where constant $C_2$ is defined in \eqref{alphabe}. This concludes the proof.
\end{IEEEproof}
	
In the following, we denote by $\omega\left(\{\vx^t\}_{t\geq1}\right)$ the set of limit points of the sequence $\{\vx^t\}_{t\geq1}$, and we define  $\text{crit}\,\mathcal{L}:=\{\vy:\vzero\in\partial\mathcal{L}\left(\vy\right)\}$ as the set of all critical points of $\mathcal{L}$.

\begin{Lemma}[Properties of limit point set]\label{lem-distance}Let $\{(\vz^t,\vu^t,\vlambda^t)\}_{t\geq1}$ be a sequence generated by Algorithm \ref{al2}. Suppose that conditions \eqref{k1_all}-\eqref{c_all} are satisfied, and let $\Omega:=\omega\left(\{(\vz^t,\vu^t,\vlambda^t,\tilde{\vz}^t,\vu^{t-1},\vz^{t-1})\}_{t\geq1}\right)$, we have the following results:
	\begin{enumerate}
		\item the set $\Omega$ is nonempty and compact;
		\item $\lim\limits_{t\to\infty}\text{dist}\Big[\left(\vz^t,\vu^t,\vlambda^t,\tilde{\vz}^t,\vu^{t-1},\vz^{t-1}\right),\Omega\Big]=0$;
		\item$\Omega\subseteq\{\left(\vz,\vu,\vlambda,\vz,\vu,\vz\right):\left(\vz,\vu,\vlambda\right)\in\text{crit}\,\mathcal{L}\}$;
		\item any critical point of $\mathcal{L}$ is a KKT point of problem \eqref{f};
		\item the potential function $\varsigma^t$ is finite and constant on $\Omega$.
	\end{enumerate}
\end{Lemma}
\begin{IEEEproof}
	We prove the results item by item below.
	\begin{enumerate}
		\item Since $\{\left(\vz^t,\vu^t,\vlambda^t,\tilde{\vz}^t\right)\}_{t\geq1}$ is bounded by Lemma \ref{lem-bounded}, and thus $\Omega$ is nonempty and bounded. By the definition of $\Omega$, it is closed and therefore compact.
		\item As a consequence of the limit point definition.
		\item Let $\left(\vz^\ast,\vu^\ast,\vlambda^\ast\right)$ be a limit point of the sequence $\{\left(\vz^t,\vu^t,\vlambda^t\right)\}_{t\geq1}$, which exists since the sequence $\{\left(\vz^t,\vu^t,\vlambda^t\right)\}_{t\geq1}$ is bounded. Consequently, there exists a subsequence $\{\left(\vz^{t_k},\vu^{t_k},\vlambda^{t_k}\right)\}_{k\geq1}$ of $\{\left(\vz^t,\vu^t,\vlambda^t\right)\}_{t\geq1}$ that converges to $\left(\vz^\ast,\vu^\ast,\vlambda^\ast\right)$. Note that Lemma \ref{lem-decrease} with Lemma \ref{var_low_bound} implies
		\begin{align*}
		\sum_{i\in\mathcal{N}}\|\tilde{\vz}^{t+1}_i-\vz^{t}_i\|^2\to0,\quad\sum_{i\in\mathcal{N}}\|\vu_i^{t+1}-\vu_i^{t}\|^2\to0,\quad\sum_{i\in\mathcal{N}}\|\vz_i^{t+1}-\vz_i^{t}\|^2\to0,\quad\text{as}~t\to\infty,
		\end{align*}
		which means $\{\left(\tilde{\vz}^{t_k},\vu^{t_k-1},\vz^{t_k-1}\right)\}_{k\geq1}$ converges to $\left(\vz^\ast,\vu^\ast,\vz^\ast\right)$. Hence, $\left(\vz^\ast,\vu^\ast,\vlambda^\ast,\vz^\ast,\vu^\ast,\vz^\ast\right)\in\Omega$, and from the continuity of $\varsigma^t$, it follows that
		\begin{align*}
		\lim\limits_{k\to\infty}\varsigma^{t_k}=\varsigma^{\ast}.
		\end{align*}
		On the other hand, from Lemma \ref{lem-decrease},  Lemma \ref{lem-partial-bound}, and Lemma \ref{var_low_bound}, we know that $\vv^{t_k}\in\partial\varsigma^{t_k}~\text{and}~\vv^{t_k}\rightarrow\vzero~\text{as}~k\to\infty.$ The closeness property of $\partial\varsigma$ (see [\cite{bolte2014proximal}, Remark 1 (ii)]) implies that
		\begin{align}\label{parvar}
		\vzero\in\partial\varsigma^{\ast}.
		\end{align}
		In addition, according to \eqref{nablaL}  we also have
		\begin{align}\label{parL}
		\partial\varsigma^{\ast}=\partial\mathcal{L}\left(\vz^\ast,\vu^\ast,\vlambda^\ast\right).
		\end{align}
		Combining \eqref{parvar} with \eqref{parL}, we conclude that  $\left(\vz^\ast,\vu^\ast,\vlambda^\ast\right)$ is a critical point of $\mathcal{L}$.
		\item Let $\left(\vz^\ast,\vu^\ast,\vlambda^\ast\right)\in\text{crit}\,\mathcal{L}$. Thanks to the separable property of $\mathcal{L}$ about node $i$, we have
		\begin{align}\label{lz}
		\vzero=&\nabla_{\vz_i}\mathcal{L}_i\left(\vz_i^\ast,\vu_i^\ast,\vlambda_i^\ast\right)=\nabla_{\vz_i}F_i(\vz_i^\ast,\vu_i^\ast)+\mA_i^T\vlambda_i^\ast+c\mA_i^T\mA_i\vz_i^\ast,\\\label{lu}\vzero\in\,&\partial{\vu_i}\mathcal{L}_i\left(\vz_i^\ast,\vu_i^\ast,\vlambda_i^\ast\right)=\nabla_{\vu_i}F_i(\vz_i^\ast,\vu_i^\ast)+\partial\delta_{\mathcal{B}^{N_i}}\left(\vu_{i}^\ast\right),\\\label{llambda}\vzero=&\nabla_{\vlambda_i}\mathcal{L}_i\left(\vz_i^\ast,\vu_i^\ast,\vlambda_i^\ast\right)=\mA_i\vz_i^\ast.
		\end{align}
		Thus, it follows from \eqref{lz} and \eqref{llambda} that
		\begin{align}\label{nablz}
		\nabla_{\vz_i}F_i(\vz_i^\ast,\vu_i^\ast)+\mA_i^T\vlambda_i^\ast=\vzero.
		\end{align}
		Since $F_i(\vz_i,\vu_i)+\langle\vlambda_i,\mA_i\vz_i\rangle$ is convex over $\mathcal{X}$ and $\mathcal{Z}$ for any fixed $\vu_i$, then the above expression \eqref{nablz} implies that $\vz_i^\ast$ also satisfies
		\begin{align}\label{KKTz}
		\vz_i^\ast\in\arg\min\lim\limits_{\substack{\vz\in\mathcal{Z}\\\vx\in\mathcal{X}}}F_i(\vz_i,\vu_i^\ast)+\langle\vlambda_i^\ast,\mA_i\vz_i\rangle.
		\end{align}
		Combining \eqref{lu}, \eqref{llambda}, and \eqref{KKTz}, we obtain $\left(\vz^\ast,\vu^\ast,\vlambda^\ast\right)$ as a KKT point of problem \eqref{f}.
		\item The sequence $\{\varsigma^t\}_{t\geq1}$ decreases by Lemma \ref{lem-decrease} and is bounded from below by Lemma \ref{var_low_bound}, thereby implying its convergence to some ﬁnite limit $\varsigma^{\ast}$. It follows that $\varsigma^t$ is constant on $\Omega$.
	\end{enumerate}
\end{IEEEproof}

To achieve our main goal to establish the global convergence of the whole sequence, we recall the following key result obtained in the existing literature \cite{bolte2014proximal,guo2017convergence,bolte2018first,boct2020proximal,yashtini2022convergence}.
	
\begin{Lemma}[Uniform KL property]\label{unif-KL}
	Let $\Omega$ be a compact set and let $\sigma:\mathbb{R}^d\to(-\infty,\infty]$ a proper and lower semicontinuous function. Assume that $\sigma$ is constant on $\Omega$, and satisfies
	the KL property at each point of $\Omega$. Then there exist $\varepsilon>0,\eta>0$ and $\varphi\in\Phi_{\eta}$ such that for all $\bar{u}$ in $\Omega$, and all $u$ in the following intersection
	\begin{align*}
	\left\{u\in\mathbb{R}^d:\operatorname{dist}(u,\Omega)<\varepsilon\right\}\cap[\sigma(\bar{u})<\sigma(u)<\sigma(\bar{u})+\eta],
	\end{align*}
	one has,	
	\begin{align}
	\varphi^{\prime}(\sigma(u)-\sigma(\bar{u})) \operatorname{dist}(0,\partial\sigma(u))\geq 1.
	\end{align}
\end{Lemma}
	
We can now conveniently summarize our convergence results.
\begin{Lemma}[Global convergence]Let $\{\vy^t=(\vz^t,\vu^t,\vlambda^t)\}_{t\geq1}$ be a sequence generated by Algorithm \ref{al2}. Suppose that conditions \eqref{k1_all}-\eqref{c_all} are satisfied, then the following statements are true:
	\begin{enumerate}
		\item the sequence $\{\vy^t\}_{t\geq1}$ is bounded and has finite length, namely,
		\begin{align*}
		\sum_{t=1}^{\infty}\|\vy^{t+1}-\vy^{t}\|<+\infty;
		\end{align*}
		\item the sequence $\{\vy^t\}_{t\geq0}$ converges to a KKT point of the problem \eqref{f}.
	\end{enumerate}
\end{Lemma}
	\begin{IEEEproof}
	As in Lemma \ref{lem-distance}, we denote by $\Omega:=\omega\left(\{(\vz^t,\vu^t,\vlambda^t,\tilde{\vz}^t,\vu^{t-1},\vz^{t-1})\}_{t\geq1}\right)$, which is a nonempty and compact set. From the proof of Lemma \ref{lem-distance} and the continuity of $\varsigma^t$, it follows that
	\begin{align}\label{lim-varsigma}
	\lim\limits_{t\to\infty}\varsigma^t=\varsigma^\ast,~\text{for all}~(\vz^{\ast},\vu^{\ast},\vlambda^{\ast},\vz^{\ast},\vu^{\ast},\vz^{\ast})\in\Omega.
	\end{align}
	Moreover, from Lemma \ref{lem-decrease}, we have
	\begin{align}\label{low}
	C_3\|\vy^{t+1}-\vy^{t}\|^2\leq\varsigma^t-\varsigma^{t+1},~\forall~t\geq1.
	\end{align}
	where  $C_3:=\min\{c,1,\frac{\rho}{4},C_0\}$. We consider two cases.
	\begin{enumerate}
		\item If there exists an integer $\bar{t}\geq0$ such that $\varsigma^{\bar{t}}=\varsigma^\ast$, then using the decreasing property obtained by Lemma \ref{lem-decrease} and \eqref{low}, we have
		\begin{align*}
		C_3\|\vy^{t+1}-\vy^{t}\|^2\leq\varsigma^{t}-\varsigma^{t+1}\leq\varsigma^{\bar{t}}-\varsigma^\ast=0,\quad t\geq\bar{t}.
		\end{align*}
		Thus, $\vy^{t+1}=\vy^{t}$ for any $t\geq\bar{t}$, and it is clear that  $\sum_{t=1}^{\infty}\|\vy^{t+1}-\vy^{t}\|<+\infty$ holds.
		\item Since $\varsigma^{t}$ is non-increasing from Lemma \ref{lem-decrease}, we have $\varsigma^{t}\geq\varsigma^{\ast}$ for all $t$. Therefore, let us assume that $\varsigma^{t}>\varsigma^{\ast}$ for all $t$. Again from \eqref{lim-varsigma}, we know that for any $\eta>0$, there exists $t_0>0$ such that
		\begin{align*}
		\varsigma^{t}<\varsigma^{\ast}+\eta,\quad t\geq t_0.
		\end{align*}
		According to Lemma \ref{lem-distance}, we know that $\lim\limits_{t\to\infty}\text{dist}\left[\left(\vz^t,\vu^t,\vlambda^t,\tilde{\vz}^t,\vu^{t-1},\vz^{t-1}\right),\Omega\right]=0$. This implies that for any $\varepsilon>0$, there exists $t_1\geq1$ such that $\text{dist}\left[\left(\vz^t,\vu^t,\vlambda^t,\tilde{\vz}^t,\vu^{t-1},\vz^{t-1}\right),\Omega\right]<\varepsilon$ for all $t\geq t_1$. Summing up all these facts, we obtain the following inequalities for any  $\eta,\varepsilon>0$
		\begin{align*}
		\text{dist}\left[\left(\vz^t,\vu^t,\vlambda^t,\tilde{\vz}^t,\vu^{t-1},\vz^{t-1}\right),\Omega\right]<\varepsilon,\quad\text{and}\quad\varsigma^{\ast}<\varsigma^{t}<\varsigma^{\ast}+\eta,\quad \text{for all}~t\geq t_2:=\max\{t_0,t_1\}.
		\end{align*}
		According to Lemma \ref{lem-distance}, it has been established that $\varsigma$ is a constant on $\Omega$. Moreover, we note that $\varsigma^t$ is semi-algebraic (a polynomial function) and satisfies the KL property (see [\cite{bolte2018first}, Theorem 6.1]). Therefore, we can apply Lemma \ref{unif-KL} with $\Omega$. Consequently, for any $t\geq t_2$, we have
		\begin{align}\label{uni-KL}
		\varphi^{\prime}\left(\varsigma^t-\varsigma^\ast\right)\cdot\text{dist}\left(\vzero,\partial\varsigma^t\right)\geq1.
		\end{align}
		Due to the concavity of $\varphi$, we get
		\begin{align}\label{concavi}
		\varphi\left(\varsigma^t-\varsigma^\ast\right)-\varphi\left(\varsigma^{t+1}-\varsigma^\ast\right)\geq\varphi^{\prime}\left(\varsigma^t-\varsigma^\ast\right)\left(\varsigma^t-\varsigma^{t+1}\right).
		\end{align}
		From \eqref{e} of Lemma \ref{lem-partial-bound}, we also have
		\begin{align*}
		\text{dist}\left(\vzero,\partial\varsigma^t\right)\leq C_2\|\vy^{t}-\vy^{t-1}\|,\quad C_2>0.
		\end{align*}
		By combining this with \eqref{concavi}, \eqref{uni-KL}, and $\varphi^{\prime}\left(\varsigma^t-\varsigma^\ast\right)>0$, we obtain
		\begin{align}\nonumber
		\varsigma^t-\varsigma^{t+1}\leq&\frac{\varphi\left(\varsigma^t-\varsigma^\ast\right)-\varphi\left(\varsigma^{t+1}-\varsigma^\ast\right)}{\varphi^{\prime}\left(\varsigma^t-\varsigma^\ast\right)}\\\label{upper}\leq&C_2\|\vy^{t}-\vy^{t-1}\|\cdot\left[\varphi\left(\varsigma^t-\varsigma^\ast\right)-\varphi\left(\varsigma^{t+1}-\varsigma^\ast\right)\right].
		\end{align}
		For convenience, we define the following for two arbitrary nonnegative integers $p$ and $q$
		\begin{align*}
		\Delta_{p,q}:=\varphi\left(\varsigma^p-\varsigma^\ast\right)-\varphi\left(\varsigma^{q}-\varsigma^\ast\right).
		\end{align*}
		Combining \eqref{upper} with \eqref{low}, we conclude that for any $t\geq t_2$, the following inequality holds
		\begin{align*}
		C_3\|\vy^{t+1}-\vy^{t}\|^2\leq\varsigma^t-\varsigma^{t+1}\leq C_2\|\vy^{t}-\vy^{t-1}\|\Delta_{t,t+1},
		\end{align*}
		Using the fact that $2\sqrt{ab}\leq a+b$ for all $a,b\geq0$, we can infer
		\begin{align}\label{sum}
		2\|\vy^{t+1}-\vy^{t}\|\leq\|\vy^{t}-\vy^{t-1}\|+\frac{C_2}{C_3}\Delta_{t,t+1}.
		\end{align}
		Let us now prove that for any $t>t_2$, the following inequality holds
		\begin{align*}
		\sum_{k=t_2+1}^{t}\|\vy^{k+1}-\vy^{k}\|\leq\|\vy^{t_2+1}-\vy^{t_2}\|+\frac{C_2}{C_3}\Delta_{t_2,t+1}.
		\end{align*}
		Summing up \eqref{sum} for $k=t_2+1,t_2+2,\ldots,t$ yields
		\begin{align*}
		2\sum_{k=t_2+1}^{t}\|\vy^{k+1}-\vy^{k}\|\leq&\sum_{k=t_2+1}^{t}\|\vy^{k}-\vy^{k-1}\|+\frac{C_2}{C_3}\sum_{k=t_2+1}^{t}\Delta_{k,k+1}\\\leq&\sum_{k=t_2+1}^{t+1}\|\vy^{k}-\vy^{k-1}\|+\frac{C_2}{C_3}\sum_{k=t_2+1}^{t}\Delta_{k,k+1}\\\leq&\sum_{k=t_2+1}^{t}\|\vy^{k+1}-\vy^{k}\|+\|\vy^{t_2+1}-\vy^{t_2}\|+\frac{C_2}{C_3}\Delta_{t_2+1,t+1},
		\end{align*}
		where the last inequality follows from the fact that $\Delta_{p,q}+\Delta_{q,r}=\Delta_{p,r}$ for arbitrary nonnegative integers $p,q$, and $r$. Since $\varphi\geq0$, we have for any $t>t_2$ that
		\begin{align*}
		\sum_{k=t_2+1}^{t}\|\vy^{k+1}-\vy^{k}\|\leq\|\vy^{t_2+1}-\vy^{t_2}\|+\frac{C_2}{C_3}\cdot\varphi\left(\varsigma^{t_2+1}-\varsigma^{\ast}\right).
		\end{align*}
		As the right hand-side of the inequality above does not depend on $t$ at all, it
		immediately follows that the sequence $\{\vy^t\}_{t\in\mathbb{N}}$ has finite length, that is
		\begin{align*}
		\sum_{t=1}^{\infty}\|\vy^{t+1}-\vy^{t}\|<\infty.
		\end{align*}
		This means that the sequence $\{\vy^t\}_{t\geq1}$ is a Cauchy sequence and hence a convergent sequence.
	\end{enumerate}
	Thanks to Lemma \ref{lem-distance}, there exists $\vy^{\ast}\in\text{crit}\mathcal{L}$ such that $\lim\limits_{t\to\infty}\vy^t=\vy^{\ast}$ and $\vy^{\ast}$ is also a KKT point of problem \eqref{f}. We complete the proof.
\end{IEEEproof}}
	
	\begin{Lemma}\label{lem1-3}
		Suppose $c\mB_i^T\mB_i$ takes the form of \eqref{cbb}, and the sequence $\left\{\left(\vz_i^t,\vu_i^t,\vlambda_i^t\right)\right\}$ is generated by the Algorithm \ref{al2}. Then we have
		\begin{align*}
		&\sum_{i\in\mathcal{N}}\|\vz_i^t-{\rm proj}_{\mathcal{X},\mathcal{Z}}\left(\vz_i^t-\left(\nabla_{\vz_i}F_i\left(\vz_i^t,\vu_i^t\right)+\mA_i^T\vlambda_i^t\right)\right)\|^{2}\\\leq
		&\sum_{i\in\mathcal{N}}\left[\sigma_{1}\|\tilde{\vz}_i^{t+1}-\vz_i^t\|^{2}+\sigma_{2}\|\vz_i^t-\vz_i^{t+1}\|^{2}+3c^{2}\left(N_{\max}+1\right)^{2}\|\mA_i\vz_i^t\|^{2}\right],
		\end{align*}
		where
		\begin{align}\nonumber
		\sigma_{1}&=3\left(2\left(c+1\right)N_{\max}-1\right)^{2}+6\max\left\{\left(1+c\right)^{2},\left(1+\frac{1}{c}\right)^{2}\right\},\\\label{sigma}
		\sigma_{2}&=\frac{3}{2}\max\left\{\left(1+c\right)^{2},\left(1+\frac{1}{c}\right)^{2}\right\}.
		\end{align}
	\end{Lemma}
	
	\begin{IEEEproof}
		Using \eqref{ztilde} and the deﬁnition of the $F_i$ given in problem \eqref{f}, we have
		\begin{align*}
		\nabla_{\vz_i}F_i\left(\vz_i^t,\vu_i^t\right)+\mA_i^T\boldsymbol{\lambda}_i^t=
		&\mQ_i^T\mQ_i\vz_i^t-\mQ_i^T\mD_i\vu_i^t+\mA_i^T\boldsymbol{\lambda}_i^t\\=
		&-\mW_i(\tilde{\vz}_i^{t+1}-\vz_i^t)+c\mA_i^T\mA_i\vz_i^t.
		\end{align*}
		Hence, with the above equality, we have
		\begin{align}\nonumber	&\|\vz_i^t-\text{proj}_{\mathcal{X},\mathcal{Z}}\left(\vz_i^t-\left(\nabla_{\vz_i}F_i\left(\vz_i^t,\vu_i^t\right)+\mA_i^T\vlambda_i^t\right)\right)\|\\\nonumber=
		&\|\vz_i^t-\text{proj}_{\mathcal{X},\mathcal{Z}}\left(\tilde{\vz}_i^{t+1}\right)+\text{proj}_{\mathcal{X},\mathcal{Z}}\left(\tilde{\vz}_i^{t+1}\right)-\text{proj}_{\mathcal{X},\mathcal{Z}}\left(\vz_i^t+\left(\mW_i(\tilde{\vz}_i^{t+1}-\vz_i^t)+c\mA_i^T\mA_i\vz_i^t\right)\right)\|\\\label{conv1}\leq
		&\|\vz_i^t-\text{proj}_{\mathcal{X},\mathcal{Z}}\left(\tilde{\vz}_i^{t+1}\right)\|+\|\left(\mI_i-\mW_i\right)\left(\tilde{\vz}_i^{t+1}-\vz_i^t\right)-c\mA_i^T\mA_i\vz_i^t\|,
		\end{align}
		where the inequality dues to the triangle inequality and the nonexpansive property of the projection operator \cite{Hong2017}. By completing the square and using the Cauchy-Schwarz inequality, we further obtain
		\begin{align}\nonumber
		&\|\vz_i^t-\text{proj}_{\mathcal{X},\mathcal{Z}}\left(\vz_i^t-\left(\nabla_{\vz_i}F_i\left(\vz_i^t,\vu_i^t\right)+\mA_i^T\vlambda_i^t\right)\right)\|^{2}\\\label{nonexp}\leq
		&3\|\vz_i^t-\text{proj}_{\mathcal{X},\mathcal{Z}}\left(\tilde{\vz}_i^{t+1}\right)\|^{2}+3c^{2}\|\mA_i^T\mA_i\vz_i^t\|^{2}+3\|\left(\mI_i-\mW_i\right)\left(\tilde{\vz}_i^{t+1}-\vz_i^t\right)\|^{2}.
		\end{align}
		Now, let us consider solving $\text{proj}_{\mathcal{X},\mathcal{Z}}\left(\tilde{\vz}_i^{t+1}\right)$. Since
		\begin{subequations}
			\begin{align*}
			\text{proj}_{\mathcal{X},\mathcal{Z}}\left(\tilde{\vz}_i^{t+1}\right)=\arg\min_{\vz}&\frac{1}{2}\sum_{i\in\mathcal{N}}\left\|\vz_i-\tilde{\vz}_i^{t+1}\right\|^{2}\\
			\text{subject to}&\,\,\vz\in\mathcal{X}, \vz\in\mathcal{Z}.
			\end{align*}
		\end{subequations}
		By a similar argument as the proof of Remark \ref{r_z_ztilde}, we can derive
		\begin{equation}\label{prox}
		\text{proj}_{\mathcal{X,Z}}\left(\tilde{\vz}_i^{t+1}\right)=\widetilde{\mW}_i\left(\tilde{\vz}_i^{t+1}-\vz_i^{t+1}\right)+\tilde{\vz}_i^{t+1},
		\end{equation}
		where $\tilde{\vx}_{i,i}^{t+1}=\va_i$ if $i\in\mathcal{A}$ and
		\begin{equation}\label{w_tilde}
		\widetilde{\mW}_i=\frac{1}{2}\cdot\textbf{D}\text{iag}\left(\left[0,\left(c+1\right)\cdot\vone_{N_i}^T,\frac{c+1}{c}\cdot\vone_{N_i}^T\right]\right)\otimes\mI_{n}.
		\end{equation}
		Then, by \eqref{prox} and Cauchy–Schwartz inequality, we have
		\begin{align}\nonumber
		\|\vz_i^t-\text{proj}_{\mathcal{X},\mathcal{Z}}\left(\tilde{\vz}_i^{t+1}\right)\|^{2}=
		&\|\vz_i^t-\tilde{\vz}_i^{t+1}-\widetilde{\mW}_i\left(\tilde{\vz}_i^{t+1}-\vz_i^t+\vz_i^t-\vz_i^{t+1}\right)\|^{2}\\\nonumber\leq
		&2\|\left(\mI_i+\widetilde{\mW}_i\right)\left(\vz_i^t-\tilde{\vz}_i^{t+1}\right)\|^{2}+2\|\widetilde{\mW}_i\left(\vz_i^t-\vz_i^{t+1}\right)\|^{2}\\\label{zproj}\leq
		&2\|\mI_i+\widetilde{\mW}_i\|^{2}\|\vz_i^t-\tilde{\vz}_i^{t+1}\|^{2}+2\|\widetilde{\mW}_i\|^{2}\|\vz_i^t-\vz_i^{t+1}\|^{2}.
		\end{align}
		For the third part of \eqref{nonexp}, by using the Cauchy–Schwartz inequality again
		\begin{equation}\label{ie2w}
		\|\left(\mI_i-\mW_i\right)(\tilde{\vz}_i^{t+1}-\vz_i^t)\|^{2}\leq\|\mI_i-\mW_i\|^{2}\|\tilde{\vz}_i^{t+1}-\vz_i^t\|^{2}.
		\end{equation}
		Substituting \eqref{zproj}-\eqref{ie2w} into \eqref{nonexp}, we have
		\begin{align}\nonumber
		&\|\vz_i^t-\text{proj}_{\mathcal{X},\mathcal{Z}}\left(\vz_i^t-\left(\nabla_{\vz_i}F_i\left(\vz_i^t,\vu_i^t\right)+\mA_i^T\vlambda_i^t\right)\right)\|^{2}\\\nonumber\leq
		&\left(6\|\mI_i+\widetilde{\mW}_i\|^{2}+3\|\mI_i-\mW_i\|^{2}\right)\|\vz_i^t-\tilde{\vz}_i^{t+1}\|^{2}+6\|\widetilde{\mW}_i\|\|\vz_i^t-\vz_i^{t+1}\|^{2}+3c^{2}\|\mA_i^T\mA_i\vz_i^t\|^{2}\\\nonumber=
		&\sigma_{1}\|\vz_i^t-\tilde{\vz}_i^{t+1}\|^{2}+\sigma_{2}\|\vz_i^t-\vz_i^{t+1}\|^{2}+3c^{2}\left(N_{\max}+1\right)^{2}\|\mA_i\vz_i^t\|^{2},
		\end{align}
		where the last equality follows from \eqref{AtA}, \eqref{w} and \eqref{w_tilde}, and $\sigma_{1},\sigma_{2}$ are defined by \eqref{sigma}. Hence we finish the proof of Lemma \ref{lem1-3}.
	\end{IEEEproof}
	
	Note that it follows from \eqref{F}, Lemma \ref{lem2-6}, and  Lemma \ref{lem1-3} that the sequence $\left\{\left(\vz_i,\vu_i,\boldsymbol{\lambda}_i\right)\right\}$ generated by SP-ADMM algorithm converges to critical point of the original problem \eqref{f} when the parameters satisfy \eqref{k1_all}-\eqref{c_all}.
	
	Finally, we prove the third part of the Theorem \ref{them2}. Let us scale the upper bound of $\mathcal{F}\left(\vz,\vu,\boldsymbol{\lambda}\right)$ even further. Using the Cauchy-Schwarz inequality, we have
	\begin{align}\label{azt}
	\|\mA_i\vz_i^t\|^{2}=\|\mA_i\left(\vz_i^t-\vz_i^{t+1}+\vz_i^{t+1}\right)\|^{2}\leq2\|\mA_i\left(\vz_i^t-\vz_i^{t+1}\right)\|^{2}+2\|\mA_i\vz_i^{t+1}\|^{2}.
	\end{align}
	Together \eqref{F} with Lemma \ref{lem1-3} and \eqref{azt}, we have
	\begin{align}\label{Fbound}
	&\mathcal{F}\left(\vz^t,\vu^t,\boldsymbol{\lambda}^t\right)\\\nonumber\leq&
	\sum_{i\in\mathcal{N}}\sigma_{1}\|\tilde{\vz}_i^{t+1}-\vz_i^t\|^{2}+\sigma_{2}\|\vz_i^t-\vz_i^{t+1}\|^{2}+\|\vu_i^t-\vu_i^{t-1}\|^{2}+2\sigma_{3}\|\mA_i\left(\vz_i^t-\vz_i^{t+1}\right)\|^{2}+2\sigma_{3}\|\mA_i\vz_i^{t+1}\|^{2}\\\nonumber\leq&\sum_{i\in\mathcal{N}}\sigma_{1}\|\tilde{\vz}_i^{t+1}-\vz_i^t\|^{2}+(\sigma_{2}+2\sigma_{3}(N_{\max}+1))\|\vz_i^t-\vz_i^{t+1}\|^{2}+\|\vu_i^t-\vu_i^{t-1}\|^{2}+\frac{2\sigma_{3}}{c^2}\|\vlambda_i^{t+1}-\vlambda_i^{t}\|^{2},
	\end{align}
	where $\sigma_{1}$ and $\sigma_{2}$ are defined in \eqref{sigma}, $\sigma_{3}=3c^{2}\left(N_{\max}+1\right)^{2}+1$. The last inequality holds due to the update step of $\vlambda_i^{t+1}$ and $\|\mA_i^T\mA_i\|=N_i+1$. Matching the bounds of Lemma \ref{lem-decrease} with those of \eqref{Fbound}, we obtain
	\begin{equation}\label{ap4}
	\mathcal{F}(\vz^t,\vu^t,\boldsymbol{\lambda}^t)\leq\epsilon\left(\varsigma^t-\varsigma^{t+1}\right),
	\end{equation}
	where $\epsilon=\frac{\min\{c,1,\frac{\rho}{4},C_0,C_1,\frac{\rho}{4}-d_{\max}^2(\kappa_{1}+\kappa_{2})\}}{\max\{\sigma_{1},\sigma_{2}+2\sigma_{3}(1+N_{\max}),1,\frac{2\sigma_{3}}{c^2}\}}$. Suppose that $\mathcal{F}(\vz^t,\vu^t,\boldsymbol{\lambda}^t)$ reaches the lower bound $\epsilon_{1}$ for the first time in step $T$, with Lemma \ref{var_low_bound} and add up the inequality \eqref{ap4} of the previous $T$ step, we can get
	\begin{align}\nonumber
	\epsilon_{1}\leq&\frac{1}{T-1}\sum_{t=1}^T\mathcal{F}(\vz^t,\vu^t,\boldsymbol{\lambda}^t)\\\nonumber
	\leq&\frac{1}{T-1}\epsilon(\varsigma^{1}-\varsigma^{T+1})\\\label{epsilon2} \leq&\frac{1}{T-1}\epsilon(\varsigma^{1}-\underline{\varsigma})=\frac{\epsilon_{2}}{T-1}.
	\end{align}
	because of $\epsilon_{2}=\epsilon(\varsigma^{1}-\underline{\varsigma})>0$ is a constant, we can deduce that $\mathcal{F}(\vz^t,\vu^t,\boldsymbol{\lambda}^t)$ converges at $\mathcal{O}(1/T)$ rate.

\end{document}